\documentclass[review]{elsarticle}
\usepackage{amsthm,amsmath,amssymb,amsopn}
\usepackage{graphicx}
\usepackage{caption}
\usepackage{subcaption}
\usepackage{booktabs}
\usepackage{siunitx}
\usepackage[colorlinks,linktocpage,linkcolor=blue]{hyperref}
\usepackage[table,xcdraw]{xcolor}
\usepackage{url,color}
\usepackage{booktabs}
\usepackage{colortbl,multirow}
\usepackage{tikz}
\usetikzlibrary{shapes.geometric}
\usepackage{fontspec}
\usepackage{algorithm,algorithmic}
\numberwithin{equation}{section}
\numberwithin{table}{section}

\theoremstyle{definition}
\newtheorem{prop}{Proposition}
\newcommand\bm{\boldsymbol}
\def\dx{\mathrm{d}x}
\def\dy{\mathrm{d}y}

\def\dv{\mathrm{d}V}   

\def\det{\text{det}}

\newcommand{\gv}[1]{\ensuremath{\boldsymbol{#1}}}
\newcommand{\grad}[1]{\gv{\nabla} #1} %
\renewcommand{\div}[1]{\gv{\nabla} \cdot #1} %
\renewcommand{\d}[2]{\frac{\partial #1}{\partial #2}} %
\definecolor{lightgray}{gray}{0.9}

\journal{Journal of Computational Physics}

\begin{document}
\begin{frontmatter}
\title{A Thermodynamically Consistent High-Order Framework for Staggered Lagrangian Hydrodynamics}
\author[add1]{Zhiyuan Sun}
\ead{zysun.math@gmail.com}

\author[add1]{Jun Liu}
\ead{caepcfd@126.com}

\author[add1]{Pei Wang}
\ead{wangpei@iapcm.ac.cn}

\address[add1]{Institute of Applied Physics and Computational Mathematics, Beijing 100094, P. R. China}

\begin{abstract}
We present a consistent high-order staggered Lagrangian hydrodynamics framework designed to reconcile a fundamental mismatch in existing curvilinear formulations: the disparity between quadrature-based "strong" mass conservation and the discrete degrees of freedom (DOFs) of thermodynamic variables. By mathematically coupling the numerical quadrature rule with the density representation, our approach ensures rigorous point-wise consistency among density, internal energy, and pressure. This synchronization eliminates the ambiguity of equation-of-state (EOS) updates inherent in previous high-order staggered methods. To stabilize the discretization, we develop a high-order generalization of the subzonal pressure method by conceptually enriching the pressure field from the $Q^{m-1}$ to the $Q^m$ finite element space. We prove that evaluating this enriched field via a high-order quadrature rule naturally generates a restorative anti-hourglass force, exactly recovering the classical $Q^1-P^0$ compatible hydrodynamics algorithm as a limiting case for $m=1$. Furthermore, we introduce a concise, algorithmic formulation of tensor artificial viscosity that streamlines implementation and significantly reduces computational overhead in high-order settings. The resulting framework yields strictly diagonal mass matrices for both the momentum and energy equations, enabling highly efficient, fully explicit time integration without global linear solves. Extensive numerical benchmarks, including smooth convergence tests and complex shock-dominated flows, demonstrate that the proposed method achieves optimal high-order accuracy while maintaining superior geometric robustness.

\vspace{0.5em}
 
\vspace{0.5em}
\end{abstract}

\end{frontmatter}
\section{Introduction}\label{sec:introduction}

Lagrangian hydrodynamics, wherein the computational mesh moves concurrently with the material, serves as a fundamental approach for multi-physics and multi-material simulations~\cite{Benson1992Computational,Barlow2016Arbitrary,Chen2023Various}. Among existing spatial discretizations, staggered grid hydrodynamics (SGH) and cell-centered hydrodynamics (CCH) emerge as the two primary strategies, distinguished by the topological placement of their kinematic and thermodynamic variables. This work focuses on high-order SGH methods, wherein these variables are discretized across distinct finite element spaces.

Since the 1940s, SGH methods have achieved considerable success in one-dimensional (1D) planar shock simulations~\cite{VonNeumann1950method,Morgan2021On}, subsequently extending to two-dimensional (2D) and three-dimensional (3D) regimes~\cite{Wilkins1963Calculation,Schulz1964Two,Wilkins1999Computer,Lee2006Computer}. However, classical SGH formulations typically rely on leapfrog time integration, updating kinematic and thermodynamic variables in a temporally staggered manner. This asynchrony precludes exact total energy conservation. To resolve this, Caramana et al.~\cite{Caramana1998construction} introduced the compatible hydrodynamics algorithm (CHA), employing a predictor-corrector scheme to synchronize updates. In CHA, internal energy evolves strictly through the work defined by nodal forces and velocities, mathematically ensuring exact total energy conservation.

Most classical SGH methods are built upon the $Q^1$-$P^0$ finite element pair~\cite{Scovazzi2008Multi}: kinematic variables are defined at mesh nodes via continuous $Q^1$ elements, while thermodynamic variables are defined at cell centers via discontinuous $P^0$ elements. For shock-dominated flows, artificial viscosity is essential to mediate discontinuities~\cite{VonNeumann1950method,Landshoff1955,Caramana1998Formulations}, while multi-dimensional simulations require hourglass control to suppress nonphysical, zero-energy deformation modes~\cite{Menchen1964tensor,Flanagan1981uniform,Caramana1998Elimination}.

Driven by growing computational power, significant attention has shifted toward next-generation, high-fidelity hydrocodes~\cite{Anderson2020Multiphysics}. Dobrev et al.~\cite{Dobrev2011Curvilinear,Dobrev2012High} proposed a generalized high-order Lagrangian finite element framework founded on three key principles: (1) kinematic variables use continuous high-order spaces, while thermodynamic variables use discontinuous spaces (typically one order lower); (2) mass conservation is enforced at individual quadrature points (the "strong mass conservation" principle); and (3) artificial viscosity is incorporated into the momentum equation to satisfy Rankine–Hugoniot conditions. This framework has since been extended to axisymmetric problems~\cite{Dobrev2013Highorder}, elastic-plastic models~\cite{DOBREV2014High}, and multi-material ALE hydrodynamics~\cite{Dobrev2016Multimaterial,Anderson2018High}, with a robust implementation available in the open-source miniapp Laghos~\cite{Laghos2024,Anderson2021MFEM}.

Despite these developments, the strong mass conservation principle remains insufficiently clarified. Enforcing mass conservation strictly at quadrature points introduces an implicit dependence of the density degrees of freedom (DOFs) on the chosen quadrature rule. This creates a fundamental mismatch between the DOFs of the density field and those of other thermodynamic variables. For example, utilizing a $Q^m$-$Q^{m-1}$ pair, any quadrature rule exceeding $m^2$ points implies the density DOFs are effectively associated with a $Q^m$ space, erroneously exceeding the bounds of the $Q^{m-1}$ thermodynamic space.

This mathematical inconsistency is evident even in the baseline $Q^1$-$P^0$ case. The framework in~\cite{Dobrev2011Curvilinear} cannot directly recover the classical CHA~\cite{Caramana1998construction} without a highly specific low-order assumption; doing so instead yields a scheme equivalent to that of Cheng et al.~\cite{Cheng2015Elimination}. In such formulations, each element contains four density DOFs but only one internal energy DOF. This structural disparity results in severe inconsistencies during equation-of-state (EOS) updates, where the computed thermodynamic state may no longer represent a physically realizable configuration.

The primary objective of this work is to definitively resolve this inconsistency, establishing a unified, high-order SGH framework that is both mathematically compatible and thermodynamically consistent. We explicitly couple the underlying quadrature rule, density representation, and thermodynamic discretization, completely eliminating the ambiguity inherent in quadrature-based mass conservation.

The main contributions of the proposed framework are summarized as follows:
\begin{enumerate}
\item[{\bf C1};] A consistent $Q^m$-$Q^{m-1}$ discretization uniquely tying the quadrature rule to the density DOFs, guaranteeing strict consistency between density and thermodynamic variables.
\item[{\bf C2};] The exact recovery of the classical compatible hydrodynamics algorithm as the $Q^1$-$P^0$ special case, without requiring additional low-order assumptions.
\item[{\bf C3};] A unified interpretation of quadrature, DOFs, and EOS consistency, realized by aligning all thermodynamic variables with the specific quadrature points.
\item[{\bf C4};] A high-order extension of the subzonal pressure method, obtained by conceptually enriching the pressure field from $Q^{m-1}$ to $Q^m$, which naturally recovers the classical anti-hourglass formulation when $m=1$.
\item[{\bf C5};] A concise, explicitly computable formulation of tensor artificial viscosity that drastically improves implementation efficiency in high-order settings.
\item[{\bf C6};] A unified treatment of artificial viscosity and hourglass control built upon a shared quadrature and interpolation architecture.
\end{enumerate}

Within this framework, both the momentum and energy mass matrices are perfectly diagonal. This mass lumping eliminates the need for expensive global linear solves, enabling highly efficient explicit time integration. Furthermore, co-locating the density, internal energy, and pressure ensures unconditionally robust EOS updates

Nevertheless, high-order SGH remains susceptible to hourglass distortion—a nonphysical geometric manifestation of discrete rank deficiencies~\cite{margolin2013arbitrary,Menchen1964tensor}. Classical controls rely on viscosity-based dampening~\cite{Flanagan1981uniform,Margolin1987method} or stiffness-based restoring forces~\cite{Caramana1998Elimination,Cheng2015Elimination,Sun2022Understanding}. While previous literature suggested strong mass conservation might naturally alleviate these modes~\cite{Anderson2020Multiphysics}, achieving strict thermodynamic consistency inevitably reintroduces this geometric instability. To resolve this, we extend the methodology of~\cite{Sun2022Understanding} into the high-order regime. By evaluating an enriched higher-order pressure field via a higher-order quadrature rule, we generate a mathematical anti-hourglass force. Consistent with the low-order case, these pressure variations are predicted using subzonal density variations and the local sound speed, elegantly generalizing the classical subzonal pressure method~\cite{Caramana1998Elimination}.

Artificial viscosity is incorporated following standard high-order tensor approaches~\cite{Kolev2009tensor,Dobrev2012High}. Rather than introducing a new viscous model, we focus on algorithmic efficiency. Inspired by classical explicit hydrocodes~\cite{Hallquist1980User},  we propose a concise formulation that balances implementation simplicity with computational cost, significantly reducing programming complexity and memory bottlenecks in high-order settings. Because both artificial viscosity and hourglass control dictate geometric robustness, we employ a shared $(m+1)^2$-point quadrature rule for both mechanisms to minimize the computational overhead.

The remainder of this paper is organized as follows. Section~\ref{sec:discretization} introduces the governing Euler equations and the chosen finite element spaces. Section~\ref{sec:gauss_quadrature} analyzes the critical role of quadrature in maintaining exact mass conservation. Sections~\ref{sec:momentum_conservation} and~\ref{sec:energy_conservation} present the explicit discretizations for the momentum conservation law and the internal energy balance equation. Section~\ref{sec:hourglass_and_viscosity} details the coupled construction of the hourglass control algorithm and the artificial viscosity tensor. Section~\ref{sec:original_discretization} introduces a mass-lumped discretization of the original high-order SGH formulation~\cite{Dobrev2012High} as a comparative baseline. Section~\ref{sec:numerical_results} briefly outlines the implicit-explicit Runge-Kutta (IMEX-RK) time-integration schemes utilized to preserve total energy (with full stage details provided in Appendix~\ref{app:imex_rk}), and subsequently presents extensive numerical validations using both smooth convergence tests and shock-dominated benchmark problems. Finally, conclusions are drawn in Section~\ref{sec:conclusion}.

\section{Compressible Euler equations and discretization space}\label{sec:discretization}

The compressible Euler equations under a Lagrangian frame are given by
\begin{equation}\label{eq:lagrangian}
    \begin{aligned}
    \frac{d \rho}{d t}  &=-\rho \nabla \cdot \vec{u}, \\
    \rho \frac{d \vec{u} }{d t}    &= - \nabla p,\\
    \rho \frac{d e}{d t}   &= - p \nabla \cdot \vec{u},\\
    p &= \text{EOS}(\rho,e).
    \end{aligned}
\end{equation}
where $d/dt$ denotes the material derivative, and $e$ is the specific internal energy. 

Here, we focus on the two-dimensional Cartesian geometry, denoting the spatial coordinates by $(x,y)$. Cylindrical geometries will be considered in future work. In this paper, we present a unified framework for high-order staggered grid hydrodynamics (SGH) employing a $Q^m-Q^{m-1}$ finite element space pair. The classical $Q^1-P^0$ SGH method is exactly recovered when $m=1$, requiring no additional assumptions.

Following the definitions in~\cite{Dobrev2012High}, $Q^m$ denotes the tensor-product continuous space of polynomials of degree not exceeding $m$ in each coordinate direction, while $Q^{m-1}$ denotes the tensor-product discontinuous space of polynomials of degree not exceeding $m-1$. We utilize the $Q^m$ space to discretize kinematic variables (velocity and position) and the $Q^{m-1}$ space for thermodynamic variables (density, internal energy, and pressure). We denote by $\bm{V}_h$ and $Q_h$ the respective FEM approximation spaces:
\begin{align*}
\bm{V}_h &\triangleq \left\{ \vec{v}_h \in [H^1(\Omega)]^2 : \vec{v}_h|_K \in [Q^m(K)]^2, \text{ continuous} \right\},\\
Q_h &\triangleq \left\{ q_h \in L^2(\Omega) : q_h|_K \in Q^{m-1}(K), \text{ discontinuous} \right\},
\end{align*}
where $K$ indexes the elements in the grid. For a given element $K$, the shape functions of the conforming $Q^m$ space are denoted by $N_i$, while those of the discontinuous $Q^{m-1}$ space are denoted by $\phi_k$. Consequently, the physical variables within an element are expressed as:
\begin{align*}
        \vec{u}_h|_K &=\sum_{i=1}^{\text{kdof}} \vec{u}_{i}N_{i},\\
        \rho_h|_K &=\sum_{k=1}^{\text{tdof}} \rho_{k}\phi_{k},\\
        p_h|_K &=\sum_{k=1}^{\text{tdof}} p_{k}\phi_{k},\\
        e_h|_K &=\sum_{k=1}^{\text{tdof}} e_{k}\phi_{k}.
\end{align*}
Let $\text{kdof}$ and $\text{tdof}$ represent the total number of degrees of freedom (DOFs) per element for the kinematic and thermodynamic spaces, respectively. The coefficients $\vec{u}_{i}$ are the velocities at the kinematic DOFs, while $\rho_{k}$, $p_{k}$, and $e_{k}$ are the respective thermodynamic values evaluated at the discontinuous DOFs. 

\subsection{Quadrature rules and nodal arrangements}
The arrangement of DOFs within this $Q^m-Q^{m-1}$ framework is intentionally non-standard. The DOFs for the kinematic variables are collocated with the $(m+1)^2$ Gauss-Lobatto quadrature points, whereas the thermodynamic variables are collocated with the $m^2$ Gauss-Legendre quadrature points. 

Before detailing specific element pairs, we briefly state the properties of the 1D Gauss quadrature rules used to construct these shape functions~\cite{Press2007Numerical}. On the reference interval $[-1,1]$, an $n$-point Gauss-Lobatto quadrature rule (which includes the boundary endpoints) achieves a maximum degree of exactness of $2n-3$. Conversely, an $n$-point Gauss-Legendre rule (which does not include the endpoints) achieves a degree of exactness of $2n-1$. This exactness ensures the rule integrates all polynomials up to the stated degree precisely. Table~\ref{tab:quadrature_combined} lists the 1D Gauss-Legendre and Gauss-Lobatto quadrature rules for various orders. In the remainder of this paper, we will use "Gauss quadrature" to denote Gauss-Legendre quadrature and "Lobatto quadrature" to denote Gauss-Lobatto quadrature for brevity.

\begin{table}[htbp]
    \centering
    \caption{Gauss-Legendre (Gauss) and Gauss-Lobatto (Lobatto) quadrature rules defined on $[-1,1]$.}
    \label{tab:quadrature_combined}
    \small
    \renewcommand{\arraystretch}{1.2}
    \begin{tabular}{c c c}
        \toprule
        \multicolumn{3}{c}{Gauss-Legendre (Degree of exactness $2n-1$)} \\
        \midrule
        $n$ & Points & Weights \\
        \midrule
        1 & $0$ & $2$ \\
        2 & $\pm\sqrt{1/3}$ & $1$ \\
        3 & $0,\ \pm\sqrt{3/5}$ & $8/9,\ 5/9$ \\
        4 & $\pm\sqrt{\frac{3+2\sqrt{6/5}}{7}},\ \pm\sqrt{\frac{3-2\sqrt{6/5}}{7}}$ & $\frac{1}{2}-\frac{1}{6\sqrt{6/5}},\ \frac{1}{2}+\frac{1}{6\sqrt{6/5}}$ \\
        \midrule
        \multicolumn{3}{c}{Gauss-Lobatto (Degree of exactness $2n-3$)} \\
        \midrule
        $n$ & Points & Weights \\
        \midrule
        2 & $\pm 1$ & $1$ \\
        3 & $0,\ \pm 1$ & $4/3,\ 1/3$ \\
        4 & $\pm 1,\ \pm\sqrt{1/5}$ & $1/6,\ 5/6$ \\
        5 & $0,\ \pm 1,\ \pm\sqrt{3/7}$ & $32/45,\ 1/10,\ 49/90$ \\
        \bottomrule
    \end{tabular}
\end{table}

For higher-dimensional geometries, such as quadrilateral or hexahedral elements, multiple integrals are evaluated as iterated integrals, applying the 1D rules across each coordinate direction independently. This formulation yields tensor-product quadrature rules. On the 2D reference element $[-1,1]^2$, an $n^2$-point tensor-product Gauss quadrature rule achieves exactness for polynomials of the form $x^p y^q$ where $p, q \leq 2n-1$. The 2D shape functions defined below are constructed directly via the tensor product of the 1D Lagrange polynomials associated with the respective quadrature points. We now provide specific examples of the DOF arrangements and basis functions for varying orders of $m$.

\subsection{Examples of $Q^m-Q^{m-1}$ space pairs}

\textbf{The $Q^1-P^0$ space pair:} Utilized in classical SGH, the kinematic variables are discretized using a bilinear conforming $Q^1$ space, while the thermodynamic variables employ a piecewise constant discontinuous $P^0$ space. Figure~\ref{fig:Q_1_P_0_reference_element} illustrates the reference element. The kinematic DOFs are located at the four vertices (identical to the $2 \times 2$ Lobatto quadrature points). The single thermodynamic DOF is located at the element's center (the $1 \times 1$ Gauss quadrature point). The corresponding $Q^1$ shape functions are:

\begin{equation}\label{eq:Q1_basis_kinematic_1D}
\begin{split}
&\xi_a =\frac{1}{2}(1\pm \xi),\ \eta_b= \frac{1}{2}(1\pm\eta), \quad a,b=1,2.\\
&N_i=\xi_a \eta_b, \quad  i=1,2,3,4.
\end{split}
\end{equation}

The DOFs are enumerated using a lexicographical (row-major) ordering, mapping $i = a + 2(b-1)$.

\begin{figure}[htbp]
    \begin{center}
    \begin{tikzpicture}[scale = 1]
    \fill[black] (-3,-3) circle (2pt);
    \fill[black] (3,-3) circle (2pt);
    \fill[black] (3,3) circle (2pt);
    \fill[black] (-3,3) circle (2pt);
    
    \draw[black] (-3,-3) -- (3,-3);
    \draw[black] (3,-3) -- (3,3);
    \draw[black] (3,3) -- (-3,3);
    \draw[black] (-3,3) -- (-3,-3);
    \node[] at (-3,-3.5) {$1,(-1,-1)$};
    \node[] at (3,-3.5)  {$2,(1,-1)$};
    \node[] at (-3,3.5)  {$3,(-1,1)$};
    \node[] at (3,3.5)   {$4,(1,1)$};

    \draw[red,fill=red] (0,0) circle (2pt);
    \node[] at (0,-.5) {\textcolor{red}{$1,(0,0)$}};
    \end{tikzpicture}
    \caption{Reference element of the $Q^1-P^0$ FEM space pair, showing the DOFs of the kinematic variable (black circles) and the thermodynamic variable (red circle).}
    \label{fig:Q_1_P_0_reference_element}
    \end{center}
\end{figure} 

\textbf{The $Q^2-Q^1$ space pair:} Figure~\ref{fig:Q_2_Q_1_reference_element} depicts the $Q^2-Q^1$ element. The $Q^2$ element aligns perfectly with the standard nine-point Lobatto quadrature nodes. The thermodynamic DOFs lie on the four interior Gauss quadrature points. Since the 1D Lagrange polynomials for $\xi$ and $\eta$ are structurally identical, we only list those for $\xi$ for brevity. The kinematic shape functions are defined as:
\begin{equation}\label{eq:Q2_basis_kinematic_1D}
    \begin{split}
    &\xi_1 =-\frac{1}{2}\xi(1-\xi),\ \xi_2 =(1-\xi)(1+\xi), \ \xi_3 =\frac{1}{2}\xi(1+\xi).\\
    &N_i=\xi_a\eta_b, \quad i= 1,...,9, \ a,b=1,2,3.
    \end{split}
\end{equation}

The discontinuous $Q^1$ shape functions for the thermodynamic variables, denoted by $\phi_j$, are constructed using the 1D polynomials $\bar{\xi}_k,\bar{\eta}_l$:
\begin{equation}\label{eq:Q1_basis_thermodynamic_1D}
    \begin{split}
    &\bar{\xi}_k =\frac{\sqrt{3}}{2}(\frac{1}{\sqrt{3}}\pm \xi), \bar{\eta}_l= \frac{\sqrt{3}}{2}(\frac{1}{\sqrt{3}}\pm \eta) , k,l=1,2.\\
    &\phi_{j} =\bar{\xi}_k \bar{\eta}_l, j=1,2,3,4.
    \end{split}
\end{equation}

\begin{figure}[htbp]
    \begin{center}
    \begin{tikzpicture}[scale = 1]
    \fill[black] (-3,-3) circle (2pt);
    \fill[black] (0,-3) circle (2pt);
    \fill[black] (3,-3) circle (2pt);

    \fill[black] (-3,0) circle (2pt);
    \fill[black] (0,0) circle (2pt);   
    \fill[black] (3,0) circle (2pt); 

    \fill[black] (-3,3) circle (2pt);
    \fill[black] (0,3) circle (2pt);                           
    \fill[black] (3,3) circle (2pt);
    
    \draw[black] (-3,-3) -- (3,-3);
    \draw[black] (3,-3) -- (3,3);
    \draw[black] (3,3) -- (-3,3);
    \draw[black] (-3,3) -- (-3,-3);
    \node[] at (-3,-3.5) {$1,(-1,-1)$};
    \node[] at (0,-3.5) {$2,(0,-1)$};
    \node[] at (3,-3.5)  {$3,(1,-1)$};

    \node[] at (-3.75,0)  {$4,(-1,0)$};
    \node[] at (0.75,0)  {$5,(0,0)$};
    \node[] at (3.75,0)  {$6,(1,0)$};

    \node[] at (-3,3.5)  {$7,(-1,1)$};
    \node[] at (0,3.5)  {$8,(0,1)$};        
    \node[] at (3,3.5)   {$9,(1,1)$};

    \draw[red,fill=red] (-1.732,-1.732) circle (2pt);
    \draw[red,fill=red] (1.732,-1.732) circle (2pt);
    \draw[red,fill=red] (-1.732,1.732) circle (2pt);
    \draw[red,fill=red] (1.732,1.732) circle (2pt);
    \node[] at (-1.732,-2.232) {\textcolor{red}{$1,(-\frac{1}{\sqrt{3}},-\frac{1}{\sqrt{3}})$}};
    \node[] at (1.732,-2.232) {\textcolor{red}{$2,(\frac{1}{\sqrt{3}},-\frac{1}{\sqrt{3}})$}};
    \node[] at (-1.732,2.232) {\textcolor{red}{$3,(-\frac{1}{\sqrt{3}},\frac{1}{\sqrt{3}})$}};
    \node[] at (1.732,2.232) {\textcolor{red}{$4,(\frac{1}{\sqrt{3}},\frac{1}{\sqrt{3}})$}};
    
    \end{tikzpicture}
    \caption{Reference element of the $Q^2-Q^1$ FEM space pair, showing the DOFs of the kinematic variables (black circles) and the thermodynamic variables (red circles).}
    \label{fig:Q_2_Q_1_reference_element}
    \end{center}
\end{figure}

\textbf{The $Q^3-Q^2$ space pair:} The $Q^3-Q^2$ element (Figure~\ref{fig:Q_3_Q_2_reference_element}) departs from standard nodal definitions. Kinematic DOFs are strictly tied to the sixteen Lobatto quadrature points, and thermodynamic DOFs are tied to the nine Gauss quadrature points. For visual clarity, the reference coordinates of the interior DOFs are omitted from the figure but are strictly defined by the tensor product of the 1D Lobatto nodes $\{-1, -\sqrt{1/5}, \sqrt{1/5}, 1\}$ and 1D Gauss nodes $\{-\sqrt{3/5}, 0, \sqrt{3/5}\}$. 

The 1D Lagrange polynomials utilized to construct the tensor-product $Q^3$ element are:
\begin{equation}\label{eq:Q3_basis_kinematic_1D}
\begin{split}
&\xi_1 =-\frac{5}{8}(\xi-1)\left(\xi-\frac{1}{\sqrt{5}}\right)\left(\xi+\frac{1}{\sqrt{5}}\right), \quad \xi_2 =\frac{5\sqrt{5}}{8}(\xi-1)\left(\xi-\frac{1}{\sqrt{5}}\right)(\xi+1),\\
&\xi_3 =-\frac{5\sqrt{5}}{8}(\xi-1)\left(\xi+\frac{1}{\sqrt{5}}\right)(\xi+1), \quad \xi_4 =\frac{5}{8}\left(\xi-\frac{1}{\sqrt{5}}\right)\left(\xi+\frac{1}{\sqrt{5}}\right)(\xi+1).
\end{split}
\end{equation}
 The 2D shape functions for the kinematic variables are defined directly by their tensor products:
\begin{equation*}\label{eq:Q3_basis_kinematic_2D}
    N_i=\xi_a\eta_b, i=1,...,16, \ a,b=1,2,3,4.
 \end{equation*} 

The non-standard $Q^2$ element shape functions for the thermodynamic variables are listed as follows:
\begin{equation}\label{eq:Q2_basis_thermodynamic_1D}
    \bar{\xi}_1 =-\frac{5}{6}\xi\left(\sqrt{\frac{3}{5}}-\xi\right), \ \bar{\xi}_2 =\frac{5}{3}\left(\sqrt{\frac{3}{5}}+\xi\right)\left(\sqrt{\frac{3}{5}}-\xi\right), \  \bar{\xi}_3 =\frac{5}{6}\xi\left(\sqrt{\frac{3}{5}}+\xi\right).
\end{equation}
\begin{equation}\label{eq:Q2_basis_thermodynamic}
    \phi_{j} =\bar{\xi}_k \bar{\eta}_l, \ j=1,...,9,\ k,l=1,2,3.
\end{equation}

\begin{figure}[htbp]
\begin{center}
\begin{tikzpicture}[scale = 1]
    \fill[black] (-3,-3) circle (2pt);
    \fill[black] (-1.3416,-3) circle (2pt);
    \fill[black] (1.3416,-3) circle (2pt);
    \fill[black] (3,-3) circle (2pt);

    \fill[black] (-3,-1.3416) circle (2pt);
    \fill[black] (-1.3416,-1.3416) circle (2pt);
    \fill[black] (1.3416,-1.3416) circle (2pt);
    \fill[black] (3,-1.3416) circle (2pt);

    \fill[black] (-3,1.3416) circle (2pt);
    \fill[black] (-1.3416,1.3416) circle (2pt);
    \fill[black] (1.3416,1.3416) circle (2pt);
    \fill[black] (3,1.3416) circle (2pt);
    
    \fill[black] (-3,3) circle (2pt);
    \fill[black] (-1.3416,3) circle (2pt);
    \fill[black] (1.3416,3) circle (2pt);
    \fill[black] (3,3) circle (2pt);
    
    \draw[black] (-3,-3) -- (3,-3);
    \draw[black] (3,-3) -- (3,3);
    \draw[black] (3,3) -- (-3,3);
    \draw[black] (-3,3) -- (-3,-3);
    
    \node[] at (-3.25,-3.5) {$1,(-1,-1)$};
    \node[] at (-1.1416,-3.5) {$2,(-\sqrt{1/5},-1)$};
    \node[] at (1.1416,-3.5) {$3,(\sqrt{1/5},-1)$};
    \node[] at (3.25,-3.5)  {$4,(1,-1)$};

    \node[] at (-4.25,-1.3416)  {$5,(-1,-\sqrt{1/5})$};
    \node[] at (-1.5416,-1.3416)  {$6$};
    \node[] at (1.5416,-1.3416)  {$7$};
    \node[] at (4.15,-1.3416)  {$8,(1,-\sqrt{1/5})$};

    \node[] at (-4.25,1.3416)  {$9,(-1,\sqrt{1/5})$};
    \node[] at (-1.6416,1.3416)  {$10$};
    \node[] at (1.6416,1.3416)  {$11$};
    \node[] at (4.15,1.3416)  {$12,(1,\sqrt{1/5})$};

    \node[] at (-3.25,3.5) {$13,(-1,1)$};
    \node[] at (-1.1416,3.5) {$14,(-\sqrt{1/5},1)$};
    \node[] at (1.1416,3.5) {$15,(\sqrt{1/5},1)$};
    \node[] at (3.25,3.5)  {$16,(1,1)$};

    \draw[red,fill=red] (-2.3238,-2.3238) circle (2pt);
    \draw[red,fill=red] (0,-2.3238) circle (2pt);
    \draw[red,fill=red] (2.3238,-2.3238) circle (2pt);
    \draw[red,fill=red] (-2.3238,0) circle (2pt);
    \draw[red,fill=red] (0,0) circle (2pt);
    \draw[red,fill=red] (2.3238,0) circle (2pt);
    \draw[red,fill=red] (-2.3238,2.3238) circle (2pt);
    \draw[red,fill=red] (0,2.3238) circle (2pt);
    \draw[red,fill=red] (2.3238,2.3238) circle (2pt);
    
    \node[] at (-2.3238,-2.6238) {\textcolor{red}{$1$}};
    \node[] at (0,-2.6238) {\textcolor{red}{$2$}};
    \node[] at (2.3238,-2.6238) {\textcolor{red}{$3$}};

    \node[] at (-2.3238,-0.3) {\textcolor{red}{$4$}};
    \node[] at (0,-0.3) {\textcolor{red}{$5$}};
    \node[] at (2.3238,-0.3) {\textcolor{red}{$6$}};
    
    \node[] at (-2.3238,2.6238) {\textcolor{red}{$7$}};
    \node[] at (0,2.6238) {\textcolor{red}{$8$}};
    \node[] at (2.3238,2.6238) {\textcolor{red}{$9$}};

\end{tikzpicture}
\caption{Reference element of the $Q^3-Q^2$ FEM space pair, showing the DOFs of the kinematic variables (black circles) and the thermodynamic variables (red circles).}
\label{fig:Q_3_Q_2_reference_element}
\end{center}
\end{figure}
\section{The Gauss quadrature rule and the mass conservation law}\label{sec:gauss_quadrature}
We now turn to the discretization of the mass conservation law in Equation~\eqref{eq:lagrangian}. Its integral form over an element $K$ is given by
\begin{equation}\label{eq:mass_integral}
\frac{d}{d t}\int_{K} \rho_h \dv =0,
\end{equation}
where $V$ denotes the volume. Given the delicate nature of numerically discretizing Equation~\eqref{eq:mass_integral}, we assume a specified Gauss quadrature rule ($\omega_q, \gv{\xi}_q$) on the reference element. Here, $\omega_q$ and $\gv{\xi}_q$ represent the weight and local coordinate of the $q$-th quadrature point, respectively, ensuring that mass conservation is satisfied locally at each individual point:
\begin{equation}\label{eq:discrete_mass_conservation}
\omega_q \rho^{n}_h(\gv{\xi}_q) \det J^{n}(\gv{\xi}_q) = \omega_q \rho^{n-1}_h(\gv{\xi}_q) \det J^{n-1}(\gv{\xi}_q) = \omega_q \rho_{0}(\gv{\xi}_q) \det J_{0}(\gv{\xi}_q),
\end{equation}
where $\det J$ is the Jacobian determinant, the superscript $n$ denotes the value at time $t^n$, and $\rho_0$ is the initial density. A more compact formulation is given by:
\begin{equation}\label{eq:discrete_mass_conservation_simple}
\rho^{n}_h(\gv{\xi}_q) \det J^{n}(\gv{\xi}_q)=\rho_{0}(\gv{\xi}_q) \det J_{0}(\gv{\xi}_q).
\end{equation}

The primary objective of this section is to mathematically establish the relationship between the chosen Gauss quadrature rule and the degrees of freedom (DOFs) of the density variable. Within our high-order framework, an under-integration strategy is intentionally adopted to maintain strict consistency among the thermodynamic variables. Concurrently, a geometrically refined density field—captured using more precise Gauss quadrature rules—can be explicitly utilized to construct the hourglass control algorithm. For the $Q^m-Q^{m-1}$ FEM space pair, the core principles are summarized as follows:

\begin{itemize}
\item An $m^2$-point Gauss quadrature rule is adopted to solve the mass conservation law~\eqref{eq:discrete_mass_conservation_simple}, resulting in exactly $m^2$ DOFs for the density variable. This perfectly matches the thermodynamic DOFs ($\text{tdof}$).
\item The refined density field is obtained by applying an $(m+1)^2$-point (or higher-order) Gauss quadrature rule, yielding $(m+1)^2$ DOFs for the density variable. This matches the kinematic DOFs ($\text{kdof}$).
\item This refined density field serves solely as an auxiliary variable for constructing the hourglass control algorithm. It cannot be used directly to update the pressure variable via the EOS, as doing so would critically violate consistency among the thermodynamic variables.
\end{itemize}

We now proceed with case-by-case proofs to demonstrate these principles.

\subsection{The relationship between the quadrature rule and the number of DOFs of the density variable}
We first analyze the relationship between the chosen Gauss quadrature rule and the available DOFs of the density variable. The general proof strategy is as follows: since mass conservation is enforced at each quadrature point, the density variation fundamentally depends on the Jacobian determinant. The Jacobian determinant, in turn, can be expressed as a linear combination of shape function values and several linearly independent geometric vectors. Consequently, the number of independent DOFs for the density variable strictly equals the rank of the matrix formed by evaluating the shape functions at the selected quadrature points. We use the 2D FEM space pairs $Q^1-P^0$ and $Q^2-Q^1$ as representative examples to illustrate these conclusions.

\textbf{The $Q^1-P^0$ case:}
For the classical $Q^1-P^0$ 2D FEM space pair, the reference element and the corresponding Gauss quadrature points for $n \times n$ rules (where $n=1,2,3$) are displayed in Figure~\ref{fig:Q_1_reference_element}. The specific quadrature rules are visualized by color: one-point (red), four-point (blue), and nine-point (green and red).

\begin{figure}[htbp]
\begin{center}
\begin{tikzpicture}[scale = 1]
\fill[black] (-3,-3) circle (2pt);
\fill[black] (3,-3) circle (2pt);
\fill[black] (3,3) circle (2pt);
\fill[black] (-3,3) circle (2pt);

\draw[black] (-3,-3) -- (3,-3);
\draw[black] (3,-3) -- (3,3);
\draw[black] (3,3) -- (-3,3);
\draw[black] (-3,3) -- (-3,-3);
\node[] at (-3,-3.5) {$1,(-1,-1)$};
\node[] at (3,-3.5)  {$2,(1,-1)$};
\node[] at (-3,3.5)  {$3,(-1,1)$};
\node[] at (3,3.5)   {$4,(1,1)$};

\draw[green,fill=red] (0,0) circle (2pt);
\draw[black,fill=blue] (-1.732508075688,-1.732508075688) circle (2pt);
\draw[black,fill=blue] (1.732508075688,-1.732508075688) circle (2pt);
\draw[black,fill=blue] (-1.732508075688,1.732508075688) circle (2pt);
\draw[black,fill=blue] (1.732508075688,1.732508075688) circle (2pt);

\draw[black,fill=green] (-2.32379,-2.32379) circle (2pt);
\draw[black,fill=green] (0,-2.32379) circle (2pt);
\draw[black,fill=green] (2.32379,-2.32379) circle (2pt);
\draw[black,fill=green] (-2.32379,0) circle (2pt);
\draw[black,fill=green] (2.32379,0) circle (2pt);
\draw[black,fill=green] (-2.32379,2.32379) circle (2pt);
\draw[black,fill=green] (0,2.32379) circle (2pt);
\draw[black,fill=green] (2.32379,2.32379) circle (2pt);
\end{tikzpicture}
\caption{Reference element $(\xi, \eta)$ of the $Q^1$ element and its 2D Gauss quadrature points.}
\label{fig:Q_1_reference_element}
\end{center}
\end{figure}

The associated Jacobian matrix is calculated as:
\begin{equation*}
J=
\begin{bmatrix}
\d{x}{\xi} &  \d{x}{\eta}\\
\d{y}{\xi} &  \d{y}{\eta}
\end{bmatrix}
=
\begin{bmatrix}
\sum_{i=1}^{4} x_i \d{N_i}{\xi} & \sum_{i=1}^{4} x_i \d{N_i}{\eta} \\
\sum_{i=1}^{4} y_i \d{N_i}{\xi} & \sum_{i=1}^{4} y_i \d{N_i}{\eta}
\end{bmatrix}.
\end{equation*}
Here, the Jacobian determinant can be explicitly expressed as a linear combination of the shape function values:
\begin{equation}\label{eq:Q1_Jacobian}
\begin{split}
\det J=& \frac{1}{4} N_1(-x_{21}y_{13}+x_{13}y_{21}) + \frac{1}{4} N_2(-x_{42}y_{21}+x_{21}y_{42})\\
+& \frac{1}{4} N_3(-x_{13}y_{34}+x_{34}y_{13}) + \frac{1}{4} N_4(-x_{34}y_{42}+x_{42}y_{34}),
\end{split}
\end{equation}
where $x_{ij}=x_i-x_j$ and $y_{ij}=y_i-y_j$.

Equation~\eqref{eq:Q1_Jacobian} carries an important geometric interpretation: each of the four bracketed terms—such as $(-x_{21}y_{13}+x_{13}y_{21})$—corresponds to exactly twice the area of a triangle defined by three adjacent vertices. This relationship can be cleanly organized in matrix form as an inner product:
\begin{equation}
\det J =\frac{1}{4}
\begin{bmatrix}
N_1 & N_2 & N_3 & N_4
\end{bmatrix}
\begin{bmatrix}
-x_{21}y_{13}+x_{13}y_{21}\\
-x_{42}y_{21}+x_{21}y_{42}\\
-x_{13}y_{34}+x_{34}y_{13}\\
-x_{34}y_{42}+x_{42}y_{34}
\end{bmatrix}.
\end{equation}

The classical $Q^1-P^0$ SGH employs a one-point Gauss quadrature rule (shown as the red point in Figure~\ref{fig:Q_1_reference_element}), for which the Jacobian determinant reduces to:
\begin{equation*}
\begin{split}
\det J(0,0) &= \frac{1}{16} \Big\{ (-x_{21}y_{13}+x_{13}y_{21}) + (-x_{42}y_{21}+x_{21}y_{42}) \\
&\quad + (-x_{13}y_{34}+x_{34}y_{13}) + (-x_{34}y_{42}+x_{42}y_{34}) \Big\}\\
&= \frac{1}{4} \left\{ \frac{1}{2}(x_{41}y_{32}-x_{32}y_{41}) \right\},
\end{split}
\end{equation*}
where $\frac{1}{2}(x_{41}y_{32}-x_{32}y_{41})$ represents the total area of the quadrilateral element. Consequently, each element mathematically supports only a single calculated density variable, perfectly aligning with the DOFs of the $P^0$ element space.

If a four-point Gauss quadrature rule is employed instead (shown as the blue points in Figure~\ref{fig:Q_1_reference_element}), the corresponding Jacobian determinant depends on multiple shape function values. We provide a complete listing of these evaluations in Table~\ref{table:Q1_shape_function_four_point_Gauss}.

\begin{table}[htbp]
\centering
\caption{Shape function values of the $Q_1$ element at the four Gauss quadrature points.}
\label{table:Q1_shape_function_four_point_Gauss}
\renewcommand{\arraystretch}{1.5}
\begin{tabular}{|c|c|c|c|c|}
\hline
Point $\backslash$ Function & $N_1$ &  $N_2$ & $N_3$ & $N_4$   \\ \hline
$(-\sqrt{1/3},-\sqrt{1/3})$ & $\frac{1}{3}+\frac{1}{2\sqrt{3}}$  & $\frac{1}{6}$ & $\frac{1}{6}$ & $\frac{1}{3}-\frac{1}{2\sqrt{3}}$   \\ \hline
$(\sqrt{1/3},-\sqrt{1/3})$ & $\frac{1}{6}$ & $\frac{1}{3}+\frac{1}{2\sqrt{3}}$ & $\frac{1}{3}-\frac{1}{2\sqrt{3}}$ & $\frac{1}{6}$\\ \hline
$(-\sqrt{1/3},\sqrt{1/3})$ & $\frac{1}{6}$  & $\frac{1}{3}-\frac{1}{2\sqrt{3}}$& $\frac{1}{3}+\frac{1}{2\sqrt{3}}$ & $\frac{1}{6}$  \\ \hline
$(\sqrt{1/3},\sqrt{1/3})$ & $\frac{1}{3}-\frac{1}{2\sqrt{3}}$ & $\frac{1}{6}$ & $\frac{1}{6}$ &$\frac{1}{3}+\frac{1}{2\sqrt{3}}$  \\ \hline
\end{tabular}
\end{table}

These shape function values construct the following evaluation matrix:
\begin{equation*}
\begin{bmatrix}
\frac{1}{3}+\frac{1}{2\sqrt{3}}  & \frac{1}{6} & \frac{1}{6} & \frac{1}{3}-\frac{1}{2\sqrt{3}}   \\
\frac{1}{6} & \frac{1}{3}+\frac{1}{2\sqrt{3}} & \frac{1}{3}-\frac{1}{2\sqrt{3}} & \frac{1}{6}\\
\frac{1}{6}  & \frac{1}{3}-\frac{1}{2\sqrt{3}} & \frac{1}{3}+\frac{1}{2\sqrt{3}} & \frac{1}{6}  \\
\frac{1}{3}-\frac{1}{2\sqrt{3}} & \frac{1}{6} & \frac{1}{6} &\frac{1}{3}+\frac{1}{2\sqrt{3}}  
\end{bmatrix}.
\end{equation*}
The rank of this matrix is 4. This implies that the four evaluated Jacobian determinants,
\begin{displaymath}
\det J\left(-\sqrt{1/3},-\sqrt{1/3}\right), \ \det J\left(\sqrt{1/3},-\sqrt{1/3}\right), \ \det J\left(-\sqrt{1/3},\sqrt{1/3}\right), \ \det J\left(\sqrt{1/3},\sqrt{1/3}\right),
\end{displaymath}
are linearly independent. Therefore, each element $K$ possesses exactly four independent DOFs for the Jacobian determinant, natively resulting in four DOFs for the calculated density variable. This number strictly exceeds the $\text{tdof}$ of the $P^0$ space and instead matches the $\text{kdof}$. However, the internal energy variable within this setup still possesses only one DOF. This critical structural mismatch leads to mathematical inconsistencies when utilizing the EOS to update the pressure: it obscures which physical point within the element should be queried for the update. Updating the pressure using four distinct density values but only a single internal energy value arbitrarily forces four DOFs onto the pressure variable. Consequently, the thermodynamic variables become improperly discretized across mismatched finite element spaces.

Using a nine-point or higher-order Gauss quadrature rule within the $Q^1-P^0$ space yields the exact identical number of valid DOFs as the four-point rule. The relevant shape function values for the nine-point rule are documented in Table~\ref{table:Q1_shape_function_nine_point_Gauss}.

\begin{table}[htbp]
\centering
\caption{Shape function values of the $Q_1$ element at the nine Gauss quadrature points.}
\label{table:Q1_shape_function_nine_point_Gauss}
\renewcommand{\arraystretch}{1.5}
\begin{tabular}{|c|c|c|c|c|}
\hline
Point $\backslash$ Function & $N_1$ &  $N_2$ & $N_3$ & $N_4$   \\ \hline
$(-\sqrt{0.6},-\sqrt{0.6})$ & $0.4+0.5\sqrt{0.6}$  & $0.1$ & $0.1$ & $0.4-0.5\sqrt{0.6}$   \\ \hline
$(0,-\sqrt{0.6})$ & $0.25(1+\sqrt{0.6})$ & $0.25(1+\sqrt{0.6})$ & $0.25(1-\sqrt{0.6})$ & $0.25(1-\sqrt{0.6})$\\ \hline
$(\sqrt{0.6},-\sqrt{0.6})$ & $0.1$  & $0.4+0.5\sqrt{0.6}$& $0.4-0.5\sqrt{0.6}$ & $0.1$  \\ \hline
$(-\sqrt{0.6},0)$ & $0.25(1+\sqrt{0.6})$ & $0.25(1-\sqrt{0.6})$ & $0.25(1+\sqrt{0.6})$ &$0.25(1-\sqrt{0.6})$  \\ \hline
$(0,0)$ & $0.25$ & $0.25$ & $0.25$ &$0.25$  \\ \hline
$(\sqrt{0.6},0)$ & $0.25(1-\sqrt{0.6})$ & $0.25(1+\sqrt{0.6})$ & $0.25(1-\sqrt{0.6})$ &$0.25(1+\sqrt{0.6})$  \\ \hline
$(-\sqrt{0.6},\sqrt{0.6})$ & $0.1$  & $0.4-0.5\sqrt{0.6}$ & $0.4+0.5\sqrt{0.6}$ & $0.1$   \\ \hline
$(0,\sqrt{0.6})$ & $0.25(1-\sqrt{0.6})$ & $0.25(1-\sqrt{0.6})$ & $0.25(1+\sqrt{0.6})$ & $0.25(1+\sqrt{0.6})$\\ \hline
$(\sqrt{0.6},\sqrt{0.6})$ & $0.4-0.5\sqrt{0.6}$  & $0.1$ & $0.1$ & $0.4+0.5\sqrt{0.6}$  \\ \hline
\end{tabular}
\end{table}

Arranging these shape function evaluations into an extended matrix:
\begin{equation*}
\begin{bmatrix}
0.4+0.5\sqrt{0.6}  & 0.1 & 0.1 & 0.4-0.5\sqrt{0.6} \\
0.25(1+\sqrt{0.6}) & 0.25(1+\sqrt{0.6}) & 0.25(1-\sqrt{0.6}) & 0.25(1-\sqrt{0.6})\\
0.1  & 0.4+0.5\sqrt{0.6} & 0.4-0.5\sqrt{0.6} & 0.1  \\
0.25(1+\sqrt{0.6}) & 0.25(1-\sqrt{0.6}) & 0.25(1+\sqrt{0.6}) &0.25(1-\sqrt{0.6})  \\
0.25 & 0.25 & 0.25 &0.25  \\
0.25(1-\sqrt{0.6}) & 0.25(1+\sqrt{0.6}) & 0.25(1-\sqrt{0.6}) &0.25(1+\sqrt{0.6})  \\
0.1  & 0.4-0.5\sqrt{0.6} & 0.4+0.5\sqrt{0.6} & 0.1   \\
0.25(1-\sqrt{0.6}) & 0.25(1-\sqrt{0.6}) & 0.25(1+\sqrt{0.6}) & 0.25(1+\sqrt{0.6})\\
0.4-0.5\sqrt{0.6}  & 0.1 & 0.1 & 0.4+0.5\sqrt{0.6}
\end{bmatrix}.
\end{equation*}
It is mathematically apparent that the rank of this matrix remains exactly 4. The generated Jacobian determinants are partially linearly dependent; for instance:
\begin{displaymath}
\frac{1}{2}\det J(0,0) = \det J(0,-\sqrt{0.6}) + \det J(0,\sqrt{0.6}) = \det J(-\sqrt{0.6},0) + \det J(\sqrt{0.6},0).
\end{displaymath}
Thus, there are still only four independent DOFs for the Jacobian determinant within each element. Although nine distinct Jacobian determinants are evaluated, five are strictly linearly dependent on the remaining four. Even for higher-order Gauss quadrature rules, the number of independent DOFs for the Jacobian determinant is hard-capped at 4. This plateau occurs because the upper bound for the DOFs of the Jacobian determinant is structurally dictated by the $\text{kdof}$, which equals 4 for the baseline $Q^1$ element.

\textbf{The $Q^2-Q^1$ case:}
For the $Q^2-Q^1$ FEM space pair, Figure~\ref{fig:Q_2_reference_element} illustrates the reference element and the corresponding Gauss quadrature points for $n \times n$ rules (where $n=2,3,4$). The discrete quadrature rules are visualized by color: four-point (blue), nine-point (green, alongside the central point 5), and sixteen-point (yellow).

\begin{figure}[htbp]
\begin{center}
\begin{tikzpicture}[scale = 1]
\fill[black] (-3,-3) circle (2pt);
\fill[black] (0,-3) circle (2pt);
\fill[black] (3,-3) circle (2pt);
\fill[black] (-3,0) circle (2pt);
\fill[black] (0,0) circle (2pt);   
\fill[black] (3,0) circle (2pt); 

\fill[black] (-3,3) circle (2pt);
\fill[black] (0,3) circle (2pt);                           
\fill[black] (3,3) circle (2pt);
    
\draw[black] (-3,-3) -- (3,-3);
\draw[black] (3,-3) -- (3,3);
\draw[black] (3,3) -- (-3,3);
\draw[black] (-3,3) -- (-3,-3);
\node[] at (-3,-3.5) {$1,(-1,-1)$};
\node[] at (0,-3.5) {$2,(0,-1)$};
\node[] at (3,-3.5)  {$3,(1,-1)$};

\node[] at (-3.75,0)  {$4,(-1,0)$};
\node[] at (0.75,0)  {$5,(0,0)$};
\node[] at (3.75,0)  {$6,(1,0)$};

\node[] at (-3,3.5)  {$7,(-1,1)$};
\node[] at (0,3.5)  {$8,(0,1)$};        
\node[] at (3,3.5)   {$9,(1,1)$};

\draw[black,fill=blue] (-1.7325,-1.7325) circle (2pt);
\draw[black,fill=blue] (1.7325,-1.7325) circle (2pt);
\draw[black,fill=blue] (-1.7325,1.7325) circle (2pt);
\draw[black,fill=blue] (1.7325,1.7325) circle (2pt);

\draw[black,fill=green] (-2.3238,-2.3238) circle (2pt);
\draw[black,fill=green] (0,-2.3238) circle (2pt);
\draw[black,fill=green] (2.3238,-2.3238) circle (2pt);
\draw[black,fill=green] (-2.3238,0) circle (2pt);
\draw[black,fill=green] (2.3238,0) circle (2pt);
\draw[black,fill=green] (-2.3238,2.3238) circle (2pt);
\draw[black,fill=green] (0,2.3238) circle (2pt);
\draw[black,fill=green] (2.3238,2.3238) circle (2pt);

\draw[black,fill=yellow] (-2.5834,-2.5834) circle (2pt);        
\draw[black,fill=yellow] (-1.0199,-2.5834) circle (2pt); 
\draw[black,fill=yellow] (1.0199,-2.5834) circle (2pt); 
\draw[black,fill=yellow] (2.5834,-2.5834) circle (2pt); 
\draw[black,fill=yellow] (-2.5834,-1.0199) circle (2pt);        
\draw[black,fill=yellow] (-1.0199,-1.0199) circle (2pt); 
\draw[black,fill=yellow] (1.0199,-1.0199) circle (2pt); 
\draw[black,fill=yellow] (2.5834,-1.0199) circle (2pt); 
\draw[black,fill=yellow] (-2.5834,1.0199) circle (2pt);        
\draw[black,fill=yellow] (-1.0199,1.0199) circle (2pt); 
\draw[black,fill=yellow] (1.0199,1.0199) circle (2pt); 
\draw[black,fill=yellow] (2.5834,1.0199) circle (2pt); 
\draw[black,fill=yellow] (-2.5834,2.5834) circle (2pt);        
\draw[black,fill=yellow] (-1.0199,2.5834) circle (2pt); 
\draw[black,fill=yellow] (1.0199,2.5834) circle (2pt); 
\draw[black,fill=yellow] (2.5834,2.5834) circle (2pt);        
\end{tikzpicture}
    \caption{Reference element $(\xi, \eta)$ of the $Q^2$ element and its 2D Gauss quadrature points.}
    \label{fig:Q_2_reference_element}
    \end{center}
\end{figure}

Similarly, the Jacobian matrix evaluates to:
\begin{equation*}
J=
\begin{bmatrix}
\d{x}{\xi} &  \d{x}{\eta}\\
\d{y}{\xi} &  \d{y}{\eta}
\end{bmatrix}
=
\begin{bmatrix}
\sum_{i=1}^{9} x_i \d{N_i}{\xi} & \sum_{i=1}^{9} x_i \d{N_i}{\eta} \\
\sum_{i=1}^{9} y_i \d{N_i}{\xi} & \sum_{i=1}^{9} y_i \d{N_i}{\eta}
\end{bmatrix},
\end{equation*}
where the expanded Jacobian determinant can again be explicitly formatted via shape function values:
\begin{equation}\label{eq:Q2_Jacobian}
\begin{split}
\det J =& \ N_1(x_{\Lambda}y_{\Delta}-x_{\Delta}y_{\Lambda}) + N_2(x_{\Lambda}y_{\Phi}-x_{\Phi}y_{\Lambda}) + N_3(x_{\Lambda}y_{\Gamma}-x_{\Gamma}y_{\Lambda})\\
+& \ N_4(x_{\Theta}y_{\Delta}-x_{\Delta}y_{\Theta}) + N_5(x_{\Theta}y_{\Phi}-x_{\Phi}y_{\Theta}) + N_6(x_{\Theta}y_{\Gamma}-x_{\Gamma}y_{\Theta})\\
+& \ N_7(x_{\Psi}y_{\Delta}-x_{\Delta}y_{\Psi}) + N_8(x_{\Psi}y_{\Phi}-x_{\Phi}y_{\Psi}) + N_9(x_{\Psi}y_{\Gamma}-x_{\Gamma}y_{\Psi}),
\end{split}
\end{equation}
with the intermediate $x$-coordinate terms constructed from the 1D Lagrange polynomial derivatives $\d{\xi_a}{\xi}$:
\begin{equation}\label{eq:Q2_J11J21}
\begin{split}
x_{\Lambda} &= -\frac{1}{2}(1-2\xi)x_1 - 2\xi x_2 + \frac{1}{2}(1+2\xi) x_3, \quad y_{\Lambda} = -\frac{1}{2}(1-2\xi)y_1 - 2\xi y_2 + \frac{1}{2}(1+2\xi) y_3,\\
x_{\Theta} &= -\frac{1}{2}(1-2\xi)x_4 - 2\xi x_5 + \frac{1}{2}(1+2\xi) x_6, \quad y_{\Theta} = -\frac{1}{2}(1-2\xi)y_4 - 2\xi y_5 + \frac{1}{2}(1+2\xi) y_6, \\
x_{\Psi} &= -\frac{1}{2}(1-2\xi)x_7 - 2\xi x_8 + \frac{1}{2}(1+2\xi) x_9, \quad y_{\Psi} = -\frac{1}{2}(1-2\xi)y_7 - 2\xi y_8 + \frac{1}{2}(1+2\xi) y_9,
\end{split}
\end{equation}
and the orthogonal geometric counterparts derived from $\d{\eta_b}{\eta}$:
\begin{equation}\label{eq:Q2_J12J22}
\begin{split}
x_{\Delta} &= -\frac{1}{2}(1-2\eta)x_1 - 2\eta x_4 + \frac{1}{2}(1+2\eta) x_7, \quad y_{\Delta} = -\frac{1}{2}(1-2\eta)y_1 - 2\eta y_4 + \frac{1}{2}(1+2\eta) y_7,\\
x_{\Phi} &= -\frac{1}{2}(1-2\eta)x_2 - 2\eta x_5 + \frac{1}{2}(1+2\eta) x_8, \quad y_{\Phi} = -\frac{1}{2}(1-2\eta)y_2 - 2\eta y_5 + \frac{1}{2}(1+2\eta) y_8, \\
x_{\Gamma} &= -\frac{1}{2}(1-2\eta)x_3 - 2\eta x_6 + \frac{1}{2}(1+2\eta) x_9, \quad y_{\Gamma} = -\frac{1}{2}(1-2\eta)y_3 - 2\eta y_6 + \frac{1}{2}(1+2\eta) y_9.
\end{split}
\end{equation}

When applying an under-resolved one-point Gauss quadrature rule to discretize mass conservation within the $Q^2-Q^1$ space, the generated Jacobian determinant collapses to a single degree of freedom. Therefore, the independent DOFs available for the density variable fall severely below the required $\text{tdof}$ of a valid $Q^1$ element. At the origin, $\det J(0,0)$ evaluates simplistically to:
\begin{displaymath}
\det J(0,0)=\frac{1}{4}(x_{64}y_{82}+x_{82}y_{64}).
\end{displaymath}

For the $Q^2-Q^1$ FEM space pair, the four-point Gauss quadrature rule (blue nodes, Figure~\ref{fig:Q_2_reference_element}) emerges as the correct minimal choice for mass conservation. It produces exactly four linearly independent Jacobian determinants, perfectly matching the required DOFs of the $Q^1$ element. Although variables like $x_{\Lambda}$ depend conditionally on the exact quadrature point, sub-groupings such as $x_{\Lambda}y_{\Delta}-x_{\Delta}y_{\Lambda}$ exhibit an elegant structured symmetry:
\begin{equation*}
\begin{split}
x_{\Lambda}y_{\Delta}-x_{\Delta}y_{\Lambda} &= (1-2\xi)(1-2\eta)(x_{21}y_{41}-x_{41}y_{21}) + (1-2\xi)(1+2\eta)(x_{21}y_{74}-x_{74}y_{21})\\
&\quad + (1+2\xi)(1-2\eta)(x_{32}y_{41}-x_{41}y_{32}) + (1+2\xi)(1+2\eta)(x_{32}y_{74}-x_{74}y_{32}),
\end{split}
\end{equation*}
where $x_{ij}=x_i-x_j$ references the nodes in Figure~\ref{fig:Q_2_reference_element}. Because all subsequent terms exhibit identical scaling structures, the true DOFs of the Jacobian determinant rely entirely upon the inherent rank of the evaluated shape function matrix. The matrix valuations at the four Gauss points for the $Q^2$ element are cataloged in Table~\ref{table:Q2_shape_function_four_point_Gauss}.

\begin{table}[htbp]
\centering
\caption{Shape function values of the $Q_2$ element evaluated at the four Gauss quadrature points.}
\label{table:Q2_shape_function_four_point_Gauss}
\renewcommand{\arraystretch}{1.5}
\resizebox{\textwidth}{!}{
\begin{tabular}{|c|c|c|c|c|c|c|c|c|c|}
\hline
Point $\backslash$ Function & $N_1$ &  $N_2$ & $N_3$ & $N_4$& $N_5$  &  $N_6$ & $N_7$& $N_8$  &  $N_9$     \\ \hline
$(-\sqrt{1/3},-\sqrt{1/3})$ & $\frac{1}{9}+\frac{1}{6\sqrt{3}}$  & $\frac{1}{9}+\frac{1}{3\sqrt{3}}$ & $-\frac{1}{18}$ & $\frac{1}{9}+\frac{1}{3\sqrt{3}}$& $\frac{4}{9}$  & $\frac{1}{9}-\frac{1}{3\sqrt{3}}$ & $-\frac{1}{18}$ & $\frac{1}{9}-\frac{1}{3\sqrt{3}}$ & $\frac{1}{9}-\frac{1}{6\sqrt{3}}$\\ \hline
$(\sqrt{1/3},-\sqrt{1/3})$ & $-\frac{1}{18}$ & $\frac{1}{9}+\frac{1}{3\sqrt{3}}$ & $\frac{1}{9}+\frac{1}{6\sqrt{3}}$ & $\frac{1}{9}-\frac{1}{3\sqrt{3}}$ & $\frac{4}{9}$ & $\frac{1}{9}+\frac{1}{3\sqrt{3}}$ & $\frac{1}{9}-\frac{1}{6\sqrt{3}}$ & $\frac{1}{9}-\frac{1}{3\sqrt{3}}$ & $-\frac{1}{18}$ \\ \hline
$(-\sqrt{1/3},\sqrt{1/3})$ & $-\frac{1}{18}$ & $\frac{1}{9}-\frac{1}{3\sqrt{3}}$ & $\frac{1}{9}-\frac{1}{6\sqrt{3}}$ & $\frac{1}{9}+\frac{1}{3\sqrt{3}}$ & $\frac{4}{9}$ & $\frac{1}{9}-\frac{1}{3\sqrt{3}}$ & $\frac{1}{9}+\frac{1}{6\sqrt{3}}$ & $\frac{1}{9}+\frac{1}{3\sqrt{3}}$ & $-\frac{1}{18}$ \\ \hline
$(\sqrt{1/3},\sqrt{1/3})$ & $\frac{1}{9}-\frac{1}{6\sqrt{3}}$ & $\frac{1}{9}-\frac{1}{3\sqrt{3}}$ & $-\frac{1}{18}$ & $\frac{1}{9}-\frac{1}{3\sqrt{3}}$  & $\frac{4}{9}$ & $\frac{1}{9}+\frac{1}{3\sqrt{3}}$ & $-\frac{1}{18}$ & $\frac{1}{9}+\frac{1}{3\sqrt{3}}$ &$\frac{1}{9}+\frac{1}{6\sqrt{3}}$\\ \hline
\end{tabular}}
\end{table}

If an overly refined nine-point (green nodes, Figure~\ref{fig:Q_2_reference_element}) or sixteen-point (yellow nodes) Gauss quadrature rule is mistakenly employed to solve the mass conservation law for the $Q^2-Q^1$ pair, the resulting evaluated Jacobian determinant yields an excessive 9 DOFs. This value exceeds the valid $\text{tdof}$ of the $Q^1$ space and incorrectly aligns with the higher DOFs of the $Q^2$ space, irreparably fracturing EOS consistency.

We now formalize this observation with the following proposition.
\begin{prop}[DOFs of the Jacobian determinant]
The DOFs of the Jacobian determinant for a $Q^m$ element are mathematically capped at $(m+1)^2$. This upper bound is strictly governed by the $\text{kdof}$ of the $Q^m$ space, which is exactly $(m+1)^2$. Consequently, any Gauss quadrature rule employing more than $(m+1)^2$ integration points will not increase the independent DOFs of the Jacobian determinant beyond this limit.
\label{prop:dofs_jacobian_determinant}
\end{prop}
\begin{proof}
Let the basis functions of the kinematic variables be $N_i$ for $i=1,\dots,\text{kdof}$. Each basis function can be expressed as a tensor product of the 1D Lagrange polynomials:
\begin{equation*}
    N_i=\xi_a\eta_b \quad a,b=1,...,h, \quad i=a+(b-1)h.
\end{equation*}
where $h = m+1$. The coordinates of a point within the element are exclusively interpolated via these kinematic shape functions:
\begin{equation*}
    x=\sum_{i=1}^{tdof}N_i x_i, \quad y=\sum_{i=1}^{tdof}N_i y_i.
\end{equation*}

The Jacobian matrix can then be expressed as:
\begin{equation}\label{eq:Jacobian_general}
J=
\begin{bmatrix}
\d{x}{\xi} &  \d{x}{\eta}\\
\d{y}{\xi} &  \d{y}{\eta}
\end{bmatrix}
=
\begin{bmatrix}
\sum_{i=1}^{tdof} x_i \d{N_i}{\xi} & \sum_{i=1}^{tdof} x_i \d{N_i}{\eta} \\
\sum_{i=1}^{tdof} y_i \d{N_i}{\xi} & \sum_{i=1}^{tdof} y_i \d{N_i}{\eta}
\end{bmatrix},
\end{equation}

The Jacobian determinant is computed as:
\begin{equation}\label{eq:Jacobian_determinant_general}
\det J=\frac{\partial x}{\partial \xi}\frac{\partial y}{\partial \eta}-\frac{\partial x}{\partial \eta}\frac{\partial y}{\partial \xi}.
\end{equation}

The partial derivatives can be directly calculated as:
\begin{equation*}
    \begin{split}
    \d{x}{\xi} &= \sum_{i=1}^{tdof} x_i \d{N_i}{\xi}=\sum_{b=1}^{h}\sum_{a=1}^{h} x_{a+(b-1)h} \d{\xi_a}{\xi}\eta_b=\sum_{b=1}^{h}\left(\sum_{a=1}^{h} x_{a+(b-1)h} \d{\xi_a}{\xi}\right)\eta_b,\\
    \d{x}{\eta} &= \sum_{i=1}^{tdof} x_i \d{N_i}{\eta}=\sum_{b=1}^{h}\sum_{a=1}^{h} x_{a+(b-1)h} \d{\eta_b}{\eta}\xi_a=\sum_{a=1}^{h} \left(\sum_{b=1}^{h}x_{a+(b-1)h} \d{\eta_b}{\eta}\right)\xi_a,\\
    \d{y}{\xi} &= \sum_{i=1}^{tdof} y_i \d{N_i}{\xi}=\sum_{b=1}^{h}\sum_{a=1}^{h} y_{a+(b-1)h} \d{\xi_a}{\xi}\eta_b=\sum_{b=1}^{h}\left(\sum_{a=1}^{h} y_{a+(b-1)h} \d{\xi_a}{\xi}\right)\eta_b,\\ 
    \d{y}{\eta} &= \sum_{i=1}^{tdof} y_i \d{N_i}{\eta}=\sum_{b=1}^{h}\sum_{a=1}^{h} y_{a+(b-1)h} \d{\eta_b}{\eta}\xi_a=\sum_{a=1}^{h} \left(\sum_{b=1}^{h}y_{a+(b-1)h} \d{\eta_b}{\eta}\right)\xi_a,   
    \end{split}
\end{equation*}

Let $\alpha,\beta,\chi,\psi$ denote the following groupings:
\begin{equation*}
    \begin{split}
    \alpha_b &=\sum_{a=1}^{h} x_{a+(b-1)h} \d{\xi_a}{\xi}, \\
    \beta_a&=\sum_{b=1}^{h}y_{a+(b-1)h} \d{\eta_b}{\eta},\\
    \chi_b&=\sum_{a=1}^{h} y_{a+(b-1)h} \d{\xi_a}{\xi}, \\
    \psi_a&=\sum_{b=1}^{h}x_{a+(b-1)h} \d{\eta_b}{\eta}.
    \end{split}
\end{equation*}
The Jacobian determinant~\eqref{eq:Jacobian_determinant_general} can thus be evaluated as:
\begin{equation*}
\det J=\left(\sum_{b=1}^{h}\alpha_b \eta_b \right) \left(\sum_{a=1}^{h}\beta_a \xi_a\right)-\left(\sum_{b=1}^{h}\chi_b \eta_b\right) \left(\sum_{a=1}^{h}\psi_a \xi_a\right).
\end{equation*}

The above equation can be further reorganized as:
\begin{equation*}
\det J=\sum_{a=1}^{h}\sum_{b=1}^{h}(\alpha_b\beta_a-\chi_b\psi_a)\xi_a\eta_b=\sum_{i=1}^{\text{tdof}} M_i N_i.
\end{equation*}
If we arrange $\alpha,\beta,\chi,\psi$ into column vectors:
\begin{equation*}
    \begin{bmatrix}
\alpha_1\\
\vdots\\
\alpha_h\\
    \end{bmatrix},\
        \begin{bmatrix}
\beta_1\\
\vdots\\
\beta_h\\
    \end{bmatrix},
        \begin{bmatrix}
\chi_1\\
\vdots\\
\chi_h\\
    \end{bmatrix},
        \begin{bmatrix}
\psi_1\\
\vdots\\
\psi_h\\
    \end{bmatrix}.
    \end{equation*}
The terms $\alpha_b\beta_a$ constitute an $h\times h$ matrix equal to $\vec{\alpha}\vec{\beta}^T$.

These vectors are computed by multiplying the element coordinate matrices by the one-dimensional shape function derivatives:
\begin{equation*}
        \begin{bmatrix}
\alpha_1\\
\alpha_2\\
\vdots\\
\alpha_h\\
    \end{bmatrix}=
    \begin{bmatrix}
        x_1 & x_2 & \cdots & x_h\\
        x_{h+1} & x_{h+2} & \cdots & x_{2h}\\
        \vdots & \vdots & \ddots & \vdots\\
        x_{(h-1)h+1} & x_{(h-1)h+2} & \cdots & x_{h^2}
        \end{bmatrix}
        \begin{bmatrix}
\d{\xi_1}{\xi}\\
\d{\xi_2}{\xi}\\
\vdots\\
\d{\xi_h}{\xi}
        \end{bmatrix}
\end{equation*} 

\begin{equation*}
        \begin{bmatrix}
\beta_1\\
\beta_2\\
\vdots\\
\beta_h\\
    \end{bmatrix}=
    \begin{bmatrix}
        y_1 & y_{h+1} & \cdots & y_{h(h-1)+1}\\
        y_{2} & y_{h+2} & \cdots & y_{h(h-1)+2}\\
        \vdots & \vdots & \ddots & \vdots\\
        y_{h} & y_{2h} & \cdots & y_{h^2}
        \end{bmatrix}
        \begin{bmatrix}
\d{\eta_1}{\eta}\\
\d{\eta_2}{\eta}\\
\vdots\\
\d{\eta_h}{\eta}
        \end{bmatrix}
\end{equation*} 

\begin{equation*}
        \begin{bmatrix}
\chi_1\\
\chi_2\\
\vdots\\
\chi_h\\
    \end{bmatrix}=
    \begin{bmatrix}
        y_1 & y_2 & \cdots & y_h\\
        y_{h+1} & y_{h+2} & \cdots & y_{2h}\\
        \vdots & \vdots & \ddots & \vdots\\
        y_{(h-1)h+1} & y_{(h-1)h+2} & \cdots & y_{h^2}
        \end{bmatrix}
        \begin{bmatrix}
\d{\xi_1}{\xi}\\
\d{\xi_2}{\xi}\\
\vdots\\
\d{\xi_h}{\xi}
        \end{bmatrix}
\end{equation*} 

\begin{equation*}
        \begin{bmatrix}
\psi_1\\
\psi_2\\
\vdots\\
\psi_h\\
    \end{bmatrix}=
    \begin{bmatrix}
        x_1 & x_{h+1} & \cdots & x_{h(h-1)+1}\\
        x_{2} & x_{h+2} & \cdots & x_{h(h-1)+2}\\
        \vdots & \vdots & \ddots & \vdots\\
        x_{h} & x_{2h} & \cdots & x_{h^2}
        \end{bmatrix}
        \begin{bmatrix}
\d{\eta_1}{\eta}\\
\d{\eta_2}{\eta}\\
\vdots\\
\d{\eta_h}{\eta}
        \end{bmatrix}
\end{equation*} 

The matrix $\vec{\alpha}\vec{\beta}^T$ can be expressed as:
\begin{equation*}
    \begin{bmatrix}
        x_1 & x_2 & \cdots & x_h\\
        x_{h+1} & x_{h+2} & \cdots & x_{2h}\\
        \vdots & \vdots & \ddots & \vdots\\
        x_{(h-1)h+1} & x_{(h-1)h+2} & \cdots & x_{h^2}
        \end{bmatrix}
        \begin{bmatrix}
\d{\xi_1}{\xi}\\
\d{\xi_2}{\xi}\\
\vdots\\
\d{\xi_h}{\xi}
        \end{bmatrix}
        \begin{bmatrix}
\d{\eta_1}{\eta} & \d{\eta_2}{\eta} & \cdots & \d{\eta_h}{\eta}
        \end{bmatrix}
    \begin{bmatrix}
        y_1 & y_2 & \cdots & y_h\\
        y_{h+1} & y_{h+2} & \cdots & y_{2h}\\
        \vdots & \vdots & \ddots & \vdots\\
        y_{(h-1)h+1} & y_{(h-1)h+2} & \cdots & y_{h^2}
        \end{bmatrix}                   
    \end{equation*}

The matrix $\vec{\chi}\vec{\psi}^T$ can be expressed as:
\begin{equation*}
        \begin{bmatrix}
        y_1 & y_2 & \cdots & y_h\\
        y_{h+1} & y_{h+2} & \cdots & y_{2h}\\
        \vdots & \vdots & \ddots & \vdots\\
        y_{(h-1)h+1} & y_{(h-1)h+2} & \cdots & y_{h^2}
        \end{bmatrix}  
        \begin{bmatrix}
\d{\xi_1}{\xi}\\
\d{\xi_2}{\xi}\\
\vdots\\
\d{\xi_h}{\xi}
        \end{bmatrix}
        \begin{bmatrix}
\d{\eta_1}{\eta} & \d{\eta_2}{\eta} & \cdots & \d{\eta_h}{\eta}
        \end{bmatrix}
    \begin{bmatrix}
        x_1 & x_2 & \cdots & x_h\\
        x_{h+1} & x_{h+2} & \cdots & x_{2h}\\
        \vdots & \vdots & \ddots & \vdots\\
        x_{(h-1)h+1} & x_{(h-1)h+2} & \cdots & x_{h^2}
        \end{bmatrix}                 
    \end{equation*}

Let matrix $M$ denote the following matrix Hadamard product:  
\begin{equation*}
M=\left[\vec{\alpha}\vec{\beta}^T-\vec{\chi}\vec{\psi}^T\right]\odot
    \begin{bmatrix}
        N_1 & N_2 & \cdots & N_h\\
        N_{h+1} & N_{h+2} & \cdots & N_{2h}\\
        \vdots & \vdots & \ddots & \vdots\\
        N_{(h-1)h+1} & N_{(h-1)h+2} & \cdots & N_{h^2}
        \end{bmatrix} 
\end{equation*}

The Jacobian determinant is then the sum of all elements in the matrix $M$:
\begin{equation*}
    \det J= \mathbf{1}^T M \mathbf{1}^T
\end{equation*} 
Since $M$ is an $h\times h$ matrix, the Jacobian determinant has at most $h^2$ independent DOFs. Recalling that $h^2 = (m+1)^2 = \text{kdof}$, this completes the proof.
\end{proof}

This fundamental conclusion generalizes universally to any $Q^m-Q^{m-1}$ FEM space pair. We reiterate the primary takeaway of this section: the $m^2$-point Gauss quadrature rule is the strictly appropriate numerical choice for solving the mass conservation law in high-order $Q^m-Q^{m-1}$ SGH formulations. This specific rule correctly caps the DOFs of the Jacobian determinant, simultaneously constraining the available DOFs of the generated density variable.

If fewer than $m^2$ quadrature points are deployed, the entire density system becomes ill-posed and underdetermined. Conversely, if an excessive number of points is used, the DOFs of the generated density variable outstrip the $\text{tdof}$ of the $Q^{m-1}$ element. Specifically, applying an $(m+1)^2$-point rule spawns exactly $(m+1)^2$ distinct density values, erroneously matching the spatial DOFs of the $Q^m$ space instead. Employing an aggressively higher-order rule, such as an $(m+2)^2$-point quadrature, does not expand the independent DOFs further; while the total pool of evaluated density variables nominally swells, the effective mathematical DOFs hit a hard ceiling at $(m+1)^2$, as the surplus variables simply devolve into linearly dependent combinations.

\subsection{The issues with exact integration}
We must emphasize the numerical hazards associated with the exact integration of the mass conservation law. Primarily, employing an $(m+1)^2$-point (or higher-order) Gauss quadrature rule spawns redundant DOFs for the discrete density variable. This geometric redundancy breaks the necessary one-to-one mapping with the internal energy variable, obscuring which physical point should be queried for pressure updates via the EOS. Furthermore, when a non-linear EOS is utilized, the emergent DOFs of the pressure variable become wildly ambiguous. Considering the $Q^m-Q^{m-1}$ pair as an example, employing either the $(m+1)^2$-point or the $(m+2)^2$-point Gauss rule to resolve mass conservation calculates the exact same underlying density field. However, attempting to update the pressure term across $(m+1)^2$ versus $(m+2)^2$ evaluation nodes using a non-linear EOS will yield irreconcilably different global pressure fields. Thus, the seemingly arbitrary choice of an exact quadrature rule profoundly and directly corrupts the deterministic numerical results in subsequent momentum conservation updates.

A secondary critical issue involves the baseline robustness of the SGH scheme, which is heavily compromised by negative Jacobian determinants arising from aggressive mesh distortion. The magnitude of the Jacobian determinant depends strongly on the coordinates of the chosen quadrature points. If this evaluation approaches zero or turns negative anywhere, the local density asymptotes toward a singularity or goes negative—a catastrophic condition that triggers an immediate abort criterion in standard SGH simulations. This mechanism dictates the absolute geometric robustness of the entire numerical scheme. Tracking the $Q^1-P^0$ FEM spaces as a test case, evaluating the four-point Gauss rule at the node $(-\sqrt{1/3},-\sqrt{1/3})$ yields a Jacobian determinant of:
\begin{equation*}
\begin{split}
\det J\left(-\sqrt{1/3},-\sqrt{1/3}\right) &= \left(\frac{1}{6}+\frac{1}{4\sqrt{3}}\right)S_{312} + \frac{1}{12} S_{124} + \frac{1}{12} S_{431} + \left(\frac{1}{6}-\frac{1}{4\sqrt{3}}\right)S_{243} \\
&= \frac{1}{4}\left(1-\frac{1}{\sqrt{3}}\right)S + \frac{1}{2\sqrt{3}}S_{312}\\
&\approx 0.10566 S + 0.28867 S_{312},
\end{split}
\end{equation*}
where $S$ represents the total scalar area of the quadrilateral element, and $S_{ijk}$ specifies the area of the subtriangle enclosed by boundary vertices $i$, $j$, and $k$.

Conversely, testing the nine-point Gauss rule at its geometric node $(-\sqrt{0.6},-\sqrt{0.6})$ evaluates to:
\begin{equation*}
\begin{split}
\det J(-\sqrt{0.6},-\sqrt{0.6}) &= (0.2+0.25\sqrt{0.6})S_{312} + 0.05 S_{124} + 0.05 S_{431} + (0.2-0.25\sqrt{0.6})S_{243} \\
&= \frac{1}{4}(1-\sqrt{0.6})S + \frac{1}{2}\sqrt{0.6}S_{312}\\
&\approx 0.05635 S + 0.38729 S_{312}.
\end{split}
\end{equation*}

The algebraic expansion of the Jacobian determinants detailed above illustrates the vulnerability thresholds at which they approach catastrophic singularities. Consider a fundamental distortion scenario where vertex $1$ in Figure~\ref{fig:Q_1_reference_element} translates inward along the diagonal line $l_{14}$. As the quadrilateral element is forced into a concave shape, the scalar area $S_{312}$ drops into negative values, even while the gross element area $S$ remains globally positive. Tracking the numerical drop-off, it becomes evident that at the exact moment $\det J(-\sqrt{0.6},-\sqrt{0.6})$ violently crosses zero, the four-point evaluation $\det J(-\sqrt{1/3},-\sqrt{1/3})$ still safely remains positive. This geometric failure mode clearly demonstrates that the structural robustness of the higher-order nine-point Gauss quadrature rule is strictly inferior to that of the four-point rule when bounded within the $Q^1-P^0$ FEM framework. To extend this analysis, if an extreme under-resolved one-point Gauss quadrature rule is utilized to bypass mass conservation tracking entirely in the $Q^1-P^0$ pair, the centralized Jacobian determinant strictly equals $\det J (0,0)=0.25 S$, representing an unbreakable, unconditionally positive state as long as the gross area of the element does not fully invert.

A general conclusion therefore emerges: for any given $Q^m-Q^{m-1}$ FEM space, applying the minimal viable number of quadrature points generally yields the most geometrically robust numerical scheme for resolving discrete mass conservation. However, this sweeping conclusion does not map cleanly onto the discretization of the momentum conservation law, which remains inherently tethered to the mandatory control of parasitic hourglass distortions. The specific localized tradeoffs for $Q^1-P^0$ FEM spaces were rigorously detailed in~\cite{Sun2022Understanding}. The generalized framework tradeoffs encompassing all extended $Q^m-Q^{m-1}$ spaces will be systematically addressed in our subsequent sections.

We conclude this section by examining the strategic spatial distribution of the generated thermodynamic DOFs. We strictly lock the thermodynamic DOFs directly to the evaluated Gauss quadrature points. Tracking the $Q^m-Q^{m-1}$ FEM space pair, the resultant DOFs spanning the density, internal energy, and pressure variables are completely defined and collocated at the $m^2$ interior Gauss quadrature points. This specific clustering approach aggressively diverges from traditional rigid FEM spaces, where core DOFs are normally hard-anchored to the geometric element vertices or exterior boundaries. Our interior clustering strategy unlocks several monumental computational advantages:

\begin{itemize}
\item When explicitly solving the dynamic mass conservation law, the constituent density variables can be natively updated and evolved strictly through isolated evaluations of the Jacobian determinant.
\item When shifting to the momentum conservation law, the formulation of the right-hand side tracking force matrix strictly requires current pressure or bulk stress variables. Because these values are already actively housed at the tracking Gauss points, all expensive, lossy interpolation procedures are completely bypassed, leading to the highly compact and efficient explicit computation algorithms demonstrated in the next section.
\item During the resolution of the generalized energy conservation law, the resultant internal energy mass matrix collapses beautifully into a purely diagonal form, completely eliminating the need for iterative global solver loops.
\item This centralized clustering approach massively streamlines the application of complex equations of state. Because all dependent and independent fluid variables structurally map to the exact same physical coordinates (the centralized quadrature point), mathematically perfect state consistency is guaranteed during all high-speed EOS updates.
\end{itemize}
\section{Momentum conservation law}\label{sec:momentum_conservation}

In this section, we detail the discretization of the momentum conservation law from Equation~\eqref{eq:lagrangian} within the proposed high-order staggered Lagrangian hydrodynamics framework. A Petrov-Galerkin technique is applied over each element $K$:
\begin{equation}\label{eq:momentum_integral}
  \int_{K} \rho_h \frac{d \vec{u}_h}{d t} N_j \dx\dy = \int_{K} - \nabla p_h N_j \dx\dy = \int_{K} p_h \nabla N_j \dx\dy.
\end{equation}
Substituting the discrete expansion for the velocity while temporarily retaining the continuum representations for density and pressure for brevity, we obtain:
\begin{equation}\label{eq:momentum_conservation}
  \int_{K} \rho_h \frac{d \sum_{i=1}^{\text{kdof}} \vec{u}_{i}N_{i}}{d t} N_j \dx\dy = \sum_{i=1}^{\text{kdof}} \frac{d  \vec{u}_{i}}{d t} \int_{K} \rho_h N_{i} N_j \dx\dy = \int_{K} p_h \nabla N_j \dx\dy.
 \end{equation}

Solving Equation~\eqref{eq:momentum_conservation} integrates three key components: quadrature rules, the mass lumping technique, and the treatment of hourglass distortion. An $m^2$-point Gauss quadrature rule is employed to integrate the right-hand side. We will demonstrate that the derivatives of the shape functions can be expressed as a product of the shape function values and their derivatives with respect to both the reference coordinates and the vertex coordinates—a vectorized technique that significantly reduces computational costs. Mass lumping is achieved via Lobatto quadrature, a strategy that requires the DOFs for the $Q^m$ space to be collocated with the quadrature points. Furthermore, we will establish mathematically that hourglass distortion is an inherent byproduct of this framework, thereby necessitating the specifichourglass control and artificial viscosity detailed in subsequent sections.

\subsection{Calculation of the right-hand side}
The $m^2$-point Gauss quadrature rule, which is exact for polynomials up to degree $2m-1$ (as noted in Section~\ref{sec:discretization}), is fully sufficient to exactly integrate the right-hand side $\int_{K} p_h \nabla N_j \dx\dy$ for the $Q^m-Q^{m-1}$ FEM space pair.

The derivation proceeds as follows: the approximating polynomial for the pressure term $p_h$ has a maximum degree of $\xi^{m-1}\eta^{m-1}$ in the reference space, while the gradient of the shape function, $\nabla N_j$, exhibits a maximum degree of $\xi^{m-1}\eta^{m}$ or $\xi^m \eta^{m-1}$. The highest-order term resulting from their product is therefore exactly $\xi^{2m-2}\eta^{2m-1}$ or $\xi^{2m-1}\eta^{2m-2}$, which falls perfectly within the exact integration limits of the $m^2$-point Gauss rule.

The discrete right-hand side is thus calculated as:
\begin{equation}\label{eq:discrete_rhs_momentum_conservation}
    \int_{K} p_h \nabla N_j \dx\dy=\sum_{q=1}^{m^2} \omega_{q} p_h(\gv{\xi}_q)\nabla N_j  \det J(\gv{\xi}_q).
\end{equation}

We now offer a brief remark on the selection of this Gauss quadrature rule. In contrast to the discretization of the mass conservation law, the use of higher-order Gauss quadrature rules here has no effect on the computed result or numerical robustness. Because the isolated Jacobian determinant in Equation~\eqref{eq:discrete_rhs_momentum_conservation} is never utilized independently, the grouped expression $\det J(\gv{\xi}_q) \nabla N_j$ is evaluated as a single entity. The pressure variables are defined exactly at the $m^2$ Gauss quadrature points to simplify calculations. Employing a higher-order quadrature rule would needlessly necessitate complex interpolation procedures and inflate computational costs without providing any physical gain.

Evaluating the grouped term $\det J(\gv{\xi}_q) \nabla N_j$ as a combined entity is fundamentally advantageous for reducing computational overhead. In 2D Cartesian geometry, the gradient of the shape function $\nabla N_j$ is calculated as:
\begin{equation}\label{eq:velocity_shape_function_gradient}
    \nabla N_j =    
    \begin{bmatrix}
        \d{N_j}{x}\\
        \d{N_j}{y} 
    \end{bmatrix}
        =
    \{J^{-1}\}^{T}
    \begin{bmatrix}
        \d{N_j}{\xi}\\
        \d{N_j}{\eta} 
    \end{bmatrix},
\end{equation}
where the Jacobian matrix is
\begin{equation*}
    J=
    \begin{bmatrix}
        \d{x}{\xi} &  \d{x}{\eta}\\
        \d{y}{\xi} &  \d{y}{\eta}
    \end{bmatrix}
    =
    \begin{bmatrix}
        J_{11} &  J_{12}\\
        J_{21} &  J_{22}
    \end{bmatrix}.
\end{equation*}

The inverse matrix is:
\begin{equation*}
J^{-1}=\frac{1}{\det J}    
\begin{bmatrix}
    J_{22} &  -J_{12}\\
    -J_{21} &  J_{11}
\end{bmatrix},
\end{equation*}
and its transpose is:
\begin{equation*}
\{J^{-1}\}^{T}=\frac{1}{\det J}    
    \begin{bmatrix}
        J_{22} &  -J_{21}\\
        -J_{12} &  J_{11}
    \end{bmatrix}.
\end{equation*}

This yields the explicit gradient components:
\begin{equation}\label{shape_function_gradient_per_direction}
\begin{split}
    (\det J) \d{N_j}{x}= & \d{N_j}{\xi}J_{22}-\d{N_j}{\eta}J_{21}, \\
    (\det J) \d{N_j}{y}= & -\d{N_j}{\xi}J_{12}+\d{N_j}{\eta}J_{11}.
\end{split}
\end{equation}

The physical coordinates $x$ and $y$ are interpolated from the nodal values using the kinematic shape functions:
\begin{equation*}
\begin{split}
x=&\sum_{i=1}^{(m+1)^2}x_i N_i,\\
y=&\sum_{i=1}^{(m+1)^2}y_i N_i.
\end{split}   
\end{equation*} 
Consequently, we have:
\begin{equation*}
J_{11}=\d{x}{\xi}=\d{\sum_{i=1}^{(m+1)^2}x_i N_i}{\xi}=\sum_{i=1}^{(m+1)^2}x_i \d{N_i}{\xi}.
\end{equation*}

The complete shape function $N_i$ is a tensor product of 1D shape functions, denoted here by $\xi_{a}$ and $\eta_{b}$. The $Q^1, Q^2,$ and $Q^3$ shape functions were detailed previously in Section~\ref{sec:discretization}.
\begin{equation*}
\begin{split}
    \xi_{a}&, \quad a=1,2,\dots,m+1,\\
    \eta_{b}&, \quad b=1,2,\dots,m+1.
\end{split}
\end{equation*}
Thus, similar as the Section~\ref{sec:gauss_quadrature} describes, $N_i$ can be strictly decomposed into the product $\xi_{a} \eta_{b}$:
\begin{equation}
    N_i=\xi_{a}\eta_{b},\quad  i=a+(b-1)(m+1).
\end{equation}

Therefore, $J_{11}$ can be expanded as:
\begin{equation*}
    J_{11}=\sum_{i=1}^{(m+1)^2}x_i \d{N_i}{\xi}=\sum_{a=1}^{m+1}\sum_{b=1}^{m+1} x_{i} \d{\xi_{a}}{\xi}\eta_{b}, \quad i=a+(b-1)(m+1),
\end{equation*}

Here, we present an efficient algorithmic approach to calculate $J_{11}$ simultaneously for all quadrature points within an element using vectorized matrix multiplication. Let $h = m+1$. We arrange the nodal coordinates $x_i$ into an $h \times h$ matrix:
\begin{equation*}
    \begin{bmatrix}
        x_1 & x_2 & \cdots & x_h\\
        x_{h+1} & x_{h+2} & \cdots & x_{2h}\\
        \vdots & \vdots & \ddots & \vdots\\
        x_{(h-1)h+1} & x_{(h-1)h+2} & \cdots & x_{h^2}
        \end{bmatrix}.
\end{equation*}

For a single evaluation point, $J_{11}$ is computed via the following matrix-vector multiplication:
\begin{equation*}
J_{11}=       
\begin{bmatrix} 
\eta_1, \eta_2, \dots, \eta_h
\end{bmatrix}
\begin{bmatrix}
x_1 & x_2 & \cdots & x_h\\
x_{h+1} & x_{h+2} & \cdots & x_{2h}\\
\vdots & \vdots & \ddots & \vdots\\
x_{(h-1)h+1} & x_{(h-1)h+2} & \cdots & x_{h^2}
\end{bmatrix}
\begin{bmatrix}
\d{\xi_1}{\xi}\\
\d{\xi_2}{\xi}\\
\vdots\\
\d{\xi_h}{\xi}
\end{bmatrix} 
\end{equation*}

By evaluating the 1D shape functions at the $m$ distinct 1D Gauss quadrature points (indexed by $k=1,2,\dots,m$), we can extend the row vector into an $m \times h$ matrix:
\begin{equation*}
\begin{bmatrix}
\eta_1(1) & \eta_2(1) & \cdots & \eta_h(1)\\
\eta_1(2) & \eta_2(2) & \cdots & \eta_h(2)\\
\vdots & \vdots & \ddots & \vdots\\
\eta_1(m) & \eta_2(m) & \cdots & \eta_h(m)
\end{bmatrix}.
\end{equation*} 

Similarly, evaluating the shape function derivatives at these $m$ quadrature points yields an $h \times m$ matrix:
\begin{equation*}
\begin{bmatrix}
\d{\xi_1}{\xi}(1) & \d{\xi_1}{\xi}(2) & \cdots & \d{\xi_1}{\xi}(m)\\
\d{\xi_2}{\xi}(1) & \d{\xi_2}{\xi}(2) & \cdots & \d{\xi_2}{\xi}(m)\\
\vdots & \vdots & \ddots & \vdots\\
\d{\xi_h}{\xi}(1) & \d{\xi_h}{\xi}(2) & \cdots & \d{\xi_h}{\xi}(m)
\end{bmatrix}.
\end{equation*}

The full evaluation is then structured as:
\begin{equation*}
\begin{bmatrix}
\eta_1(1)  & \cdots & \eta_h(1)\\
\eta_1(2) & \cdots & \eta_h(2)\\
\vdots & \ddots & \vdots\\
\eta_1(m) & \cdots & \eta_h(m)
\end{bmatrix}
\begin{bmatrix}
x_1 & \cdots & x_h\\
x_{h+1}  & \cdots & x_{2h}\\
\vdots & \ddots & \vdots\\
x_{(h-1)h+1}  & \cdots & x_{h^2}
\end{bmatrix}
\begin{bmatrix}
\d{\xi_1}{\xi}(1)  & \cdots & \d{\xi_1}{\xi}(m)\\
\d{\xi_2}{\xi}(1) & \cdots & \d{\xi_2}{\xi}(m)\\
\vdots & \ddots & \vdots\\
\d{\xi_h}{\xi}(1)  & \cdots & \d{\xi_h}{\xi}(m)
\end{bmatrix}
\end{equation*}
This operation yields an $m\times m$ matrix containing the exact values of $J_{11}$ at all $m^2$ interior quadrature points simultaneously.

Following an identical methodology, the remaining Jacobian components expand as:
\begin{equation*}
\begin{split}
    J_{12}&=\sum_{a=1}^{m+1}\sum_{b=1}^{m+1} x_{i} \xi_{a}\d{\eta_{b}}{\eta},\\
    J_{21}&=\sum_{a=1}^{m+1}\sum_{b=1}^{m+1} y_{i} \d{\xi_{a}}{\xi}\eta_{b},\\
    J_{22}&=\sum_{a=1}^{m+1}\sum_{b=1}^{m+1} y_{i} \xi_{a}\d{\eta_{b}}{\eta}, \quad i=a+(b-1)(m+1).
\end{split}
\end{equation*}

Expressed in matrix form for a generic quadrature point:
\begin{equation*}
J_{12}=       
\begin{bmatrix} 
\d{\eta_1}{\xi}, \d{\eta_2}{\xi}, \dots, \d{\eta_h}{\xi}
\end{bmatrix}
\begin{bmatrix}
x_1 & x_2 & \cdots & x_h\\
x_{h+1} & x_{h+2} & \cdots & x_{2h}\\
\vdots & \vdots & \ddots & \vdots\\
x_{(h-1)h+1} & x_{(h-1)h+2} & \cdots & x_{h^2}
\end{bmatrix}
\begin{bmatrix}
\xi_1\\
\xi_2\\
\vdots\\
\xi_h
\end{bmatrix} 
\end{equation*}

\begin{equation*}
J_{21}=       
\begin{bmatrix} 
\eta_1, \eta_2, \dots, \eta_h
\end{bmatrix}
\begin{bmatrix}
y_1 & y_2 & \cdots & y_h\\
y_{h+1} & y_{h+2} & \cdots & y_{2h}\\
\vdots & \vdots & \ddots & \vdots\\
y_{(h-1)h+1} & y_{(h-1)h+2} & \cdots & y_{h^2}
\end{bmatrix}
\begin{bmatrix}
\d{\xi_1}{\xi}\\
\d{\xi_2}{\xi}\\
\vdots\\
\d{\xi_h}{\xi}
\end{bmatrix} 
\end{equation*}

\begin{equation*}
J_{22}=       
\begin{bmatrix} 
\d{\eta_1}{\xi}, \d{\eta_2}{\xi}, \dots, \d{\eta_h}{\xi}
\end{bmatrix}
\begin{bmatrix}
y_1 & y_2 & \cdots & y_h\\
y_{h+1} & y_{h+2} & \cdots & y_{2h}\\
\vdots & \vdots & \ddots & \vdots\\
y_{(h-1)h+1} & y_{(h-1)h+2} & \cdots & y_{h^2}
\end{bmatrix}
\begin{bmatrix}
\xi_1\\
\xi_2\\
\vdots\\
\xi_h
\end{bmatrix} 
\end{equation*}
In these formulations, the pre-multiplying row vectors correspond to functions of $\eta$, while the post-multiplying column vectors correspond to functions of $\xi$, ensuring strict dimensional consistency. With $J_{11}, J_{12}, J_{21},$ and $J_{22}$ efficiently computed, the global gradient of any shape function is readily resolved via Equation~\eqref{shape_function_gradient_per_direction}.

\subsection*{$Q^1-P^0$ case}
For the $Q^1-P^0$ case, the shape functions possess a significantly simpler representation. Because the subscript notation can become dense, we illustrate this formulation explicitly using $N_1 = \xi_1\eta_1$ as an example:
\begin{equation}\label{shape_function_gradient_simplified}
\begin{split}
    (\det J) \d{N_1}{x}= & \d{N_1}{\xi}J_{22}-\d{N_1}{\eta}J_{21} \\
    =& \d{\xi_1}{\xi}\eta_1\sum_{a=1}^{2}\sum_{b=1}^{2} x_{i} \xi_{a}\d{\eta_{b}}{\eta} - \xi_1\d{\eta_1}{\eta}\sum_{a=1}^{2}\sum_{b=1}^{2} y_{i} \d{\xi_{a}}{\xi}\eta_{b} \\
    =&  \d{\xi_1}{\xi} \sum_{a=1}^{2}\sum_{b=1}^{2} x_i  \d{\eta_{b}}{\eta} \xi_{a} \eta_1 - \d{\eta_1}{\eta} \sum_{a=1}^{2}\sum_{b=1}^{2} y_{i} \d{\xi_{a}}{\xi} \xi_1\eta_{b} \\
    =&  \d{\xi_1}{\xi} \sum_{a=1}^{2} \left\{ \sum_{b=1}^{2} x_i \d{\eta_{b}}{\eta} \right\} \xi_{a} \eta_1 -  \d{\eta_1}{\eta}  \sum_{b=1}^{2} \left\{ \sum_{a=1}^{2} y_{i} \d{\xi_{a}}{\xi} \right\} \xi_1\eta_{b},
\end{split}
\end{equation}
where
\begin{equation*}
\begin{split}
    \xi_{a} \eta_1 =& N_{1},N_{2} ,\quad a=1,2,\\
    \xi_1\eta_{b} =& N_{1},N_{3} ,\quad b=1,2.
\end{split}
\end{equation*}

In the $Q^1-P^0$ case, ordering the vertices as shown previously in Figure~\ref{fig:Q_1_reference_element}, Equation~\eqref{shape_function_gradient_simplified} evaluates cleanly to:
\begin{equation}\label{eq:Q1_velocity_shape_function_gradient_1}
   (\det J) \nabla N_1 =    \det J
    \begin{bmatrix}
        \d{N_1}{x}\\
        \d{N_1}{y}
    \end{bmatrix}
        =
    \begin{bmatrix}
        -N_1 y_{32}-N_2 y_{42}-N_3 y_{34}\\
        N_1 x_{32}+N_2 x_{42}+N_3 x_{34} 
    \end{bmatrix},
\end{equation}
and correspondingly for the remaining nodes:
\begin{equation}\label{eq:Q1_velocity_shape_function_gradient_234}
    \begin{split}
    (\det J) \nabla N_2&
         =
     \begin{bmatrix}
         -N_1 y_{13}+N_2 y_{41}-N_4 y_{34}\\
         N_1 x_{13}-N_2 x_{41}+N_4 x_{34} 
     \end{bmatrix},\\
     (\det J) \nabla N_3&
     =
    \begin{bmatrix}
     -N_1 y_{21}-N_3 y_{41}-N_4 y_{42}\\
     N_1 x_{21}+N_3 x_{41}+N_4 x_{42} 
    \end{bmatrix}, \\      
     (\det J) \nabla N_4&
     =
 \begin{bmatrix}
    -N_2 y_{21}-N_3 y_{13}+N_4 y_{32}\\
    N_2 x_{21}+N_3 x_{13}-N_4 x_{32}
 \end{bmatrix}.
    \end{split}
 \end{equation}

Applying the one-point Gauss quadrature rule collapses these expressions into the highly compact forms widely utilized in established legacy hydrocodes (Note: The subscript numbering here deviates slightly from standard conventions due to the specific nodal sequencing defined in Figure~\ref{fig:Q_1_reference_element}):
\begin{equation}\label{eq:Q1_velocity_shape_function_gradient_compact}
    \begin{split}
    &(\det J) \nabla N_1 (0,0) =
         \frac{1}{2}
     \begin{bmatrix}
         -y_{32}\\
         x_{32}
     \end{bmatrix}, \ 
     (\det J) \nabla N_2 (0,0) =
         \frac{1}{2}
     \begin{bmatrix}
         y_{41}\\
         -x_{41}
     \end{bmatrix}, \\
     &(\det J) \nabla N_3 (0,0) =
     \frac{1}{2}
 \begin{bmatrix}
     -y_{41}\\
     x_{41}
 \end{bmatrix}, \
 (\det J) \nabla N_4 (0,0) =
     \frac{1}{2}
 \begin{bmatrix}
     y_{32}\\
     -x_{32}
 \end{bmatrix}.
    \end{split}
 \end{equation}

\subsection*{$Q^2-Q^1$ case}
The 1D shape functions of the $Q^2$ element are documented in Equation~\eqref{eq:Q2_basis_kinematic_1D}. Recalling the intermediate variables $x_{\Lambda}$, $x_{\Theta}$, $x_{\Psi}$, $x_{\Delta}$, $x_{\Phi}$, and $x_{\Gamma}$ defined previously in Section~\ref{sec:gauss_quadrature}, these represent the explicit geometric combinations of the following summations:
\begin{equation*}
     \left\{ \sum_{a=1}^{3} x_{i} \d{\xi_{a}}{\xi} \right\} ,\quad \left\{ \sum_{a=1}^{3} x_{i} \d{\eta_{a}}{\eta} \right\},
\end{equation*}
with the orthogonal $y$ variables ($y_{\Lambda}$, $y_{\Theta}$, etc.) defined identically. Based on the preceding derivations, the condensed vector $[x_{\Lambda}, x_{\Theta}, x_{\Psi}]^T$ is cleanly generated as:
\begin{equation*}
 \begin{bmatrix}
    x_{\Lambda} \\
    x_{\Theta}\\
    x_{\Psi} 
\end{bmatrix}
=
\begin{bmatrix}
    x_1 & x_2 & x_3\\
    x_4 & x_5 & x_6\\
    x_7 & x_8 & x_9   
\end{bmatrix}
\begin{bmatrix}
    \d{\xi_1}{\xi}\\
    \d{\xi_2}{\xi}\\
    \d{\xi_3}{\xi}
\end{bmatrix}
\end{equation*}

\begin{equation*}
 \begin{bmatrix}
    x_{\Delta}&   x_{\Phi}&  x_{\Gamma} 
\end{bmatrix}
=
\begin{bmatrix}
    \d{\eta_1}{\eta}& \d{\eta_2}{\eta} & \d{\eta_3}{\eta}
\end{bmatrix}
\begin{bmatrix}
    x_1 & x_2 & x_3\\
    x_4 & x_5 & x_6\\
    x_7 & x_8 & x_9   
\end{bmatrix}
\end{equation*}

Consequently, the scalar components of the Jacobian matrix can be expressed directly in terms of these intermediate geometric variables and the known shape functions:
\begin{equation}\label{eq:Q2_Jacobian_element}
    \begin{split}
        J_{11}&=\eta_1 x_{\Lambda} + \eta_2 x_{\Theta} + \eta_3 x_{\Psi}, \\
        J_{12}&=\xi_1 x_{\Delta} +\xi_2 x_{\Phi} +\xi_3 x_{\Gamma} ,\\
        J_{21}&=\eta_1 y_{\Lambda} + \eta_2 y_{\Theta} + \eta_3 y_{\Psi}, \\
        J_{22}&=\xi_1 y_{\Delta} +\xi_2 y_{\Phi} +\xi_3 y_{\Gamma}. 
    \end{split}
\end{equation}

Expanded into full matrix formulation:
\begin{equation*}
J_{11}=
\begin{bmatrix}
    \eta_1 & \eta_2 & \eta_3
\end{bmatrix}
\begin{bmatrix}
    x_1 & x_2 & x_3\\
    x_4 & x_5 & x_6\\
    x_7 & x_8 & x_9   
\end{bmatrix}
\begin{bmatrix}
    \d{\xi_1}{\xi}\\
    \d{\xi_2}{\xi}\\
    \d{\xi_3}{\xi}
\end{bmatrix}
\end{equation*}

\begin{equation*}
J_{12}=
\begin{bmatrix}
   \d{\eta_1}{\eta}& \d{\eta_2}{\eta} & \d{\eta_3}{\eta}
\end{bmatrix}
\begin{bmatrix}
    x_1 & x_2 & x_3\\
    x_4 & x_5 & x_6\\
    x_7 & x_8 & x_9   
\end{bmatrix}
\begin{bmatrix}
    \xi_1\\
    \xi_2\\
    \xi_3
\end{bmatrix}
\end{equation*}

\begin{equation*}
J_{21}=
\begin{bmatrix}
    \eta_1 & \eta_2 & \eta_3
\end{bmatrix}
\begin{bmatrix}
    y_1 & y_2 & y_3\\
    y_4 & y_5 & y_6\\
    y_7 & y_8 & y_9
\end{bmatrix}
\begin{bmatrix}
    \d{\xi_1}{\xi}\\
    \d{\xi_2}{\xi}\\
    \d{\xi_3}{\xi}
\end{bmatrix}
\end{equation*}

\begin{equation*}
J_{22}=
\begin{bmatrix}
   \d{\eta_1}{\eta}& \d{\eta_2}{\eta} & \d{\eta_3}{\eta}
\end{bmatrix}
\begin{bmatrix}
    y_1 & y_2 & y_3\\
    y_4 & y_5 & y_6\\
    y_7 & y_8 & y_9
\end{bmatrix}
\begin{bmatrix}
    \xi_1\\
    \xi_2\\
    \xi_3
\end{bmatrix}
\end{equation*}
The coordinate matrix is kept the same for implementation simplicity.

\subsection*{$Q^3-Q^2$ case}
For the $Q^3-Q^2$ case, the 1D shape functions of the $Q^3$ element are provided in Equation~\eqref{eq:Q3_basis_kinematic_1D}. Omitting the explicit auxiliary variables for brevity, the exact matrix forms utilized to calculate the elements of the Jacobian matrix are seamlessly scaled as follows:
\begin{equation*}
J_{11}=
\begin{bmatrix}
    \eta_1 & \eta_2 & \eta_3& \eta_4
\end{bmatrix}
\begin{bmatrix}
    x_1 & x_2 & x_3 & x_4\\
    x_5 & x_6 & x_7 & x_8\\
    x_9 & x_{10} & x_{11} & x_{12}\\
    x_{13} & x_{14} & x_{15} & x_{16}
\end{bmatrix}
\begin{bmatrix}
    \d{\xi_1}{\xi}\\
    \d{\xi_2}{\xi}\\
    \d{\xi_3}{\xi}\\
    \d{\xi_4}{\xi}
\end{bmatrix}
\end{equation*}

\begin{equation*}
J_{12}=
\begin{bmatrix}
   \d{\eta_1}{\eta}& \d{\eta_2}{\eta} & \d{\eta_3}{\eta} & \d{\eta_4}{\eta}
\end{bmatrix}
\begin{bmatrix}
    x_1 & x_2 & x_3 & x_4\\
    x_5 & x_6 & x_7 & x_8\\
    x_9 & x_{10} & x_{11} & x_{12}\\
    x_{13} & x_{14} & x_{15} & x_{16}
\end{bmatrix}
\begin{bmatrix}
    \xi_1\\
    \xi_2\\
    \xi_3\\
    \xi_4
\end{bmatrix}
\end{equation*}

\begin{equation*}
J_{21}=
\begin{bmatrix}
    \eta_1 & \eta_2 & \eta_3& \eta_4
\end{bmatrix}
\begin{bmatrix}
    y_1 & y_2 & y_3 & y_4\\
    y_5 & y_6 & y_7 & y_8\\
    y_9 & y_{10} & y_{11} & y_{12}\\
    y_{13} & y_{14} & y_{15} & y_{16}
\end{bmatrix}
\begin{bmatrix}
    \d{\xi_1}{\xi}\\
    \d{\xi_2}{\xi}\\
    \d{\xi_3}{\xi}\\
    \d{\xi_4}{\xi}
\end{bmatrix}
\end{equation*}

\begin{equation*}
J_{22}=
\begin{bmatrix}
   \d{\eta_1}{\eta}& \d{\eta_2}{\eta} & \d{\eta_3}{\eta} & \d{\eta_4}{\eta}
\end{bmatrix}
\begin{bmatrix}
    y_1 & y_2 & y_3 & y_4\\
    y_5 & y_6 & y_7 & y_8\\
    y_9 & y_{10} & y_{11} & y_{12}\\
    y_{13} & y_{14} & y_{15} & y_{16}
\end{bmatrix}
\begin{bmatrix}
    \xi_1\\
    \xi_2\\
    \xi_3\\
    \xi_4
\end{bmatrix}
\end{equation*}

All the Jacobican matrix at each Gauss point can be computed using the above formulations with just two matrix multiplications.

\subsection{Mass matrix discretization}
We now address the discretization of the mass matrix originating from the left-hand side of the momentum equation (Equation~\eqref{eq:momentum_conservation}):
\begin{equation}\label{momentum_mass_matrix}
    \sum_{i=1}^{\text{kdof}} \frac{d  \vec{u}_{i}}{d t} \int_{K} \rho_h N_{i} N_j \dx\dy.
\end{equation}
The global discretization of the mass matrix can be approached in two primary ways. The first is to calculate the local mass matrix exactly and subsequently assemble and solve the resulting global linear system at each time step. This rigorous approach permits the unrestricted use of any polynomial shape function and ensures that the momentum conservation law is discretized with exact mathematical precision.

The second approach is to employ the mass lumping technique, which natively diagonalizes the local (and consequently global) mass matrix. This transformation enables the momentum equation to be solved completely explicitly, bypassing the prohibitive cost of solving a coupled linear system entirely. To achieve this, we adopt the Lobatto quadrature rule~\cite{Duczek2019Mass}, systematically utilizing its discrete quadrature points as the governing physical nodes to construct the shape functions. When kinematic shape functions are anchored precisely to these integration nodes, the Lobatto quadrature rule evaluates the local mass matrix as purely diagonal. While this under-integrated quadrature rule is technically inexact for integrating these continuous functions—and numerous references~\cite{Zienkiewicz2025Finite} document its adverse influence on formal convergence orders—it remains a highly advantageous and standard trade-off for speed in explicit hydrocodes.

We use the $Q^1-P^0$ and $Q^2-Q^1$ cases as fundamental examples to illustrate the stark difference between these two approaches. For the $Q^1-P^0$ case, the density is assumed to be constant within each element. For an element $K$, the generalized local mass matrix is represented as $M^K$:

\begin{equation}\label{Q1_mass_matrix}
M^K =
\rho_h(K)
\begin{bmatrix}
\int_{K} N_1 N_1 \dx\dy & \int_{K} N_1 N_2 \dx\dy  & \int_{K} N_1 N_3 \dx\dy & \int_{K} N_1 N_4 \dx\dy \\
\int_{K} N_2 N_1 \dx\dy & \int_{K} N_2 N_2 \dx\dy  & \int_{K} N_2 N_3 \dx\dy & \int_{K} N_2 N_4 \dx\dy \\
\int_{K} N_3 N_1 \dx\dy & \int_{K} N_3 N_2 \dx\dy  & \int_{K} N_3 N_3 \dx\dy & \int_{K} N_3 N_4 \dx\dy \\
\int_{K} N_4 N_1 \dx\dy & \int_{K} N_4 N_2 \dx\dy  & \int_{K} N_4 N_3 \dx\dy & \int_{K} N_4 N_4 \dx\dy
\end{bmatrix}.
\end{equation}

If exact integration via $2 \times 2$ Gauss quadrature is employed, the local mass matrix $M^K$ is exactly computed, but every resulting matrix entry remains non-zero. The first row of this exactly integrated local mass matrix evaluates explicitly to:
\begin{equation}\label{Q1_mass_matrix_1}
\begin{split}
 & (\frac{7}{36}+\frac{1}{3\sqrt{3}})\det J_{1} + \frac{1}{36} \det J_{2} + \frac{1}{36} \det J_{3} + \frac{1}{36} \det J_{4} \\
 & (\frac{1}{18}+\frac{1}{12\sqrt{3}})\det J_{1} + (\frac{1}{18}+\frac{1}{12\sqrt{3}})\det J_{2} + (\frac{1}{18}-\frac{1}{12\sqrt{3}}) \det J_{3} + (\frac{1}{18}-\frac{1}{12\sqrt{3}}) \det J_{4} \\
 & (\frac{1}{18}+\frac{1}{12\sqrt{3}})\det J_{1} + (\frac{1}{18}-\frac{1}{12\sqrt{3}})\det J_{2} + (\frac{1}{18}+\frac{1}{12\sqrt{3}}) \det J_{3} + (\frac{1}{18}-\frac{1}{12\sqrt{3}}) \det J_{4} \\  
 & \frac{1}{36} \det J_{1} + \frac{1}{36} \det J_{2} + \frac{1}{36} \det J_{3} + \frac{1}{36} \det J_{4} \\
\end{split}
\end{equation}
Here, $J_i$ represents the Jacobian determinant evaluated at the $i$-th interior Gauss quadrature point, and the entire expression is scaled globally by the uniform element density $\rho_h(K)$ to complete the row definition.

Conversely, the mass-lumping technique employs the discrete $2\times 2$ Lobatto quadrature points to approximate these spatial integrals. Because the $Q^1$ basis functions are precisely defined to be $1$ at their own node and $0$ at all others, this approach natively collapses $M^K$ into a perfectly diagonal matrix:
\begin{equation}\label{Q1_mass_matrix_lumped}
    M^K =
    \rho_h(K)
    \begin{bmatrix}
     \det J_{1}  & & & \\
     & \det J_{2}   & & \\
     & & \det J_{3}  & \\
     & & & \det J_{4}
    \end{bmatrix}
\end{equation}
where $J_i$ is now the Jacobian determinant evaluated strictly at the $i$-th Lobatto quadrature point (the geometric vertices).

Scaling up to the $Q^2-Q^1$ case, the local mass matrix expands to a $9\times 9$ tensor:
\begin{equation}\label{Q2_mass_matrix}
    M^K =
    \begin{bmatrix}
     \int_{K} \rho_h N_1 N_1 \dx\dy & \int_{K} \rho_h N_1 N_2 \dx\dy  & \cdots & \int_{K} \rho_h N_1 N_9 \dx\dy \\
     \int_{K} \rho_h N_2 N_1 \dx\dy & \int_{K} \rho_h N_2 N_2 \dx\dy  & \cdots & \int_{K} \rho_h N_2 N_9 \dx\dy \\
     \vdots & \vdots & \ddots & \vdots\\
     \int_{K} \rho_h N_9 N_1 \dx\dy & \int_{K} \rho_h N_9 N_2 \dx\dy  & \cdots & \int_{K} \rho_h N_9 N_9 \dx\dy
    \end{bmatrix},
\end{equation}
If exact calculation of $M^K$ is desired, a $3\times 3$ Gauss quadrature rule must be utilized. This strictly produces a dense, fully coupled local mass matrix—a calculation easily performed utilizing the shape functions and Jacobian determinants stored at the interior quadrature points. In aggressive contrast, employing the $3\times 3$ Lobatto quadrature completely zeroes the off-diagonal terms, yielding a strictly diagonal local mass matrix that eliminates subsequent solver coupling and drastically reduces computational assembly overhead.

\subsection{Nodal force assembly}
The conventional implicit approach to resolving the momentum conservation law involves assembling the local mass matrices and local right-hand side force vectors into their large-scale global equivalents. The momentum equation is then expressed in the monolithic matrix-vector form:
\begin{equation}
M \frac{d\vec{u}}{dt} = \vec{f},
\end{equation}
where $M$ is the global mass matrix, $\vec{u}$ is the unified nodal velocity vector, and $\vec{f}$ represents the global force vector. If the local mass matrices are non-diagonal (unlumped), simulating a time step requires the direct inversion or iterative solution of this massive, sparse global system. Subsequently, the velocity's time evolution is advanced utilizing standard ODE solvers. However, continuously assembling and solving this dense global mass matrix is prohibitively expensive for high-resolution transient shock dynamics, compelling the development of explicit bypass techniques.

Within the proposed explicit framework, the mass-lumping technique elegantly eliminates the need to ever assemble a coupled global mass matrix. Assuming the local lumped mass matrices remain constant across the immediate time increment, the strictly diagonal entries of each local mass matrix are simply accumulated at their shared nodes, producing a 1D discrete nodal mass vector that functions as a perfectly conserved variable across the timestep.

The right-hand side tracking force vector, however, must still be assembled globally by distributing local element contributions to the shared nodes. Nodal acceleration is then computed entirely explicitly by dividing this assembled global force vector by the lumped global nodal mass vector. The time integration of nodal velocities and physical positions is finally advanced using standard explicit ODE solvers.

To conclude this section, we add an important technical note regarding the localized right-hand side force components. The intermediate pressure force evaluated at each individual Gauss quadrature point must be explicitly cached in memory after the global nodal force vector is completed. These isolated, quadrature-point forces are strictly required for the subsequent correct discretization of the internal energy balance equation, which will be introduced in Section~\ref{sec:energy_conservation}. Furthermore, while the discrete momentum equations presented here form the theoretical skeleton of the dynamic framework, they are mathematically insufficient on their own to stably simulate highly deformed or shock-dominated physical flows. The mandatory numerical stabilization mechanisms—specifically, advanced hourglass control and artificial viscosity—will be systematically established in Section~\ref{sec:hourglass_and_viscosity}.
\section{Internal energy balance equation}\label{sec:energy_conservation}

This section addresses the discretization of the internal energy balance equation~\eqref{eq:lagrangian}. For a given element $K$, the integral form of the internal energy equation is:
\begin{equation}\label{eq:energy_integral}
  \int_{K} \rho_h \frac{d e_h}{d t} \phi_l \dx\dy =\int_{K} -  p_h \div{\vec{u}_h}  \phi_l \dx\dy.
\end{equation}
By discretizing the internal energy term while temporarily retaining density, pressure, and velocity in their continuum forms, we obtain the semi-discrete formulation:
\begin{equation}\label{eq:energy_conservation}
  \sum_{k=1}^{\text{tdof}} \frac{d e_k}{d t}\int_{K} \rho_h  \phi_k \phi_l \dx\dy  =\int_{K} -  p_h \div{\vec{u}_h}  \phi_l \dx\dy.
 \end{equation}

As with the momentum conservation law, both sides of Equation~\eqref{eq:energy_conservation} must be evaluated. However, because the internal energy is discretized using a discontinuous FEM space, assembling a global mass matrix or right-hand side vector is unnecessary. The temporal evolution of the internal energy is computed strictly at the element level.

The efficiency of this element-wise computation relies heavily on the nodal definition of the thermodynamic variables. For the $Q^m-Q^{m-1}$ FEM space pair, these variables are collocated with the $m^2$ Gauss quadrature points. Consequently, the shape functions $\phi_k$ naturally satisfy the Kronecker delta property at these specific points:
\begin{equation*}
    \phi_k(\gv{\xi}_q)=\delta_{k,q}, \quad \text{where} \quad
\delta_{k,q} = 
\begin{cases} 
    1, & k = q, \\ 
    0, & k \neq q.
\end{cases}
\end{equation*}

Applying this property to the left-hand side of Equation~\eqref{eq:energy_conservation}, the local energy mass matrix $M_T$ takes the form:
\begin{equation}\label{energy_mass_matrix}
    M_T = 
    \begin{bmatrix}
\int_{K} \rho_h  \phi_1 \phi_1 \dx\dy & \int_{K} \rho_h  \phi_2 \phi_1 \dx\dy & \cdots & \int_{K} \rho_h  \phi_{\text{tdof}} \phi_1 \dx\dy \\
\int_{K} \rho_h  \phi_1 \phi_2 \dx\dy & \int_{K} \rho_h  \phi_2 \phi_2 \dx\dy & \cdots & \int_{K} \rho_h  \phi_{\text{tdof}} \phi_2 \dx\dy \\
\vdots & \vdots & \ddots & \vdots \\
\int_{K} \rho_h  \phi_1 \phi_{\text{tdof}} \dx\dy & \int_{K} \rho_h  \phi_2 \phi_{\text{tdof}} \dx\dy & \cdots & \int_{K} \rho_h  \phi_{\text{tdof}} \phi_{\text{tdof}} \dx\dy 
    \end{bmatrix}.
\end{equation}
Because of the Kronecker delta property, all off-diagonal terms inherently vanish when integrated using the $m^2$-point Gauss quadrature rule. It is crucial to note that this yields a perfectly diagonal matrix without resorting to the artificial "mass lumping" required in the momentum equation; the mathematical evaluation remains exact. Utilizing the principle of mass conservation ($\rho \det J = \rho_0 \det J_0$), the $i$-th diagonal entry of $M_T$ simplifies strictly to:
\begin{equation}\label{element_mass}
    \begin{split}
    M_T(i,i) &= \int_{K} \rho_h  \phi_i \phi_i \dx\dy \\
    &= \sum_{q=1}^{m^2} \rho_h(\gv{\xi}_q) \phi_i(\gv{\xi}_q) \phi_i(\gv{\xi}_q) \omega_{q}\det J(\gv{\xi}_q)\\
    &=\rho_0 (\gv{\xi}_{i})\omega_{i}\det J_{0}(\gv{\xi}_{i}).
    \end{split}
\end{equation}

Similarly, for the right-hand side of Equation~\eqref{eq:energy_conservation}, the continuum formulation can be expanded as:
\begin{equation}\label{rhs_energy_discrete}
    \begin{split}
        \int_{K} -  p_h \div{\vec{u}_h}  \phi_l \dx\dy &=\int_{K} \sum_{k=1}^{\text{tdof}} -p_k\phi_k \left(\sum_{i=1}^{\text{kdof}}u_i \d{N_i}{x} +\sum_{i=1}^{\text{kdof}}v_i \d{N_i}{y}  \right) \phi_l \dx\dy, 
    \end{split}
\end{equation}
where $\vec{u}=[u,v]^T$. Applying the $m^2$-point Gauss quadrature rule and the Kronecker delta property effectively collapses the internal summations, elegantly yielding:
\begin{equation}\label{rhs_energy_discrete_2}
    \begin{split}
        \int_{K} -  p_h \div{\vec{u}_h}  \phi_l \dx\dy &=\sum_{q=1}^{m^2}\left\{\sum_{i=1}^{\text{kdof}}u_i \left(-p_h(\gv{\xi}_q)\d{N_i}{x} \det J(\gv{\xi}_q) \right) \right.\\
        &\quad \left. +\sum_{i=1}^{\text{kdof}}v_i \left(-p_h(\gv{\xi}_q)\d{N_i}{y} \det J(\gv{\xi}_q) \right)\right\} \phi_l(\gv{\xi}_q)\\
        &=\sum_{i=1}^{\text{kdof}}u_i \left(-p_h(\gv{\xi}_l)\d{N_i}{x} \det J(\gv{\xi}_l) \right) \\
        &\quad +\sum_{i=1}^{\text{kdof}}v_i \left(-p_h(\gv{\xi}_l)\d{N_i}{y} \det J(\gv{\xi}_l) \right).
    \end{split}
\end{equation}

The bracketed terms within Equation~\eqref{rhs_energy_discrete_2} are structurally identical to the discrete right-hand side of the momentum conservation law evaluated at the specified Gauss quadrature point $\gv{\xi}_l$. This structural reuse makes the internal energy evolution exceptionally efficient.

Ultimately, the variation in internal energy caused by macroscopic mechanical work is evaluated compactly point-by-point. The temporal evolution of the internal energy is obtained by simply dividing the local mechanical work~\eqref{rhs_energy_discrete_2} by the exact local mass scalar $M_T(i,i)$. Paired with a suitable ODE solver, Equation~\eqref{eq:energy_conservation} is solved explicitly and independently at each quadrature point.

Note that this formulation accounts exclusively for the reversible mechanical work performed by the thermodynamic fluid pressure. In high-order Lagrangian simulations, two distinct non-reversible mechanisms must be separately addressed: artificial viscosity, which provides the necessary dissipation for shock-capturing, and hourglass control, which is required to geometrically constrain non-physical zero-energy modes. The mechanical work contributed by both the artificial viscous forces and the anti-hourglass forces will be explicitly incorporated into this exact energy update in Section~\ref{sec:hourglass_and_viscosity}.
\section{Hourglass control and artificial viscosity}\label{sec:hourglass_and_viscosity}

The discrete momentum and energy equations derived in the previous sections establish the ideal, reversible core of the proposed high-order staggered Lagrangian framework. However, robust simulations require two additional critical components, which are addressed independently in this section: hourglass control, which geometrically constrains non-physical zero-energy modes to maintain thermodynamic consistency, and artificial viscosity, which provides the necessary mathematical dissipation to effectively capture shock waves.

While physically and functionally distinct, these two mechanisms are presented sequentially here because they natively share a unified quadrature architecture. Notably, within our framework, evaluating the artificial viscosity using an $(m+1)^2$-point Gauss-Legendre quadrature rule for the $Q^m-Q^{m-1}$ space pair inherently suppresses hourglass motion. This synergistic effect extends the findings from our previous study on lower-order $Q^1-P^0$ elements~\cite{Sun2022On}. In the following, we first analyze the fundamental origins of hourglass distortion in high-order spaces, followed by the derivation of the corresponding corrective algorithms and the explicit tensor viscosity formulations.

\subsection{Hourglass distortion and rank deficiency analysis}
Hourglass distortion is a major challenge that has plagued the staggered Lagrangian hydrodynamics community for many years. It is well-established that employing a one-point Gauss quadrature rule with a $Q^1-P^0$ FEM space induces zero-energy distortion. In this subsection, we apply a standard structural rank analysis~\cite{Flanagan1981uniform} to demonstrate that the $Q^2-Q^1$ and $Q^3-Q^2$ spaces within our proposed framework are also inherently susceptible to this geometric instability. 

Let the continuous velocity field in two dimensions be denoted by $\vec{u} = [u, v]^T$, discretized as:
\begin{equation*}
    u=\sum_{i=1}^{(m+1)^2} u_i N_i , \quad  v=\sum_{i=1}^{(m+1)^2} v_i N_i.
\end{equation*}

The strain-rate vector, denoted by $\gv{\dot{\epsilon}}$, is defined as:
\begin{equation}\label{eq:strain_rate}
    \gv{\dot{\epsilon}}=
    \begin{bmatrix}
        \d{u}{x} \\
        \d{v}{y} \\
        \d{u}{y}+\d{v}{x}
    \end{bmatrix}.
\end{equation}

By substituting the discretized velocity field, the discrete form of the strain rate at a given point is expressed as:
\begin{equation}\label{eq:strain_rate_discrete}
    \gv{\dot{\epsilon}}= B U = 
    \begin{bmatrix}
        \d{N_1}{x}&\dots&\d{N_{(m+1)^2}}{x} & 0 &\dots&0 \\
        0&\dots&0 &\d{N_1}{y}&\dots&\d{N_{(m+1)^2}}{y} \\
        \d{N_1}{y}&\dots&\d{N_{(m+1)^2}}{y}&\d{N_1}{x}&\dots&\d{N_{(m+1)^2}}{x}
    \end{bmatrix}
    \begin{bmatrix}
        u_1\\
        \vdots\\
        u_{(m+1)^2}\\
        v_1\\
        \vdots\\
        v_{(m+1)^2}
    \end{bmatrix},
\end{equation}
where $B$ is the discrete gradient operator matrix (with $3$ rows) and $U$ is the global velocity DOF vector for the element. When integrating over an element using an $m^2$-point Gauss quadrature rule, the matrices $B$ evaluated at each quadrature point are stacked vertically, resulting in a global constraint matrix of size $3m^2 \times 2(m+1)^2$.

We begin by examining the baseline $Q^1-P^0$ FEM space pair ($m=1$). The velocity field has 8 DOFs. Subtracting the 3 DOFs corresponding to rigid body motion (two translations, one rotation), the required full rank of the stacked $B$ matrix to prevent zero-energy modes is 5. When a one-point Gauss quadrature rule is employed ($m^2=1$), the stacked matrix contains only 3 rows, capping its maximum possible rank at 3. This strict lack of 2 constraint DOFs is the fundamental mathematical cause of hourglass distortion. In contrast, applying a four-point Gauss quadrature rule yields a stacked matrix with 12 rows, which easily accommodates the required rank of 5, formally demonstrating that the four-point rule naturally suppresses hourglass distortion.

For the $Q^2-Q^1$ element ($m=2$), the velocity field possesses 18 DOFs. After subtracting the 3 rigid body DOFs, the required rank is 15. However, the four-point Gauss quadrature rule ($m^2=4$) yields a stacked matrix with only 12 rows. This rank deficiency of 3 inevitably leaves zero-energy modes unconstrained, dictating that explicit hourglass control must be introduced.

For the $Q^3-Q^2$ element ($m=3$), the velocity field has 32 DOFs. Excluding the 3 rigid body DOFs, the required rank is 29. Employing a nine-point Gauss quadrature rule ($m^2=9$) produces a stacked matrix with 27 rows. The resulting rank deficiency of 2 dictates the existence of two unconstrained zero-energy modes, again confirming the mathematical necessity of hourglass control.

Interestingly, for cases where $m \ge 4$, this specific matrix row-counting technique becomes theoretically inconclusive. Because the number of constraint rows scales as $3m^2$, it automatically equals or exceeds the required structural rank of $2(m+1)^2-3$ for $m \ge 4$. Therefore, while the theoretical rank bounds for $m \ge 4$ require alternative spectral analyses deferred to future work, the practical implementation and validation of this manuscript deliberately focus on the $Q^1-P^0$ through $Q^3-Q^2$ spaces, covering the most highly utilized regimes in production hydrocodes.

\subsection{The high-order hourglass control algorithm}
To constrain the zero-energy modes demonstrated above, we recall the right-hand side of the semi-discrete momentum equation, which represents the internal force driven by the fluid pressure:
\begin{equation*}\label{eq:momentum_conservation_2}
    \vec{F}^{\text{int}}_j = \int_{K} p_h \nabla N_j \dx\dy.
\end{equation*}

Because the actual pressure field $p_h$ resides in the $Q^{m-1}$ space, evaluating this integral using an $m^2$-point Gauss quadrature rule yields an exact result but provides an insufficient number of geometric constraints to prevent hourglassing:
\begin{equation}\label{p_h_RHS}
    \int_{K} p_h \nabla N_j \dx\dy=  \sum_{q=1}^{m^2} \omega_{q} p_h(\gv{\xi}_q) \nabla N_j(\gv{\xi}_q) \det J(\gv{\xi}_q).
\end{equation}

Conceptually, if we assume the existence of an enriched pressure polynomial $\tilde{p} \in Q^{m}$ possessing additional DOFs, the corresponding internal force would successfully eliminate the hourglass deficiencies. Because $\tilde{p}$ is of higher order, this enriched term must be integrated using an $(m+1)^2$-point Gauss quadrature rule:
\begin{equation}\label{p_tilde_RHS}
    \int_{K} \tilde{p}\nabla N_j \dx\dy  = \sum_{q=1}^{(m+1)^2} \omega_{q}\tilde{p}(\gv{\xi}_q) \nabla N_j(\gv{\xi}_q) \det J(\gv{\xi}_q).
\end{equation}

Inspired by the subzonal pressure hourglass control methods developed for the $Q^1-P^0$ case~\cite{Caramana1998Elimination,Sun2022Understanding}, we construct the restoring anti-hourglass force, $\vec{f}^{\text{hg}}_j$, directly from the mathematical difference between these two formulations. Crucially, because $p_h \in Q^{m-1}$, its integral is also evaluated exactly by the $(m+1)^2$-point rule. This allows us to combine the integrals point-by-point over the enriched quadrature points:
\begin{equation}\label{eq:hg_force}
    \begin{split}
        \vec{f}^{\text{hg}}_{j} &= \int_{K} \tilde{p} \nabla N_j \dx\dy - \int_{K} p_h \nabla N_j \dx\dy \\
        &=\sum_{q=1}^{(m+1)^2} \omega_{q} \delta p(\gv{\xi}_q) \nabla N_j(\gv{\xi}_q) \det J(\gv{\xi}_q),
    \end{split}
\end{equation}
where $\delta p = \tilde{p}-p_h$ represents the subzonal pressure variation. Following standard practice~\cite{Wilkins1963Calculation,Caramana1998Elimination}, this variation is approximated using the local acoustic relation:
\begin{equation}\label{pressure_variation}
    \delta p \approx c_s^2 \delta \rho. 
\end{equation}

To evaluate Equation~\eqref{pressure_variation} at the $(m+1)^2$ points, we must project the baseline sound speed and density from their native $m^2$ points using an interpolation matrix $M_{\text{interp}}$ defined by the $Q^{m-1}$ shape functions:
\begin{equation}\label{eq:interpolation}
    \tilde{\mathbf{c}}_s = M_{\text{interp}} \mathbf{c}_{s,h}, \quad 
    \delta \gv{\rho} = \tilde{\gv{\rho}} - M_{\text{interp}} \gv{\rho}_h,
\end{equation}
where $\tilde{\mathbf{c}}_s$, $\tilde{\gv{\rho}}$, and $\delta \gv{\rho}$ are vectors of length $(m+1)^2$, while $\mathbf{c}_{s,h}$ and $\gv{\rho}_h$ are vectors of length $m^2$. The enriched subzonal density $\tilde{\rho}$ is obtained point-by-point by enforcing strong mass conservation at the $(m+1)^2$ points:
\begin{equation}\label{discrete_mass_conservation_mp1}
    \tilde{\rho}^{n}(\gv{\xi}_{q}) \det J^{n}(\gv{\xi}_{q})=\rho_{0}(\gv{\xi}_{q}) \det J_{0}(\gv{\xi}_{q}), \quad \text{for } q=1,\dots,(m+1)^2.
\end{equation}

\subsubsection*{Case-by-case implementation}

\textbf{The $Q^1-P^0$ case:} 
The interpolation matrix $M_{\text{interp}}$ reduces to a simple $4 \times 1$ column vector $[1,1,1,1]^T$. The high-order density variable $\tilde{\rho}$ resides in the $Q^1$ space, and the discrete mass conservation~\eqref{discrete_mass_conservation_mp1} is independently satisfied at the four standard Gauss points $\gv{\xi}_q \in \{ (\pm 1/\sqrt{3}, \pm 1/\sqrt{3}) \}$.

\textbf{The $Q^2-Q^1$ case:} 
The interpolation matrix $M_{\text{interp}}$ is of size $9 \times 4$. It is constructed by evaluating the four $Q^1$ shape functions (from Equation~\eqref{eq:Q1_basis_thermodynamic_1D}) at the nine $Q^2$ Gauss quadrature points, yielding:
\begin{equation*}
    M_{\text{interp}} = 
    \begin{bmatrix}
        0.7+1.5\sqrt{1/5} & -0.2 & -0.2 & 0.7-1.5\sqrt{1/5} \\
        0.25+0.75\sqrt{1/5} & 0.25+0.75\sqrt{1/5} & 0.25-0.75\sqrt{1/5} & 0.25-0.75\sqrt{1/5} \\
        -0.2 & 0.7+1.5\sqrt{1/5} & 0.7-1.5\sqrt{1/5} & -0.2\\
        0.25+0.75\sqrt{1/5} & 0.25-0.75\sqrt{1/5} & 0.25+0.75\sqrt{1/5} & 0.25-0.75\sqrt{1/5} \\
        0.25 & 0.25 & 0.25 & 0.25 \\
        0.25-0.75\sqrt{1/5} & 0.25+0.75\sqrt{1/5} & 0.25-0.75\sqrt{1/5} & 0.25+0.75\sqrt{1/5} \\
        -0.2 & 0.7-1.5\sqrt{1/5} &  0.7+1.5\sqrt{1/5} & -0.2 \\
        0.25-0.75\sqrt{1/5} & 0.25-0.75\sqrt{1/5} & 0.25+0.75\sqrt{1/5} & 0.25+0.75\sqrt{1/5} \\
        0.7-1.5\sqrt{1/5} & -0.2 & -0.2 & 0.7+1.5\sqrt{1/5}
    \end{bmatrix}.  
\end{equation*}
The enriched density field $\tilde{\rho}$ is subsequently computed directly from Equation~\eqref{discrete_mass_conservation_mp1} at these $3\times 3$ integration points.

\textbf{The $Q^3-Q^2$ case:}
Similarly, the $16 \times 9$ interpolation matrix is constructed by evaluating the nine $Q^2$ shape functions(Equation~\eqref{eq:Q2_basis_thermodynamic_1D}~\eqref{eq:Q2_basis_thermodynamic}) at the sixteen Gauss quadrature points. The explicit matrix form is omitted here for brevity; within our hydrocode, this matrix is populated algorithmically via standard shape function subroutines rather than being hard-coded. The enriched density field $\tilde{\rho}$ is then obtained at these $4\times 4$ points by Equation~\eqref{discrete_mass_conservation_mp1}.

\subsection{The tensor artificial viscosity algorithm}
To capture shock waves and prevent unphysical oscillations, artificial viscosity is introduced into the momentum equation, yielding $\rho \frac{d \vec{u} }{d t} = \div (\mu \gv{\sigma}_{a})$. Following the approach proposed in~\cite{Dobrev2012High}, we define $\gv{\sigma}_{a}$ using the symmetrized velocity gradient tensor:
\begin{equation}\label{viscosity_form_symmetry}
    \gv{\sigma}_{a}= \epsilon(\vec{u}) = \frac{1}{2} (\nabla \vec{u}+ \vec{u}\nabla)=
    \begin{bmatrix}
       \frac{\partial u }{\partial x }  & \frac{1}{2}\left(\frac{\partial v }{\partial x }+\frac{\partial u }{\partial y }\right)\\
       \frac{1}{2}\left(\frac{\partial v }{\partial x }+\frac{\partial u }{\partial y }\right)  & \frac{\partial v }{\partial y }
    \end{bmatrix}.
\end{equation}

Applying the Petrov-Galerkin method, the discrete viscous force vector for node $j$ (shown here for the $x$-direction) evaluates to:
\begin{equation}\label{eq:viscosity_rhs_discretization}
    \begin{split}
    f^{\text{visc}}_{j,R} = \int_{K} \mu \nabla u_h \cdot \nabla N_j \dx\dy &= \sum_{q=1}^{(m+1)^2}  \mu_q \omega_{q} \frac{1}{\det J(\gv{\xi}_q)} \left ( \sum_{i=1}^{\text{kdof}}u_i \nabla  N_i(\gv{\xi}_q)\det J(\gv{\xi}_q) \right) \cdot \nabla N_j(\gv{\xi}_q)\det J(\gv{\xi}_q).
    \end{split} 
\end{equation}
Crucially, the geometric gradient term $\nabla N_j(\gv{\xi}_q) \det J(\gv{\xi}_q)$ is precomputed at the $(m+1)^2$ Gauss quadrature points. This allows for direct reuse between the hourglass control and artificial viscosity calculations, maximizing computational efficiency.

\subsubsection*{Example: $Q^1-P^0$ Tensor Viscosity Force}
To illustrate the difference between classical explicit unrolling and our algorithmic approach, consider the baseline $Q^1$ element with nodes sequenced as in Figure~\ref{fig:Q_1_reference_element}. The resulting nodal viscous force for node 1 is:
\begin{equation}\label{Q_1_viscosity_force_1R}
    \begin{split}
        \int_{K} \mu \nabla u_h \cdot \nabla N_1 \dx\dy = \sum_{q=1}^{4} \mu_q \frac{1}{\det J(\gv{\xi}_q)} \left(\nabla u_h(\gv{\xi}_q)\det J(\gv{\xi}_q)\right) \cdot \left(\nabla N_1(\gv{\xi}_q) \det J(\gv{\xi}_q)\right).
    \end{split}
\end{equation}

For the $Q^1$ element, the nodes are collocated with the four Lobatto quadrature points, meaning the shape functions inherently satisfy the Kronecker delta property ($N_i(\gv{\xi}_q) = \delta_{i,q}$). Consequently, the gradient terms $\nabla N_1 \det J$ evaluated at each integration point can be represented explicitly by the nodal coordinates:
\begin{equation*}
\nabla N_1 (\gv{\xi}_1) \det J(\gv{\xi}_1) = \begin{bmatrix} -y_{32} \\ x_{32} \end{bmatrix}, \quad
\nabla N_1 (\gv{\xi}_2) \det J(\gv{\xi}_2) = \begin{bmatrix} -y_{42} \\ x_{42} \end{bmatrix}, \quad
\nabla N_1 (\gv{\xi}_3) \det J(\gv{\xi}_3) = \begin{bmatrix} -y_{34} \\ x_{34} \end{bmatrix},
\end{equation*}
while $\nabla N_1 (\gv{\xi}_4) \det J(\gv{\xi}_4) = [0, 0]^T$. The gradients for the remaining shape functions are obtained analogously.

The discrete velocity gradients at the four integration points then expand into coordinate-based expressions:
\begin{equation*}
 \begin{split}
 \nabla u_h(\gv{\xi}_1) \det J(\gv{\xi}_1) &=
 u_1
 \begin{bmatrix}
     -y_{32}\\
    x_{32} 
 \end{bmatrix}
 +
 u_2
 \begin{bmatrix}
     -y_{13}\\
    x_{13} 
 \end{bmatrix}
 +
 u_3
 \begin{bmatrix}
     -y_{21}\\
    x_{21} 
 \end{bmatrix}, \\
 \nabla u_h(\gv{\xi}_2) \det J(\gv{\xi}_2)  &=
 u_1
 \begin{bmatrix}
     -y_{42}\\
    x_{42} 
 \end{bmatrix}
 +
 u_2
 \begin{bmatrix}
     y_{41}\\
    -x_{41} 
 \end{bmatrix}
 +
 u_4
 \begin{bmatrix}
     -y_{21}\\
    x_{21} 
 \end{bmatrix}, \\
 \nabla u_h(\gv{\xi}_3) \det J(\gv{\xi}_3)  &=
 u_1
 \begin{bmatrix}
     -y_{34}\\
    x_{34} 
 \end{bmatrix}
 +
 u_3
 \begin{bmatrix}
     -y_{41}\\
    x_{41} 
 \end{bmatrix}
 +
 u_4
 \begin{bmatrix}
     -y_{13}\\
    x_{13} 
 \end{bmatrix}, \\
 \nabla u_h(\gv{\xi}_4) \det J(\gv{\xi}_4)  &=
 u_2
 \begin{bmatrix}
     -y_{34}\\
    x_{34} 
 \end{bmatrix}
 +
 u_3
 \begin{bmatrix}
     -y_{42}\\
    x_{42} 
 \end{bmatrix}
 +
 u_4
 \begin{bmatrix}
     y_{32}\\
    -x_{32} 
 \end{bmatrix}.
 \end{split}
 \end{equation*}

Substituting these into the weak form yields an explicit algebraic representation of the tensor viscosity force for node $1$:
\begin{equation}\label{Q_1_viscosity_force_1Explicit}
\begin{split}
\int_{K} \mu \nabla u_h \cdot \nabla N_1 \dx\dy &= \sum_{q=1}^{4} \mu_q \frac{1}{\det J(\gv{\xi}_q)} \left( \nabla u_h(\gv{\xi}_q)\det J(\gv{\xi}_q) \right) \cdot \left( \nabla N_1(\gv{\xi}_q) \det J(\gv{\xi}_q) \right) \\
& = \frac{\mu_1}{\det J (\gv{\xi}_1)} \left[ u_1(x_{32}^2+y_{32}^2) + u_2(x_{13} x_{32} + y_{13} y_{32}) + u_3 (x_{21}x_{32}+y_{21} y_{32}) \right] \\
& + \frac{\mu_2}{\det J (\gv{\xi}_2)} \left[ u_1(x_{42}^2+y_{42}^2) + u_2 (-x_{41}x_{42}-y_{41}y_{42}) + u_4(x_{21}x_{42}+y_{21}y_{42}) \right] \\
& + \frac{\mu_3}{\det J (\gv{\xi}_3)} \left[ u_1(x_{34}^2+y_{34}^2) + u_3(x_{41}x_{34}+y_{41}y_{34}) + u_4(x_{13}x_{34}+y_{13}y_{34}) \right].
\end{split}
\end{equation}

Equation~\eqref{Q_1_viscosity_force_1Explicit} is consistent with classical formulations~\cite{Campbell2001tensor,Kolev2009tensor,Sun2022On}. However, while this explicit unrolling allows for deep variable reuse in the $Q^1$ case, it scales disastrously to higher-order discretizations. Hard-coding such expansions for $Q^2$ or $Q^3$ elements demands an excessive number of temporary variables, obscuring code logic and complicating maintenance. 

Therefore, for the general $Q^m-Q^{m-1}$ framework, we advocate for the algorithmic summation form presented in Equation~\eqref{eq:viscosity_rhs_discretization}. By leveraging $(m+1)^2$-point Gauss quadrature and precomputed structured shape function derivatives ($\nabla N_i \det J$), the velocity gradients and symmetrized viscous stresses $\gv{\sigma}_{a}$ are reconstructed dynamically via nodal summation. This approach significantly enhances code clarity and extensibility with only a negligible increase in floating-point operations.

\subsection{Viscosity parameter and time step control}
With the spatial discretization established, the scalar artificial viscosity parameter $\mu$ must be evaluated at the $(m+1)^2$ active quadrature points:
\begin{equation}\label{viscosity_vorticity_parameter}
    \mu = \rho \left( c_1 c_{\text{vor}} c_{s} l_{c} + c_2 l_{c}^2 |\lambda_m| \right),    
\end{equation}
where $\lambda_m$ is the minimum eigenvalue of the symmetrized velocity gradient $\epsilon(\vec{u})$, serving as an active compression indicator:
\begin{equation}\label{eigenvalue}
    \lambda_{m}= \frac{1}{2}\left(\frac{\partial u }{\partial x } + \frac{\partial v }{\partial y }\right) - \frac{1}{2}\sqrt{\left(\frac{\partial u }{\partial x } - \frac{\partial v }{\partial y }\right)^2 + \left(\frac{\partial u }{\partial y }+\frac{\partial v}{\partial x }\right)^2}.
\end{equation}
The corresponding unnormalized eigenvector $\vec{e}$ is given by:
\begin{equation}
    \vec{e}=
    \begin{bmatrix}
        \frac{1}{2}\left(\frac{\partial v }{\partial x }+\frac{\partial u }{\partial y }\right) \\[1ex]
        \lambda_m - \frac{\partial u }{\partial x }    
    \end{bmatrix}.   
\end{equation}

The characteristic element length $l_c$ along the direction of maximum compression is computed by mapping the normalized eigenvector $\vec{n}_e = \vec{e}/|\vec{e}|$ through the deformation Jacobian:
\begin{equation}\label{characteristic_length}
 l_c = l_0 \frac{|(J^0)^{-1}J \vec{e}|}{|\vec{e}|} = l_0 \sqrt{\frac{1}{\det J^0} \left[ M_{11}^2 + M_{12}^2 + M_{21}^2 + M_{22}^2 + 2(M_{11}M_{12} + M_{21}M_{22})ne_{1}ne_{2} \right]},
\end{equation}
where $M_{ij}$ are the standard components of the mapped matrix $(J^0)^{-1} J$. 

The baseline characteristic length $l_0$ requires careful treatment on unstructured meshes. A strictly global average, such as $l_{\text{global}} = \frac{1}{m}\sqrt{Vol / nele}$, fails to account for severe local size variations. Thus, we allow $l_0$ to vary element-wise within bounded limits $l_{min} \le l_0 \le l_{max}$, defined as:
\begin{equation*}
\begin{aligned}
    l_{min} &= \frac{1}{m} \min\left(\sqrt{Vol_{K}}, \sqrt{\frac{Vol}{nele}}\right), \\
    l_{max} &= \frac{1}{m} \max\left(\sqrt{Vol_{K}}, \sqrt{\frac{Vol}{nele}}\right).
\end{aligned}    
\end{equation*}

Following~\cite{Dobrev2012High}, a dimensionless vorticity modifier $c_{\text{vor}}$ is multiplied with the linear bulk term to suppress spurious viscosity in shear-dominated flows. It is computed using previously evaluated spatial derivatives:
\begin{displaymath}
    c_{\text{vor}} = \frac{|\div \vec{u}|}{ \|\grad \vec{u}\|} = \frac{\left| \frac{\partial u }{\partial x } + \frac{\partial v }{\partial y } \right|}{\sqrt{ \left(\frac{\partial u }{\partial x }\right)^2 + \left(\frac{\partial u }{\partial y }\right)^2 + \left(\frac{\partial v }{\partial x }\right)^2 + \left(\frac{\partial v }{\partial y }\right)^2 }}.
\end{displaymath}

Because thermodynamic variables (sound speed $c_s$ and density $\rho$) are natively defined at the $m^2$ quadrature points, they must be projected to the $(m+1)^2$ evaluation points. We seamlessly reuse the interpolation matrix $M_{\text{interp}}$ and the enriched subzonal density $\tilde{\rho}$ derived during the hourglass control formulation (Equation~\ref{eq:interpolation}) for this task.

Finally, to guarantee numerical stability, the global operational time-step $\tau$ is dictated by the most restrictive condition across all $(m+1)^2$ active quadrature points in the domain:
\begin{displaymath}
    \tau = \text{CFL} \min_{\text{all } K, \gv{\xi}_q}  \left(\frac{c_s}{l_{\tau}} + \frac{\mu}{\rho l_{\tau}^2}\right)^{-1},
\end{displaymath}
where $l_{\tau}^2$ is derived analytically from the minimum singular value of the Jacobian matrix $J$:
\begin{displaymath}
    l_{\tau}^2 = \frac{1}{2}(J_{11}^2+J_{12}^2+J_{21}^2+J_{22}^2) - \frac{1}{2}\sqrt{(J_{11}^2+J_{21}^2-J_{12}^2-J_{22}^2)^2 + 4(J_{11}J_{12}+J_{21}J_{22})^2}.
\end{displaymath}

\subsection{Consistent Internal Energy Update}
To maintain strict thermodynamic consistency, the mechanical work done by both the hourglass forces and the artificial viscous force must be perfectly captured and converted into internal energy. The total external mathematical force acting on node $j$ at the $(m+1)^2$ quadrature points $\gv{\xi}_q$ is simply their sum: $\vec{f}^{\text{mp}}_{j}(\gv{\xi}_q) = \vec{f}^{\text{hg}}_{j}(\gv{\xi}_q) + \vec{f}^{\text{visc}}_{j}(\gv{\xi}_q)$. 

The discrete internal energy update is formulated to seamlessly absorb this composite mechanical work:
\begin{equation}\label{eq:energy_hg_av_discrete}
\sum_{k=1}^{\text{tdof}} \frac{d e_k}{d t}\int_{K} \rho_h  \phi_k \phi_l \dx\dy  =  \sum_{q=1}^{(m+1)^2} \omega_{q} \left(\sum_{j=1}^{\text{kdof}}\vec{u}_j \cdot \vec{f}^{\text{mp}}_{j} (\gv{\xi}_q)\right)  \phi_l(\gv{\xi}_q).
\end{equation}
By computing the dissipative energy addition via the exact dot product of the nodal velocities and the derived total forces, we guarantee exact total energy conservation. Furthermore, this approach achieves maximum computational efficiency by ensuring that the expensive shape function gradients are evaluated exactly once during the momentum sweep and reused directly for the point-by-point energy update.
\section{A mass-lumped discretization of the original high-order formulation}\label{sec:original_discretization}

In this section, we present a mass-lumped discretization of the high-order staggered Lagrangian formulation originally proposed by Dobrev et al.~\cite{Dobrev2011Curvilinear,Dobrev2012High}, which will serve as a rigorous comparative baseline for our numerical experiments in Section~\ref{sec:numerical_results}. While the finite element approximation spaces in this baseline are identical to those defined in Section~\ref{sec:discretization}, this formulation differs fundamentally in its treatment of the momentum mass matrix—which is approximated by a diagonal lumped matrix rather than evaluated exactly—and in its handling of the discrete mass conservation law. Conversely, the internal energy mass matrix remains exact, benefiting from the localized collocation of the thermodynamic variables at the standard $m^2$ Gauss quadrature points.

For the $Q^{m}-Q^{m-1}$ finite element pair, the continuous mass conservation law~\eqref{eq:mass_integral} in this baseline approach is enforced strictly at the enriched $(m+1)^2$ Gauss quadrature points within each element. As analyzed in Section~\ref{sec:gauss_quadrature}, evaluating density at these $(m+1)^2$ points fully engages the capacity of the $Q^m$ space, aligning the density degrees of freedom (DOFs) with the kinematic DOFs ($\text{kdof}$). Consequently, the density field possesses 9 DOFs for the $Q^2-Q^1$ case and 16 DOFs for the $Q^3-Q^2$ case. In practice, this means we directly adopt the enriched subzonal density $\tilde{\rho}$ from Equation~\eqref{discrete_mass_conservation_mp1} as the primary density field $\rho_h$. 

To avoid any potential ambiguity for the reader, we explicitly specify the notations and definitions for the thermodynamic variables utilized in this section:
\begin{equation}
\begin{array}{r@{\quad : \quad}l}
\rho_h & \text{density evaluated at the $(m+1)^2$ quadrature points}, \\
e_h & \text{internal energy defined natively at the $m^2$ quadrature points}, \\
\tilde{e} & \text{internal energy interpolated to the $(m+1)^2$ quadrature points}, \\
p_h & \text{pressure evaluated at the $(m+1)^2$ quadrature points}.
\end{array}
\end{equation}

Regarding the momentum equation, the right-hand side integral from Equation~\eqref{eq:momentum_integral} is integrated using the $(m+1)^2$-point Gauss quadrature rule rather than the under-integrated $m^2$-point rule. Thus, Equation~\eqref{eq:discrete_rhs_momentum_conservation} is replaced by:
\begin{equation}\label{eq:discrete_rhs_momentum_conservation_original}
    \int_{K} p_h \nabla N_j \dx\dy = \sum_{q=1}^{(m+1)^2} \omega_{q} p_h(\gv{\xi}_q)\nabla N_j(\gv{\xi}_q)  \det J(\gv{\xi}_q).
\end{equation}

To evaluate this enriched integral, the pressure $p_h$ must be determined at all $(m+1)^2$ quadrature points. Because the internal energy $e_h$ natively resides in the $m^2$ points of the $Q^{m-1}$ space, it must first be projected from its native points to the enriched $(m+1)^2$ evaluation nodes using the interpolation matrix $M_{\text{interp}}$ defined in Section~\ref{sec:hourglass_and_viscosity}:
\begin{equation}\label{eq:internal_energy_interpolation}
  \begin{bmatrix} \tilde{e}_1 \\ \vdots \\ \tilde{e}_{(m+1)^2} \end{bmatrix} = M_{\text{interp}}  \begin{bmatrix} e_{1} \\ \vdots \\ e_{m^2} \end{bmatrix}.
\end{equation}
With both the density $\rho_h$ and the interpolated internal energy $\tilde{e}$ available at the $(m+1)^2$ points, the pressure is updated pointwise via the equation of state (EOS):
\begin{equation}
p_h(\gv{\xi}_q) = \text{EOS}(\rho_h(\gv{\xi}_q), \tilde{e}(\gv{\xi}_q)), \quad \text{for } q=1,2,\dots,(m+1)^2.    
\end{equation}

A notable characteristic of this enriched integration scheme is that it naturally suppresses zero-energy hourglass modes by utilizing the higher-order $(m+1)^2$ quadrature rule for internal force evaluation, rendering explicit subzonal anti-hourglass control unnecessary. However, tensor artificial viscosity is still required to capture shock waves and maintain stability, which is computed exactly as detailed in Section~\ref{sec:hourglass_and_viscosity}. 

To maintain strict thermodynamic consistency, the discretization of the internal energy equation mirrors that presented in Section~\ref{sec:energy_conservation}. The total internal force vector $\vec{f}_{j}(\gv{\xi}_q)$ at each enriched quadrature point is formulated as the sum of the physical pressure force and the artificial viscous force:
\begin{equation}\label{eq:force_original}
 \vec{f}_{j}(\gv{\xi}_q) = \vec{f}^{\text{pre}}_{j}(\gv{\xi}_q) + \vec{f}^{\text{visc}}_{j}(\gv{\xi}_q).
\end{equation}
The temporal evolution of the internal energy is then driven by the mechanical work done by these combined forces, evaluated consistently over the $(m+1)^2$ quadrature points:
\begin{equation}\label{eq:energy_discrete_original}
\sum_{k=1}^{\text{tdof}} \frac{d e_k}{d t}\int_{K} \rho_h  \phi_k \phi_l \dx\dy  =  \sum_{q=1}^{(m+1)^2} \omega_{q} \left(\sum_{j=1}^{\text{kdof}}\vec{u}_j \cdot \vec{f}_{j} (\gv{\xi}_q)\right)  \phi_l(\gv{\xi}_q).
\end{equation}

Because both the momentum and energy mass matrices are diagonal in this formulation, their respective temporal evolution equations can be solved explicitly and efficiently through simple scalar division by the lumped nodal masses and local quadrature point masses. Time integration is performed using conservative implicit-explicit Runge-Kutta (IMEX-RK) schemes~\cite{Sandu2021Conservative}, ensuring that total energy remains exactly conserved. 

Numerical results generated using this mass-lumped discretization scheme will be presented in Section~\ref{sec:numerical_results} to serve as a direct comparative baseline against our newly proposed framework, specifically highlighting the advantages of our tunable hourglass control over the rigid stiffness of the original formulation.
\section{Numerical examples}\label{sec:numerical_results}
The numerical performance of our proposed framework is evaluated through a comprehensive set of benchmarks, with a primary focus on the $Q^2-Q^1$ and $Q^3-Q^2$ element pairs. The evaluation is structurally divided into two distinct phases based on their verification objectives:
\begin{itemize}
    \item \textbf{Accuracy and Convergence:} Using smooth analytical solutions (the Taylor-Green vortex and the Kidder problem) to rigorously verify the optimal convergence rates of the high-order spatial discretizations without the dissipative influence of artificial viscosity.
    \item \textbf{Robustness and Stabilization:} Using strong shock-driven problems (Sod, Noh, Sedov, etc.) to validate the robustness of the framework, the necessity of the proposed hourglass control mechanism, and the shock-capturing capability via artificial viscosity.
\end{itemize}

\subsection{Time integration scheme}\label{subsec:time_integration_summary}
Before presenting the numerical results, we briefly outline the time integration strategy employed across all simulations. To rigorously guarantee total energy conservation, we utilize implicit-explicit Runge-Kutta (IMEX-RK) schemes~\cite{Sandu2021Conservative}. The formal order of the temporal discretization is strictly chosen to be one order higher than the spatial discretization to render temporal integration errors negligible; i.e., the $Q^2-Q^1$ pair is coupled with a 3-stage IMEX-RK3 method, and the $Q^3-Q^2$ pair with a 4-stage IMEX-RK4 method.

As an illustrative example, we consider the second-order, energy-conserving RK2-average scheme. In this approach, the internal energy $e$ is evaluated strictly through the specific work done by the explicitly advanced nodal forces $\vec{f}$ over the velocity field $\vec{u}$. The scheme is formulated as follows:

The first stage of the update is given by:
\begin{equation}
\begin{bmatrix}
\vec{u}^{n+\frac{1}{2}}\\
\vec{x}^{n+\frac{1}{2}}\\
e^{n+\frac{1}{2}}
\end{bmatrix}
=
\begin{bmatrix}
\vec{u}^{n}\\
\vec{x}^{n}\\
e^{n}
\end{bmatrix}
+
\frac{\Delta t}{2}
\begin{bmatrix}
\frac{\vec{f}^n}{M_K}\\
\vec{u}^{n+\frac{1}{2}}\\
\frac{\vec{f}^n\cdot \vec{\tilde{u}}^{n+\frac{1}{2}} }{M_T}
\end{bmatrix}.
\end{equation}
Here, $M_K$ and $M_T$ denote the diagonal mass matrices for the kinematic and thermodynamic variables, respectively. The division is done at the corresponding kinematic and thermodynamic DOFs level, eliminating the need to solve the full algebraic linear system. With the updated position $\vec{x}^{n+\frac{1}{2}}$ and internal energy $e^{n+\frac{1}{2}}$, the equation of state (EOS) is evaluated to obtain the intermediate pressure $p^{n+\frac{1}{2}}$, which subsequently yields the intermediate force $\vec{f}^{n+\frac{1}{2}}$. The second stage then proceeds as follows:
\begin{equation}
\begin{bmatrix}
\vec{u}^{n+1}\\
\vec{x}^{n+1}\\
e^{n+1}
\end{bmatrix}
=
\begin{bmatrix}
\vec{u}^{n}\\
\vec{x}^{n}\\
e^{n}
\end{bmatrix}
+
\Delta t
\begin{bmatrix}
\frac{\vec{f}^{n+\frac{1}{2}}}{M_K}\\
\frac{\vec{u}^{n}+\vec{u}^{n+1}}{2}\\
\frac{\vec{f}^{n+\frac{1}{2}}\cdot \frac{\vec{u}^{n}+\vec{u}^{n+1}}{2} }{M_T}
\end{bmatrix}.
\end{equation}

More generally, the unified $s$-stage IMEX-RK scheme advancing from $t^{n}$ to $t^{n+1}$ can be formulated compactly, mapping intermediate stages via implicit weights $b_i$ and explicit coefficient matrices $a_{i,j}^{\{E\}}$. 

As the explicit stage-by-stage formulas and Butcher tableaux for the higher-order IMEX-RK3 and IMEX-RK4 schemes are quite extensive, they are omitted here for brevity. The complete formulations and step-by-step algorithmic implementations for these higher-order temporal updates are detailed in Appendix~\ref{app:imex_rk}.

\subsection{Accuracy and convergence}
The Taylor-Green vortex and the Kidder problem are employed to evaluate the formal convergence rates of the proposed method. Because both problems possess smooth analytical solutions, we can rigorously assess the spatial accuracy of the finite element spaces without the dissipative influence of artificial viscosity. The results demonstrate that the proposed method achieves optimal convergence rates for both high-order element pairs.

\subsubsection{2D Taylor-Green vortex} 
The 2D Taylor-Green vortex is a smooth, shock-free problem. The computational domain is the unit square $[0,1]^2$, governed by no-penetration wall boundary conditions $\vec{u}\cdot \vec{n} =0$. Although the problem describes an incompressible flow, it is solved here using the fully compressible Euler equations. The domain is filled with an ideal gas characterized by a constant adiabatic index $\gamma=\frac{5}{3}$. The simulation is advanced to a terminal time of $t=0.75$ with hourglass control enabled and the scaling parameter set to $1$. The initial conditions are defined as follows:
\begin{equation}
    \begin{split}
        \vec{u}_{(t=0)}&=
        \begin{bmatrix}
            \sin(\pi x)\cos(\pi y)\\
            -\cos(\pi x)\sin(\pi y)
        \end{bmatrix},\\
        p_{(t=0)}&=\frac{\rho}{4}(\cos(2\pi x)+\cos(2\pi y))+1.
    \end{split}
\end{equation}  
Because this flow is physically incompressible, the density field remains strictly constant in both space and time; we therefore set $\rho \equiv 1$. Consequently, an external energy source term must be actively added to the energy equation to maintain thermodynamic consistency:
\begin{equation}
    e_{\text{source}}=\frac{3}{8}\pi \left(\cos(3\pi x)\cos(\pi y)-\cos(\pi x)\cos(3\pi y)\right).
\end{equation}

\begin{figure}[htbp]
  \centering
  \begin{subfigure}[b]{0.49\textwidth}
    \includegraphics[width=\textwidth]{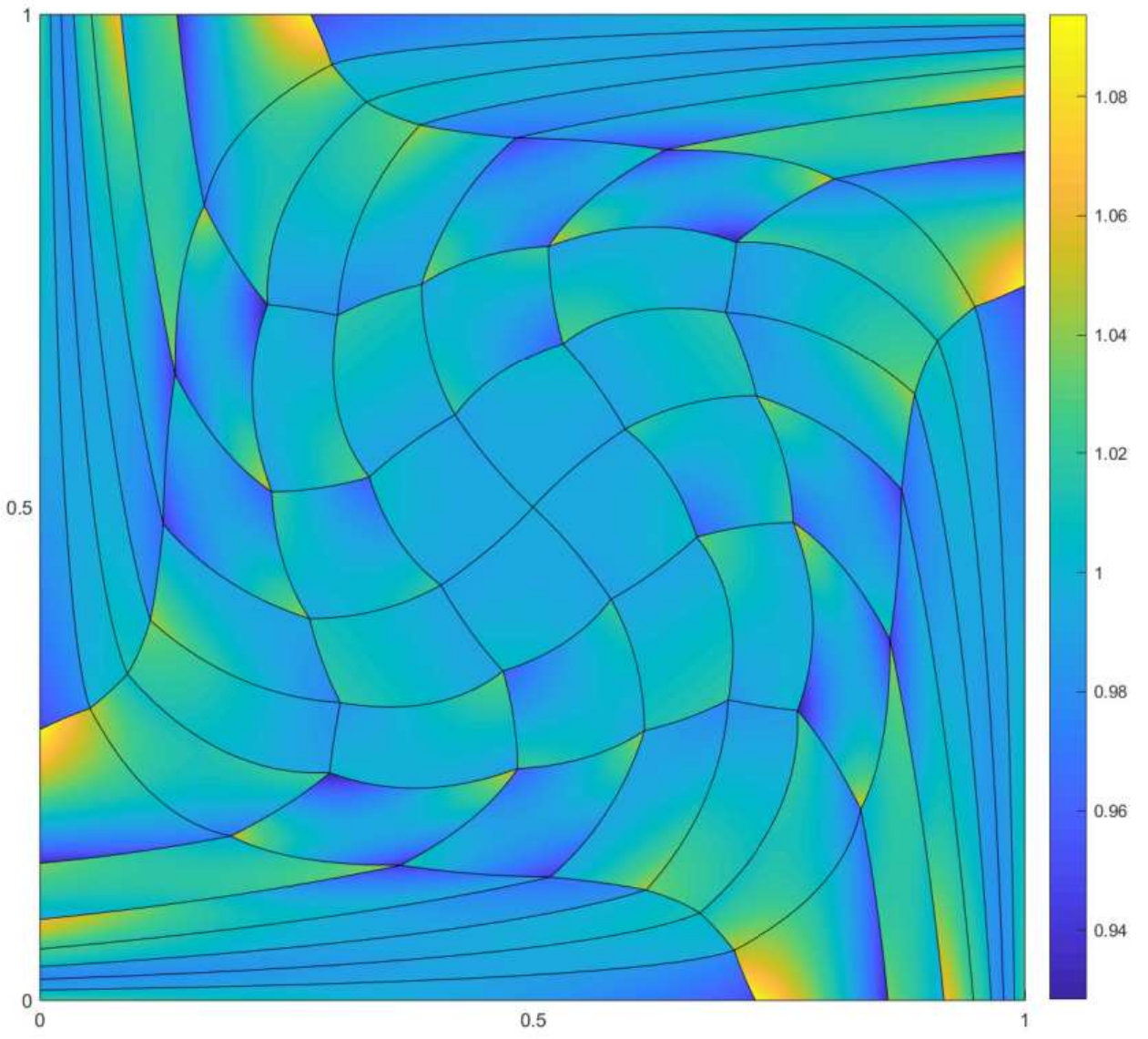}
    \caption{Density field ($Q^2-Q^1$).}
    \label{fig:Taylor_Green_vortex_field_a}
  \end{subfigure}
  \hfill
  \begin{subfigure}[b]{0.49\textwidth}
    \includegraphics[width=\textwidth]{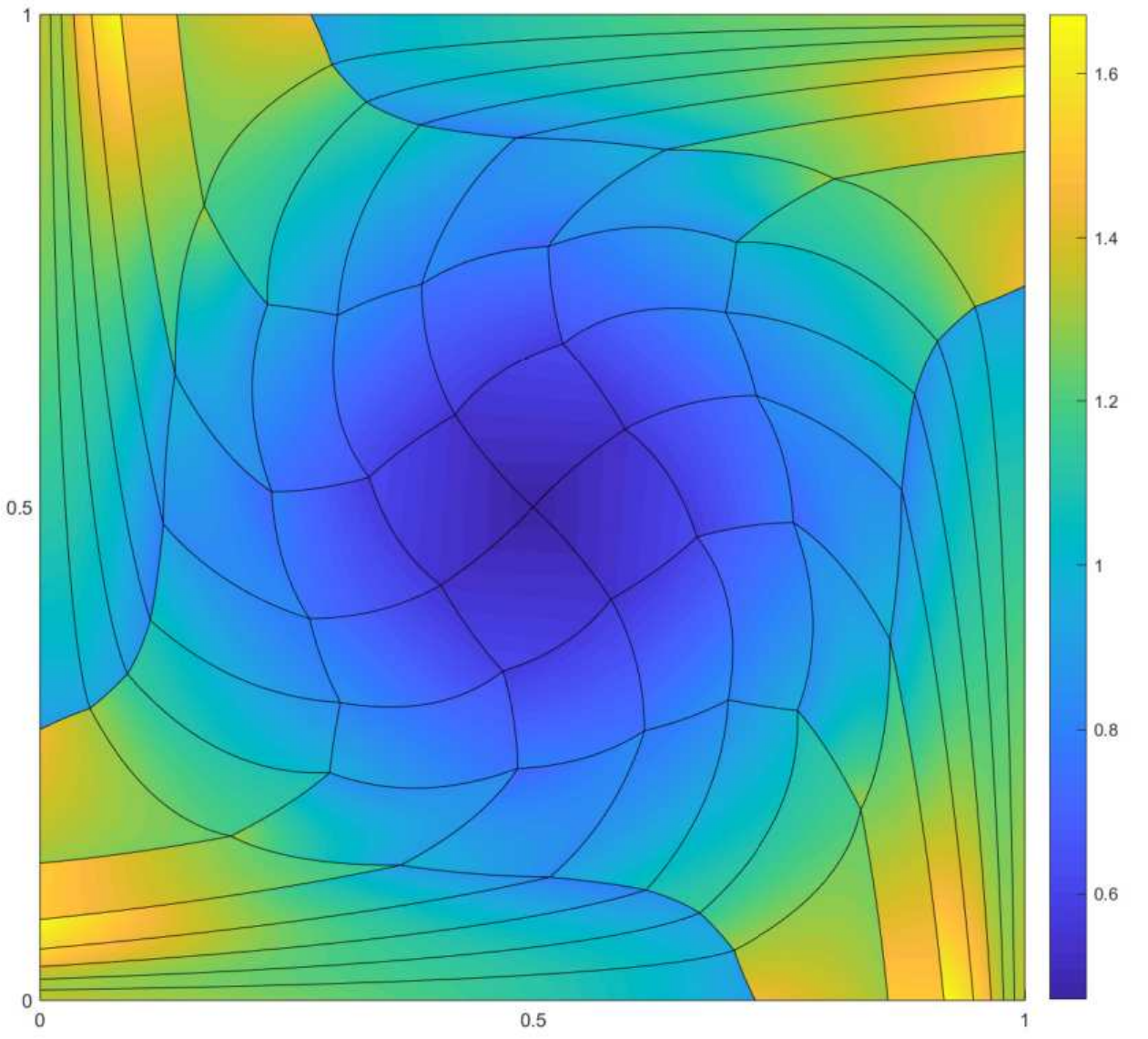}
    \caption{Pressure field ($Q^2-Q^1$).}
    \label{fig:Taylor_Green_vortex_field_b}
  \end{subfigure}
  
  \vspace{0.1em}
  
  \begin{subfigure}[b]{0.49\textwidth}
    \includegraphics[width=\textwidth]{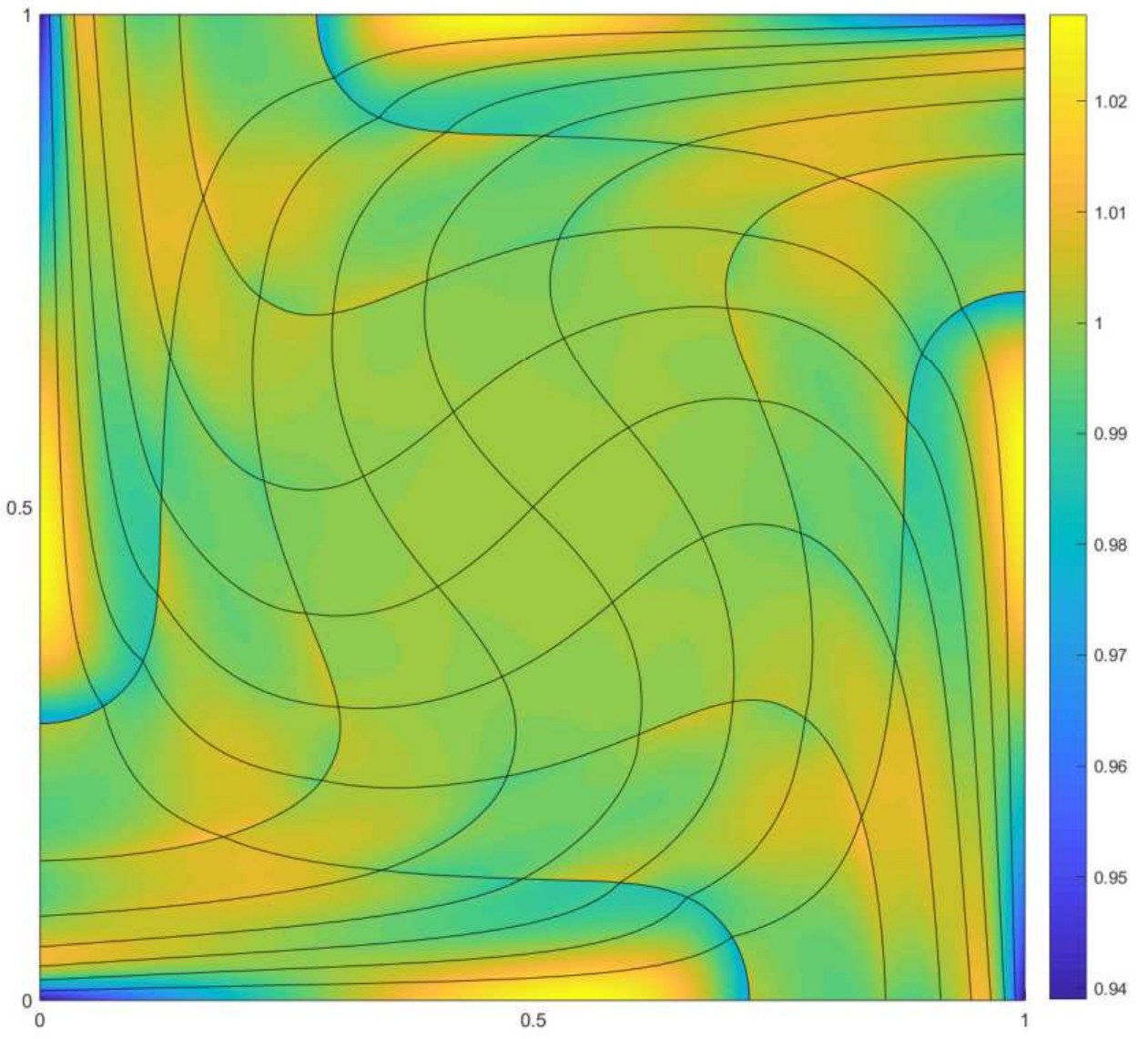}
    \caption{Density field ($Q^3-Q^2$).}
    \label{fig:Taylor_Green_vortex_field_c}
  \end{subfigure}
  \hfill
  \begin{subfigure}[b]{0.49\textwidth}
    \includegraphics[width=\textwidth]{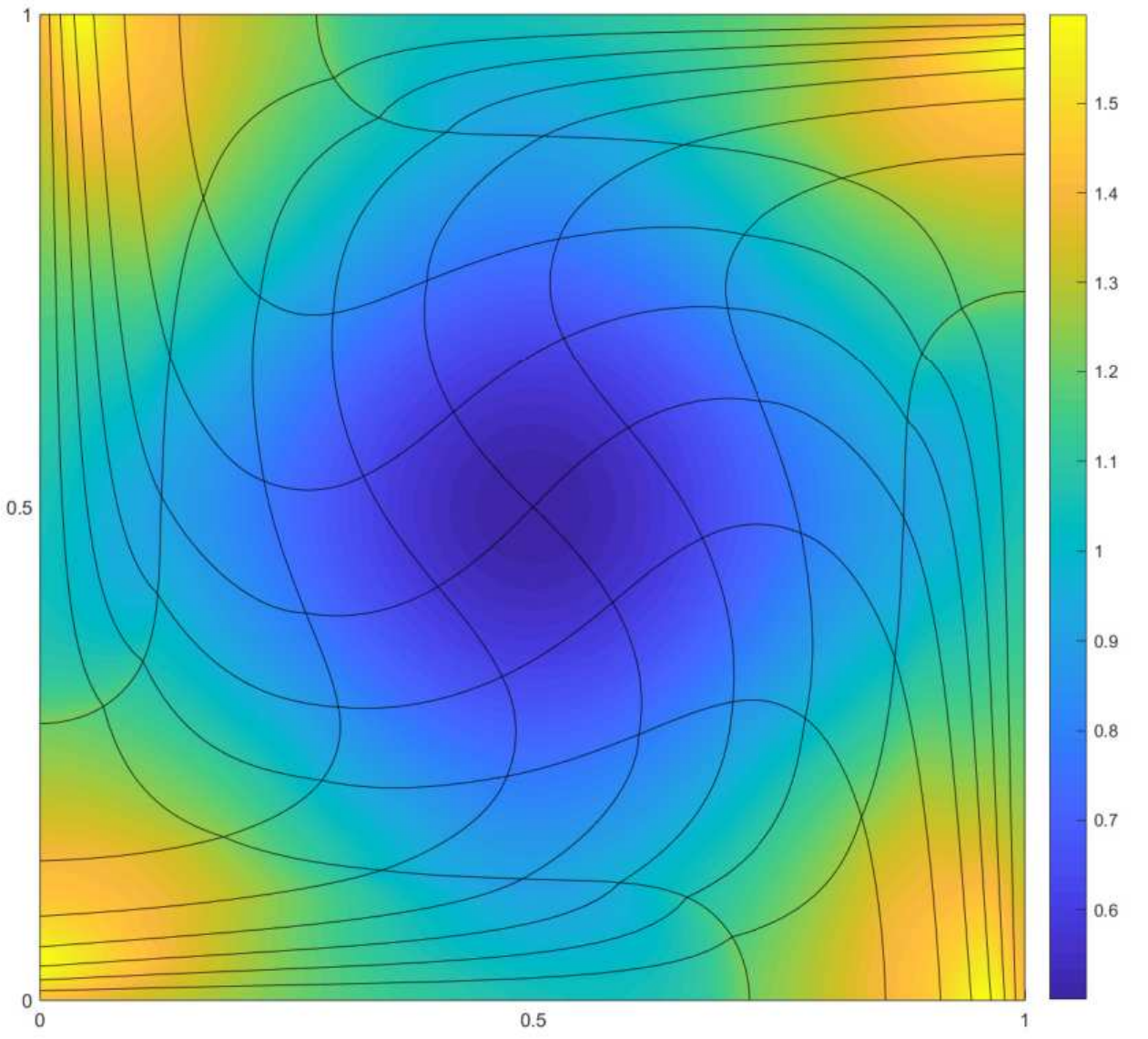}
    \caption{Pressure field ($Q^3-Q^2$).}
    \label{fig:Taylor_Green_vortex_field_d}
  \end{subfigure}
  \caption{Density and pressure fields for the Taylor-Green vortex problem on a coarse mesh ($h = 1/8$) at $t = 0.75$.}
  \label{Taylor_Green_vortex_field}
\end{figure}

\begin{figure}[htbp]
  \centering
  \begin{subfigure}[b]{0.49\textwidth}
    \includegraphics[width=\textwidth]{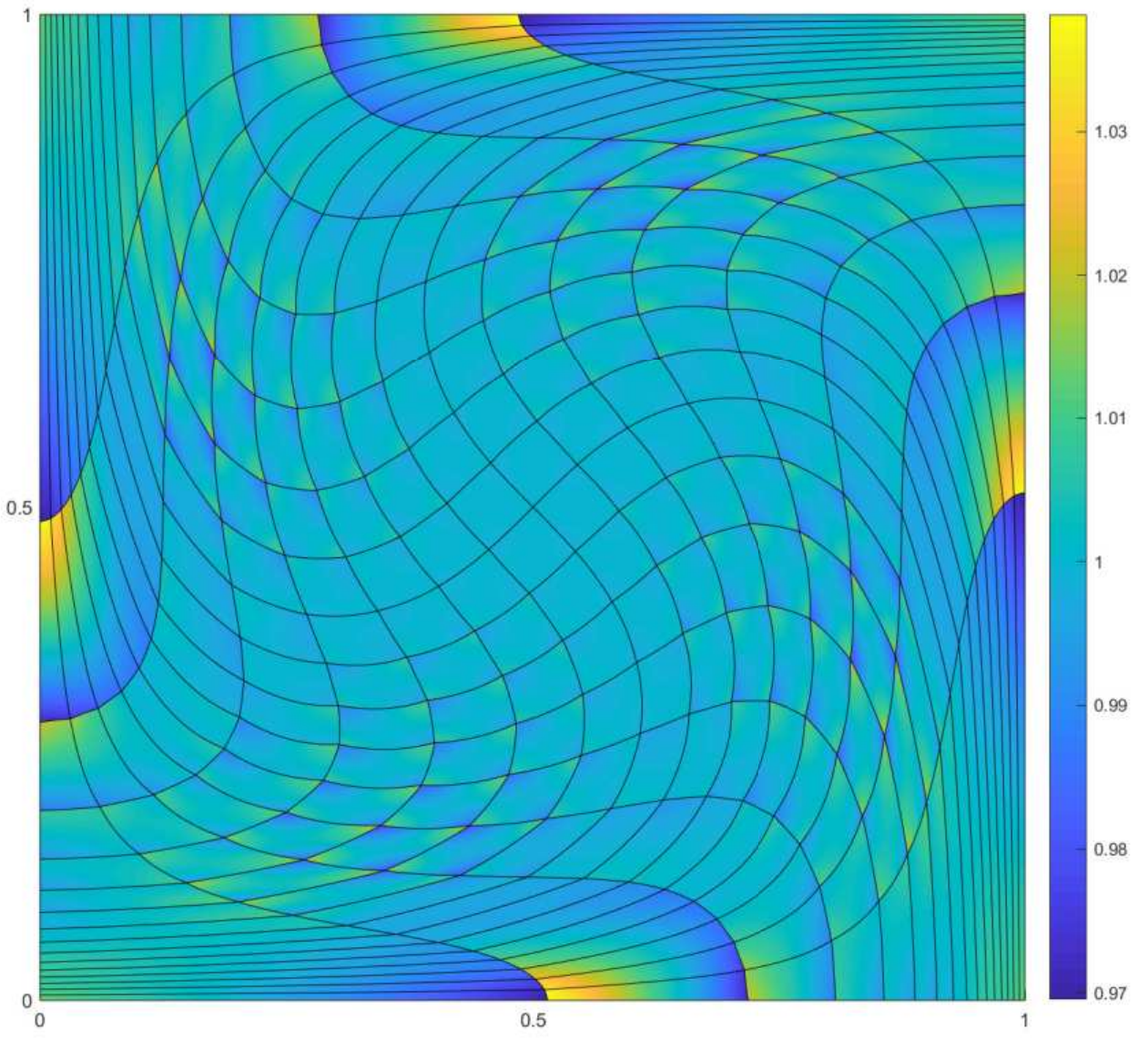}
    \caption{Density field ($Q^2-Q^1$).}
    \label{fig:Taylor_Green_vortex_field_fine_a}
  \end{subfigure}
  \hfill
  \begin{subfigure}[b]{0.49\textwidth}
    \includegraphics[width=\textwidth]{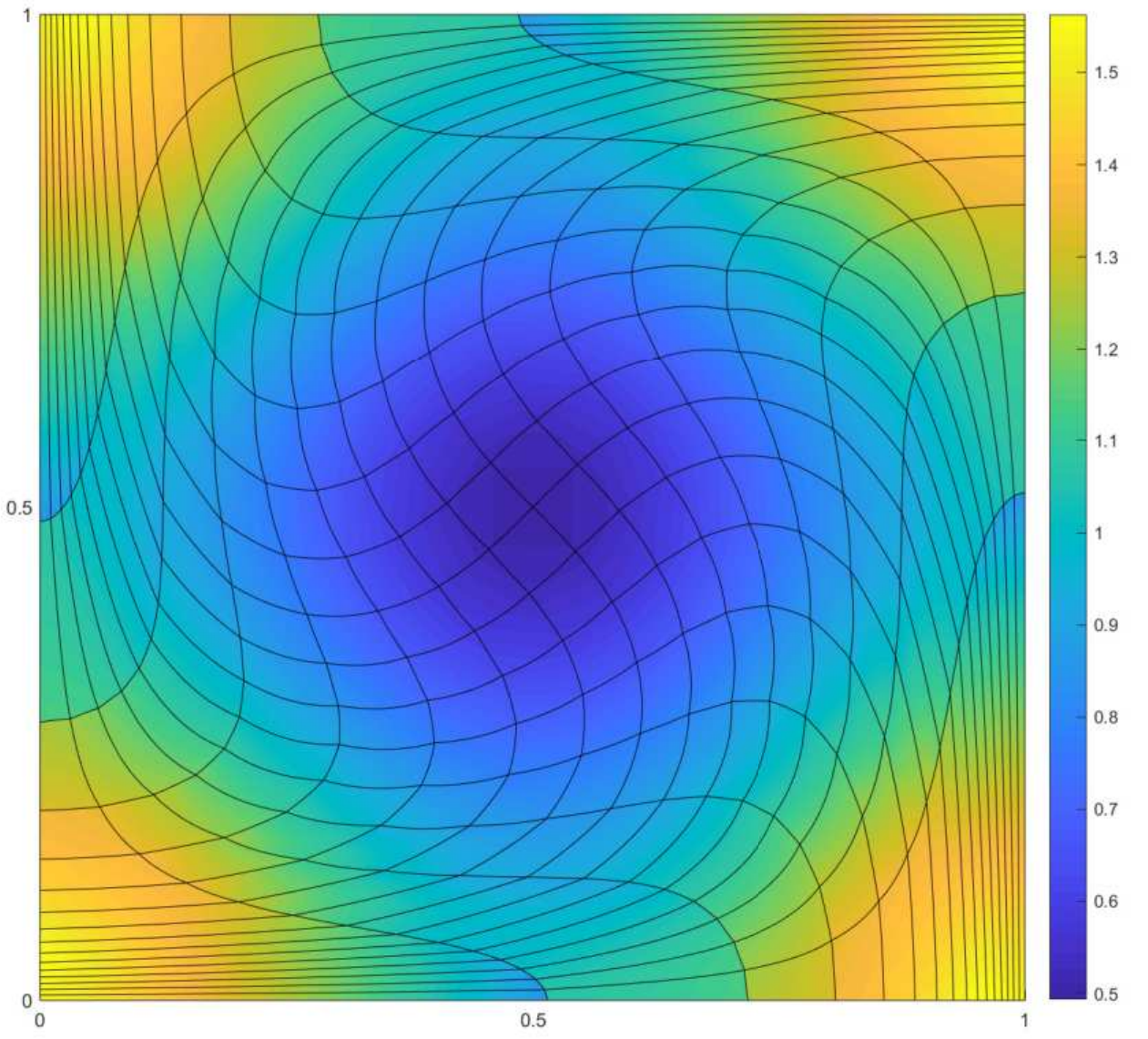}
    \caption{Pressure field ($Q^2-Q^1$).}
    \label{fig:Taylor_Green_vortex_field_fine_b}
  \end{subfigure}
  
  \vspace{0.5em}
  
  \begin{subfigure}[b]{0.49\textwidth}
    \includegraphics[width=\textwidth]{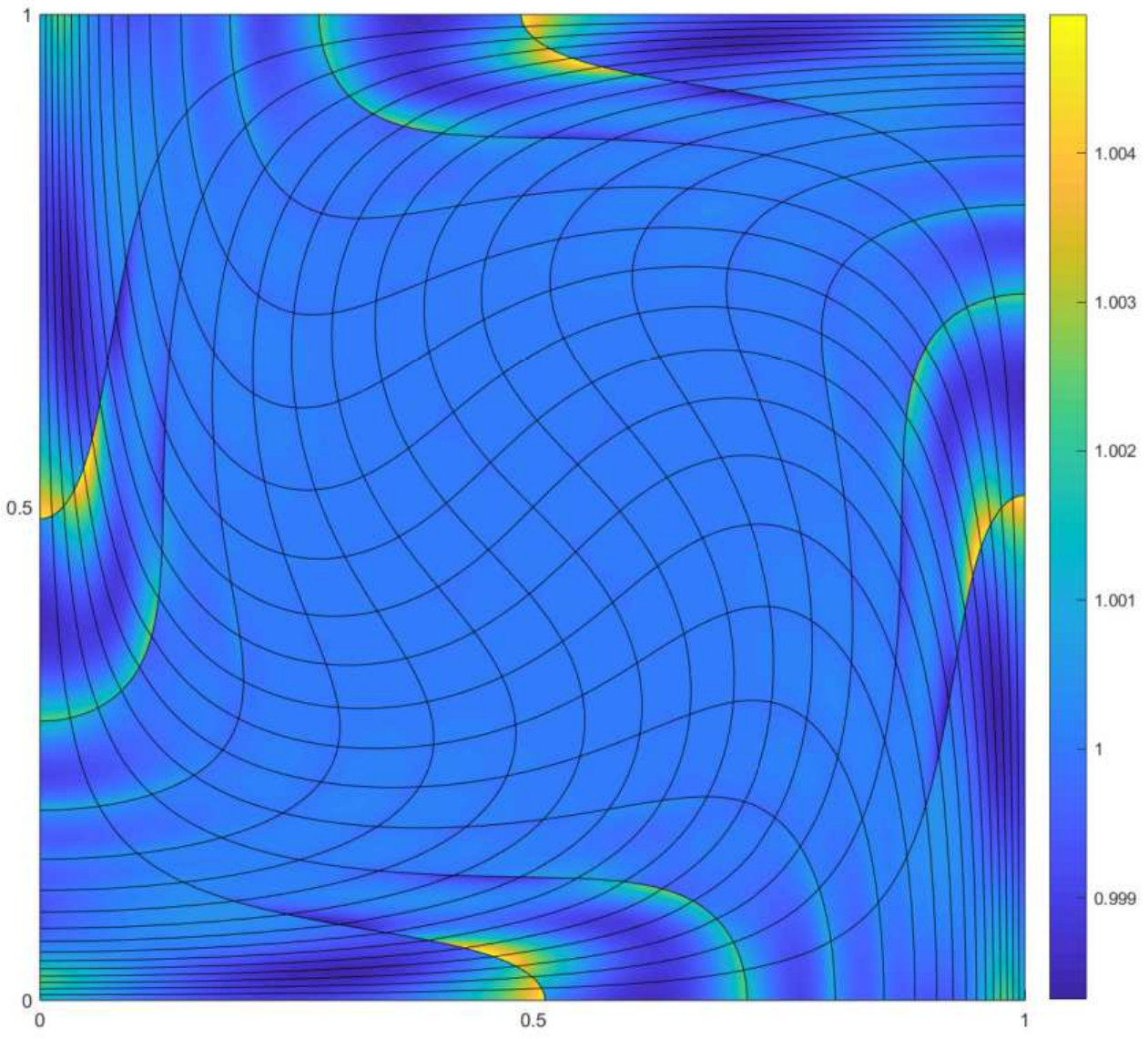}
    \caption{Density field ($Q^3-Q^2$).}
    \label{fig:Taylor_Green_vortex_field_fine_c}
  \end{subfigure}
  \hfill
  \begin{subfigure}[b]{0.49\textwidth}
    \includegraphics[width=\textwidth]{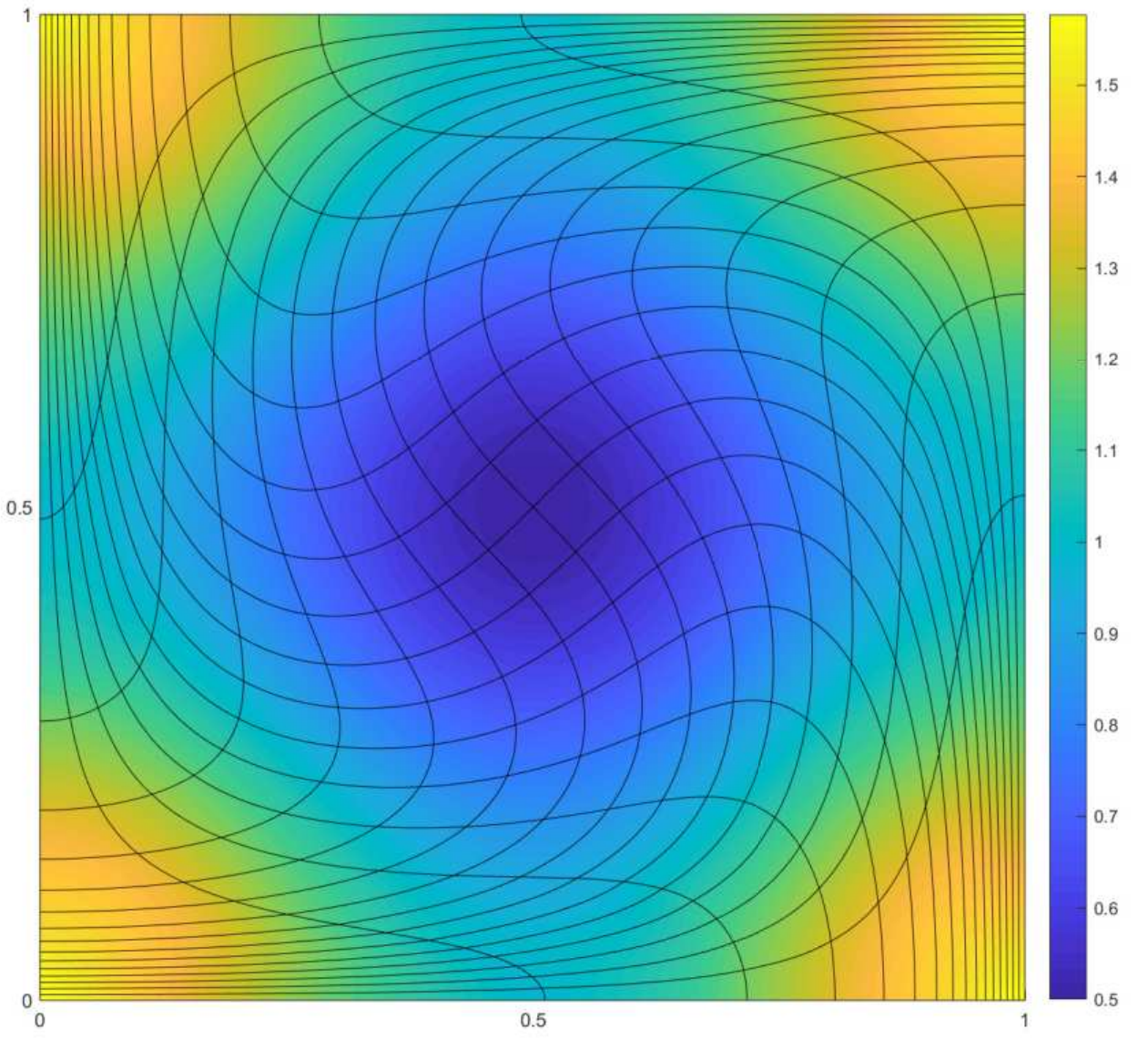}
    \caption{Pressure field ($Q^3-Q^2$).}
    \label{fig:Taylor_Green_vortex_field_fine_d}
  \end{subfigure}
  \caption{Density and pressure fields for the Taylor-Green vortex problem on a fine mesh ($h = 1/16$) at $t = 0.75$.}
  \label{Taylor_Green_vortex_field_fine}
\end{figure}

\begin{figure}[htbp]
  \centering
  \begin{subfigure}[b]{0.32\textwidth}
    \includegraphics[width=\textwidth]{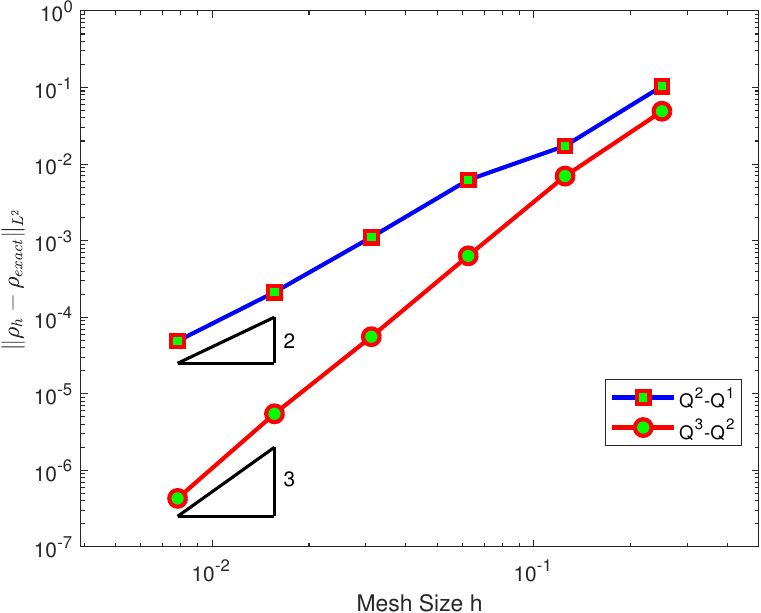}
    \caption{Density convergence rate.}
    \label{Taylor_Green_vortex_density_error}
  \end{subfigure}
  \hfill
  \begin{subfigure}[b]{0.32\textwidth}
    \includegraphics[width=\textwidth]{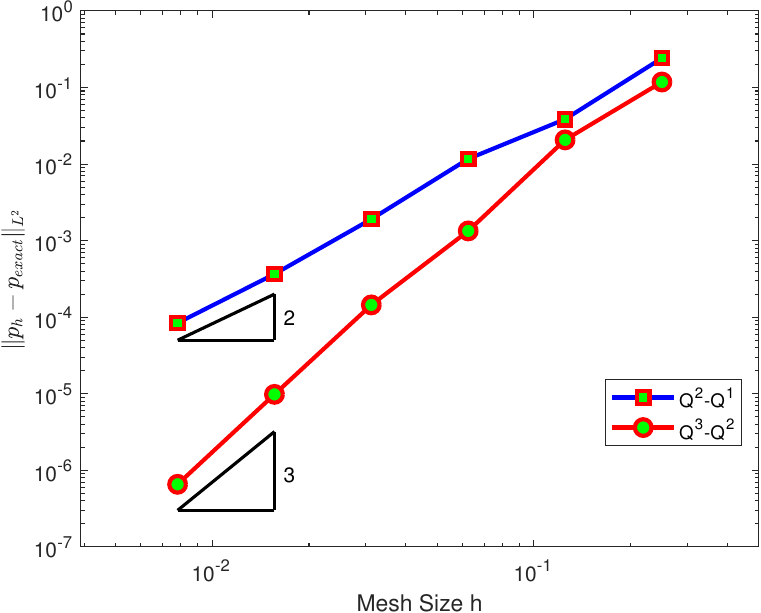}
    \caption{Pressure convergence rate.}
    \label{Taylor_Green_vortex_pressure_error}
  \end{subfigure}
  \hfill
  \begin{subfigure}[b]{0.32\textwidth}
    \includegraphics[width=\textwidth]{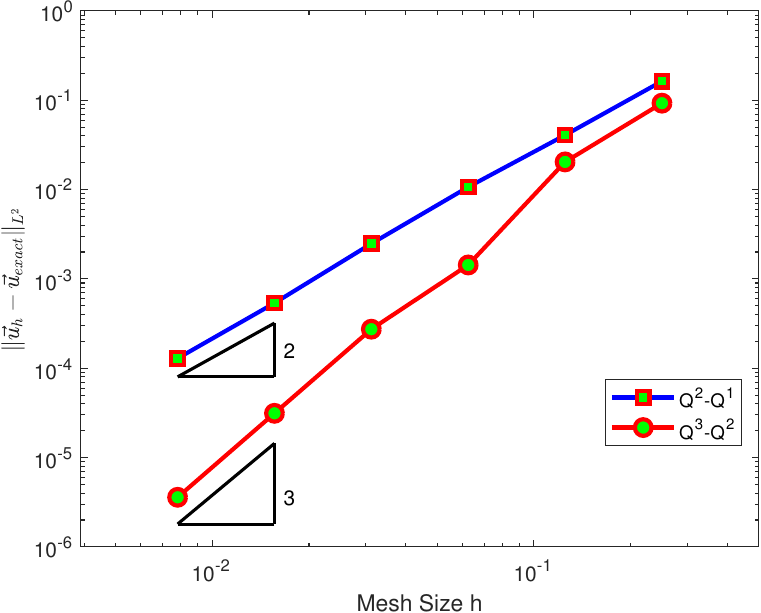}
    \caption{Velocity convergence rate.}
    \label{Taylor_Green_vortex_velocity_error}
  \end{subfigure}

  \vspace{0.5em}

  \begin{subfigure}[b]{0.32\textwidth}
    \includegraphics[width=\textwidth]{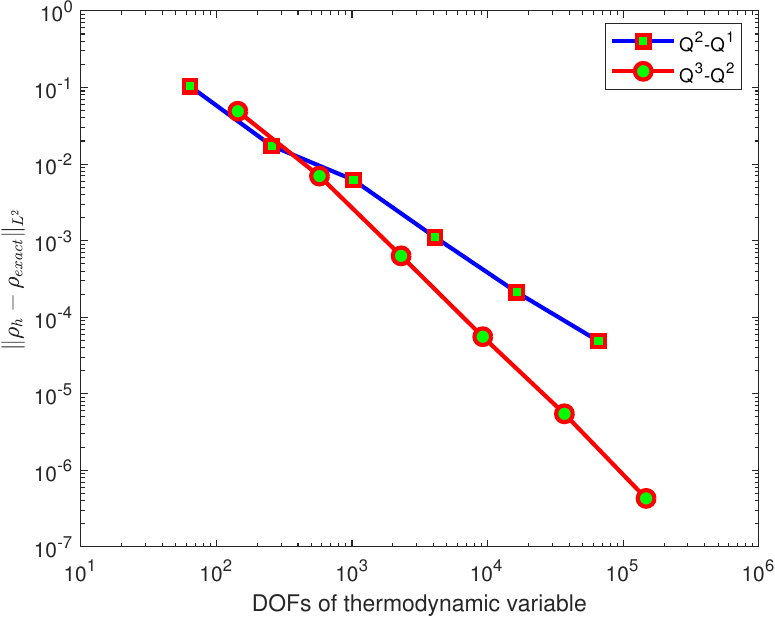}
    \caption{Density error vs. thermodynamic DOFs.}
    \label{Taylor_Green_vortex_error_rho_DOFs}
  \end{subfigure}
  \hfill
  \begin{subfigure}[b]{0.32\textwidth}
    \includegraphics[width=\textwidth]{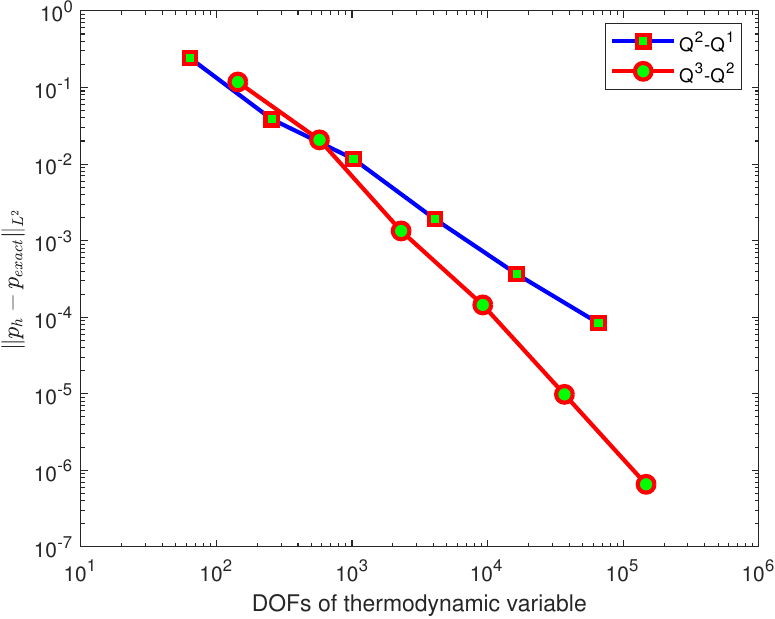}
    \caption{Pressure error vs. thermodynamic DOFs.}
    \label{Taylor_Green_vortex_error_pre_DOFs}
  \end{subfigure}
  \hfill
  \begin{subfigure}[b]{0.32\textwidth}
    \includegraphics[width=\textwidth]{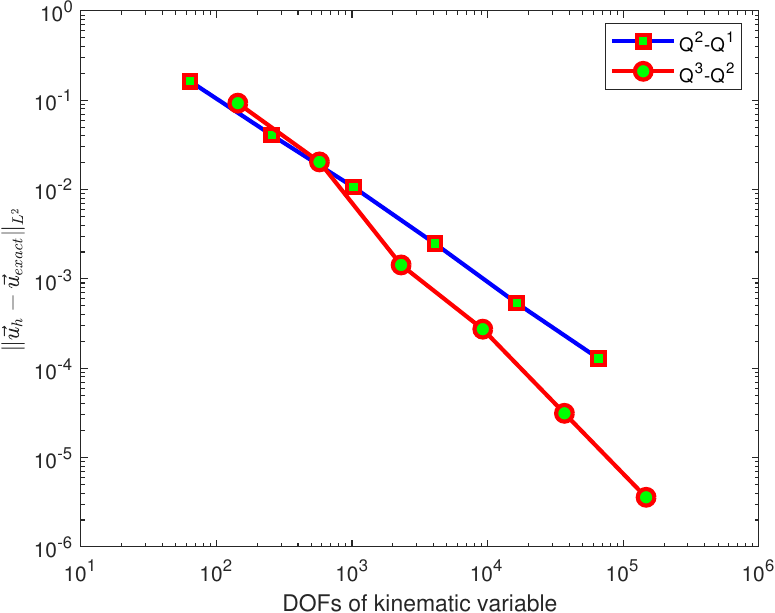}
    \caption{Velocity error vs. kinematic DOFs.}
    \label{Taylor_Green_vortex_error_vel_DOFs}
  \end{subfigure}

  \caption{Convergence rates of $L^2$ error norms (top row) and error norms versus DOFs (bottom row) for the Taylor-Green vortex problem at $t = 0.75$.}
  \label{Taylor_Green_vortex_error}
\end{figure}

Figures~\ref{Taylor_Green_vortex_field} and~\ref{Taylor_Green_vortex_field_fine} show the density and pressure fields at $t=0.75$ for the $Q^2-Q^1$ and $Q^3-Q^2$ element pairs, respectively. With the hourglass control active, unphysical grid distortion is effectively suppressed. Without this stabilization, severe mesh tangling inherently prevents the computation from reaching the terminal time. The $L^2$ error norms of the density, pressure, and velocity fields are detailed in Figure~\ref{Taylor_Green_vortex_error} and Tables~\ref{tab:Q_2_Q_1_error_t.75} and~\ref{tab:Q_3_Q_2_error_t.75}, closely matching theoretical predictions.

Figure~\ref{Taylor_Green_vortex_error} and Tables~\ref{tab:Q_2_Q_1_error_t.75} and~\ref{tab:Q_3_Q_2_error_t.75} also list the corresponding DOFs of the thermodynamic and kinematic variables for specific mesh sizes. This data is crucial for demonstrating the core advantage of high-order methods: achieving superior accuracy with significantly fewer DOFs. For instance, at $h = 1/64$, the $Q^2-Q^1$ space pair utilizes 16,641 kinematic and 16,384 thermodynamic DOFs, yet it yields higher errors than the $Q^3-Q^2$ pair at $h = 1/32$ (which uses only 9,409 and 9,216 DOFs, respectively). Furthermore, the $Q^3-Q^2$ results at $h = 1/32$ are comparable in accuracy to the $Q^2-Q^1$ results at $h = 1/128$, despite the latter requiring roughly seven times more DOFs. The bottom panels of Figure~\ref{Taylor_Green_vortex_error} visualize this advantage clearly, showing the error norm of the $Q^3-Q^2$ pair shifted diagonally downward relative to the $Q^2-Q^1$ pair.

\begin{table}[htbp]
\centering
\caption{Convergence rates of the $Q^2-Q^1$ element pair for the Taylor-Green vortex problem at $t = 0.75$.}
\label{tab:Q_2_Q_1_error_t.75}
\begin{tabular}{c c S[table-format=1.4e-1] c S[table-format=1.4e-1] c S[table-format=1.4e-1] c}
\toprule
{$h$} & {DOFs} & {$\rho$ Error} & {Order} & {$p$ Error} & {Order} & {$\vec{u}$ Error} & {Order} \\
\midrule
$\frac{1}{4}$  & (81,64)    & 1.0376E-1 & --     & 2.4207E-1 & --     & 1.6326E-1 & -- \\
$\frac{1}{8}$  & (289,256)  & 1.7286E-2 & 2.5855 & 3.8559E-2 & 2.6503 & 4.0727E-2 & 2.003 \\
$\frac{1}{16}$ & (1089,1024)& 6.2032E-3 & 1.4785 & 1.1675E-2 & 1.7236 & 1.0683E-2 & 1.9306 \\
$\frac{1}{32}$ & (4225,4096)& 1.1088E-3 & 2.4839 & 1.9083E-3 & 2.6130 & 2.4808E-3 & 2.1064 \\
$\frac{1}{64}$ & (16641,16384)&2.1195E-4& 2.3872 & 3.6603E-4 & 2.3822 & 5.3623E-4 & 2.2099 \\
$\frac{1}{128}$& (66049,65536)&4.8818E-5& 2.1182 & 8.3855E-5 & 2.1260 & 1.2891E-4 & 2.0564 \\
\bottomrule
\end{tabular}
\end{table}

\begin{table}[htbp]
\centering
\caption{Convergence rates of the $Q^3-Q^2$ element pair for the Taylor-Green vortex problem at $t = 0.75$.}
\label{tab:Q_3_Q_2_error_t.75}
\begin{tabular}{c c S[table-format=1.4e-1] c S[table-format=1.4e-1] c S[table-format=1.4e-1] c}
\toprule
{$h$} & {DOFs} & {$\rho$ Error} & {Order} & {$p$ Error} & {Order} & {$\vec{u}$ Error} & {Order} \\
\midrule
$\frac{1}{4}$  & (169,144)    & 4.8974E-2 & --     & 1.1796E-1 & --     & 9.2793E-2 & -- \\
$\frac{1}{8}$  & (625,576)    & 6.9662E-3 & 2.8135 & 2.0691E-2 & 2.5111 & 2.0403E-2 & 2.1852 \\
$\frac{1}{16}$ & (2401,2304)  & 6.3266E-4 & 3.4608 & 1.3380E-3 & 3.9508 & 1.4328E-3 & 3.8318 \\
$\frac{1}{32}$ & (9409,9216)  & 5.5592E-5 & 3.5084 & 1.4461E-4 & 3.2098 & 2.7356E-4 & 2.3889 \\
$\frac{1}{64}$ & (37249,36864)& 5.4709E-6 & 3.3450 & 9.8134E-6 & 3.8813 & 3.1284E-5 & 3.1283 \\
$\frac{1}{128}$& (148225,147456)&4.2979E-7& 3.6700 & 6.5620E-7 & 3.9025 & 3.5892E-6 & 3.1236 \\
\bottomrule
\end{tabular}
\end{table}

\subsubsection{2D Kidder Problem} 
The Kidder problem is a well-known isentropic compression benchmark featuring a self-similar analytical solution, originally proposed by Kidder in~\cite{Kidder1974Theory,Kidder1976Laser}. It models an ideal gas being isentropically compressed toward the origin. The initial computational domain is a concentric ring defined by $r_1 \le r \le r_2$. The initial density, pressure, and internal energy distributions are given by:
\begin{equation*}
\begin{split}
\rho(r,0) &= \left( \frac{\rho_2^{\gamma-1}(r^2-r_1^2)+\rho_1^{\gamma-1}(r_2^2-r^2)}{r_2^2-r_1^2} \right)^{\frac{1}{\gamma-1}},\\
p(r,0) &= \left( \frac{\rho_2^{\gamma-1}(r^2-r_1^2)+\rho_1^{\gamma-1}(r_2^2-r^2)}{r_2^2-r_1^2} \right)^{\frac{\gamma}{\gamma-1}},\\
e(r,0) &= \frac{1}{\gamma-1}  \frac{\rho_2^{\gamma-1}(r^2-r_1^2)+\rho_1^{\gamma-1}(r_2^2-r^2)}{r_2^2-r_1^2},  
\end{split}
\end{equation*}
where $\rho_1=1$ and $\rho_2=2$ are constants. The ideal gas equation of state is $p = (\gamma-1)\rho e$, with $\gamma = 2$ for this two-dimensional case. The initial velocity field is governed by $u(r,0) = r \dot{h}(t)$, where $h(t)^2 = 1-\frac{t^2}{\tau^2}$, and $\tau$ is the characteristic time scale defined as $\tau=\sqrt{\frac{r_2^2-r_1^2}{2(\gamma-1)(\rho_2^{\gamma-1}-\rho_1^{\gamma-1})}}$.

\begin{figure}[htbp]
  \centering
  \begin{subfigure}[b]{0.485\textwidth}
    \includegraphics[width=\textwidth]{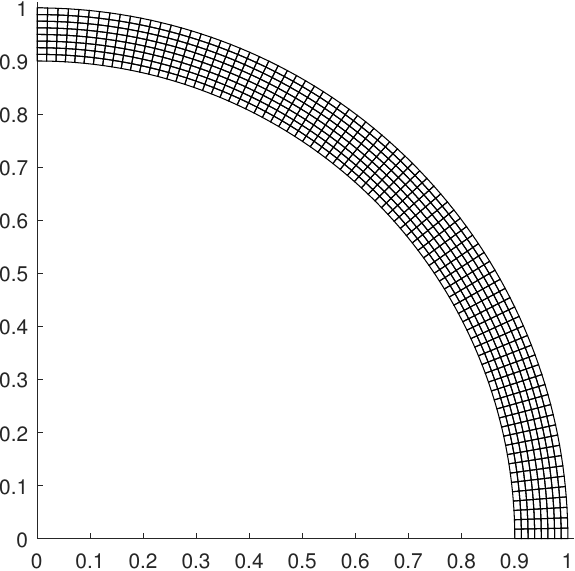}
    \caption{Initial computational mesh.}
    \label{kidder_mesh}
  \end{subfigure}
  \hfill
  \begin{subfigure}[b]{0.485\textwidth}
    \includegraphics[width=\textwidth]{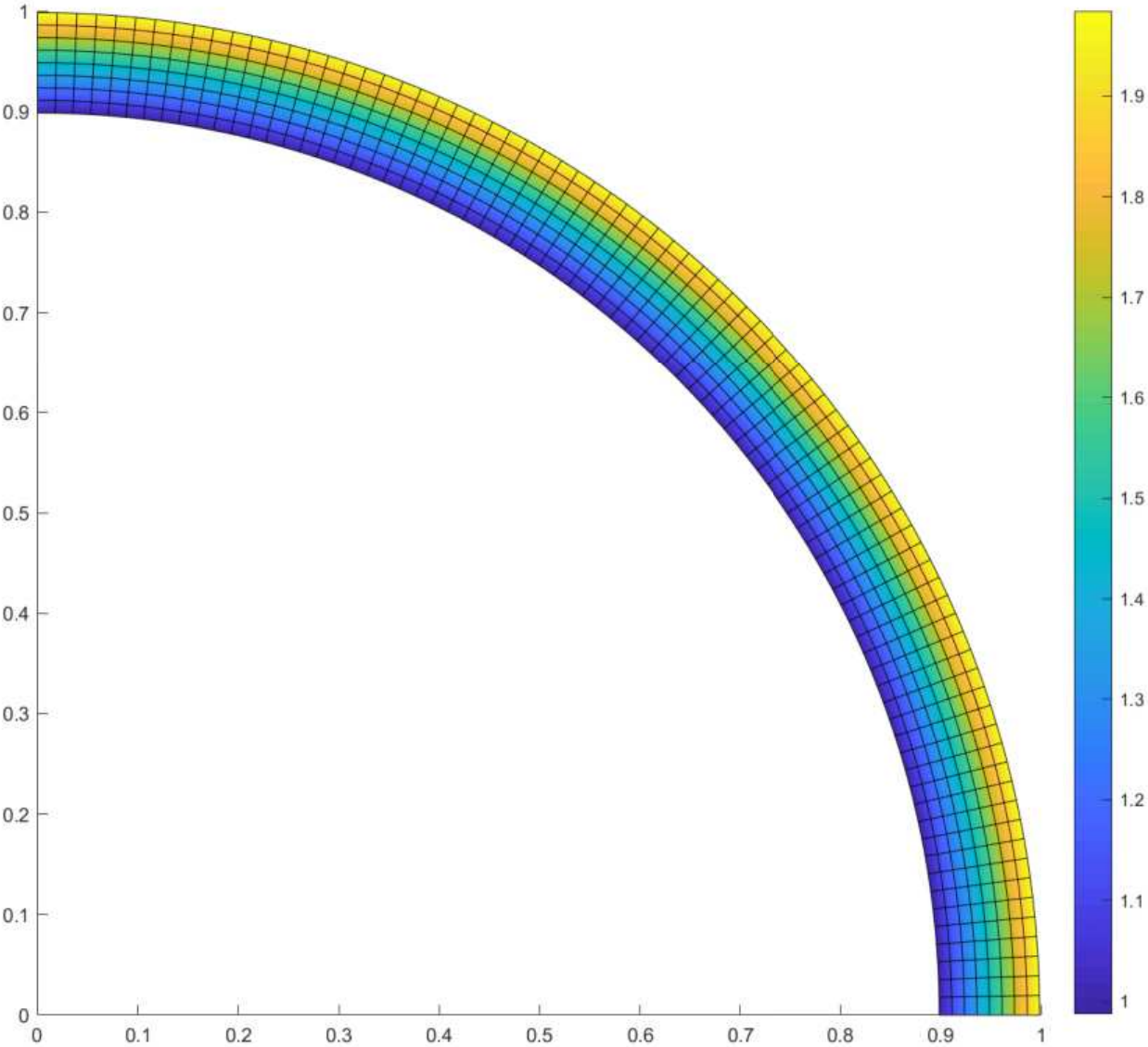}
    \caption{Initial density distribution.}
    \label{kidder_rho0}
  \end{subfigure}
  \caption{Initial computational mesh and density distribution for the Kidder problem.}
\label{Kidder_initial}
\end{figure}

For this test, we set $r_1 = 0.9$ and $r_2 = 1$. The initial mesh is generated by uniformly dividing the domain into quadrilateral elements with mesh sizes $h \in \{\frac{1}{20}, \frac{1}{40}, \frac{1}{80}, \frac{1}{160}, \frac{1}{320}, \frac{1}{640}\}$. Figure~\ref{Kidder_initial} displays the initial mesh and density distribution. 

\begin{figure}[htbp]
  \centering
  \begin{subfigure}[b]{0.485\textwidth}
    \includegraphics[width=\textwidth]{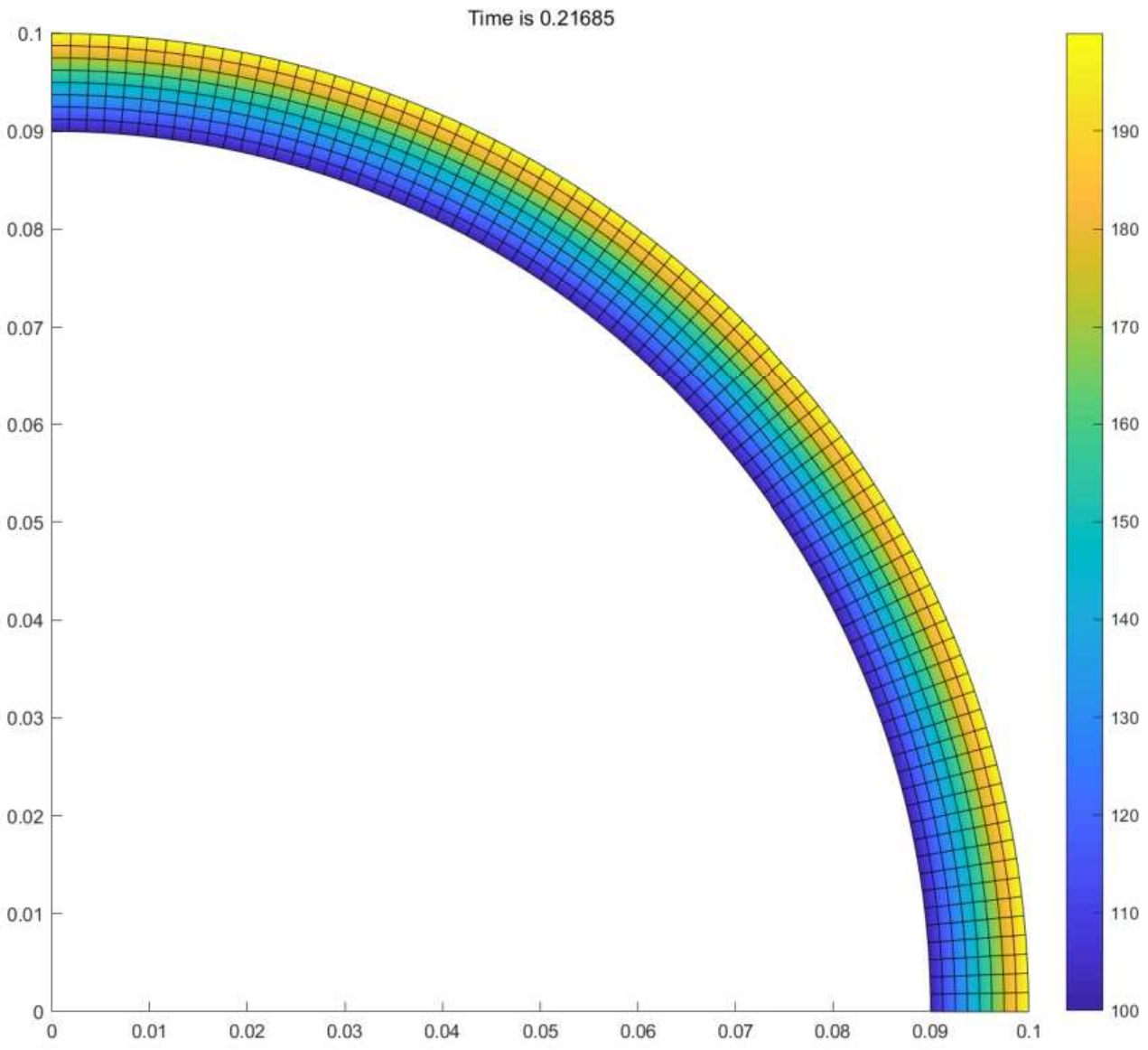}
    \caption{$Q^2-Q^1$ pair at $t = \frac{\sqrt{99}}{10}\tau$.}
    \label{kidder_Q2_0.1}
  \end{subfigure}
  \hfill
  \begin{subfigure}[b]{0.485\textwidth}
    \includegraphics[width=\textwidth]{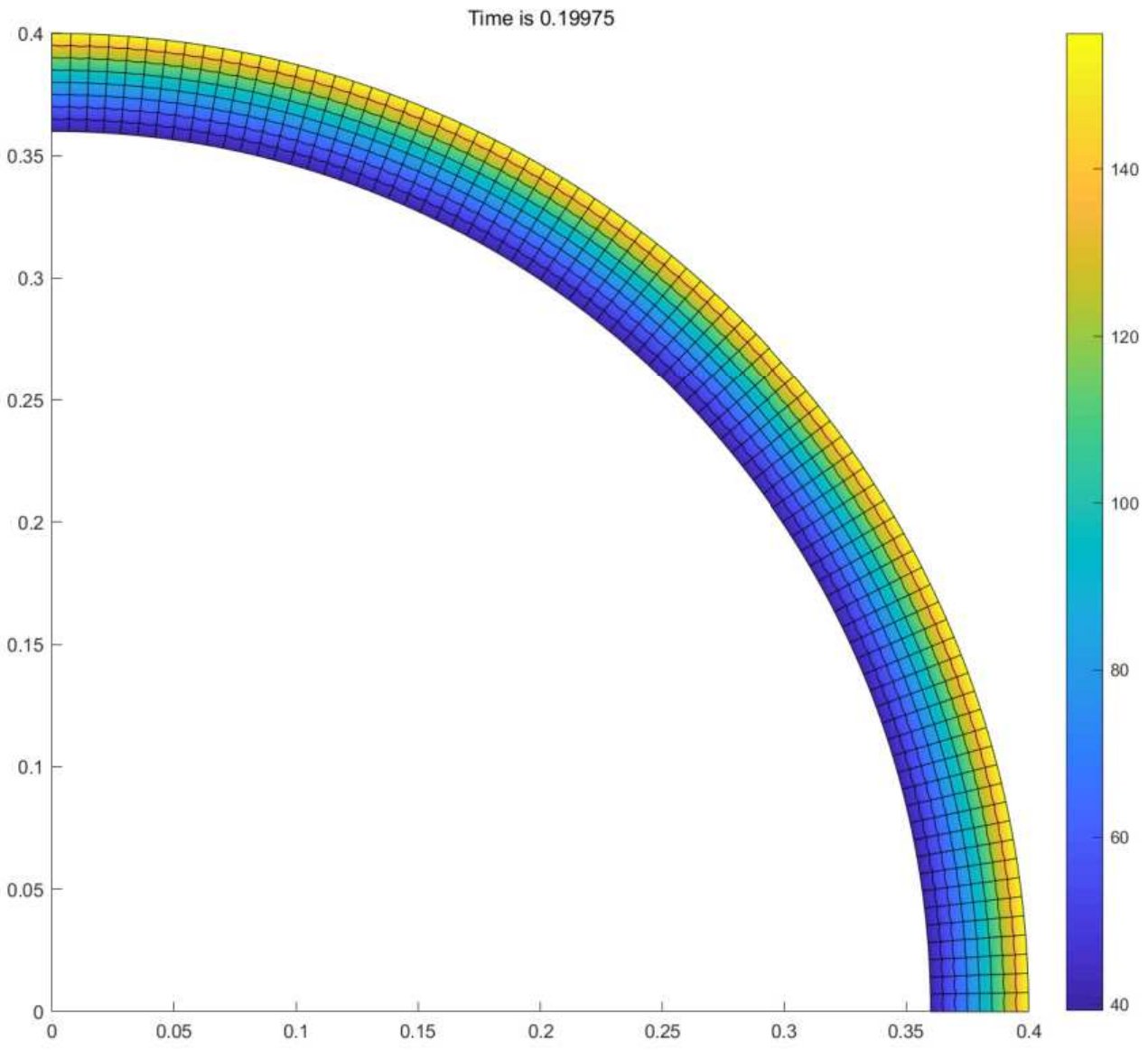}
    \caption{$Q^3-Q^2$ pair at $t = \frac{\sqrt{5.25}}{2.5}\tau$.}
    \label{kidder_Q3_0.4}
  \end{subfigure}
  \caption{Density contours for the Kidder problem without hourglass control, illustrating varying sensitivities to zero-energy modes.}
\label{Kidder_no_hg}
\end{figure}

Beyond initial convergence, the Kidder problem also serves as an excellent stress test for mesh robustness under extreme deformation. Before applying our hourglass control algorithm, we examine the baseline behavior of the elements. Figure~\ref{Kidder_no_hg} compares the numerical results for the $Q^2-Q^1$ and $Q^3-Q^2$ space pairs without the hourglass control algorithm. The macroscopic behaviors of the two cases are strikingly different. For the $Q^2-Q^1$ case, gross hourglass distortion is not immediately triggered in this highly symmetric configuration, allowing the simulation to proceed even when the ring is compressed to $\frac{1}{10}$ of its initial radius. Meanwhile, the $Q^3-Q^2$ case is extraordinarily susceptible to zero-energy modes. The simulation crashes when the ring is compressed to only about $\frac{1}{3}$ of its initial radius, as severe grid distortion leads to inverted elements and negative Jacobian determinants. Consequently, all subsequent numerical results for the Kidder problem are computed with the hourglass control explicitly activated, using a stabilization scaling parameter of $1$. 

\begin{figure}[htbp]
  \centering
  \begin{subfigure}[b]{0.32\textwidth}
    \includegraphics[width=\textwidth]{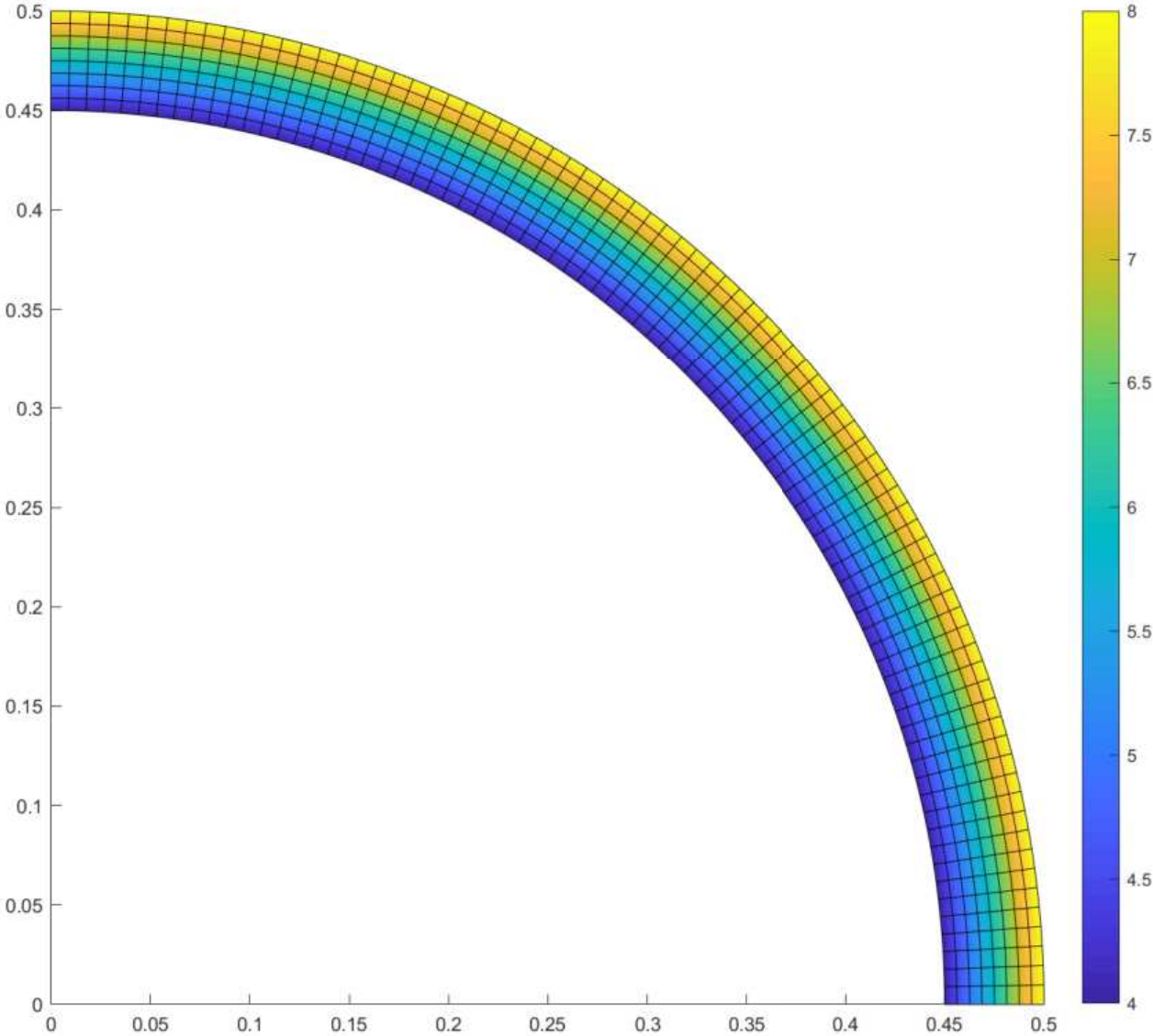}
    \caption{$t = \frac{\sqrt{3}}{2}\tau$ ($r/r_0 = 1/2$).}
    \label{kidder_half}
  \end{subfigure}
  \hfill
  \begin{subfigure}[b]{0.32\textwidth}
    \includegraphics[width=\textwidth]{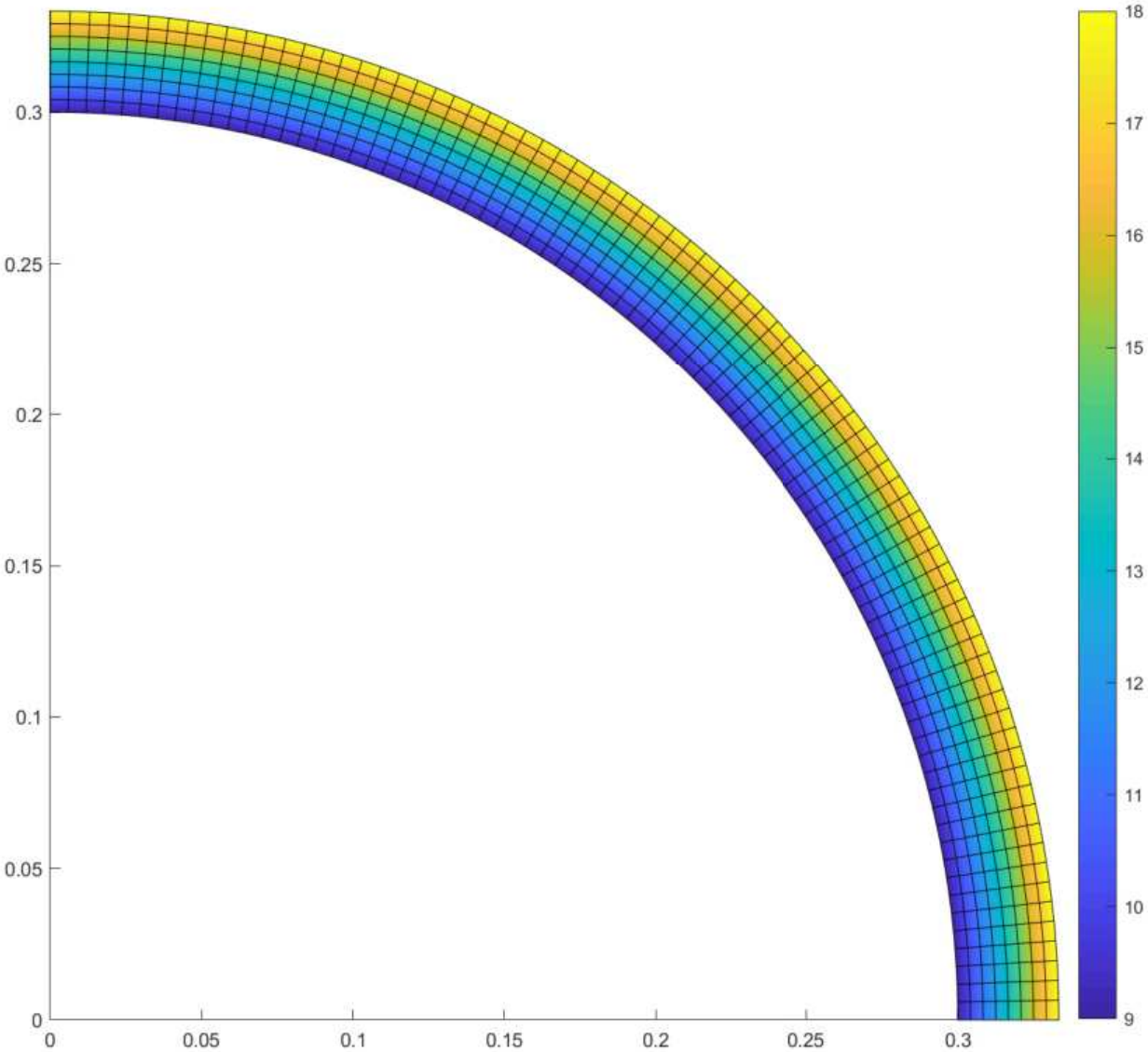}
    \caption{$t = \frac{\sqrt{8}}{3}\tau$ ($r/r_0 = 1/3$).}
    \label{kidder_third}
  \end{subfigure}
  \hfill
  \begin{subfigure}[b]{0.32\textwidth}
    \includegraphics[width=\textwidth]{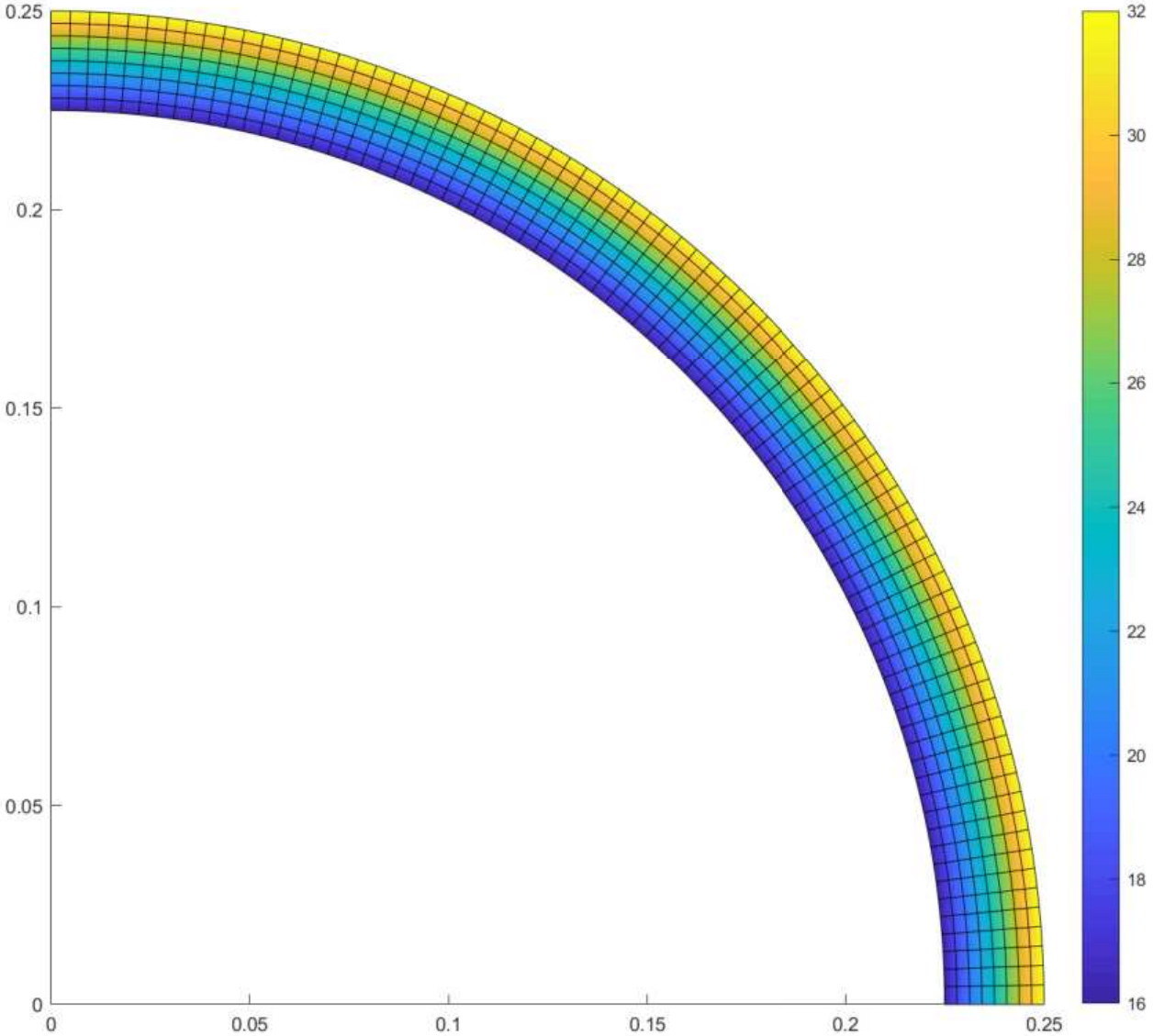}
    \caption{$t = \frac{\sqrt{15}}{4}\tau$ ($r/r_0 = 1/4$).}
    \label{kidder_tquarter}
  \end{subfigure}

  \vspace{0.5em}

  \begin{subfigure}[b]{0.32\textwidth}
    \includegraphics[width=\textwidth]{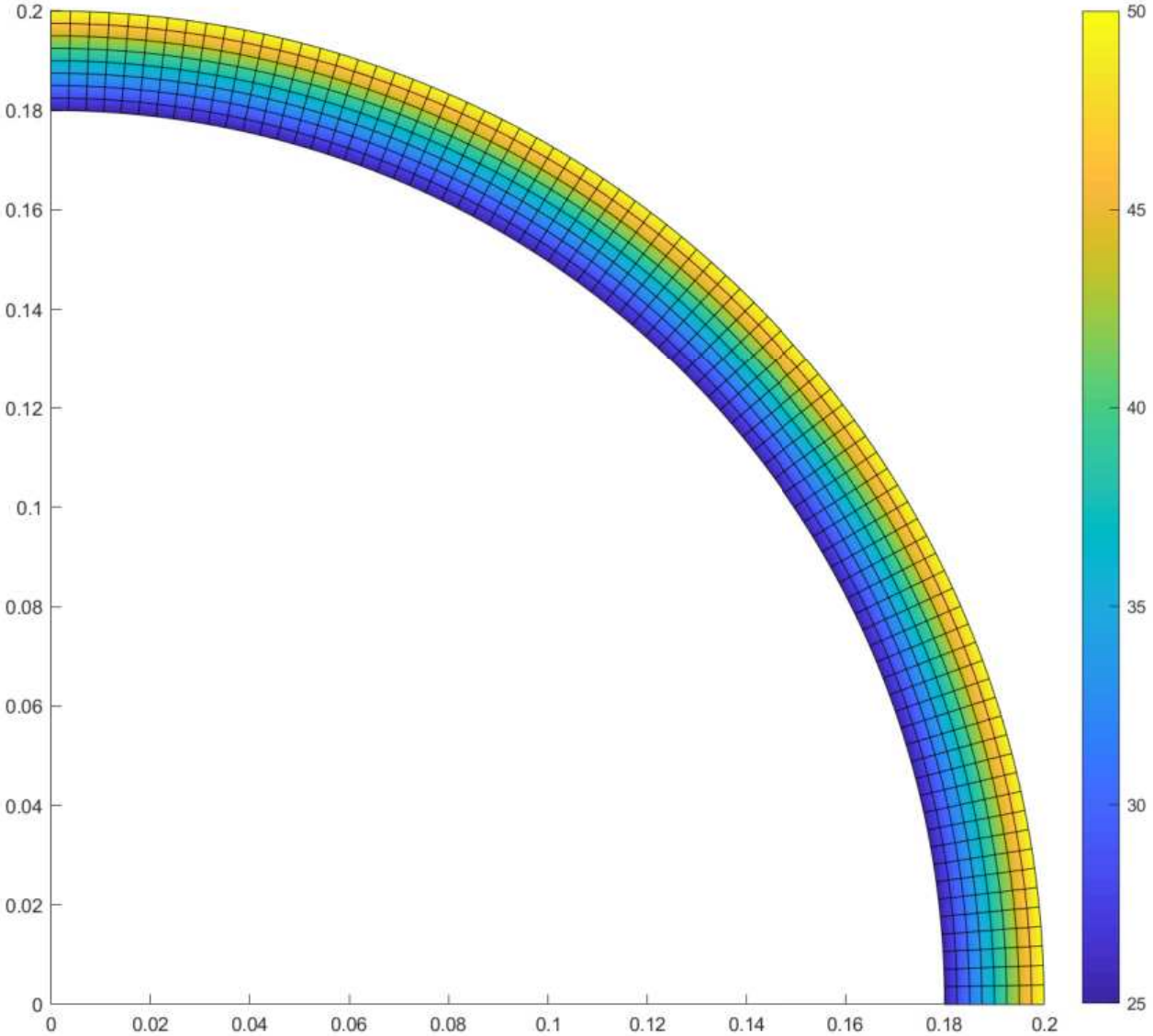}
    \caption{$t = \frac{\sqrt{24}}{5}\tau$ ($r/r_0 = 1/5$).}
    \label{kidder_fifth}
  \end{subfigure}
  \hfill
  \begin{subfigure}[b]{0.32\textwidth}
    \includegraphics[width=\textwidth]{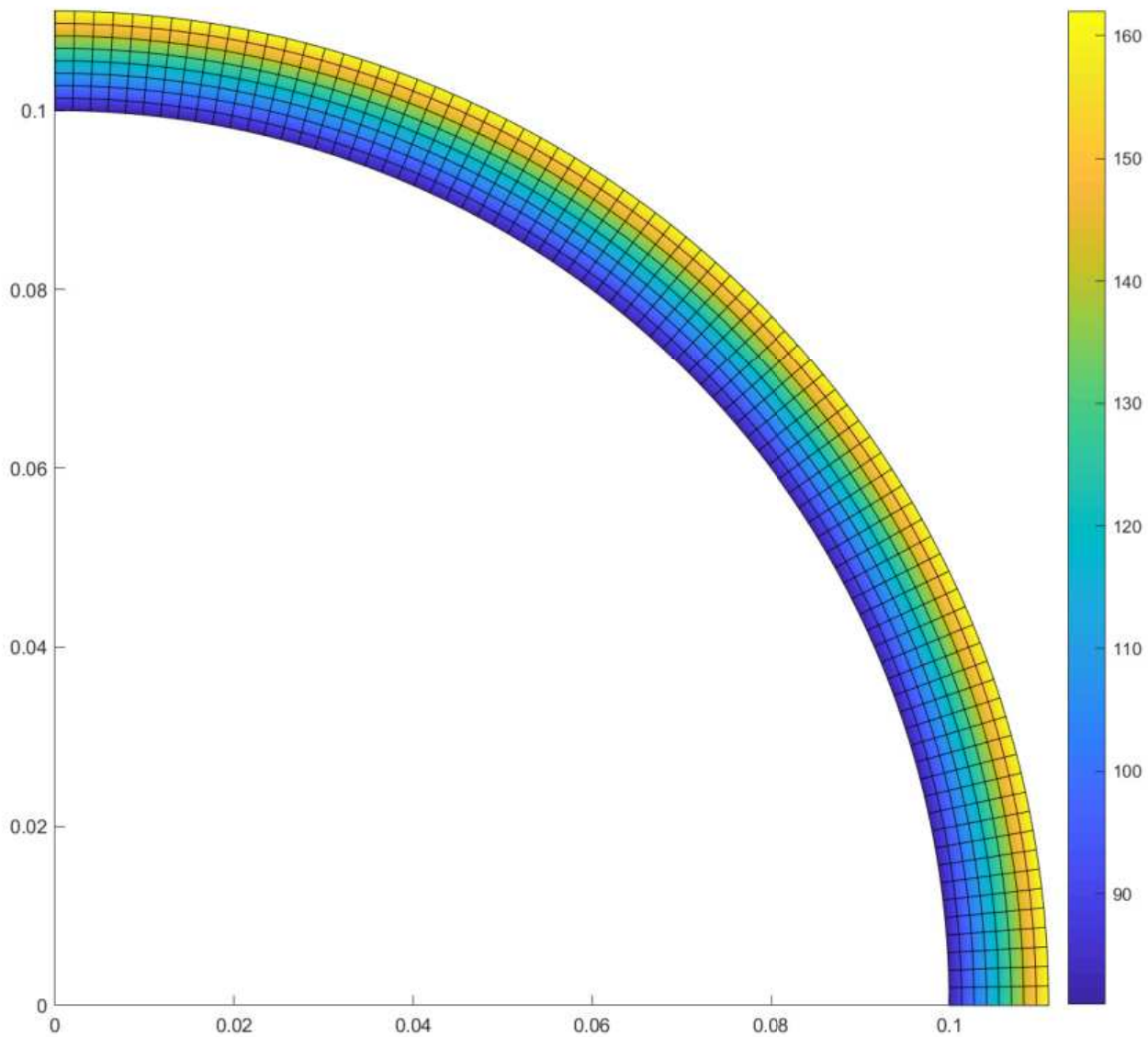}
    \caption{$t = \frac{\sqrt{80}}{9}\tau$ ($r/r_0 = 1/9$).}
    \label{kidder_ninth}
  \end{subfigure}
  \hfill
  \begin{subfigure}[b]{0.32\textwidth}
    \includegraphics[width=\textwidth]{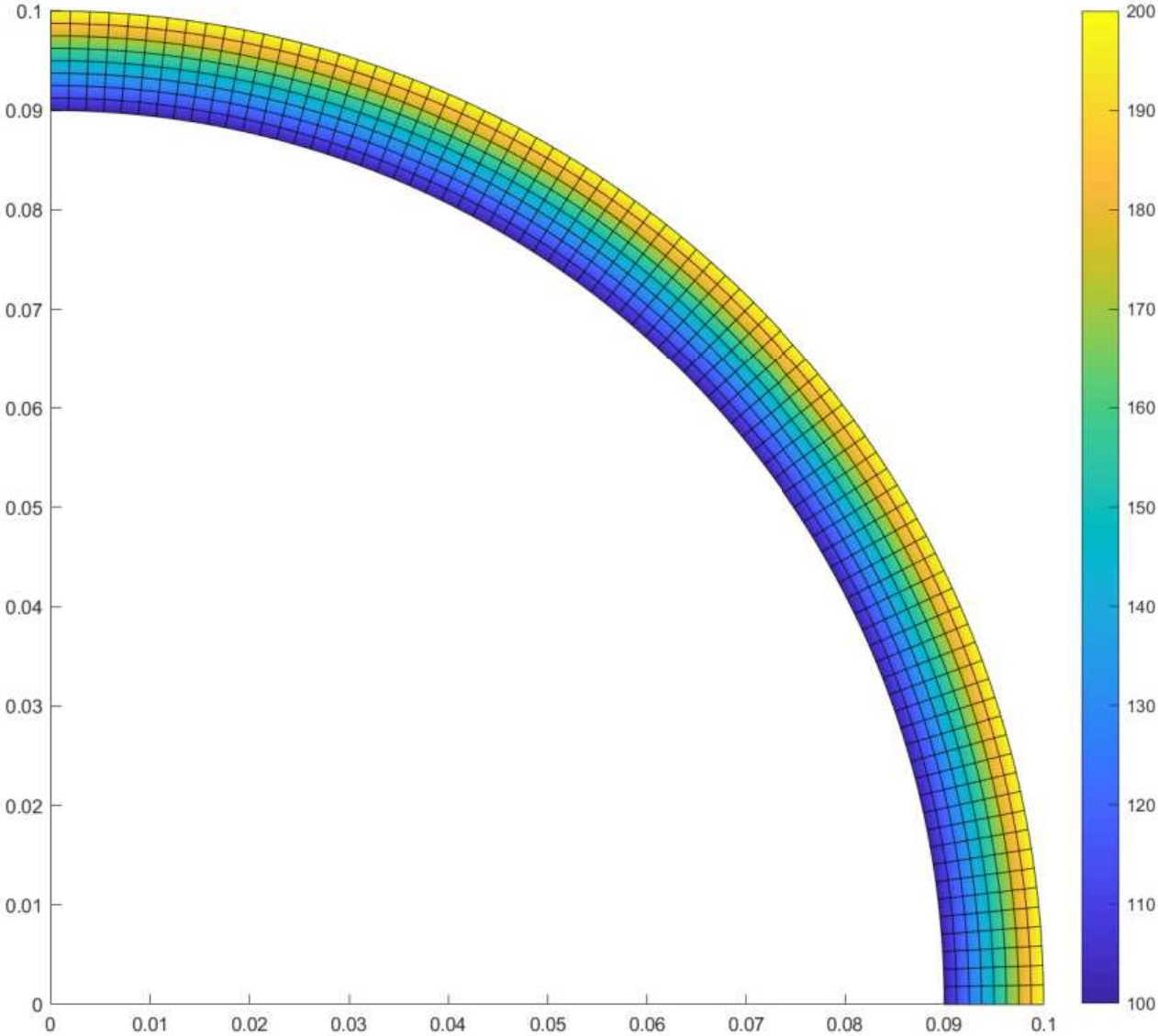}
    \caption{$t = \frac{\sqrt{99}}{10}\tau$ ($r/r_0 = 1/10$).}
    \label{kidder_tenth}
  \end{subfigure}

  \caption{Density contours for the Kidder problem at progressive compression stages using the stabilized $Q^3-Q^2$ element pair with hourglass control.}
  \label{Kidder_density}
\end{figure}

Figure~\ref{Kidder_density} shows the density fields at progressive timestamps using the stabilized $Q^3-Q^2$ space pair, advanced to specific evaluation times of $t \in \{\frac{\sqrt{3}}{2}\tau, \frac{\sqrt{8}}{3}\tau, \frac{\sqrt{15}}{4}\tau, \frac{\sqrt{24}}{5}\tau, \frac{\sqrt{80}}{9}\tau, \frac{\sqrt{99}}{10}\tau\}$. At these times, the ring's radius is symmetrically compressed to $\frac{1}{2}, \frac{1}{3}, \frac{1}{4}, \frac{1}{5}, \frac{1}{9}$, and $\frac{1}{10}$ of its initial radius, respectively. The simulated density fields maintain excellent symmetry and closely match the analytical self-similar solution. The proposed high-order SGH method exhibits highly robust performance even under extreme compression (1/10 of the initial radius)—a notoriously challenging condition for Lagrangian hydrodynamics that often triggers mesh entanglement.

\begin{figure}[htbp]
  \centering
  \begin{subfigure}[b]{0.485\textwidth}
    \includegraphics[width=\textwidth]{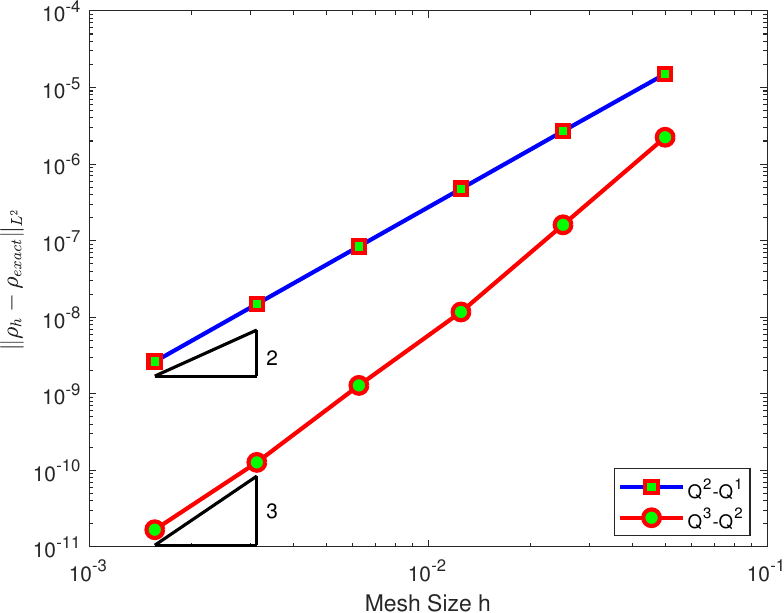}
    \caption{$L^2$ density error convergence rates.}
    \label{kidder_convergence_rate}
  \end{subfigure}
  \hfill
  \begin{subfigure}[b]{0.485\textwidth}
    \includegraphics[width=\textwidth]{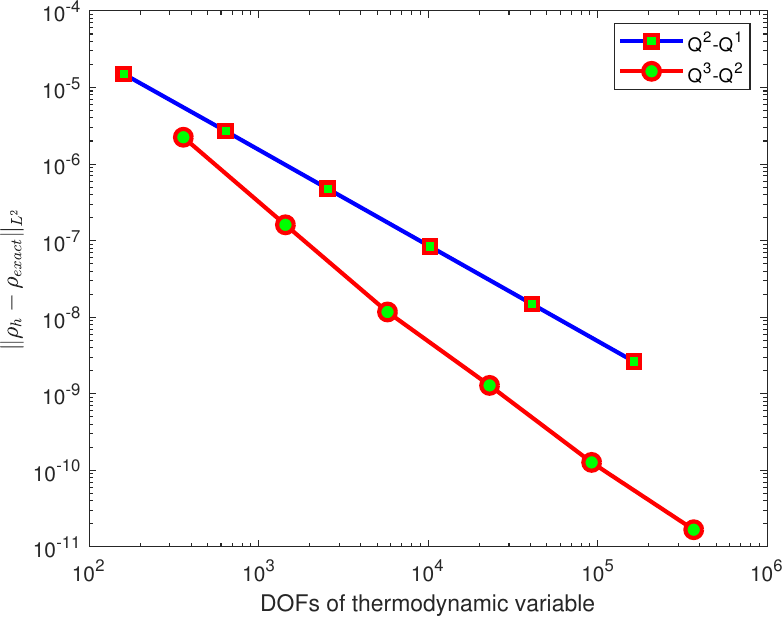}
    \caption{$L^2$ density error vs. thermodynamic DOFs.}
    \label{kidder_dofs}
  \end{subfigure}
  \caption{Error analysis and convergence properties for the Kidder problem at $t = \frac{\sqrt{3}}{2}\tau$.}
\label{Kidder_rate}
\end{figure}

Figure~\ref{Kidder_rate} and Table~\ref{tab:density_error} present the $L^2$ error norms of the density field evaluated at $t = \frac{\sqrt{3}}{2}\tau$ for both the $Q^2-Q^1$ and $Q^3-Q^2$ space pairs. The error norms are plotted against the corresponding DOFs of the thermodynamic variables in Figure~\ref{Kidder_rate}(b). Although the observed convergence rate for the $Q^3-Q^2$ case exhibits slight degradation from the perfect theoretical order—an expected consequence of introducing the subzonal anti-hourglass restoring forces—the fundamental advantage of the high-order formulation remains vividly evident. Specifically, the $Q^3-Q^2$ space pair achieves superior numerical accuracy while utilizing substantially fewer DOFs compared to the $Q^2-Q^1$ pair.

\begin{table}[htbp]
\centering
\caption{Density convergence rates for the Kidder problem evaluated at $t=\frac{\sqrt{3}}{2}\tau$.}
\label{tab:density_error}
\begin{tabular}{c c S[table-format=1.4e-1] c c S[table-format=1.4e-1] c}
\toprule
\multirow{2}{*}{$h$} & \multicolumn{3}{c}{$Q^2-Q^1$} & \multicolumn{3}{c}{$Q^3-Q^2$} \\
\cmidrule{2-4} \cmidrule{5-7}
& {DOFs} & {Error} & {Order} & {DOFs} & {Error} & {Order} \\
\midrule
$\frac{1}{20}$  & (205,160)    & 1.5029E-5 & --     & (427,360)    & 2.2386E-6 & -- \\
$\frac{1}{40}$  & (729,640)    & 2.6921E-6 & 2.4809 & (1573,1440)  & 1.6088E-7 & 3.7985 \\
$\frac{1}{80}$  & (2737,2560)  & 4.7532E-7 & 2.5017 & (6025,5760)  & 1.1694E-8 & 3.7821 \\
$\frac{1}{160}$ & (10593,10240)& 8.4030E-8 & 2.4999 & (23569,23040)& 1.2847E-9 & 3.1862 \\
$\frac{1}{320}$ & (41665,40960)& 1.4861E-8 & 2.4993 & (93217,92160)& 1.2683E-10& 3.3405 \\
$\frac{1}{640}$ & (165249,163840)&2.6292E-9& 2.4988 & (370753,368640)&1.6755E-11& 2.9202 \\
\bottomrule
\end{tabular}
\end{table}

\subsection{Robustness and shock-capturing capability}
This subsection presents numerical results for strong shock-driven benchmark problems. All cases here employ artificial viscosity. We emphasize that artificial viscosity natively dampens hourglass motion when it is locally active (i.e., near the shock front), which is why mesh quality often remains intact near the shock even without explicit hourglass control. However, in regions propagating away from the shock—where the artificial viscosity correctly switches off—severe hourglass distortions emerge. If left unchecked, these unphysical zero-energy modes lead to instability and chaotic geometric tangling in subsequent time steps. Therefore, explicit hourglass control remains absolutely necessary for shock wave problems, especially in long-term transient simulations.

In the following shock-driven tests, we systematically evaluate our proposed high-order framework. To specifically demonstrate the efficacy of our tunable hourglass control and to highlight the advancements over existing approaches, the Noh and Sedov problems will feature a three-way comparison: (1) our framework without hourglass control, (2) our framework with the proposed tunable hourglass control, and (3) the original discretization baseline (hereafter referred to as the "Original method") detailed in Section~\ref{sec:original_discretization}.

\subsubsection{Sod shock tube: 1D equivalent and 2D extension} 
The Sod problem is typically a one-dimensional benchmark problem. Here, we consider solving the Sod problem on two-dimensional structured and unstructured meshes, presenting two distinct geometrical extensions.

\begin{figure}[htbp]
  \begin{center}
    \includegraphics[width=0.99\textwidth]{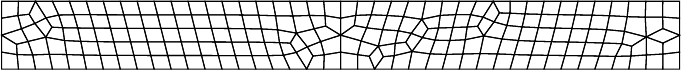}
    \caption{Initial unstructured mesh configuration for the 1D-like Sod shock tube problem.}
    \label{sod_unstructure_grid}
  \end{center}
\end{figure}

For the 1D equivalent case, the computational domain is $[0,1]\times[0,0.1]$, which is divided into two regions. The left region, $[0,0.5]\times[0,0.1]$, is occupied by an ideal gas with an initial state of $\rho=1$, $e=2.5$, and $\gamma=1.4$. Conversely, the right region, $[0.5,1]\times[0,0.1]$, has an initial state of $\rho=0.125$, $e=2$, and $\gamma=1.4$. We investigate two types of meshes for this configuration. The structured uniform initial mesh is simulated without any hourglass control. The initial unstructured mesh, generated using Gmsh~\cite{Geuzaine2009Gmsh}, is shown in Figure~\ref{sod_unstructure_grid}.

\begin{figure}[htbp]
  \centering
  \begin{subfigure}[b]{0.29\textwidth}
    \includegraphics[width=\textwidth]{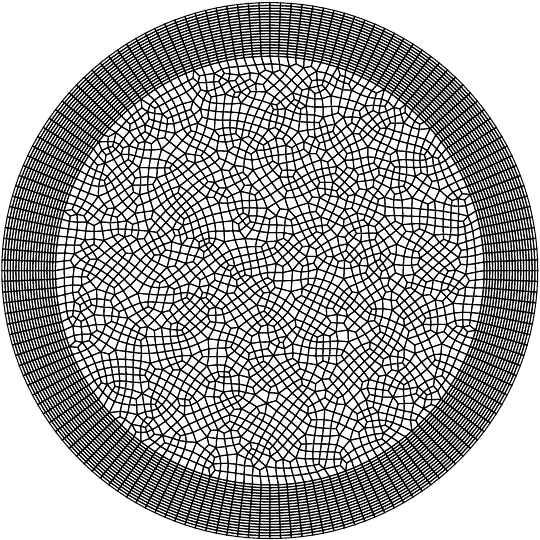}
    \caption{Initial mesh.}
    \label{Sod2D_initial_grid}
  \end{subfigure}
  \hfill
  \begin{subfigure}[b]{0.345\textwidth}
    \includegraphics[width=\textwidth]{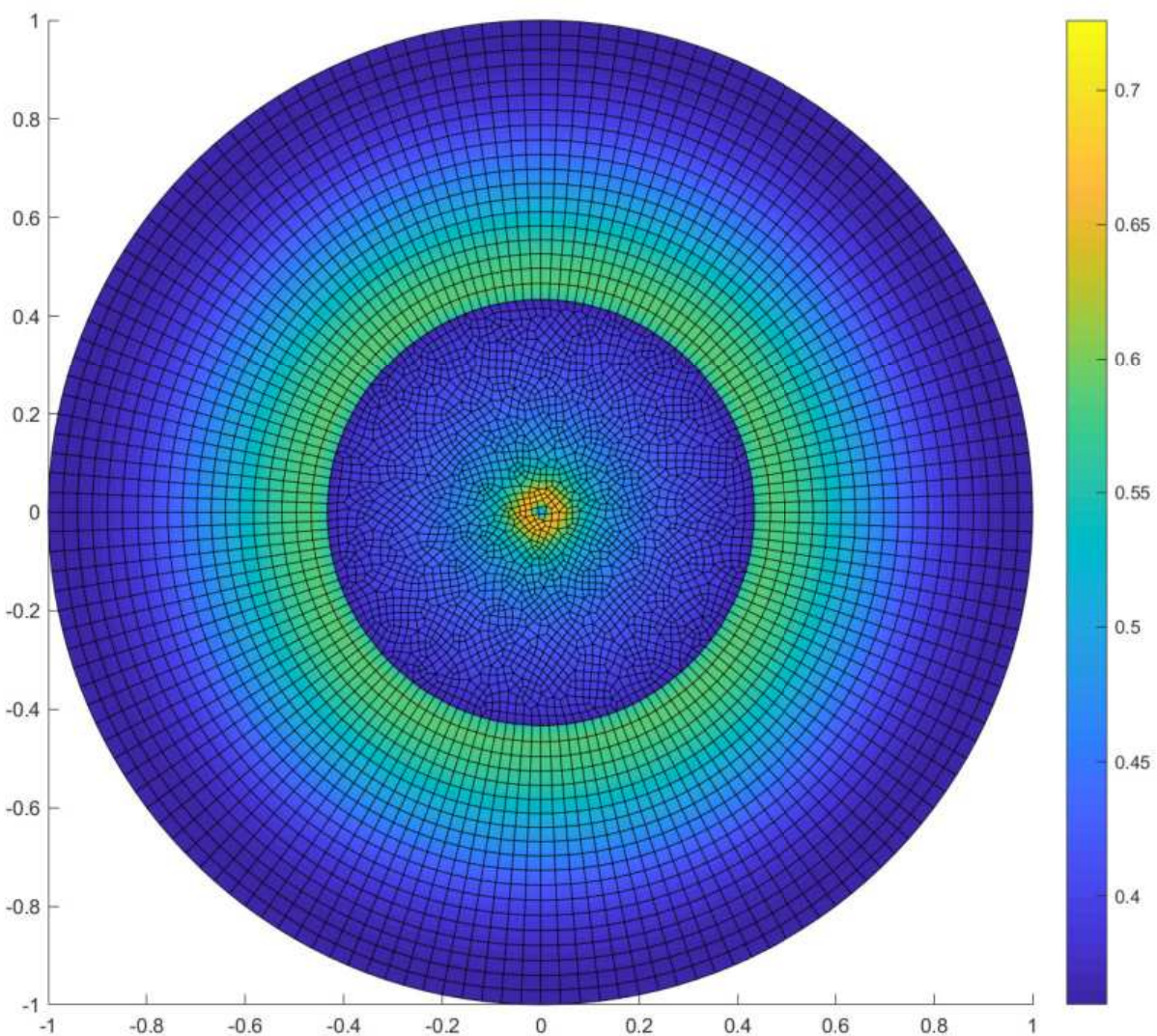}
    \caption{$Q^2-Q^1$ element pair.}
    \label{Sod2D_Q2_Q1}
  \end{subfigure}
  \hfill
  \begin{subfigure}[b]{0.345\textwidth}
    \includegraphics[width=\textwidth]{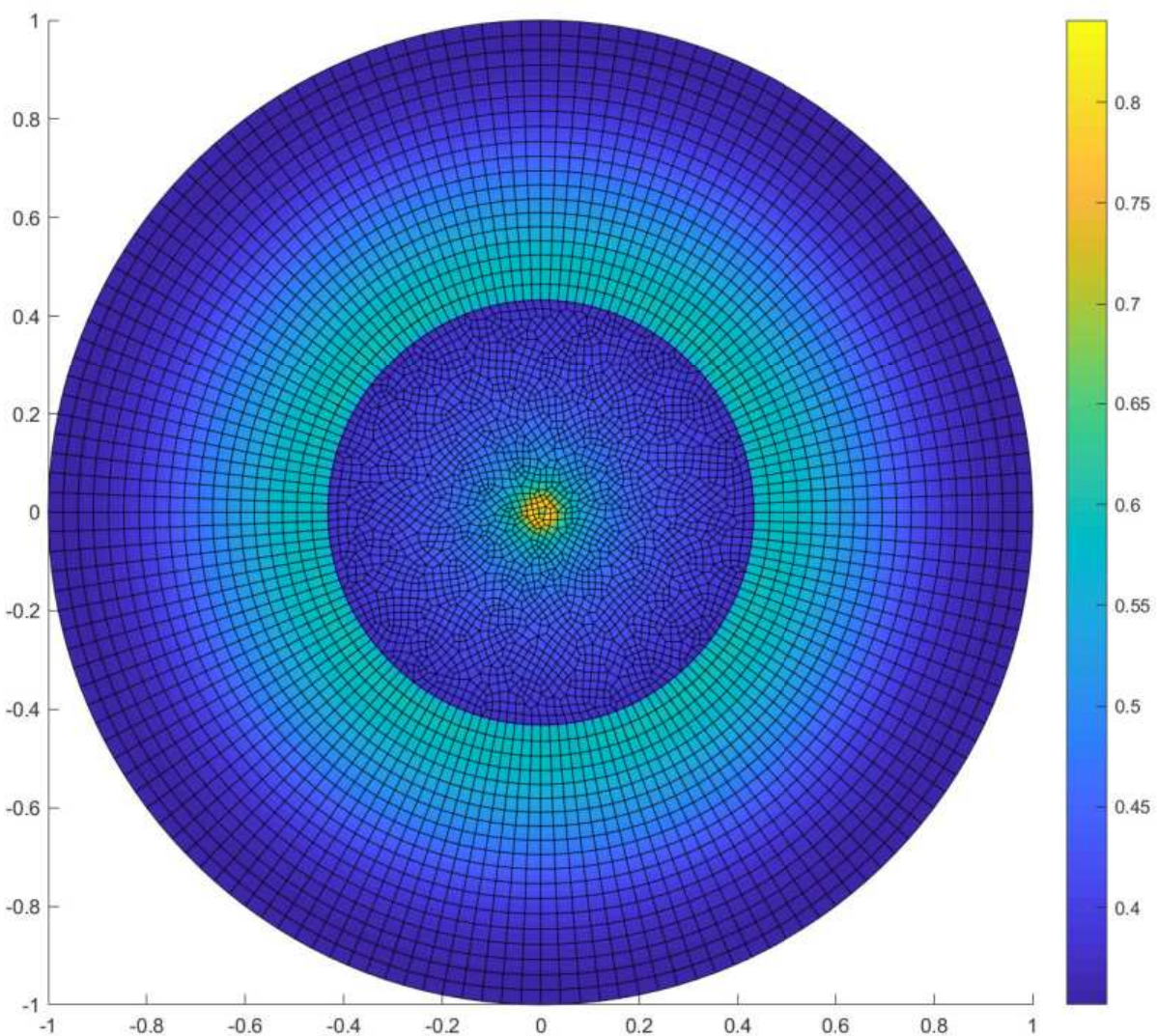}
    \caption{$Q^3-Q^2$ element pair.}
    \label{Sod2D_Q3_Q2}
  \end{subfigure}  
  \caption{Initial mesh and 2D density contours for the circular 2D Sod problem at $t=0.4$.}
\label{Sod2D_plane}
\end{figure}

For the 2D extension, the setup is identical to the "2D Sod" problem presented in~\cite{Sun2022On}. The computational domain is a unit circle centered at the origin, divided into two concentric parts, and the initial mesh is demonstrated in Figure~\ref{Sod2D_plane}a. The outer ring region is filled with a high-pressure ideal gas with an initial state of $\rho=1$, $e=2.5$, and $\gamma=1.4$, while the inner circular region contains a low-pressure gas with an initial state of $\rho=0.125$, $e=2$, and $\gamma=1.4$.

\begin{figure}[htbp]
  \centering
  \begin{subfigure}[b]{0.99\textwidth}
    \includegraphics[width=\textwidth]{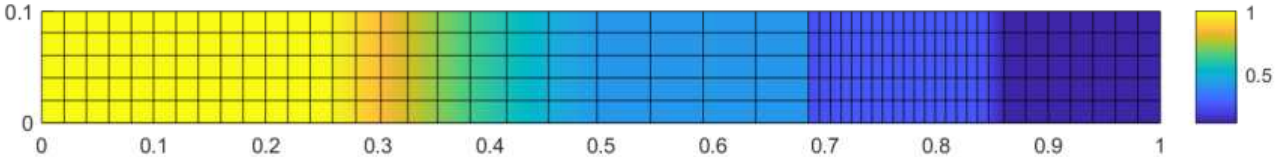}
    \caption{$Q^2-Q^1$ element pair.}
    \label{Sod_Q2_Q1_structure}
  \end{subfigure}
  
  \vspace{0.5em}
  
  \begin{subfigure}[b]{0.99\textwidth}
    \includegraphics[width=\textwidth]{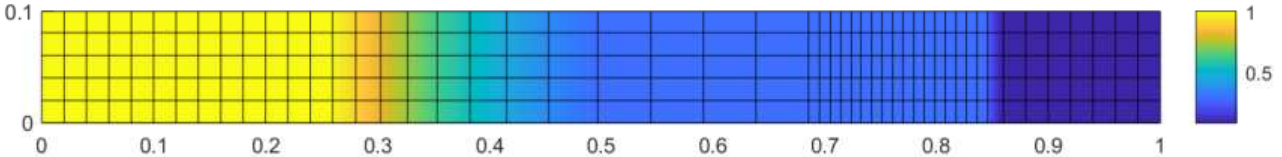}
    \caption{$Q^3-Q^2$ element pair.}
    \label{Sod_Q3_Q2_structure}
  \end{subfigure}
  \caption{Density contours and deformed meshes for the Sod problem on a structured uniform mesh without hourglass control at $t=0.2$.}
\end{figure}

\begin{figure}[htbp]
  \centering
  \begin{subfigure}[b]{0.485\textwidth}
    \includegraphics[width=\textwidth]{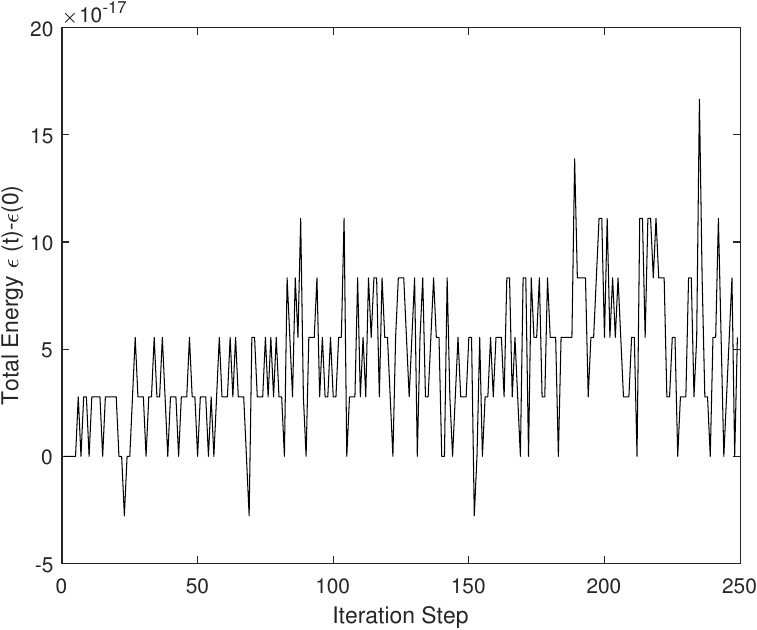}
    \caption{$Q^2-Q^1$ element pair.}
    \label{Sod_TEV_Q2_Q1}
  \end{subfigure}
  \hfill
  \begin{subfigure}[b]{0.485\textwidth}
    \includegraphics[width=\textwidth]{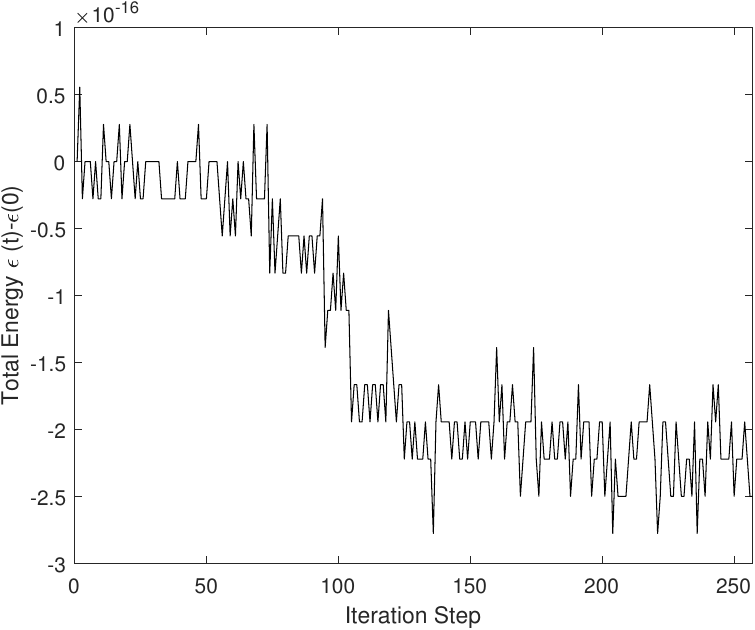}
    \caption{$Q^3-Q^2$ element pair.}
    \label{Sod_TEV_Q3_Q2}
  \end{subfigure}
  \caption{Time histories of total energy variation for the Sod problem across all time steps.}
\label{Sod_TEV}
\end{figure}

First, we present the numerical results for the structured uniform mesh without hourglass control. Figures~\ref{Sod_Q2_Q1_structure} and~\ref{Sod_Q3_Q2_structure} show the density distributions at $t=0.2$ for the $Q^2-Q^1$ and $Q^3-Q^2$ elements, respectively. The artificial viscosity parameters are set to $c_1=0.2$ and $c_2=0.2$. The numerical results clearly demonstrate the absence of hourglass motion on the structured uniform initial mesh. Both the contact discontinuity and the shock waves are well captured without any spurious oscillations. Figure~\ref{Sod_TEV} displays the total energy variation over time, which remains at machine precision throughout the simulation. 

\begin{figure}[htbp]
  \centering
  \begin{subfigure}[b]{0.99\textwidth}
    \includegraphics[width=\textwidth]{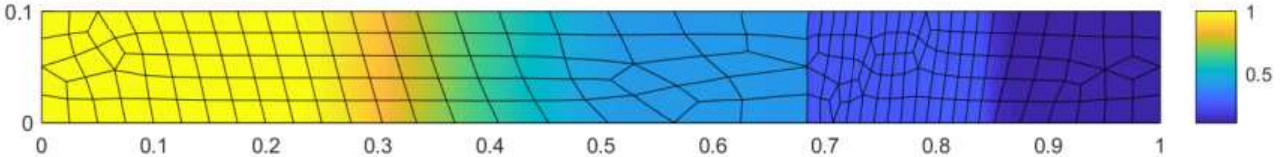}
    \caption{$Q^2-Q^1$ element pair.}
    \label{Sod_Q2_Q1_unstructure_nohg}
  \end{subfigure}
  
  \vspace{0.5em}
  
  \begin{subfigure}[b]{0.99\textwidth}
    \includegraphics[width=\textwidth]{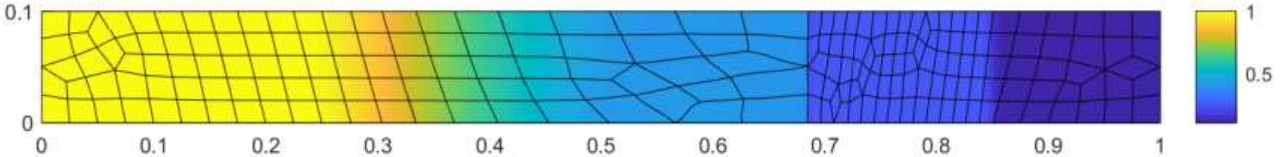}
    \caption{$Q^3-Q^2$ element pair.}
    \label{Sod_Q3_Q2_structure_nohg}
  \end{subfigure}
  
  \caption{Density contours and deformed meshes for the Sod problem on an unstructured mesh without hourglass control at $t=0.2$.}
  \label{Sod_unstructure_nohg}
\end{figure}

Figure~\ref{Sod_unstructure_nohg} shows the density distribution of the 1D equivalent Sod problem on the unstructured mesh without hourglass control at $t=0.2$. The artificial viscosity parameters are again set to $c_1=0.2$ and $c_2=0.2$. In this scenario, the grid near the contact discontinuity tends to deform, indicating the presence of hourglass motion. However, the shock waves are still captured correctly without significant spurious oscillations.

\begin{figure}[htbp]
  \centering
  \begin{subfigure}[b]{0.99\textwidth}
    \includegraphics[width=\textwidth]{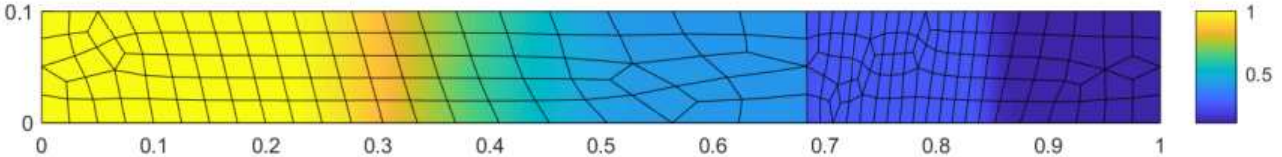}
    \caption{$Q^2-Q^1$ element pair.}
    \label{Sod_Q2_Q1_unstructure_hg}
  \end{subfigure}
  
  \vspace{0.5em}
  
  \begin{subfigure}[b]{0.99\textwidth}
    \includegraphics[width=\textwidth]{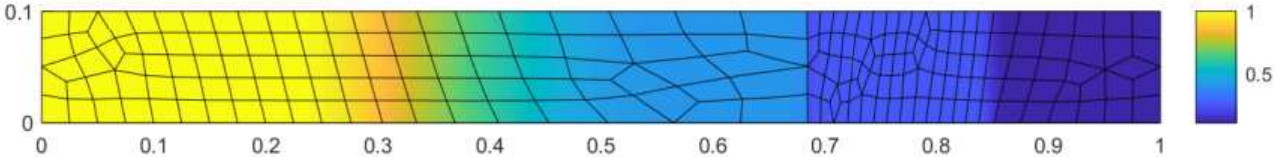}
    \caption{$Q^3-Q^2$ element pair.}
    \label{Sod_Q3_Q2_structure_hg}
  \end{subfigure}
  
  \caption{Density contours and deformed meshes for the Sod problem on an unstructured mesh with hourglass control at $t=0.2$.}
  \label{Sod_unstructure_hg}
\end{figure}

Figure~\ref{Sod_unstructure_hg} presents the density distribution on the unstructured mesh with hourglass control activated. The viscosity parameters remain unchanged, while the hourglass control parameter is set to $1$. With hourglass control applied, the grid quality is significantly improved, and the numerical results are highly satisfactory.

Figures~\ref{Sod2D_Q2_Q1} and~\ref{Sod2D_Q3_Q2} show the density distributions for the Sod extension at $t=0.4$. The artificial viscosity parameters are set to $c_1=0.5$ and $c_2=0.5$, and hourglass control is applied with the parameter set to $1$. The results demonstrate that both the $Q^2-Q^1$ and $Q^3-Q^2$ discretizations maintain grid quality well and accurately handle the complex flow features and the unstructured geometry.

\begin{figure}[htbp]
  \centering
  \begin{subfigure}[b]{0.485\textwidth}
    \includegraphics[width=\textwidth]{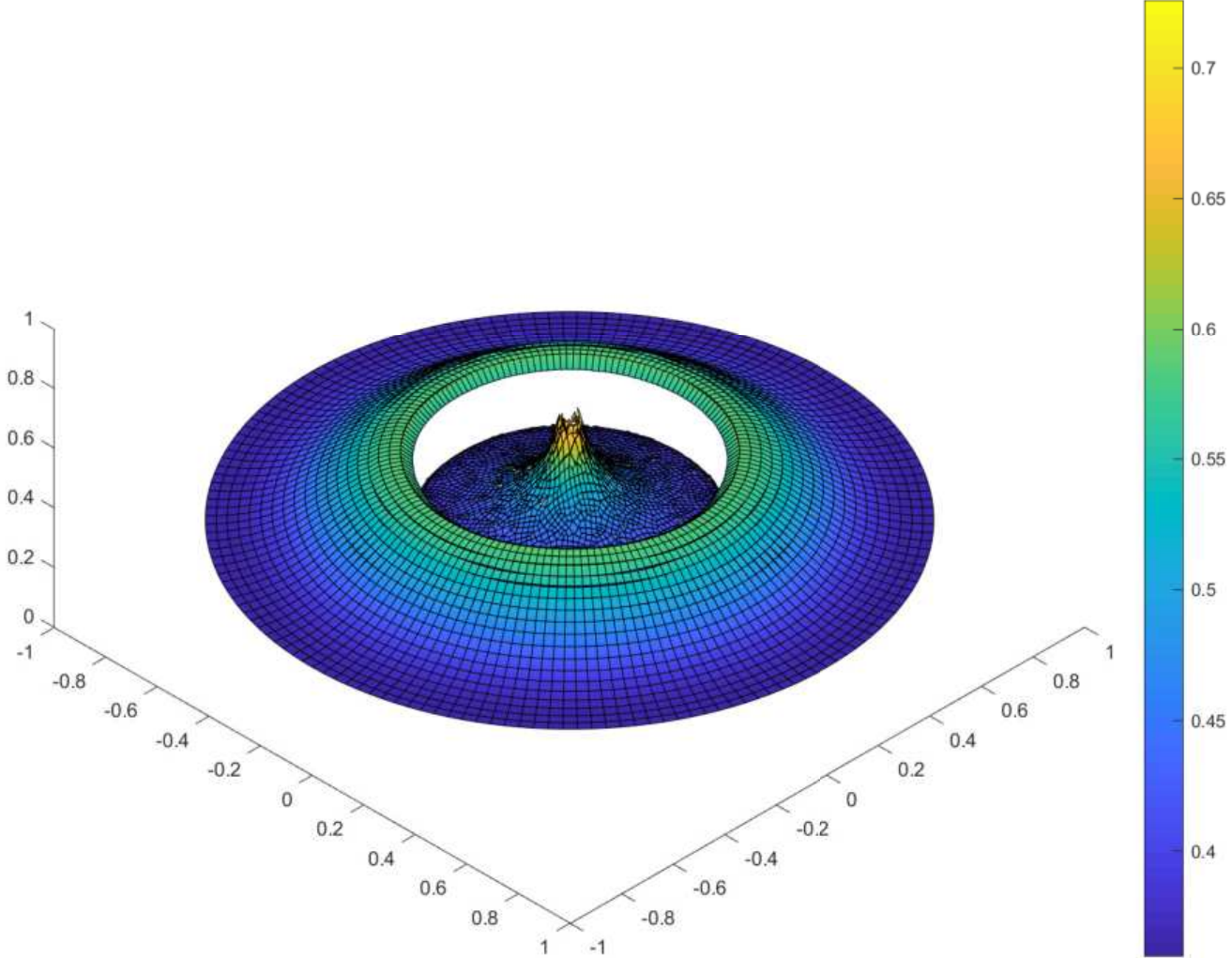}
    \caption{$Q^2-Q^1$ element pair.}
    \label{Sod2D_Q2_Q1_3D}
  \end{subfigure}
  \hfill
  \begin{subfigure}[b]{0.485\textwidth}
    \includegraphics[width=\textwidth]{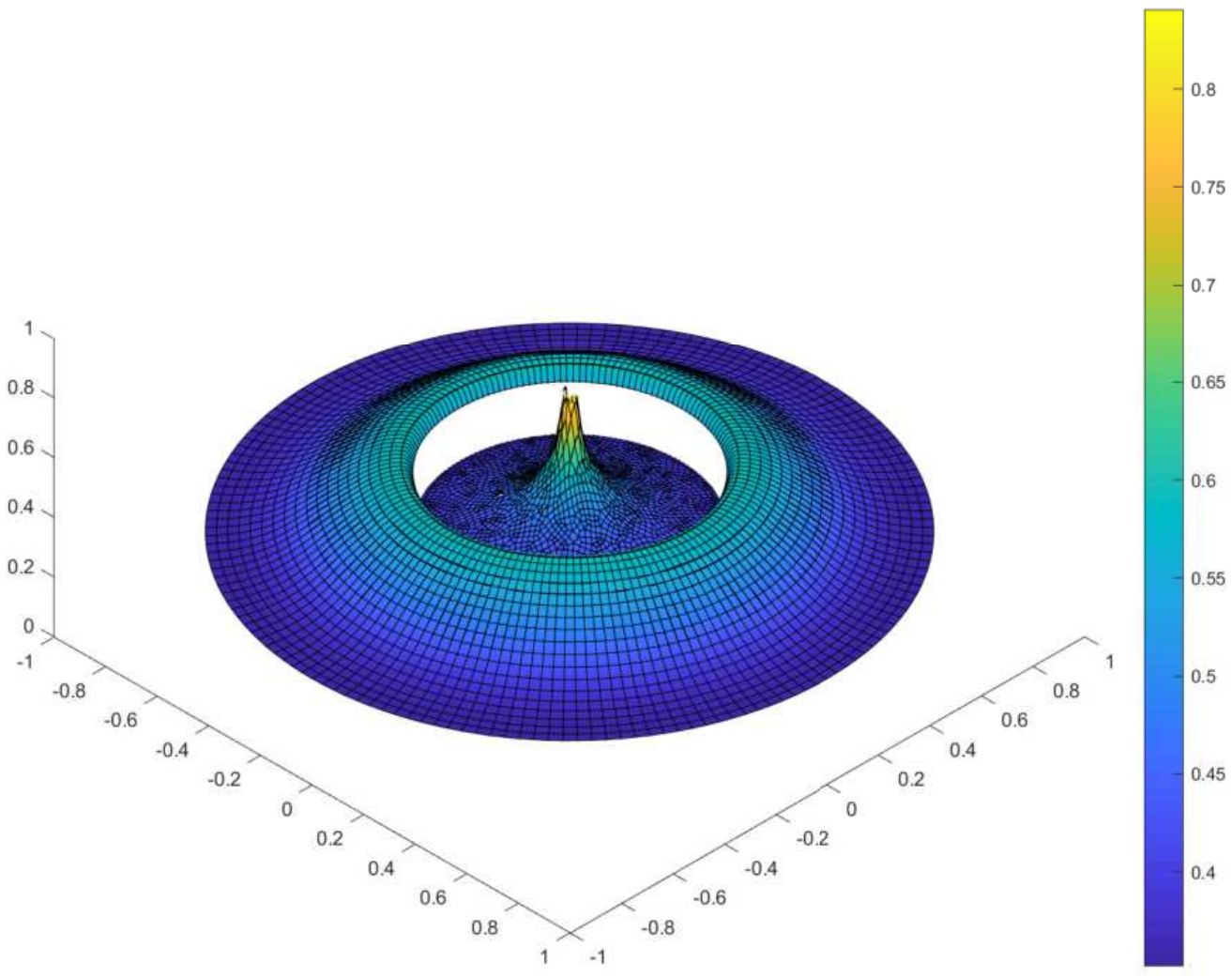}
    \caption{$Q^3-Q^2$ element pair.}
    \label{Sod2D_Q3_Q2_3D}
  \end{subfigure}  
  \caption{3D perspective views of the density distributions for the circular 2D Sod problem at $t=0.4$.}
\label{Sod2D_3D}
\end{figure}

Figure~\ref{Sod2D_3D} presents a 3D visualization of the density profiles for this circular domain. The maximum density value achieved with the $Q^3-Q^2$ element is higher than that of the $Q^2-Q^1$ element, indicating the improved resolution capabilities of the higher-order discretization using identical initial mesh.

\subsubsection{Noh implosion problem} 
An ideal gas moves toward the origin with a unit velocity, generating a perfectly symmetric strong shock at the origin that propagates outward. The initial computational domain is $[0,1]^2$, with the specific heat ratio $\gamma = 5/3$, an initial density $\rho = 1$, and an initial internal energy of $1 \times 10^{-10}$. This problem has an exact analytical solution: the shock reaches a radius of $0.2$ at $t = 0.6$, with post-shock density and pressure limits of $16$ and $5.33$, respectively.

\begin{figure}[htbp]
  \centering
  \begin{subfigure}[b]{0.485\textwidth}
    \includegraphics[width=\textwidth]{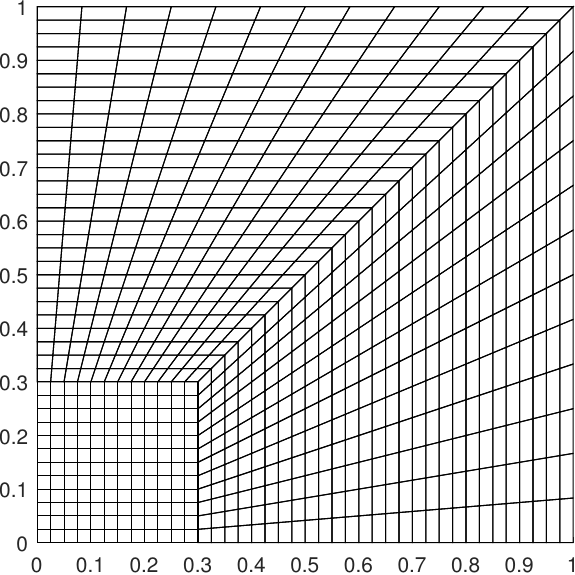}
    \caption{Butterfly mesh.}
    \label{Noh_butterfly_mesh}
  \end{subfigure}
  \hfill
  \begin{subfigure}[b]{0.485\textwidth}
    \includegraphics[width=\textwidth]{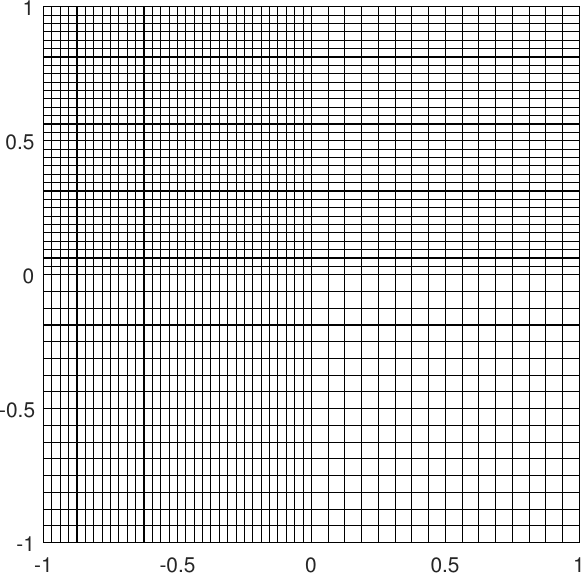}
    \caption{Asymmetric mesh.}
    \label{Noh_asym_mesh}
  \end{subfigure}
  \caption{Initial non-uniform mesh configurations for the Noh problem.}
\label{Noh_mesh}
\end{figure}

We consider three types of initial meshes for the Noh problem: a uniform structured mesh ($h=0.05$), a butterfly mesh (Figure~\ref{Noh_butterfly_mesh}), and an asymmetric mesh (Figure~\ref{Noh_asym_mesh}). For the Noh problem, we observed that the numerical performance improves when setting $c_{vor}=1$; therefore, we maintain this setting for all Noh simulations. When hourglass (HG) control is activated, the hourglass parameter is set to $1$.

\begin{figure}[htbp]
  \centering
  \begin{subfigure}[b]{0.32\textwidth}
    \includegraphics[width=\textwidth]{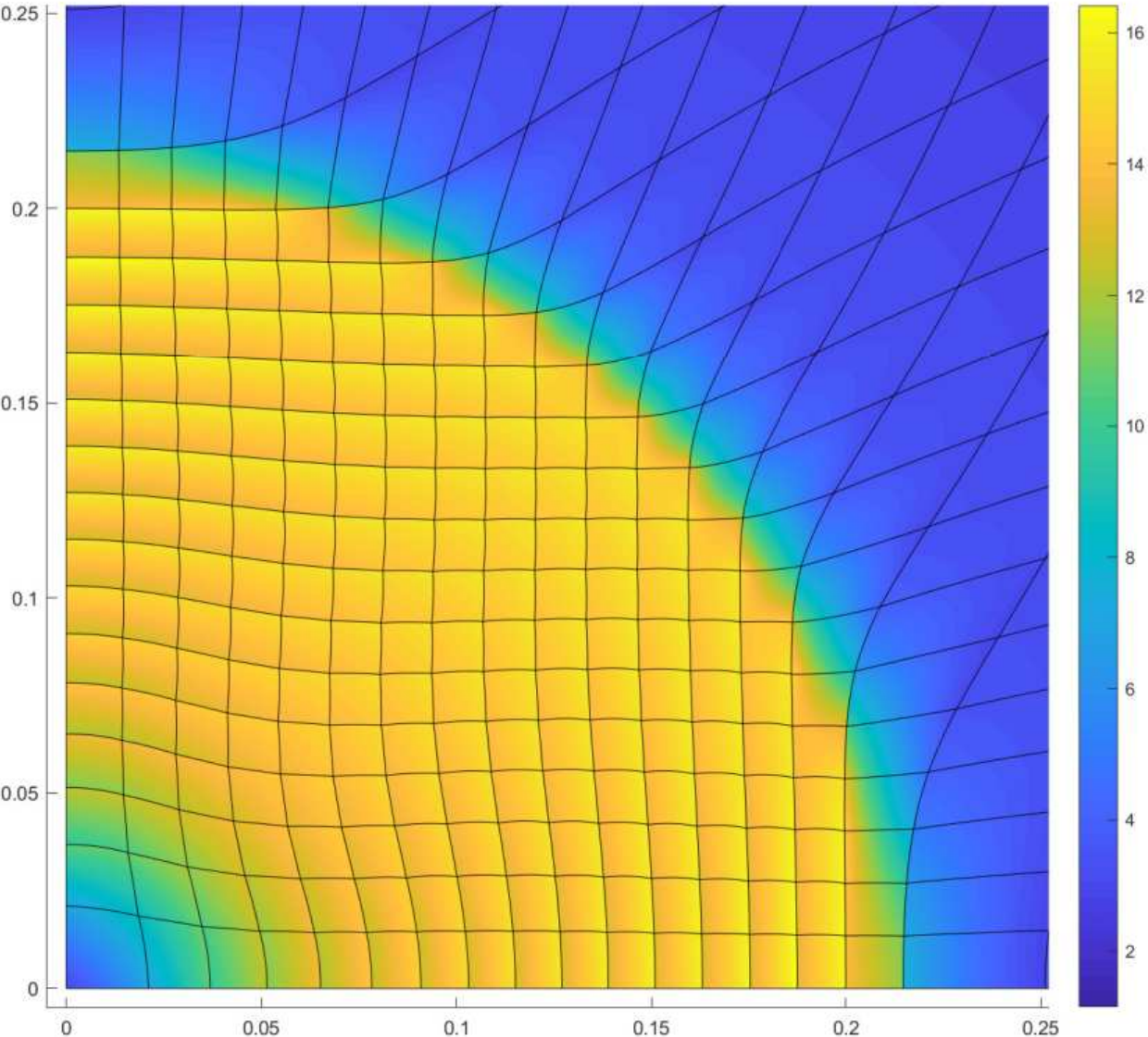}
    \caption{Density contour (without HG).}
    \label{Noh_Q2_nohg}
  \end{subfigure}
  \hfill
  \begin{subfigure}[b]{0.32\textwidth}
    \includegraphics[width=\textwidth]{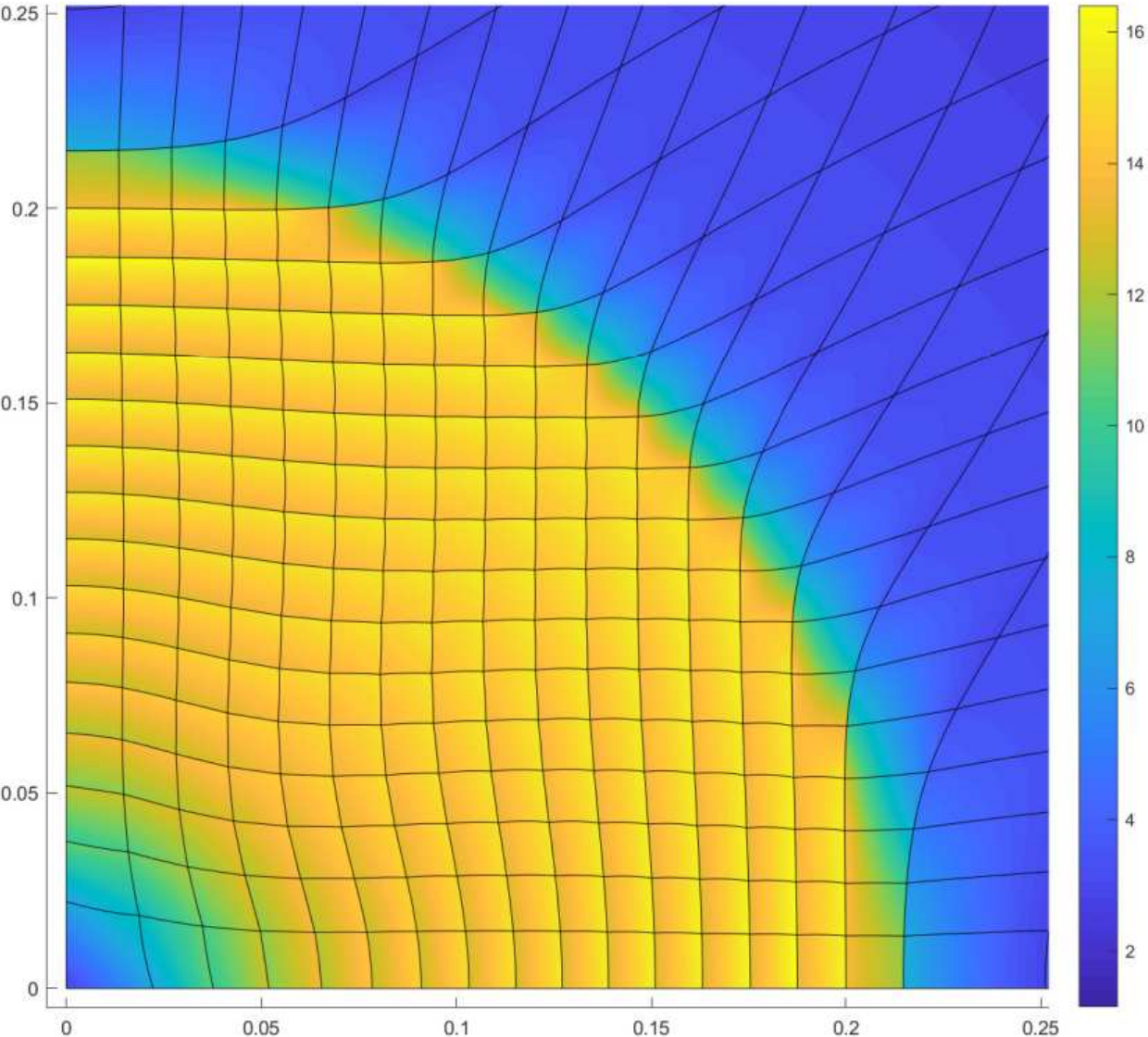}
    \caption{Density contour (with HG).}
    \label{Noh_Q2_hg}
  \end{subfigure}
  \hfill
  \begin{subfigure}[b]{0.32\textwidth}
    \includegraphics[width=\textwidth]{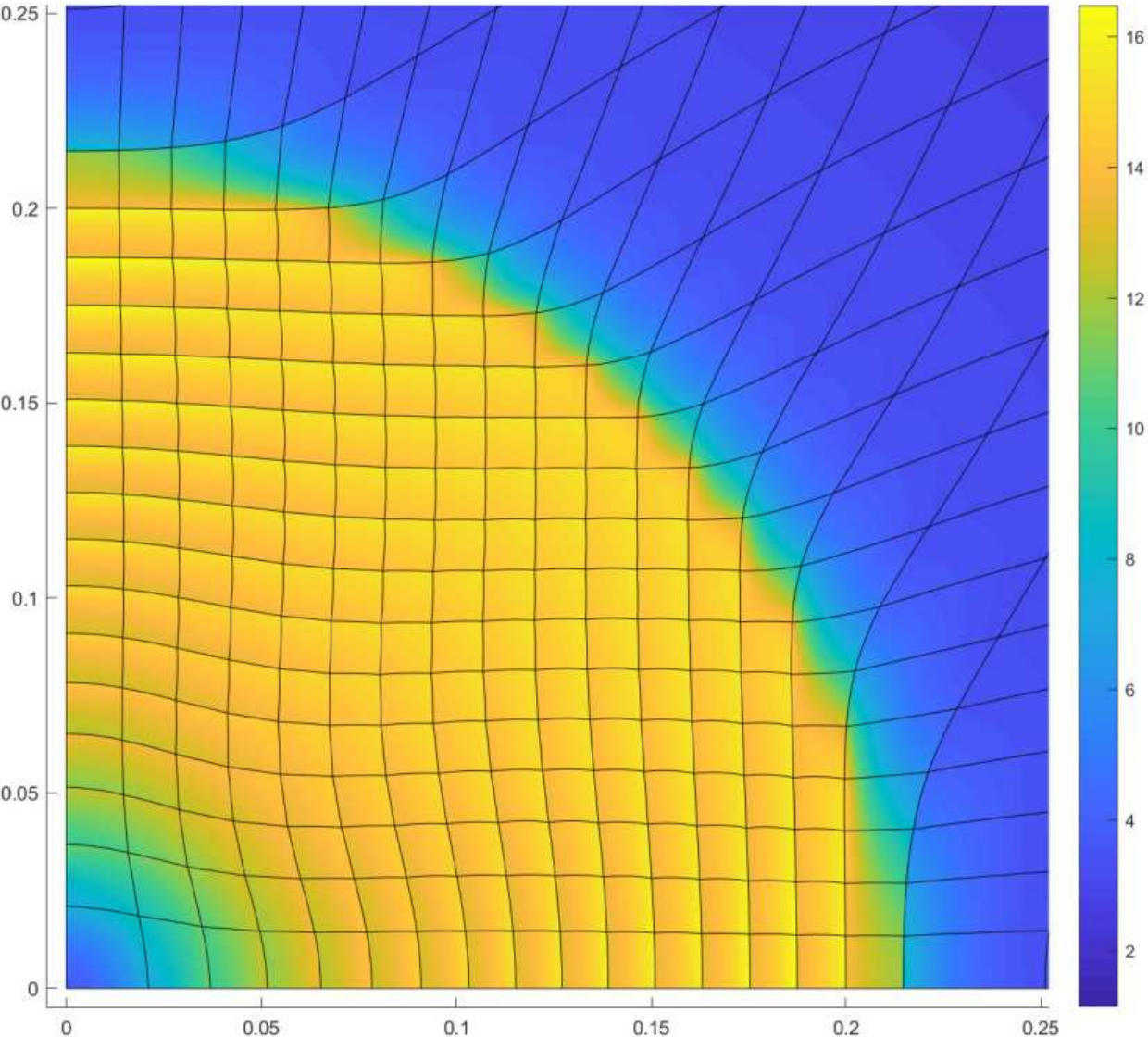}
    \caption{Density contour (Original method).}
    \label{Original_Noh_Q2}
  \end{subfigure}

  \vspace{0.5em}

  \begin{subfigure}[b]{0.32\textwidth}
    \includegraphics[width=\textwidth]{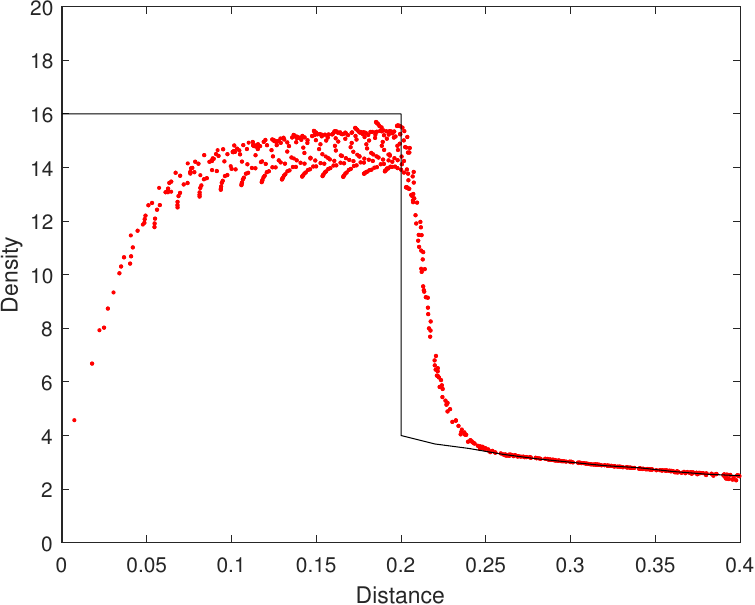}
    \caption{Pointwise density (without HG).}
    \label{Noh_Q2_point_nohg}
  \end{subfigure}
  \hfill
  \begin{subfigure}[b]{0.32\textwidth}
    \includegraphics[width=\textwidth]{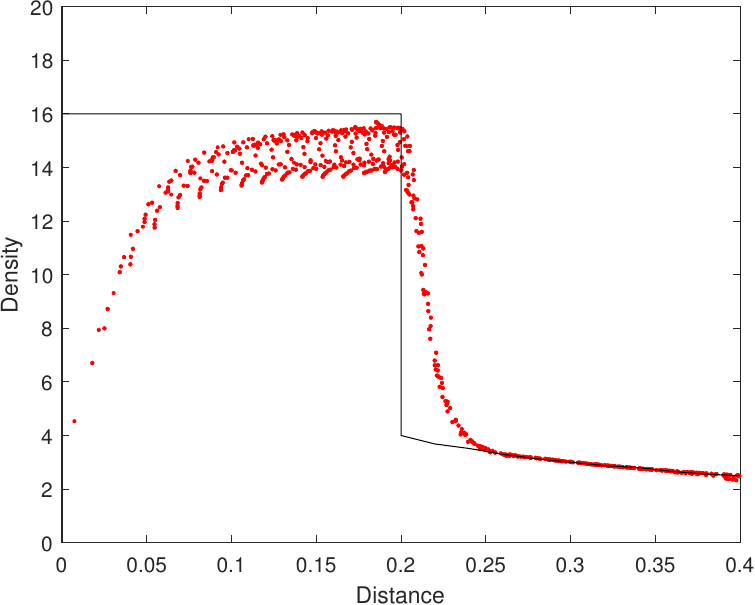}
    \caption{Pointwise density (with HG).}
    \label{Noh_Q2_point.pdf}
  \end{subfigure}
  \hfill
  \begin{subfigure}[b]{0.32\textwidth}
    \includegraphics[width=\textwidth]{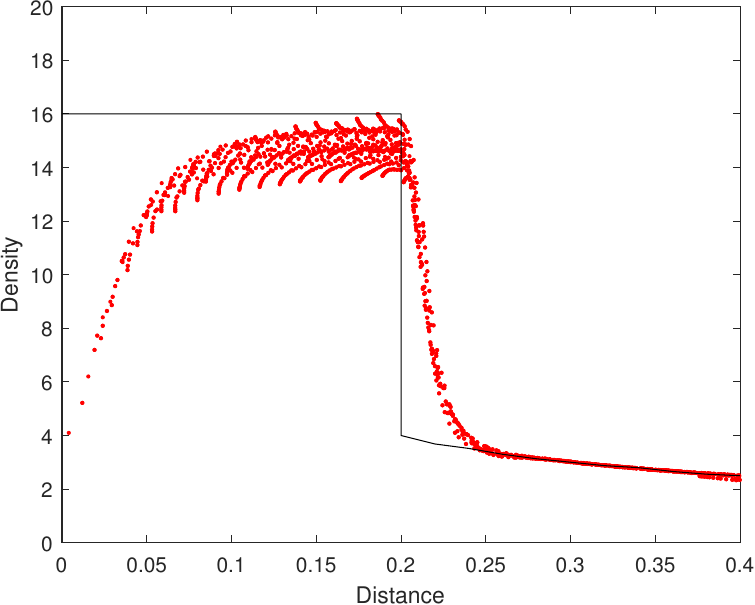}
    \caption{Pointwise density (Original method).}
    \label{Original_Noh_Q2_point}
  \end{subfigure}
  \caption{Density contours and 1D pointwise density profiles along the radius for the Noh problem using the $Q^2-Q^1$ element pair on a uniform mesh at $t = 0.6$.}
  \label{Noh_Q2_Q1}
\end{figure}

Figure~\ref{Noh_Q2_Q1} displays the density profiles and pointwise density values for the Noh problem using the $Q^2-Q^1$ element pair on a uniform structured mesh at $t=0.6$. We compare the three configurations mentioned earlier. The artificial viscosity parameters are set to $c_1=0.5$ and $c_2=0.5$. Figures~\ref{Noh_Q2_point_nohg}, \ref{Noh_Q2_point.pdf}, and \ref{Original_Noh_Q2_point} plot the pointwise density distributions against the exact solution. The Original method visibly displays a denser point distribution because it evaluates nine density values within each $Q^2$ element, whereas our proposed framework evaluates only four. Importantly, while our proposed framework does not strictly outperform the Original method, the overall numerical performance and accuracy of the three approaches remain highly comparable. Furthermore, even without hourglass control, our method maintains satisfactory grid quality.

\begin{figure}[htbp]
  \centering
  \begin{subfigure}[b]{0.32\textwidth}
    \includegraphics[width=\textwidth]{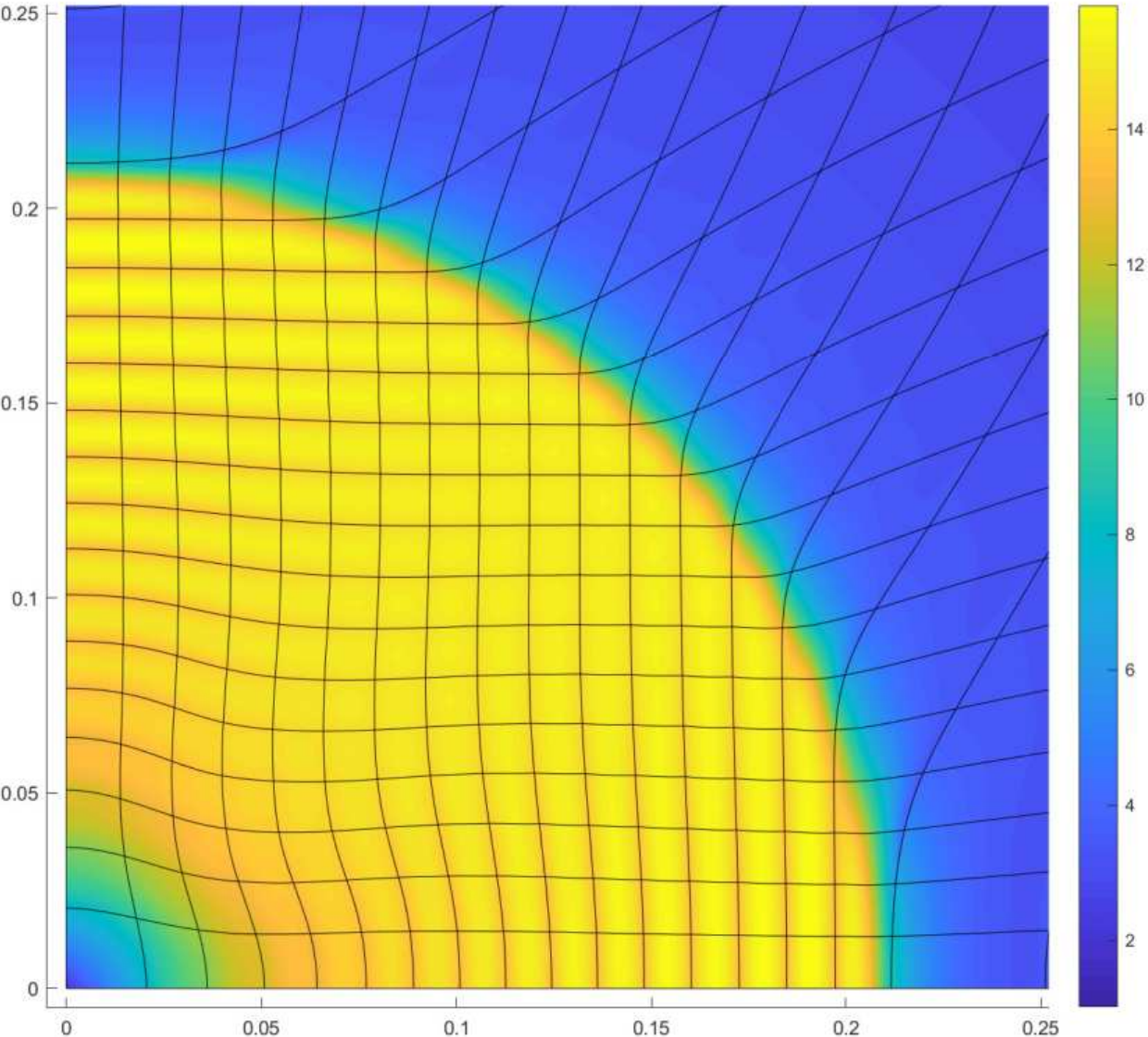}
    \caption{Density contour (without HG).}
    \label{Noh_Q3_nohg}
  \end{subfigure}
  \hfill
  \begin{subfigure}[b]{0.32\textwidth}
    \includegraphics[width=\textwidth]{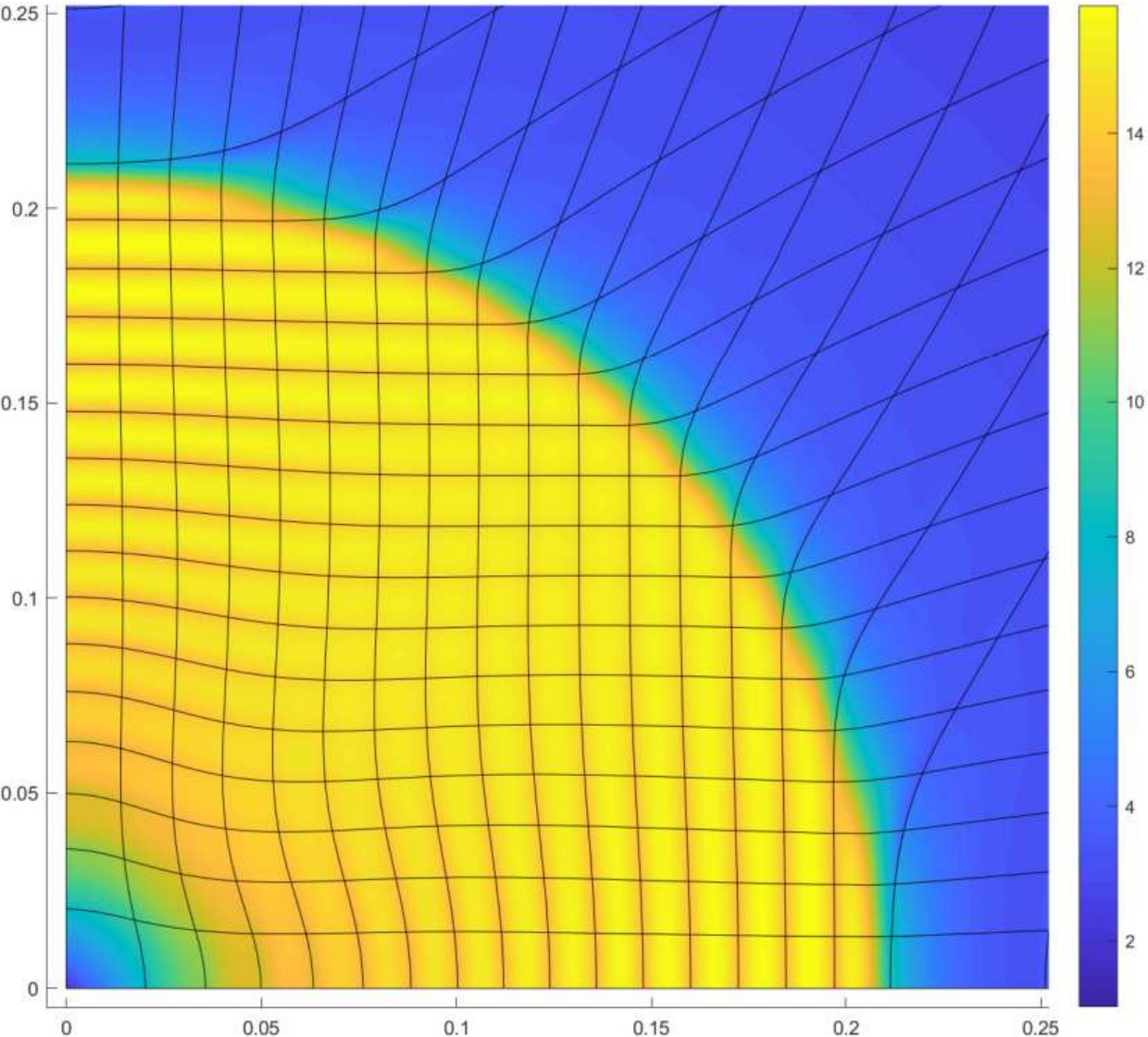}
    \caption{Density contour (with HG).}
    \label{Noh_Q3_hg}
  \end{subfigure}
  \hfill
  \begin{subfigure}[b]{0.32\textwidth}
    \includegraphics[width=\textwidth]{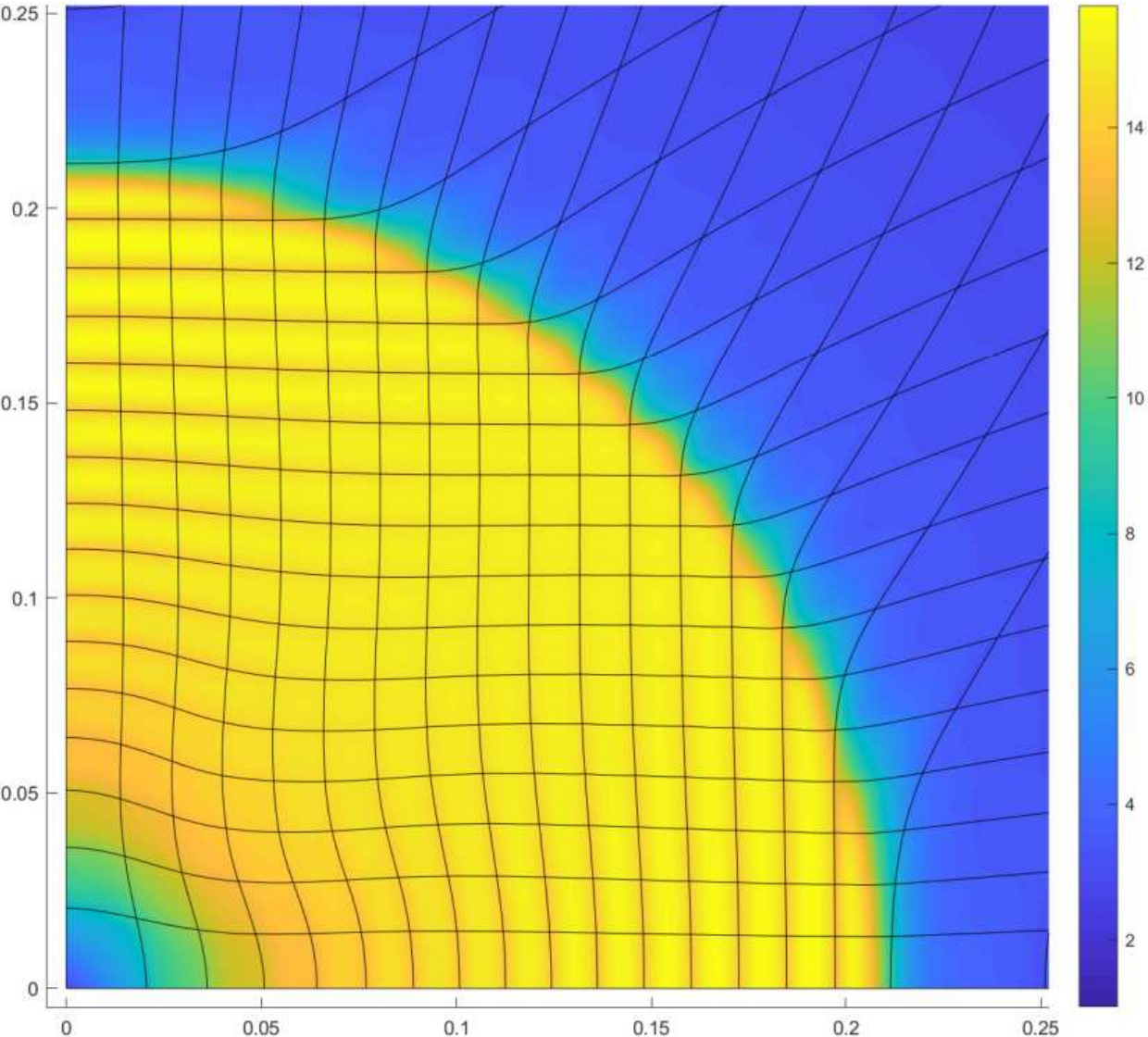}
    \caption{Density contour (Original method).}
    \label{Original_Noh_Q3}
  \end{subfigure}

  \vspace{0.5em}

  \begin{subfigure}[b]{0.32\textwidth}
    \includegraphics[width=\textwidth]{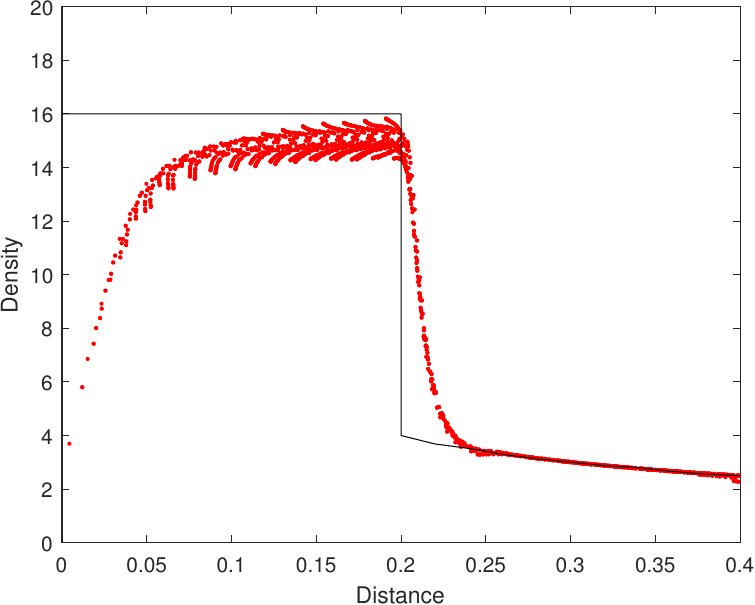}
    \caption{Pointwise density (without HG).}
    \label{Noh_Q3_point_nohg}
  \end{subfigure}
  \hfill
  \begin{subfigure}[b]{0.32\textwidth}
    \includegraphics[width=\textwidth]{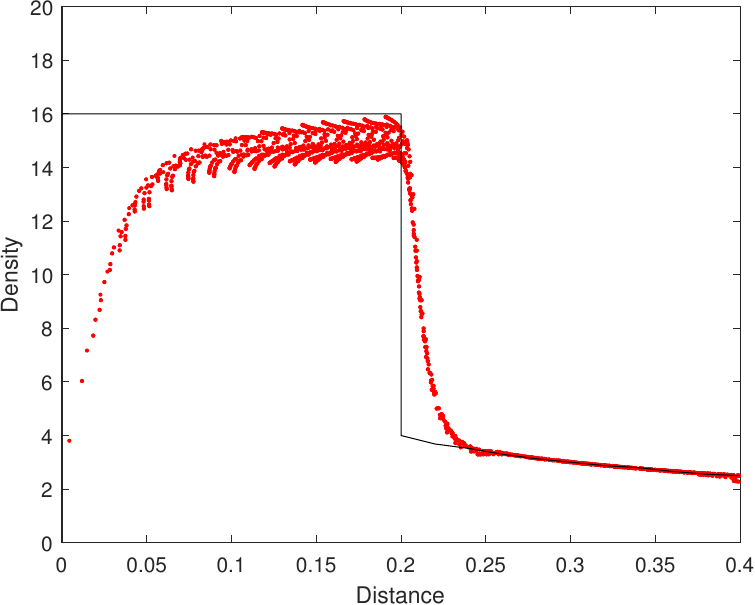}
    \caption{Pointwise density (with HG).}
    \label{Noh_Q3_point.pdf}
  \end{subfigure}
  \hfill
  \begin{subfigure}[b]{0.32\textwidth}
    \includegraphics[width=\textwidth]{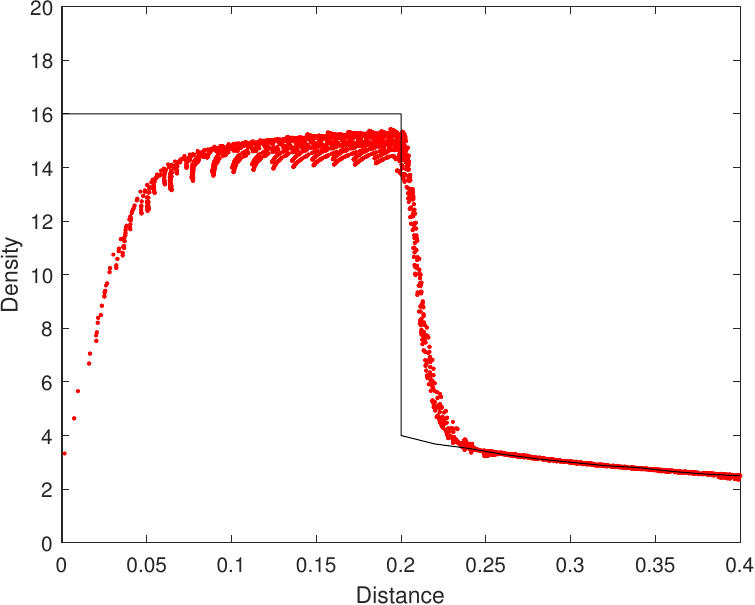}
    \caption{Pointwise density (Original method).}
    \label{Original_Noh_Q3_point}
  \end{subfigure}
  \caption{Density contours and 1D pointwise density profiles along the radius for the Noh problem using the $Q^3-Q^2$ element pair on a uniform mesh at $t = 0.6$.}
  \label{Noh_Q3_Q2}
\end{figure}

Figure~\ref{Noh_Q3_Q2} shows the corresponding results on the uniform mesh using the $Q^3-Q^2$ element pair, keeping all other parameters identical. In this scenario, the Original method exhibits slightly smaller density oscillations than our proposed framework. However, it is important to note that the maximum density values captured by our methods at the shock front are actually closer to the exact analytical solution. While reduced oscillations are often favored in Lagrangian simulations for overall stability, both approaches provide highly comparable overall solutions. Furthermore, all three strategies yield very similar grid distributions, demonstrating that our proposed methods remain robust and fully competitive.

\begin{figure}[htbp]
  \centering
  \begin{subfigure}[b]{0.32\textwidth}
    \includegraphics[width=\textwidth]{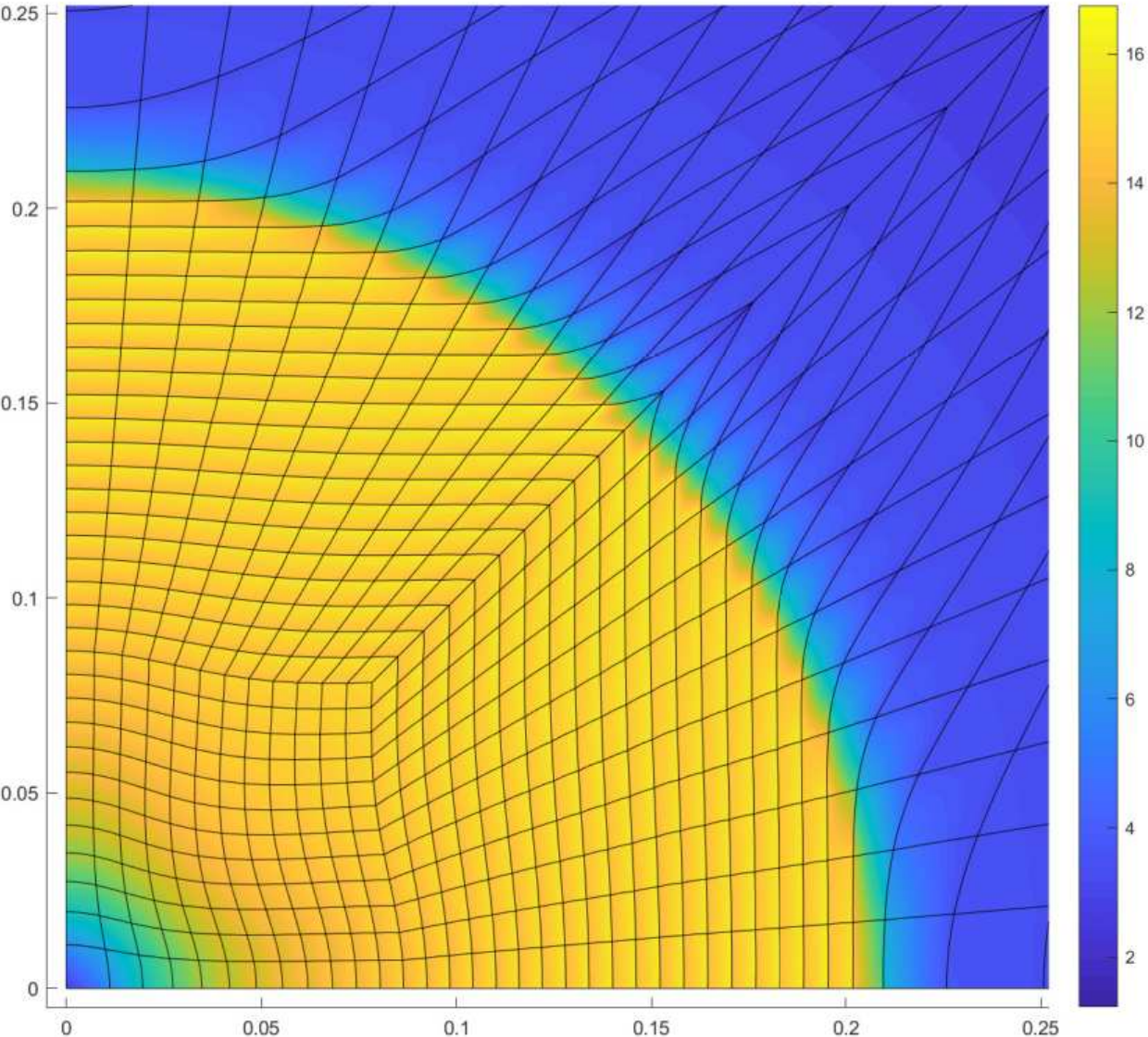}
    \caption{Density contour (without HG).}
    \label{Noh_Q2_butterfly_nohg}
  \end{subfigure}
  \hfill
  \begin{subfigure}[b]{0.32\textwidth}
    \includegraphics[width=\textwidth]{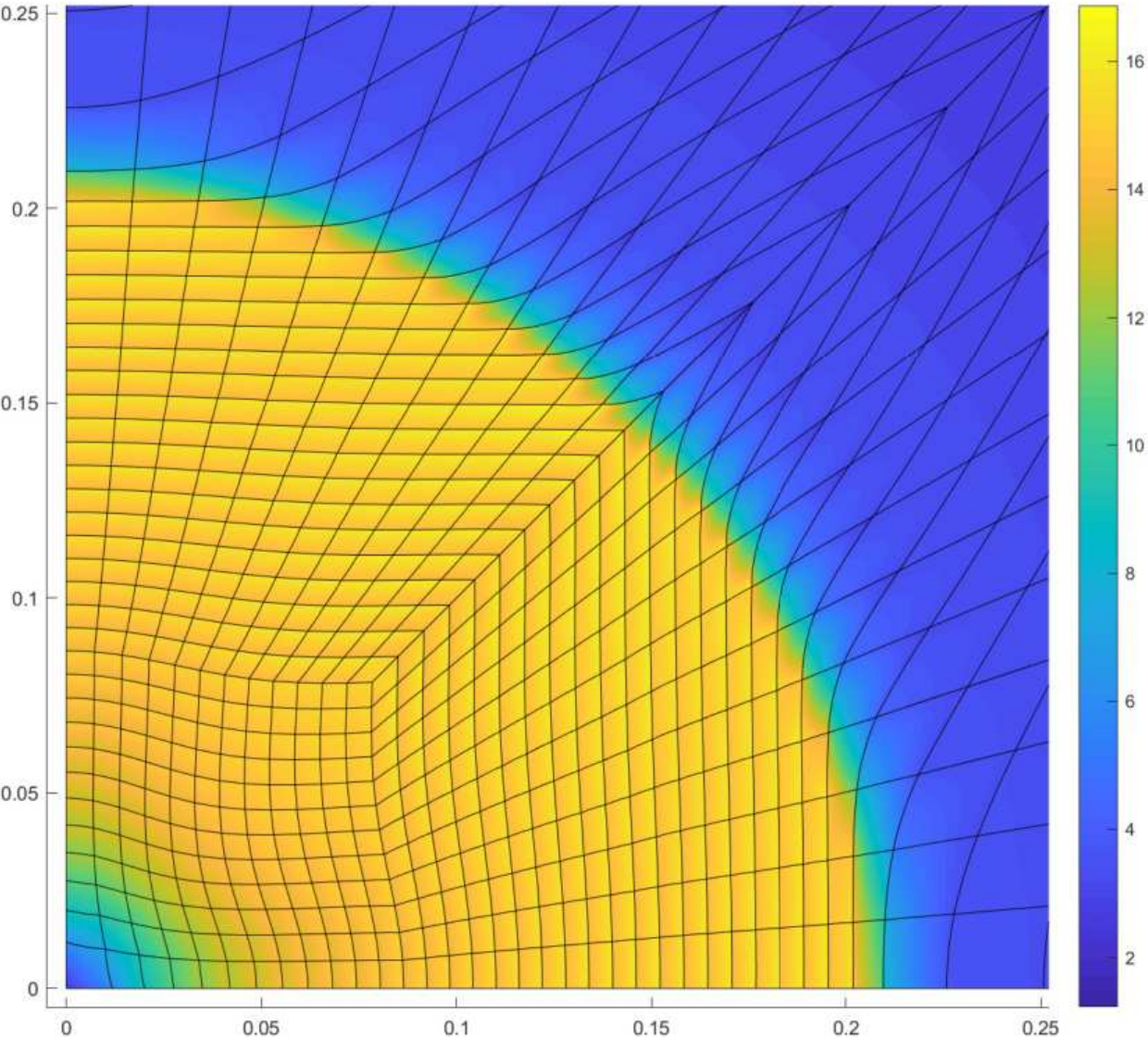}
    \caption{Density contour (with HG).}
    \label{Noh_Q2_butterfly_hg}
  \end{subfigure}
  \hfill
  \begin{subfigure}[b]{0.32\textwidth}
    \includegraphics[width=\textwidth]{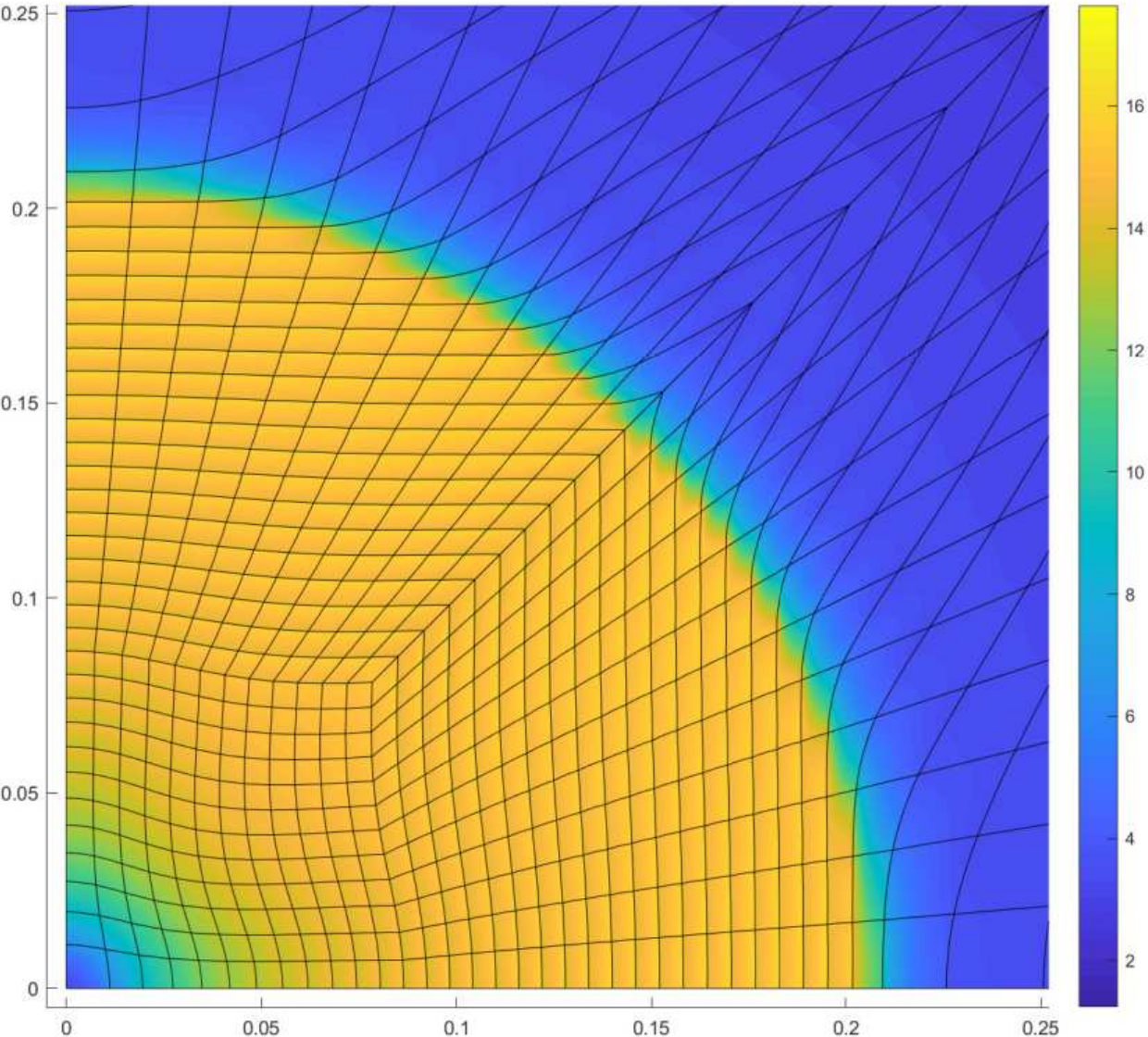}
    \caption{Density contour (Original method).}
    \label{Original_Noh_Q2_butterfly}
  \end{subfigure}

  \vspace{0.5em}

  \begin{subfigure}[b]{0.32\textwidth}
    \includegraphics[width=\textwidth]{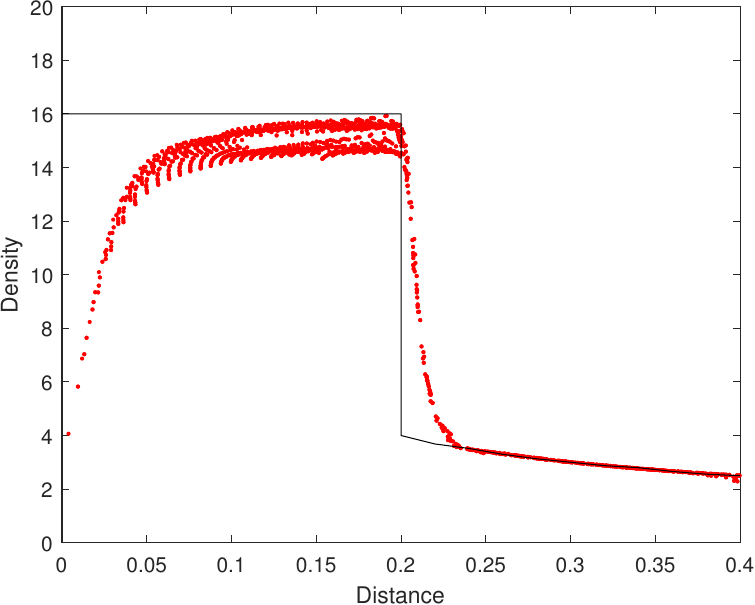}
    \caption{Pointwise density (without HG).}
    \label{Noh_Q2_butterfly_point_nohg}
  \end{subfigure}
  \hfill
  \begin{subfigure}[b]{0.32\textwidth}
    \includegraphics[width=\textwidth]{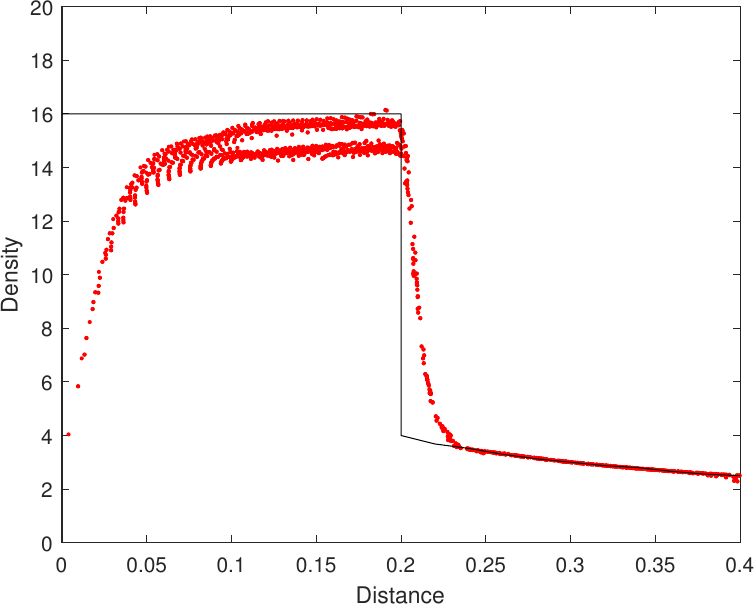}
    \caption{Pointwise density (with HG).}
    \label{Noh_Q2_butterfly_point.pdf}
  \end{subfigure}
  \hfill
  \begin{subfigure}[b]{0.32\textwidth}
    \includegraphics[width=\textwidth]{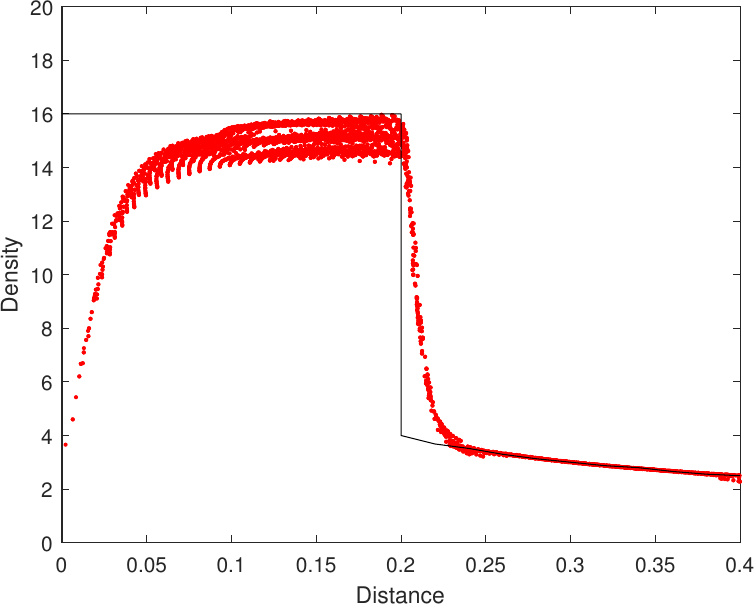}
    \caption{Pointwise density (Original method).}
    \label{Original_Noh_Q2_butterfly_point}
  \end{subfigure}
  \caption{Density contours and 1D pointwise density profiles along the radius for the Noh problem using the $Q^2-Q^1$ element pair on a butterfly mesh at $t = 0.6$.}
  \label{Noh_butterfly_Q2_Q1}
\end{figure}

Figure~\ref{Noh_butterfly_Q2_Q1} presents the density profiles for the Noh problem using the $Q^2-Q^1$ element pair on a butterfly mesh. Again, we compare the method without HG, with HG, and the Original method, using artificial viscosity parameters $c_1=0.5$ and $c_2=0.5$. Similar to the uniform mesh case, the three strategies show only marginal differences in their final grid distributions. The pointwise density profile of the Original method (Figure~\ref{Original_Noh_Q2_butterfly_point}) appears denser due to the higher number of evaluation points per $Q^2$ element, but the overall shock resolution remains highly comparable across all three approaches.

\begin{figure}[htbp]
  \centering
  \begin{subfigure}[b]{0.32\textwidth}
    \includegraphics[width=\textwidth]{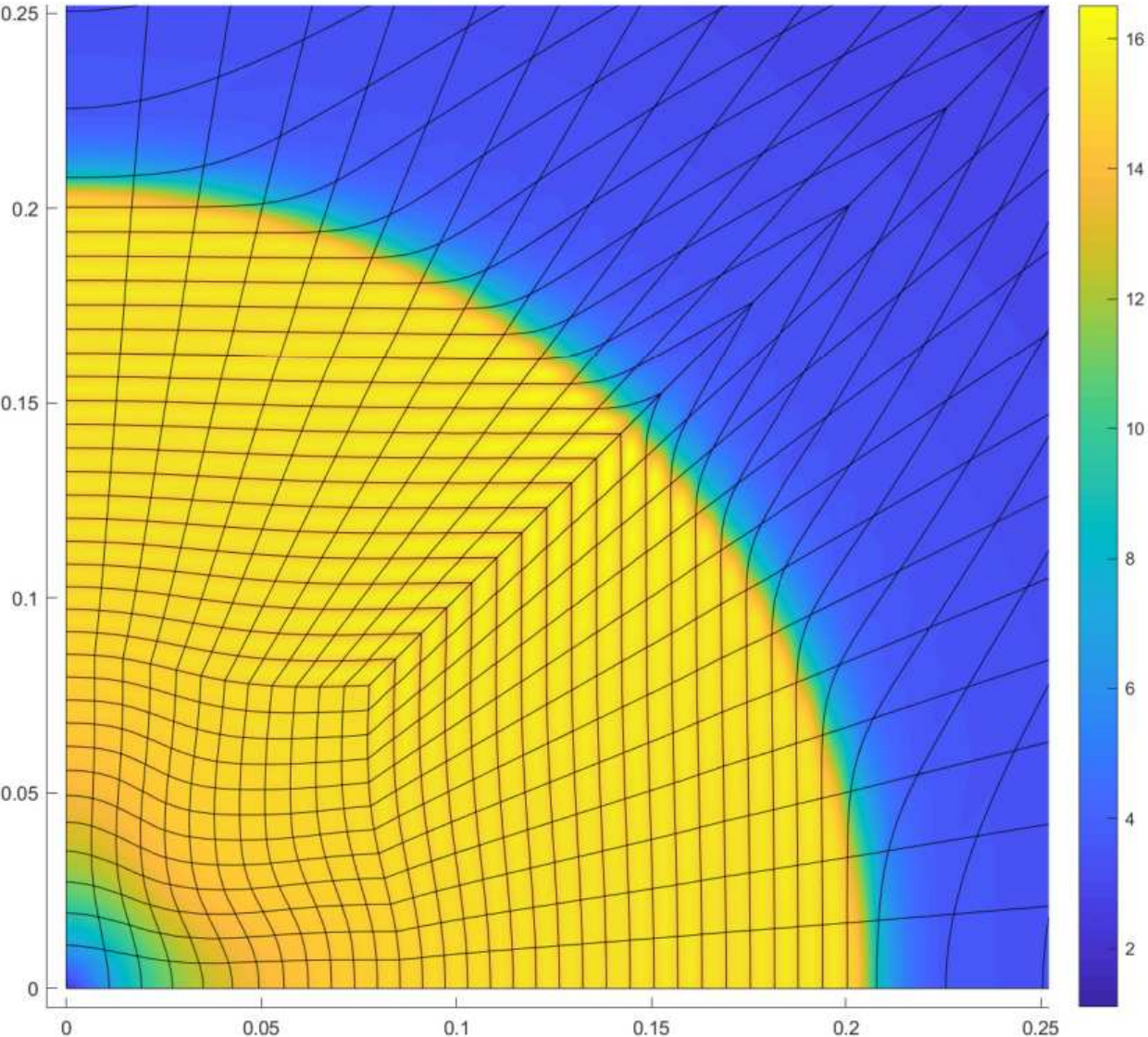}
    \caption{Density contour (without HG).}
    \label{Noh_Q3_butterfly_nohg}
  \end{subfigure}
  \hfill
  \begin{subfigure}[b]{0.32\textwidth}
    \includegraphics[width=\textwidth]{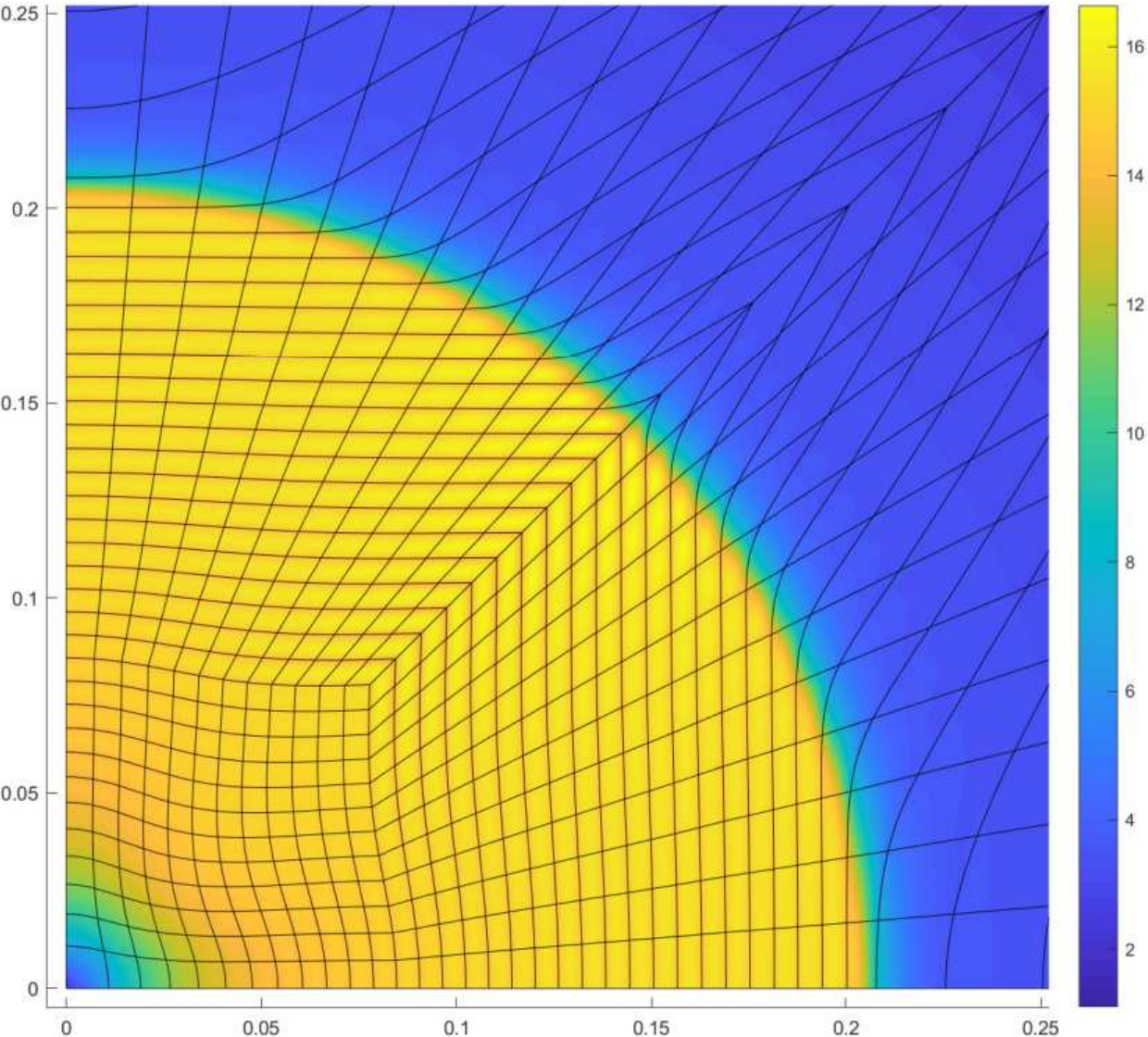}
    \caption{Density contour (with HG).}
    \label{Noh_Q3_butterfly_hg}
  \end{subfigure}
  \hfill
  \begin{subfigure}[b]{0.32\textwidth}
    \includegraphics[width=\textwidth]{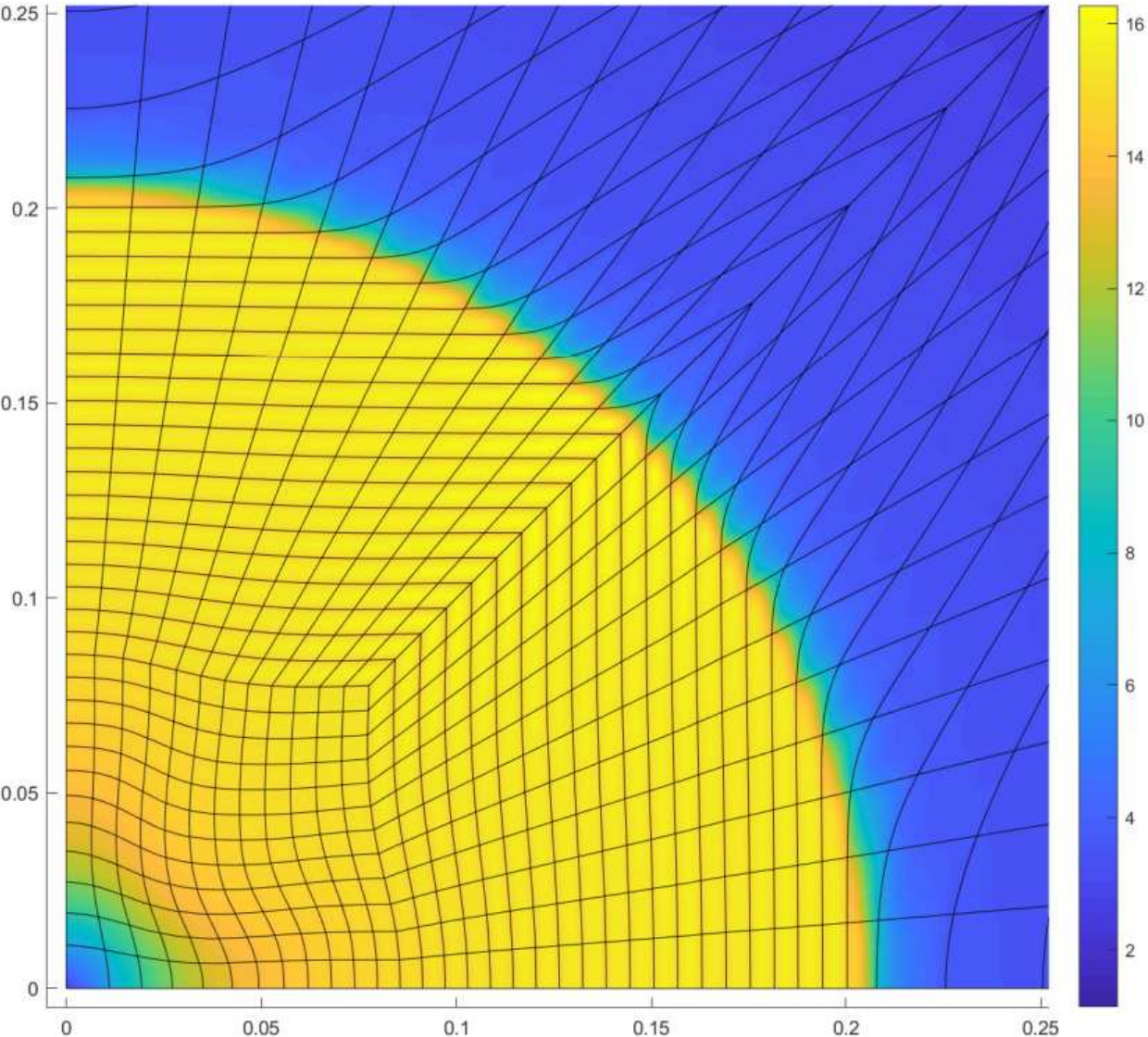}
    \caption{Density contour (Original method).}
    \label{Original_Noh_Q3_butterfly}
  \end{subfigure}

  \vspace{0.5em}

  \begin{subfigure}[b]{0.32\textwidth}
    \includegraphics[width=\textwidth]{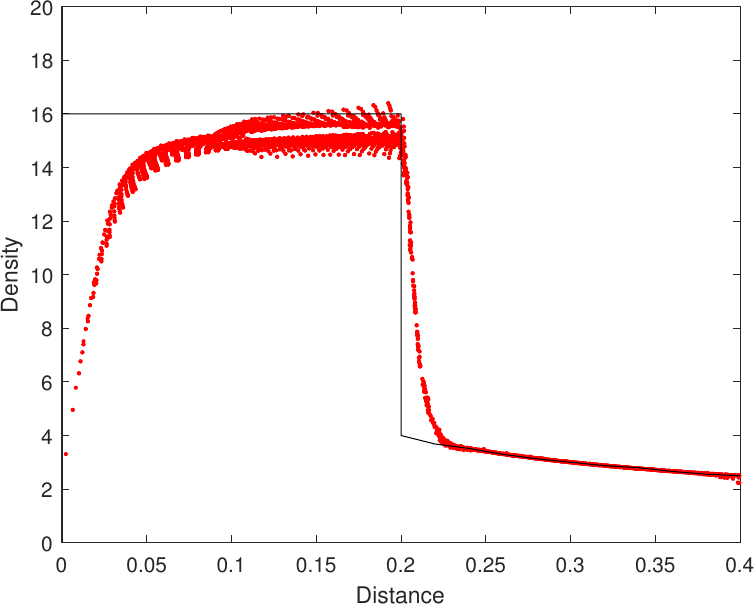}
    \caption{Pointwise density (without HG).}
    \label{Noh_Q3_butterfly_point_nohg}
  \end{subfigure}
  \hfill
  \begin{subfigure}[b]{0.32\textwidth}
    \includegraphics[width=\textwidth]{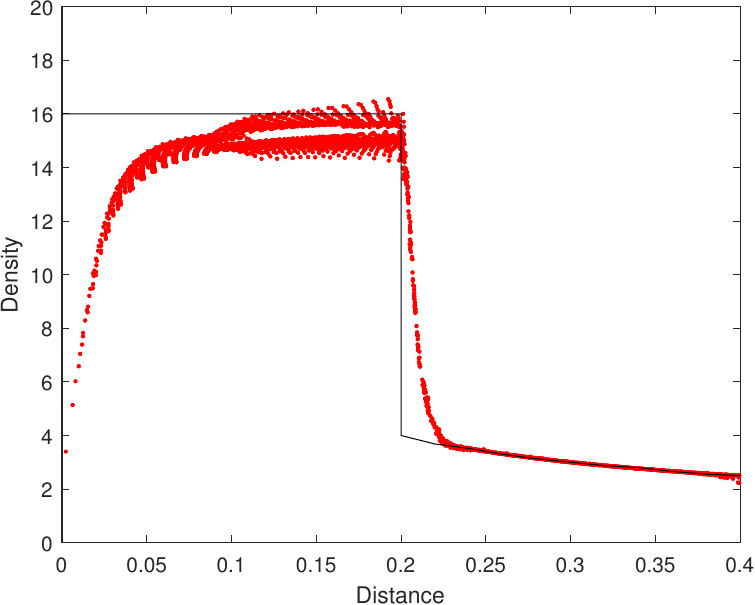}
    \caption{Pointwise density (with HG).}
    \label{Noh_Q3_butterfly_point.pdf}
  \end{subfigure}
  \hfill
  \begin{subfigure}[b]{0.32\textwidth}
    \includegraphics[width=\textwidth]{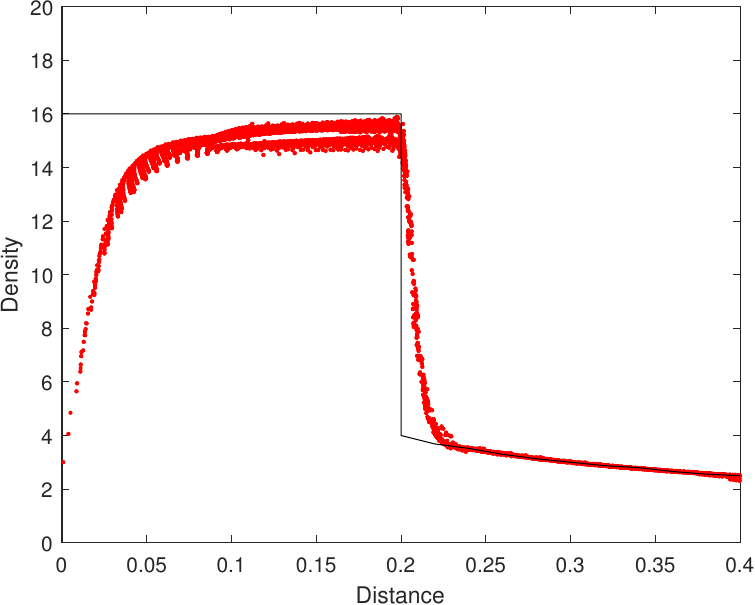}
    \caption{Pointwise density (Original method).}
    \label{Original_Noh_Q3_butterfly_point}
  \end{subfigure}
  \caption{Density contours and 1D pointwise density profiles along the radius for the Noh problem using the $Q^3-Q^2$ element pair on a butterfly mesh at $t = 0.6$.}
  \label{Noh_butterfly_Q3_Q2}
\end{figure}

Figure~\ref{Noh_butterfly_Q3_Q2} illustrates the results on the butterfly mesh using the $Q^3-Q^2$ element pair. The pointwise density profiles reveal that the Original method suppresses oscillations more effectively than our proposed framework, which is consistent with the observations from the uniform mesh case. Despite this slight difference in oscillation control, the grid evolution and structural integrity of the solution are remarkably similar for all three strategies, reaffirming that our methods perform at a highly comparable level.

\begin{figure}[htbp]
  \centering
  \begin{subfigure}[b]{0.32\textwidth}
    \includegraphics[width=\textwidth]{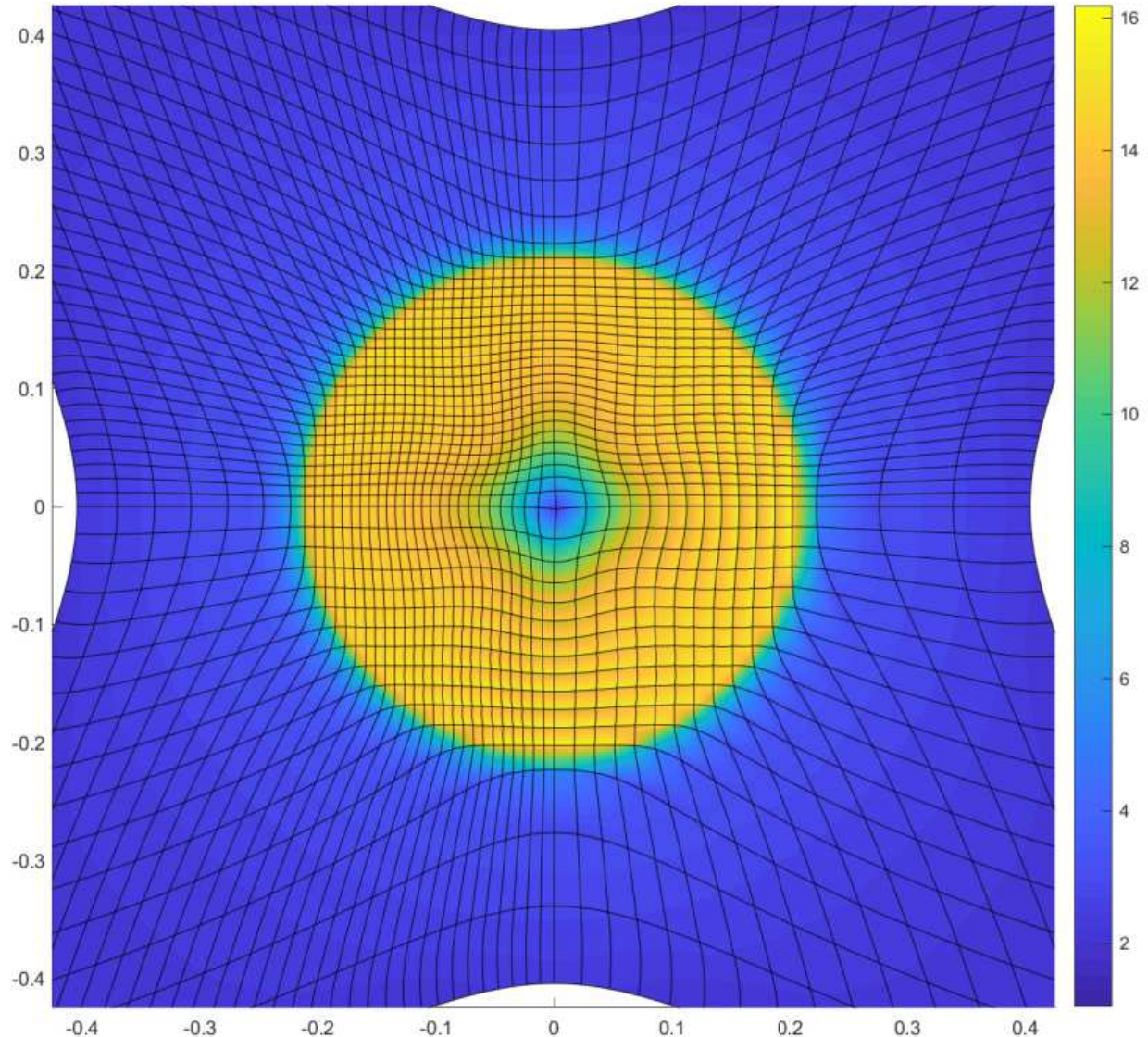}
    \caption{Density contour (without HG).}
    \label{Noh_Q2_asym_nohg}
  \end{subfigure}
  \hfill
  \begin{subfigure}[b]{0.32\textwidth}
    \includegraphics[width=\textwidth]{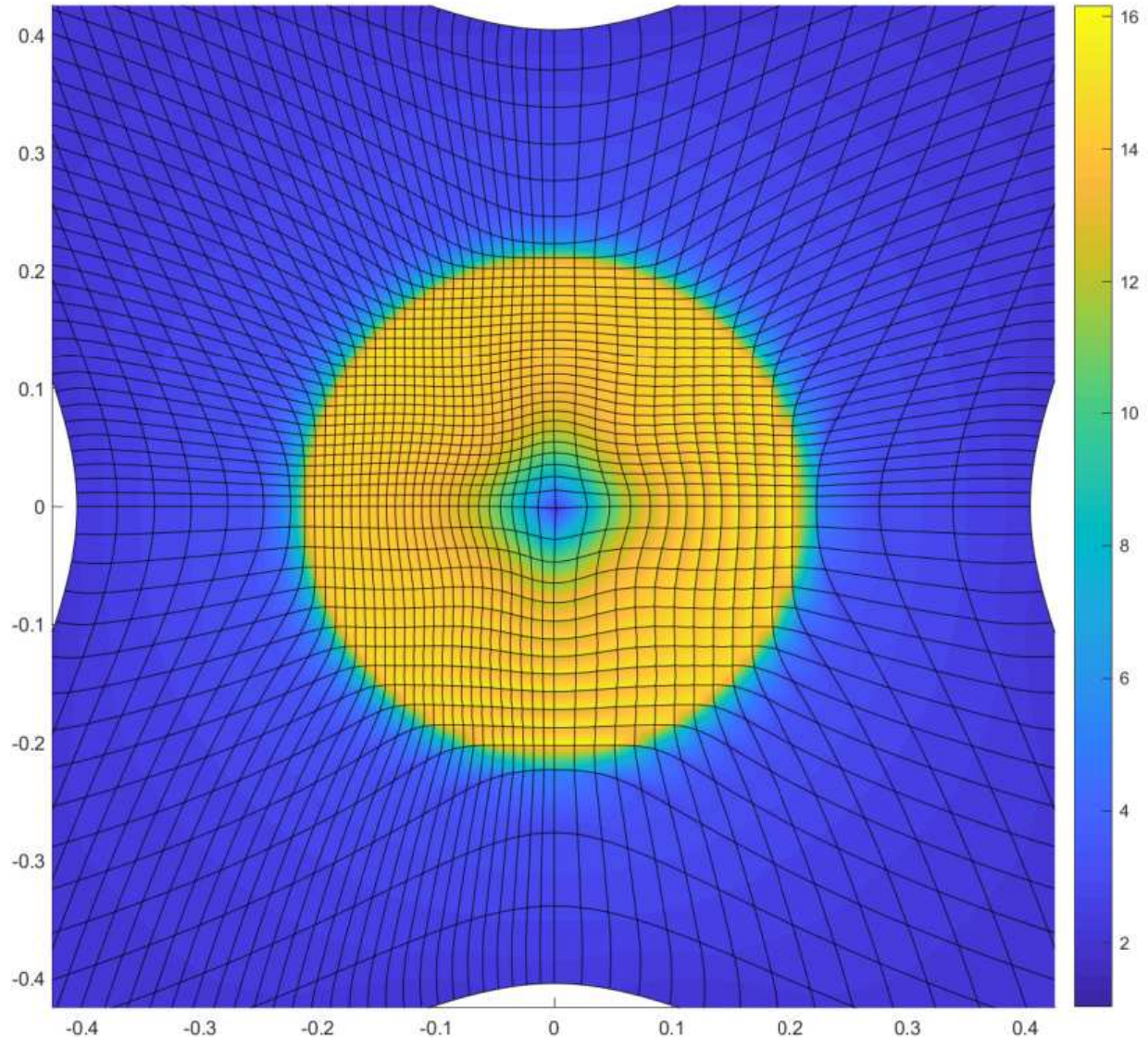}
    \caption{Density contour (with HG).}
    \label{Noh_Q2_asym_hg}
  \end{subfigure}
  \hfill
  \begin{subfigure}[b]{0.32\textwidth}
    \includegraphics[width=\textwidth]{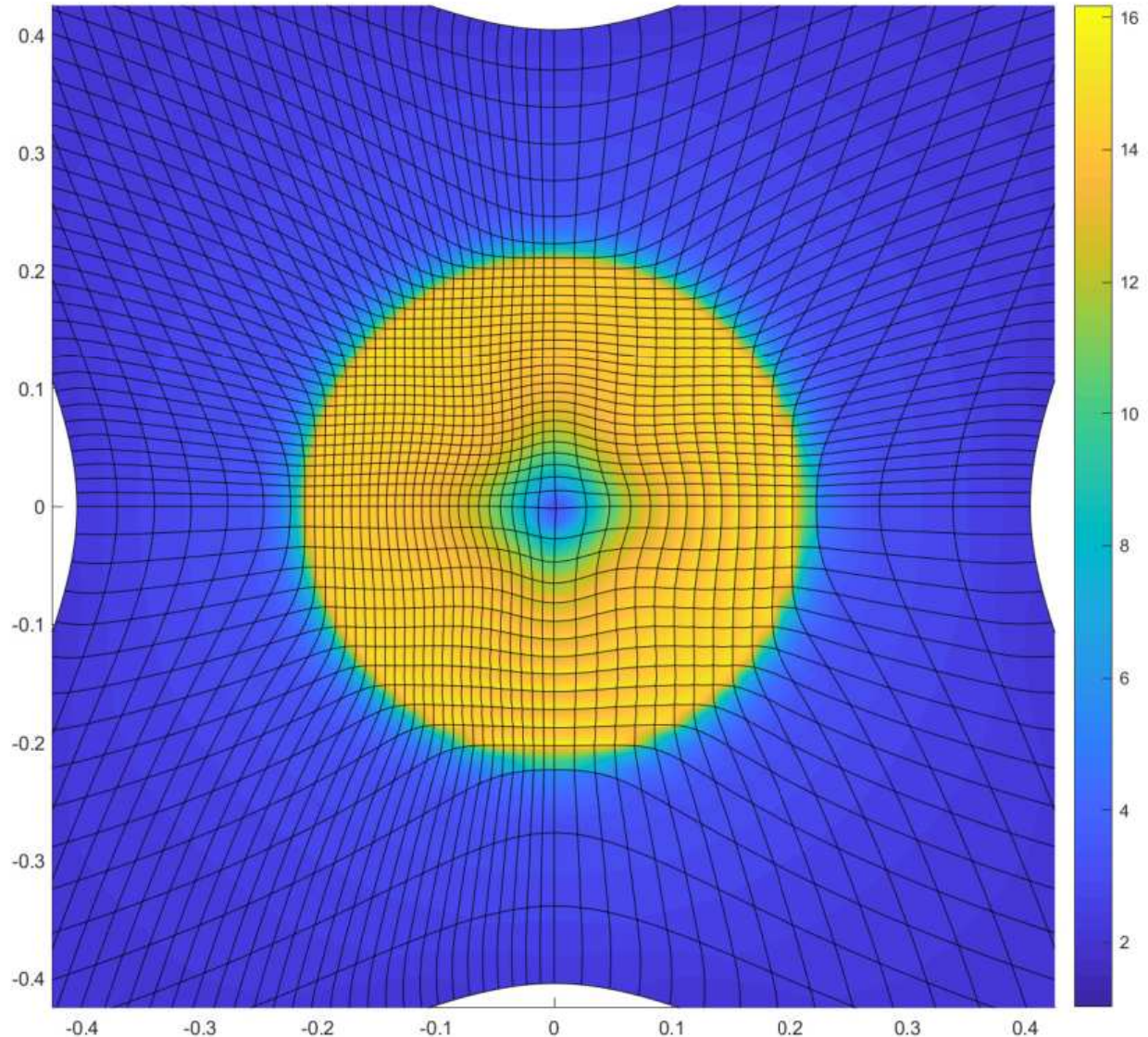}
    \caption{Density contour (Original method).}
    \label{Original_Noh_Q2_asym}
  \end{subfigure}

  \vspace{0.5em}

  \begin{subfigure}[b]{0.32\textwidth}
    \includegraphics[width=\textwidth]{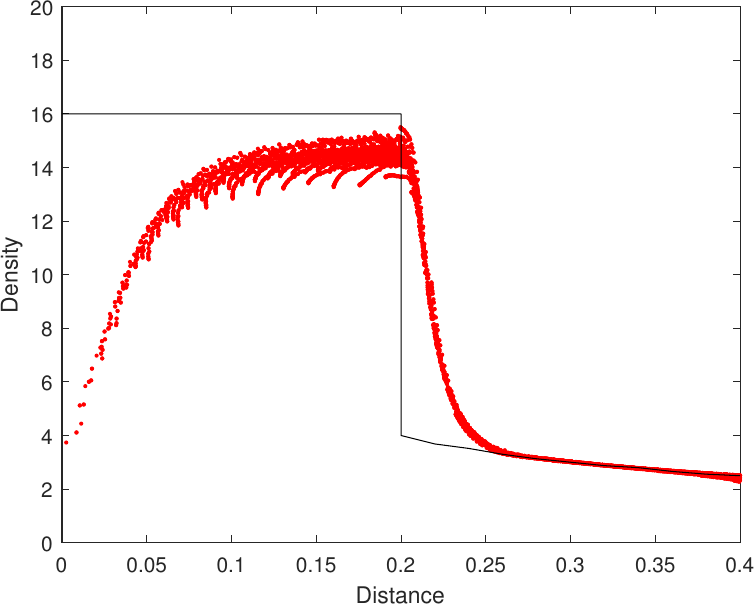}
    \caption{Pointwise density (without HG).}
    \label{Noh_Q2_asym_point_nohg}
  \end{subfigure}
  \hfill
  \begin{subfigure}[b]{0.32\textwidth}
    \includegraphics[width=\textwidth]{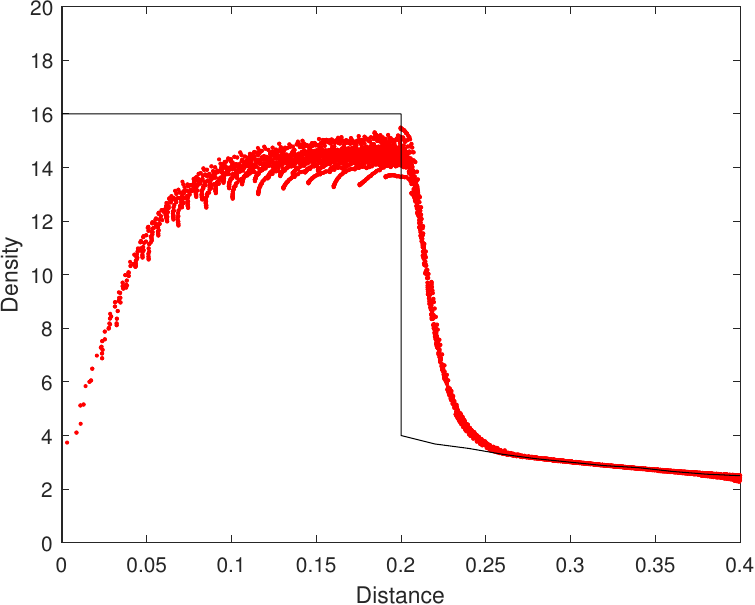}
    \caption{Pointwise density (with HG).}
    \label{Noh_Q2_asym_point.pdf}
  \end{subfigure}
  \hfill
  \begin{subfigure}[b]{0.32\textwidth}
    \includegraphics[width=\textwidth]{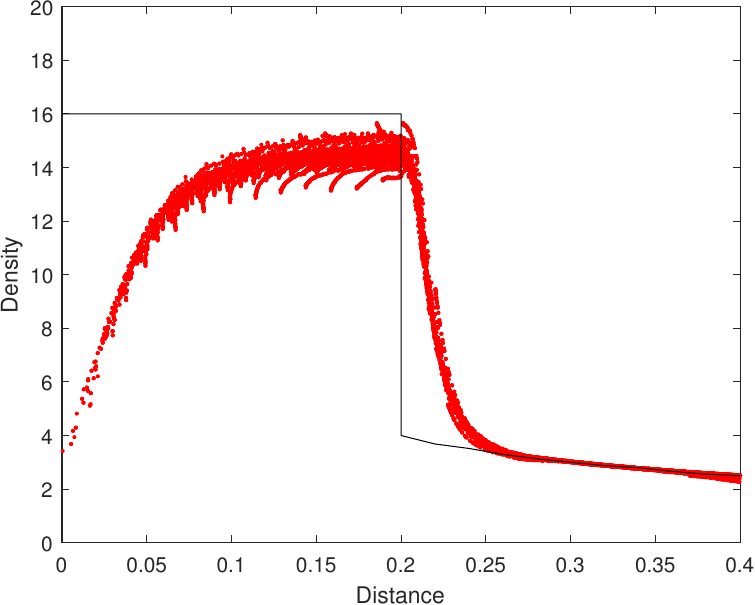}
    \caption{Pointwise density (Original method).}
    \label{Original_Noh_Q2_asym_point}
  \end{subfigure}
  \caption{Density contours and 1D pointwise density profiles along the radius for the Noh problem using the $Q^2-Q^1$ element pair on an asymmetric mesh at $t = 0.6$.}
  \label{Noh_asym_Q2_Q1}
\end{figure}

Figure~\ref{Noh_asym_Q2_Q1} displays the numerical results for the $Q^2-Q^1$ element pair on an asymmetric mesh with varying aspect ratios. For this demanding setup, we increased the artificial viscosity coefficients ($c_1=0.8$ and $c_2=0.5$) to adequately suppress oscillations. Consequently, the pointwise density profiles reflect the smearing effects of this increased viscosity. In this case, the Original method does not demonstrate any distinct advantage over our proposed framework; the overall numerical performance, oscillation suppression, and shock-capturing ability of all three methods are entirely comparable, and the resulting grid quality remains highly satisfactory.

\begin{figure}[htbp]
  \centering
  \begin{subfigure}[b]{0.32\textwidth}
    \includegraphics[width=\textwidth]{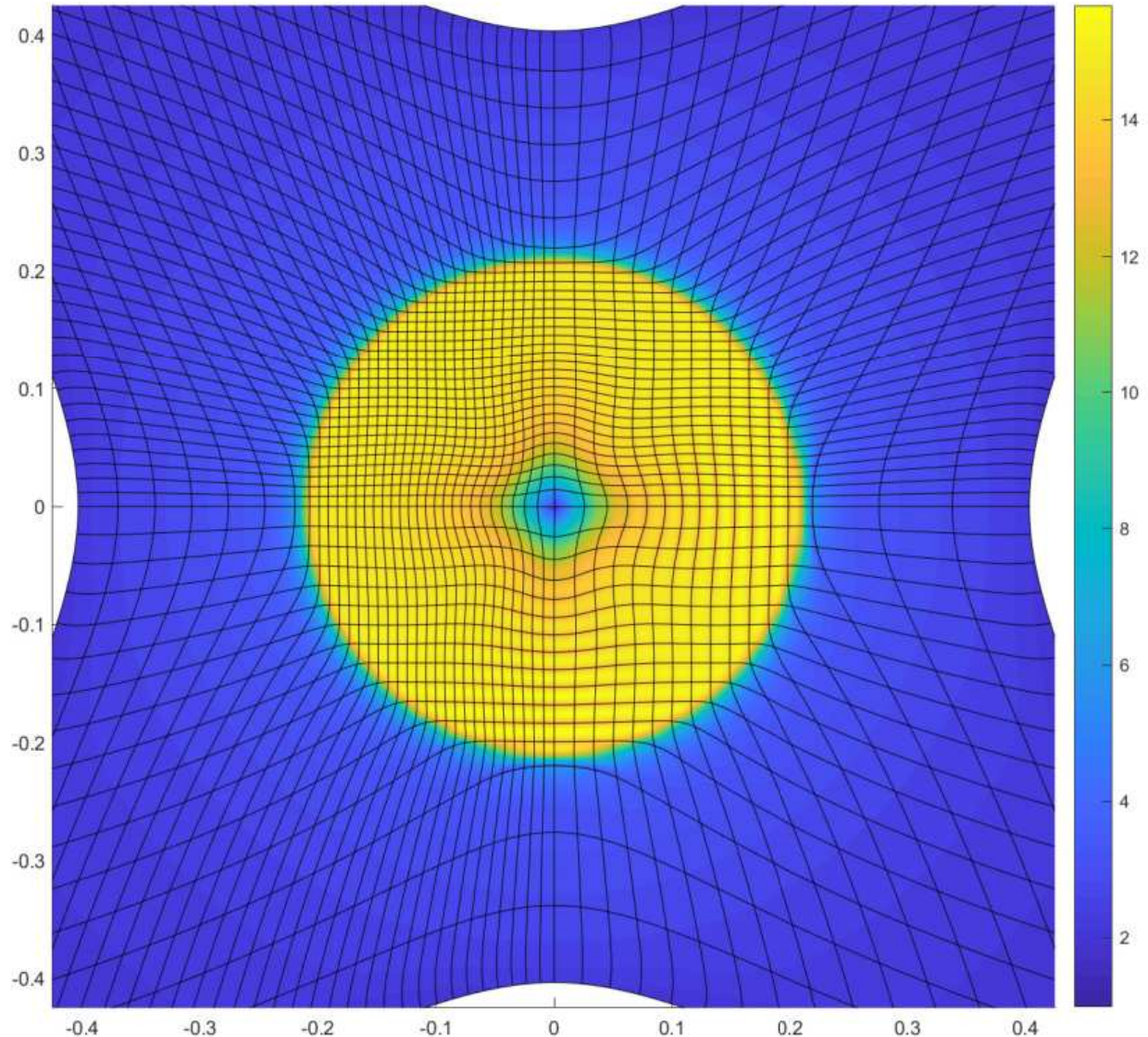}
    \caption{Density contour (without HG).}
    \label{Noh_Q3_asym_nohg}
  \end{subfigure}
  \hfill
  \begin{subfigure}[b]{0.32\textwidth}
    \includegraphics[width=\textwidth]{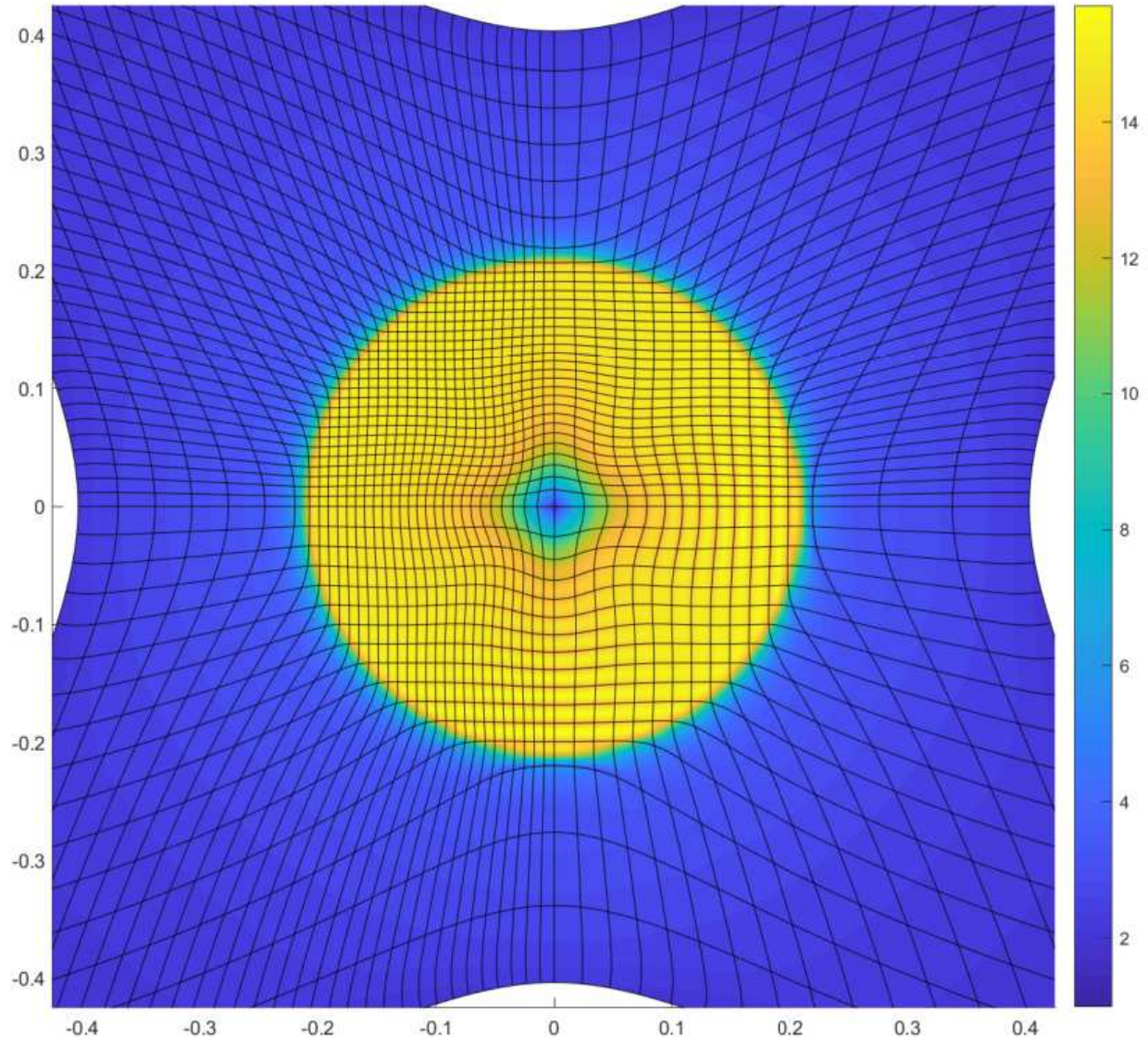}
    \caption{Density contour (with HG).}
    \label{Noh_Q3_asym_hg}
  \end{subfigure}
  \hfill
  \begin{subfigure}[b]{0.32\textwidth}
    \includegraphics[width=\textwidth]{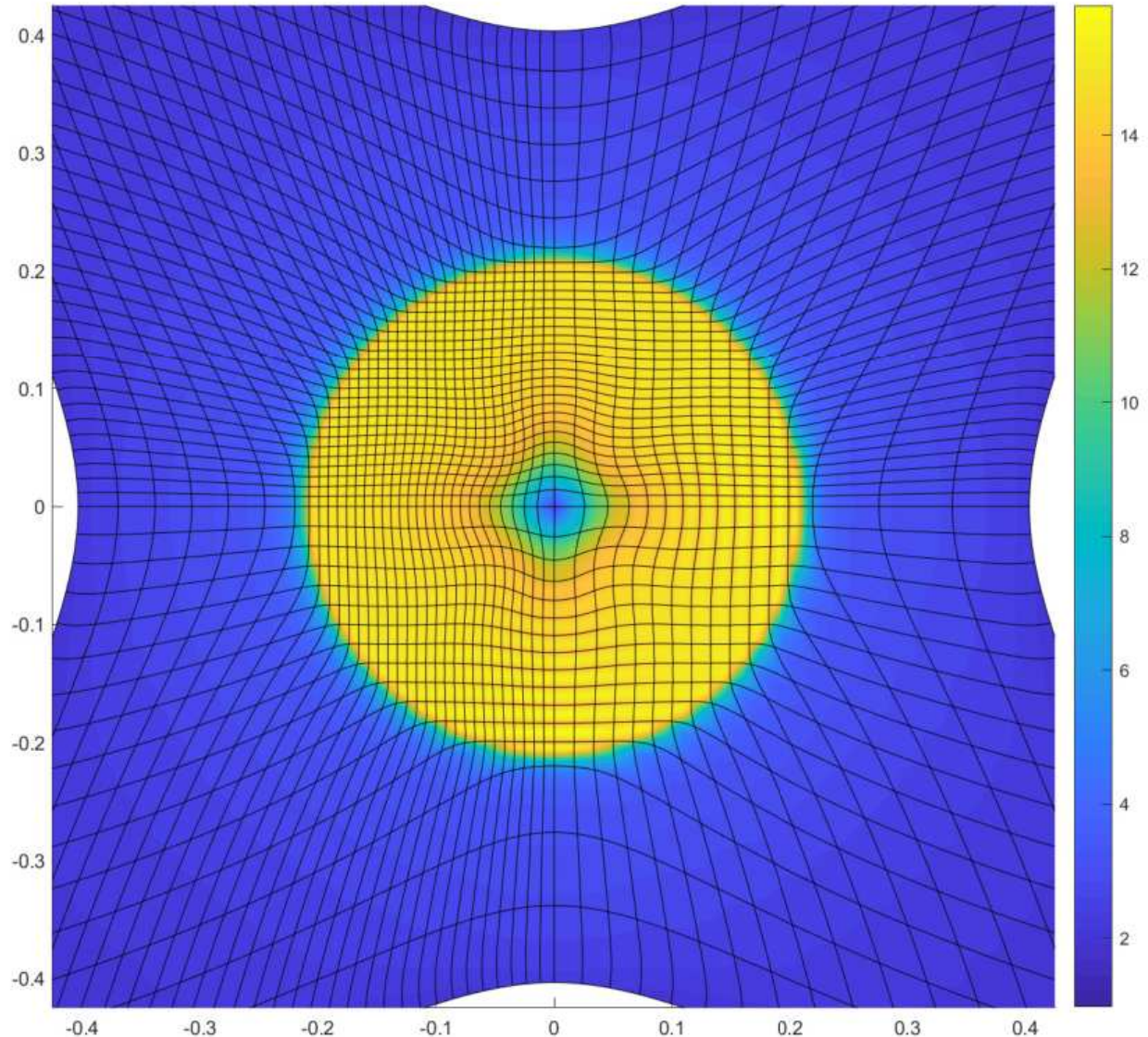}
    \caption{Density contour (Original method).}
    \label{Original_Noh_Q3_asym}
  \end{subfigure}

  \vspace{0.5em}

  \begin{subfigure}[b]{0.32\textwidth}
    \includegraphics[width=\textwidth]{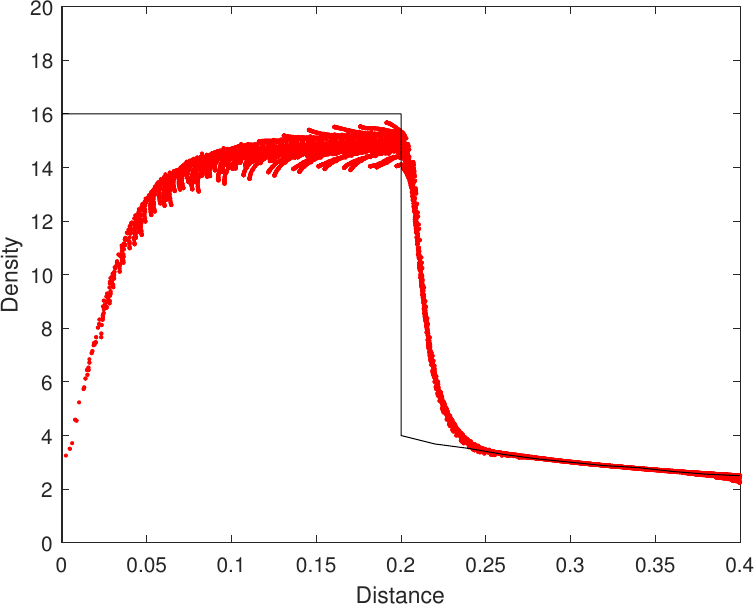}
    \caption{Pointwise density (without HG).}
    \label{Noh_Q3_asym_point_nohg}
  \end{subfigure}
  \hfill
  \begin{subfigure}[b]{0.32\textwidth}
    \includegraphics[width=\textwidth]{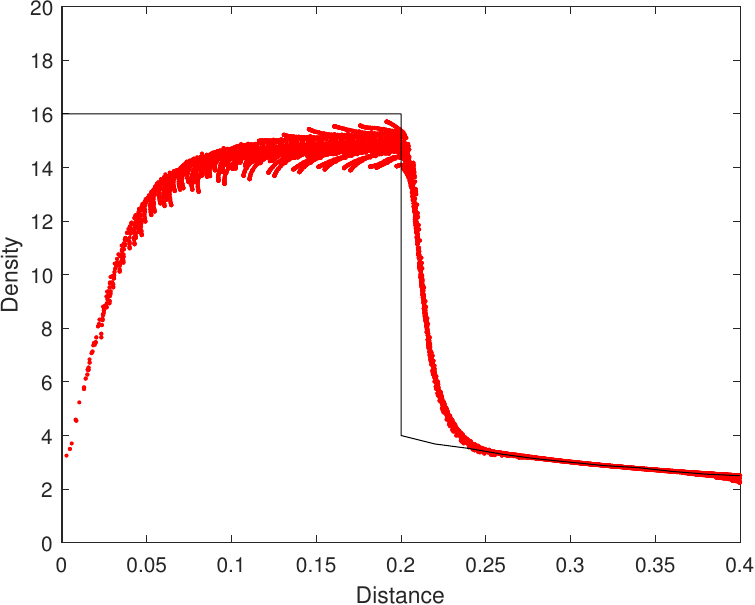}
    \caption{Pointwise density (with HG).}
    \label{Noh_Q3_asym_point.pdf}
  \end{subfigure}
  \hfill
  \begin{subfigure}[b]{0.32\textwidth}
    \includegraphics[width=\textwidth]{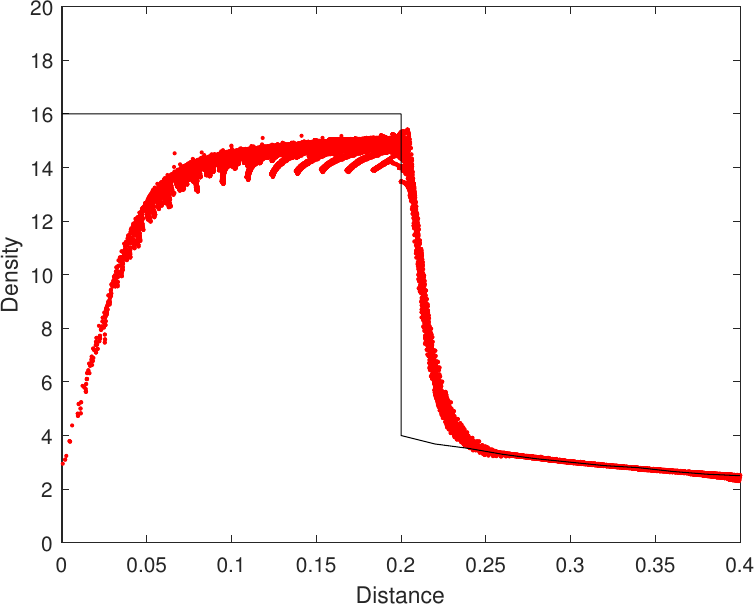}
    \caption{Pointwise density (Original method).}
    \label{Original_Noh_Q3_asym_point}
  \end{subfigure}
  \caption{Density contours and 1D pointwise density profiles along the radius for the Noh problem using the $Q^3-Q^2$ element pair on an asymmetric mesh at $t = 0.6$.}
  \label{Noh_asym_Q3_Q2}
\end{figure}

Finally, Figure~\ref{Noh_asym_Q3_Q2} demonstrates the results for the $Q^3-Q^2$ element pair on the asymmetric mesh. As with the $Q^2-Q^1$ case, the artificial viscosity coefficients were increased ($c_1=0.8$ and $c_2=0.5$). The pointwise density profiles follow the same trend seen in the previous configurations: while the Original method is slightly more effective at smoothing out unphysical oscillations, our proposed framework successfully resolves the shock wave without breaking down. Across all three strategies, the final grid distributions remain robust and fully comparable, proving the viability of our approach even under skewed initial conditions.

\subsubsection{Sedov point-blast problem} 
The Sedov problem models a strong blast wave in planar geometry, driven by a point source of internal energy deposited at the origin. The blast wave reaches a radius of $1$ at $t = 1$. The background internal energy is set to $1 \times 10^{-10}$, except within the source elements. The computational domain is either $[0,1.2]^2$ (for a corner source) or $[-1.2,1.2]^2$ (for a central source), with $\gamma = 5/3$ and an initial density $\rho = 1$. The mesh size is $h = 0.05$, and both the linear and quadratic artificial viscosity parameters are set to $0.5$.

Similar to the Noh problem, we evaluate three sets of initial mesh configurations for the Sedov problem: uniform, butterfly, and asymmetric meshes. The initialization for the asymmetric mesh differs slightly from~\cite{Dobrev2012High}. To preserve the perfect symmetry of the Sedov blast wave on an asymmetric mesh, the internal energy is initialized such that the four quadrants surrounding the origin receive an identical amount of initial internal energy. Specifically, in the first quadrant, two elements along the $y$-axis are assigned non-zero energy; in the third quadrant, two elements along the $x$-axis are set similarly. In the second quadrant, four elements near the origin are initialized with non-zero energy, while in the fourth quadrant, only the single element containing the origin is excited. We present numerical results comparing our proposed framework (with and without hourglass control) against the Original method. For the Sedov problem, the physical vorticity coefficient $c_{vor}$ is activated and calculated by definition. The artificial viscosity coefficients are set to $c_1=0.3$ and $c_2=0.5$, and the hourglass parameter is tuned case by case to achieve optimal results.

\begin{figure}[htbp]
  \centering
  \begin{subfigure}[b]{0.32\textwidth}
    \includegraphics[width=\textwidth]{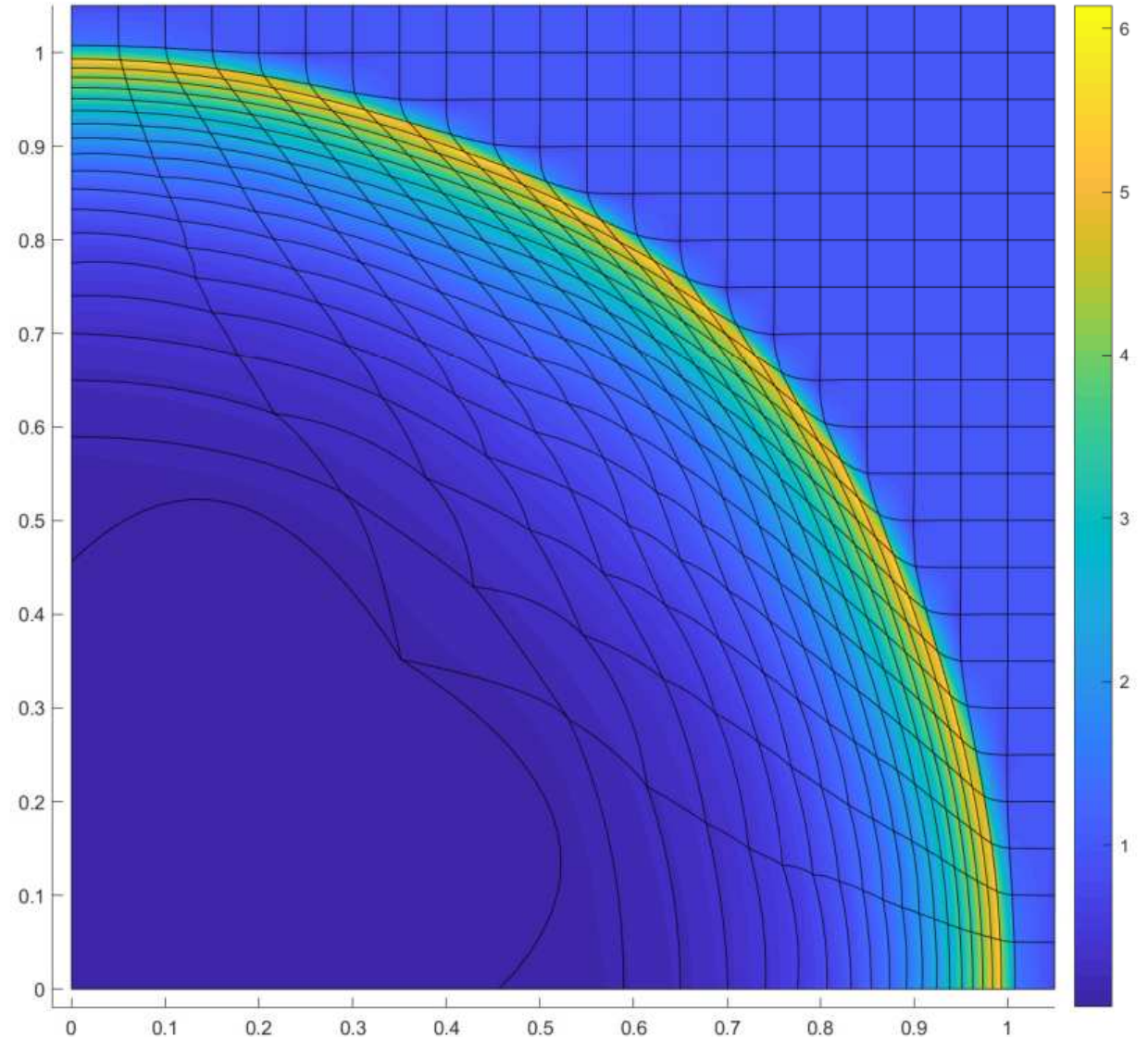}
    \caption{Density contour (without HG).}
    \label{Sedov_Q2_nohg}
  \end{subfigure}
  \hfill
  \begin{subfigure}[b]{0.32\textwidth}
    \includegraphics[width=\textwidth]{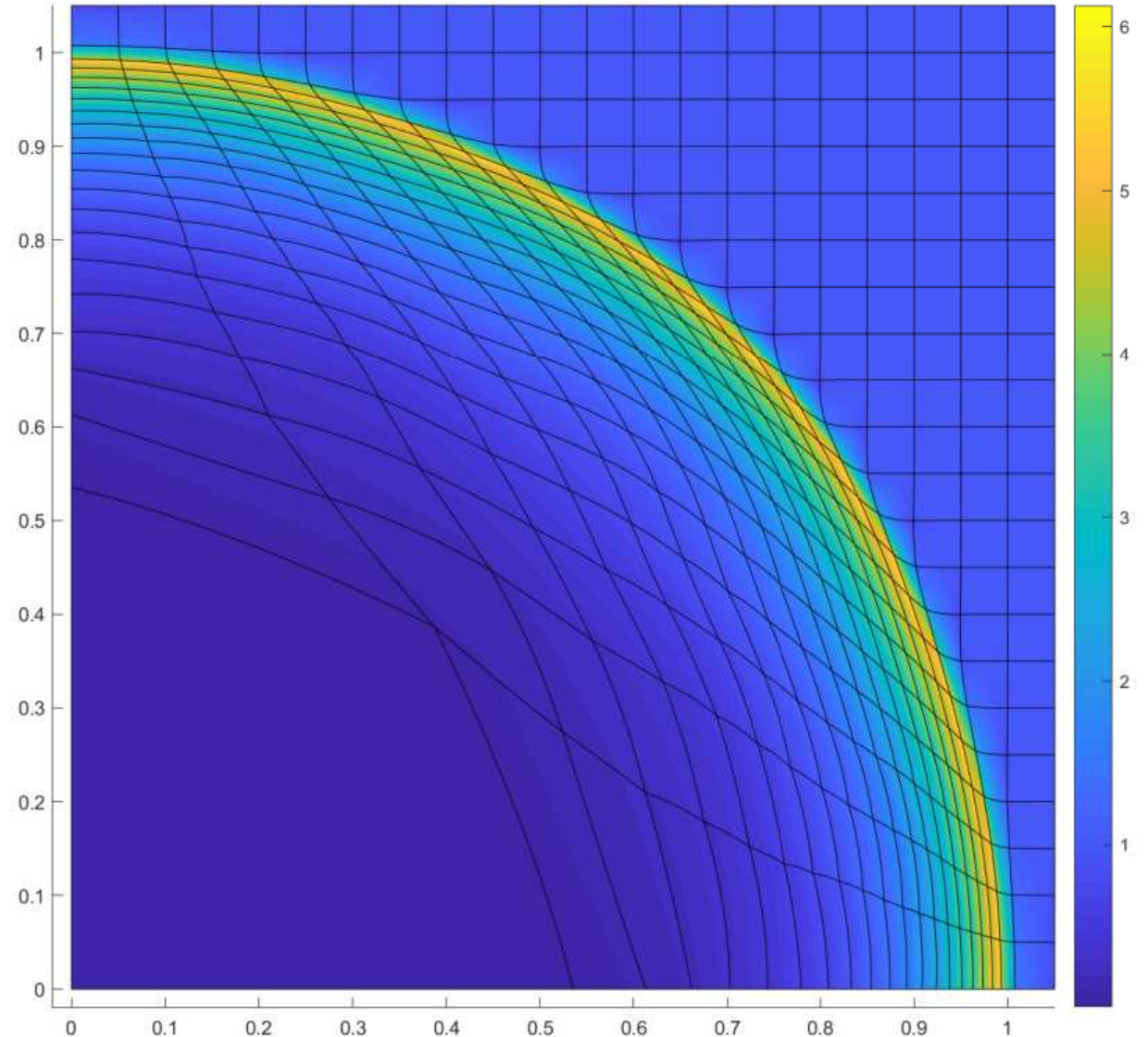}
    \caption{Density contour (with HG).}
    \label{Sedov_Q2_hg}
  \end{subfigure}
  \hfill
  \begin{subfigure}[b]{0.32\textwidth}
    \includegraphics[width=\textwidth]{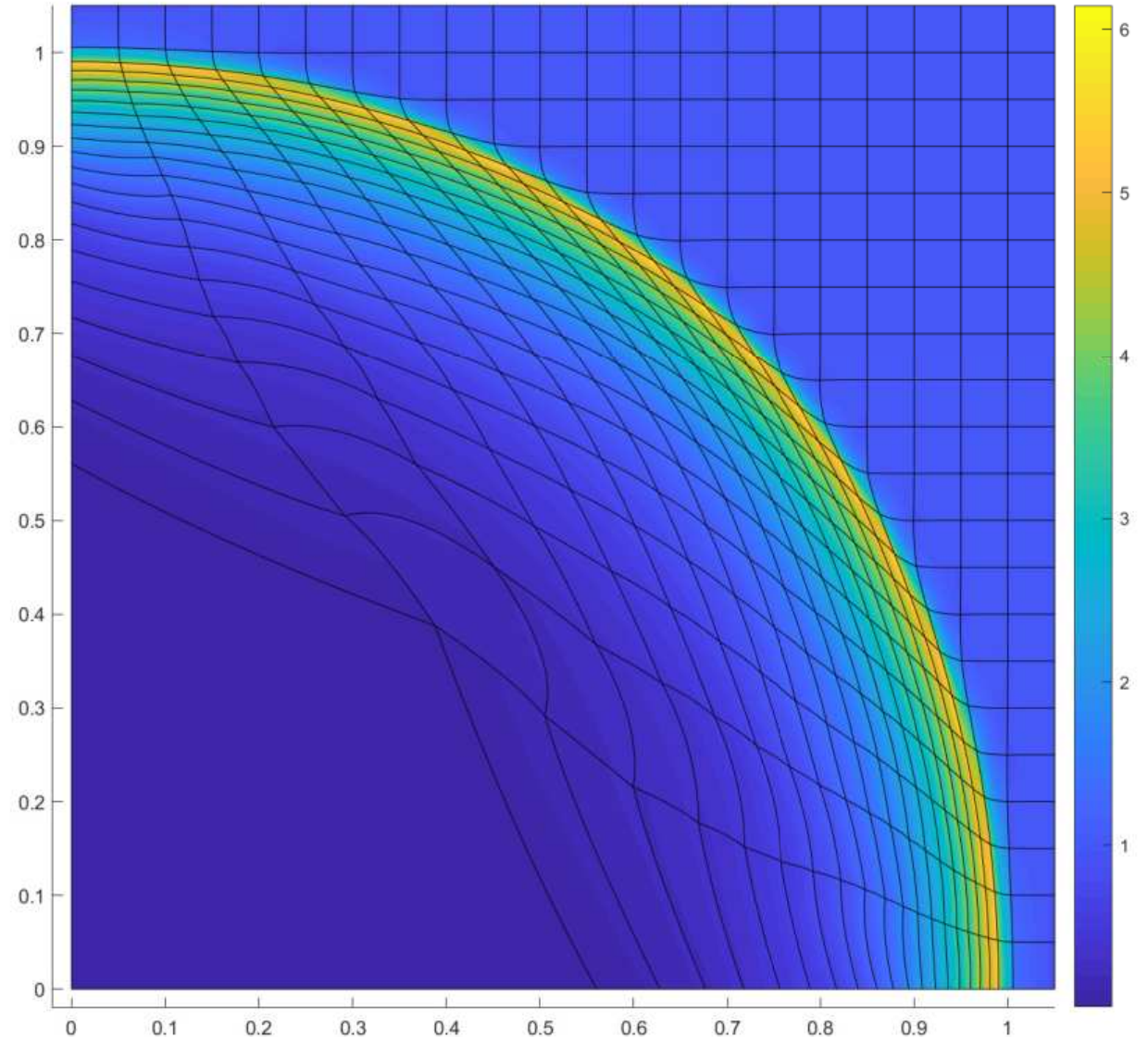}
    \caption{Density contour (Original method).}
    \label{Original_Sedov_Q2}
  \end{subfigure}

  \vspace{0.5em}

  \begin{subfigure}[b]{0.32\textwidth}
    \includegraphics[width=\textwidth]{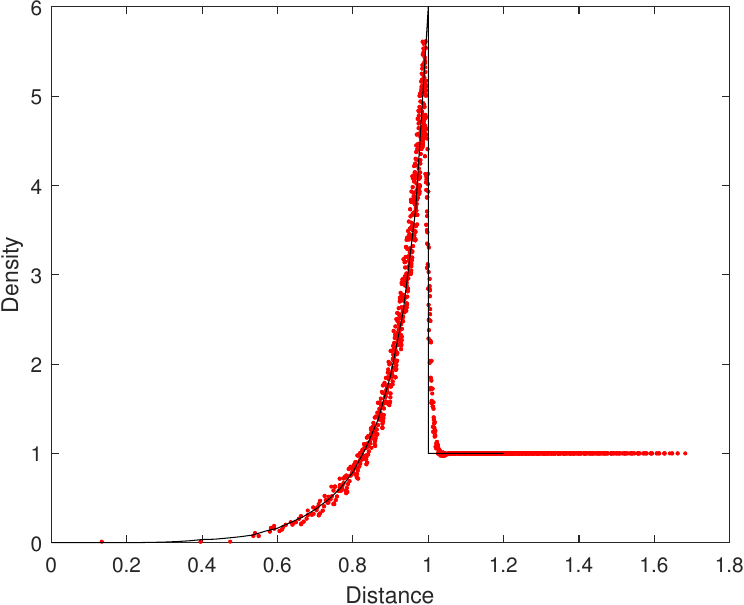}
    \caption{Pointwise density (without HG).}
    \label{Sedov_Q2_point_nohg}
  \end{subfigure}
  \hfill
  \begin{subfigure}[b]{0.32\textwidth}
    \includegraphics[width=\textwidth]{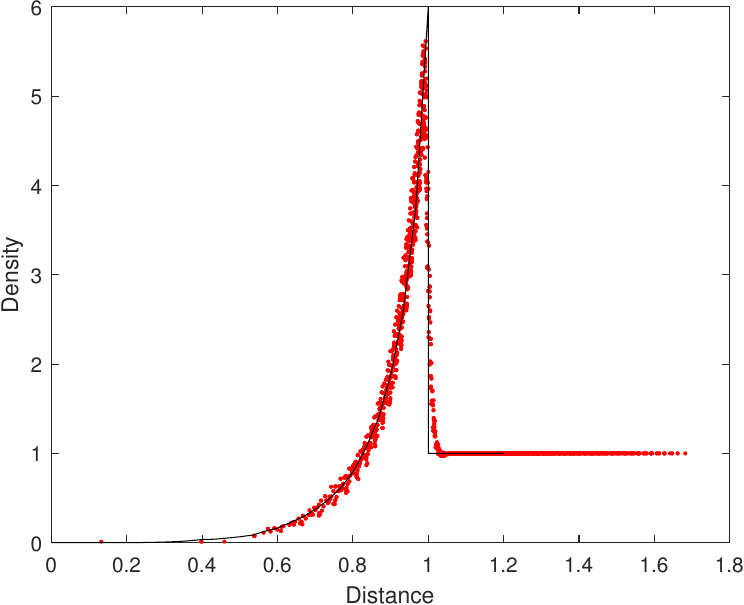}
    \caption{Pointwise density (with HG).}
    \label{Sedov_Q2_point}
  \end{subfigure}
  \hfill
  \begin{subfigure}[b]{0.32\textwidth}
    \includegraphics[width=\textwidth]{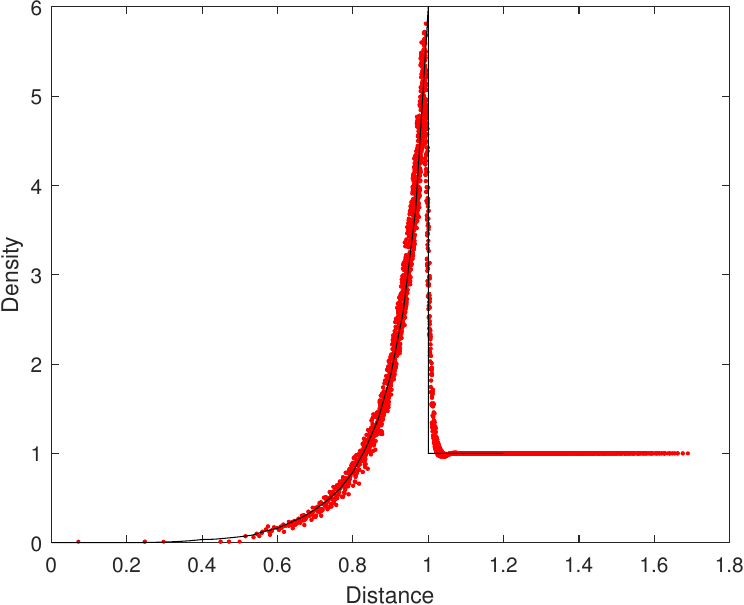}
    \caption{Pointwise density (Original method).}
    \label{Original_Sedov_Q2_point}
  \end{subfigure}
  \caption{Density contours and 1D pointwise density profiles along the radius for the Sedov problem using the $Q^2-Q^1$ element pair on a uniform mesh at $t = 1$.}
  \label{Sedov_Q2_Q1}
\end{figure}

Figure~\ref{Sedov_Q2_Q1} demonstrates the numerical results of the Sedov problem using the $Q^2-Q^1$ element pair on a uniform mesh at $t = 1$. The hourglass parameter is set to $0.1$ for Figures~\ref{Sedov_Q2_hg} and~\ref{Sedov_Q2_point}. The three strategies produce visually distinct grid distributions. Without hourglass control, Figure~\ref{Sedov_Q2_nohg} shows severe grid distortion at the origin. This distortion is effectively mitigated when our hourglass control is applied, as seen in Figure~\ref{Sedov_Q2_hg}. Conversely, the Original method (Figure~\ref{Original_Sedov_Q2}) exhibits unsatisfactory grid quality because the numerical scheme is overly stiff. This excessive stiffness is also evident in the pointwise density profile (Figure~\ref{Original_Sedov_Q2_point}). As explained in the Noh problem, controlling oscillations without introducing excessive stiffness is critical in Lagrangian simulations. Although the pointwise density profile of the Original method remains close to the analytical solution, it forces the grid to become overly rigid at the shock discontinuity. In contrast, by choosing an appropriate hourglass parameter, our proposed framework successfully balances the suppression of hourglass deformation against numerical stiffness, providing an optimal solution. (Note that if a larger hourglass parameter is chosen, our grid distribution becomes identically stiff to that of the Original method, highlighting the importance of having a tunable parameter.)

\begin{figure}[htbp]
  \centering
  \begin{subfigure}[b]{0.32\textwidth}
    \includegraphics[width=\textwidth]{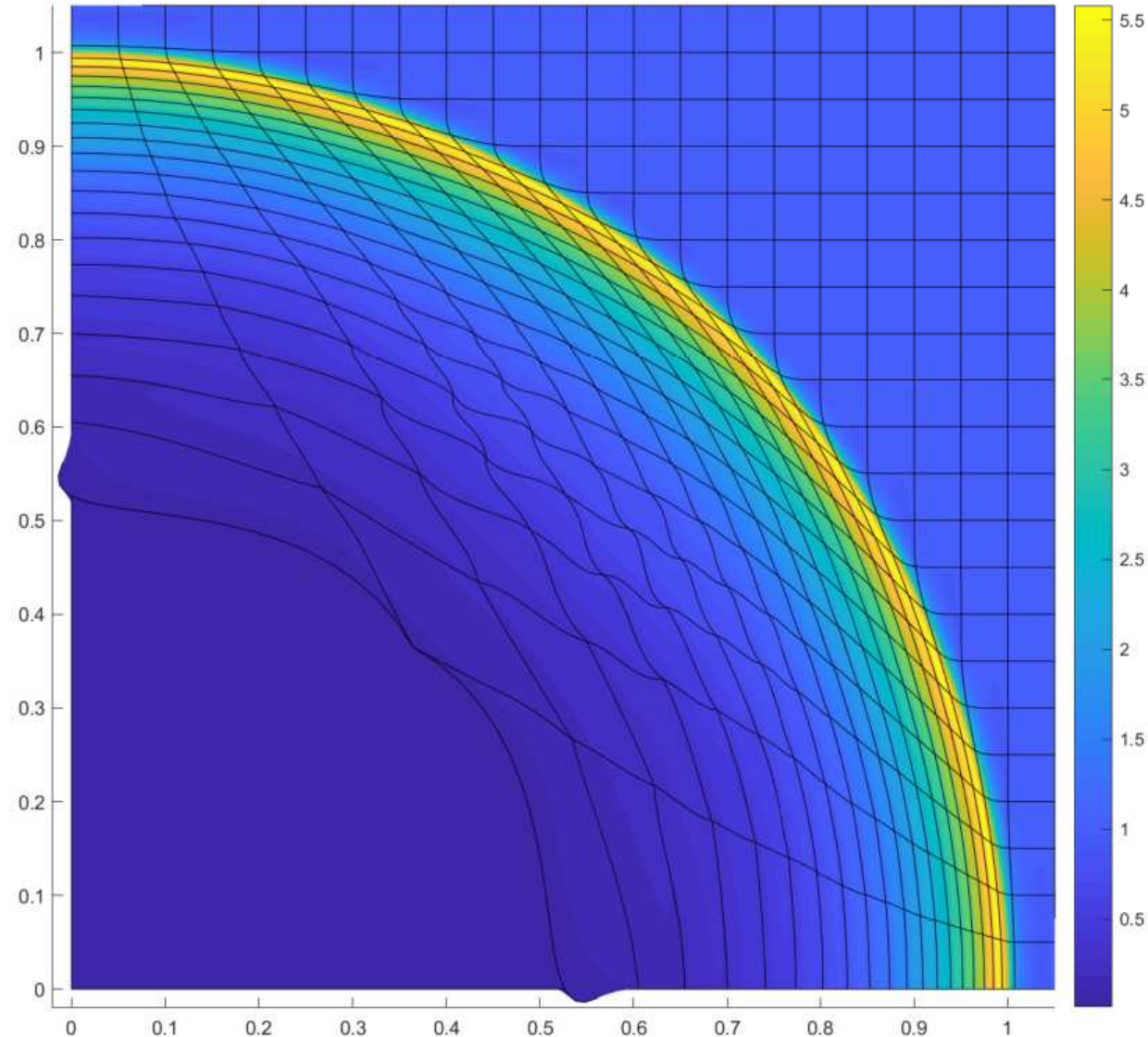}
    \caption{Density contour (without HG).}
    \label{Sedov_Q3_nohg}
  \end{subfigure}
  \hfill
  \begin{subfigure}[b]{0.32\textwidth}
    \includegraphics[width=\textwidth]{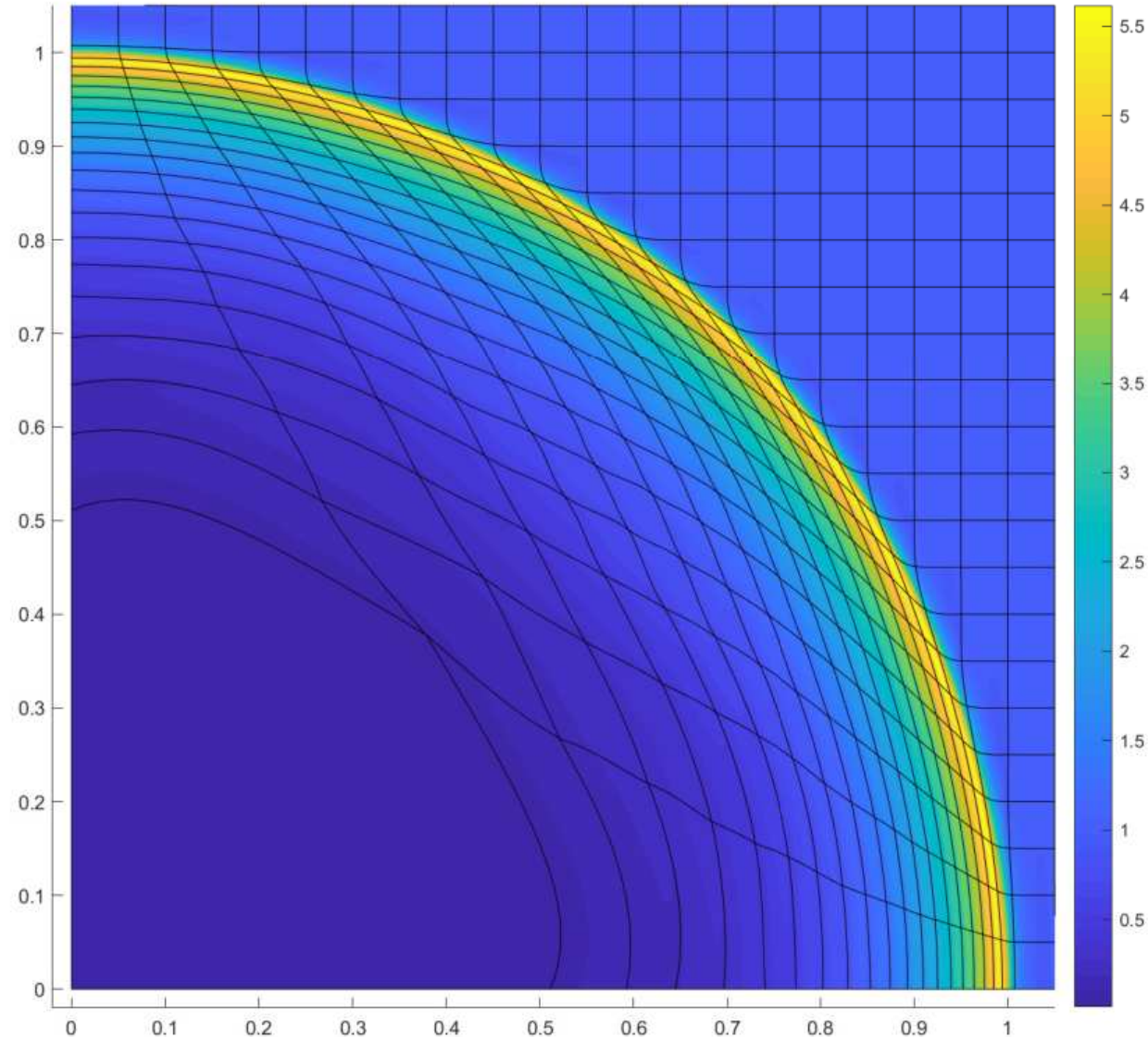}
    \caption{Density contour (with HG).}
    \label{Sedov_Q3_hg}
  \end{subfigure}
  \hfill
  \begin{subfigure}[b]{0.32\textwidth}
    \includegraphics[width=\textwidth]{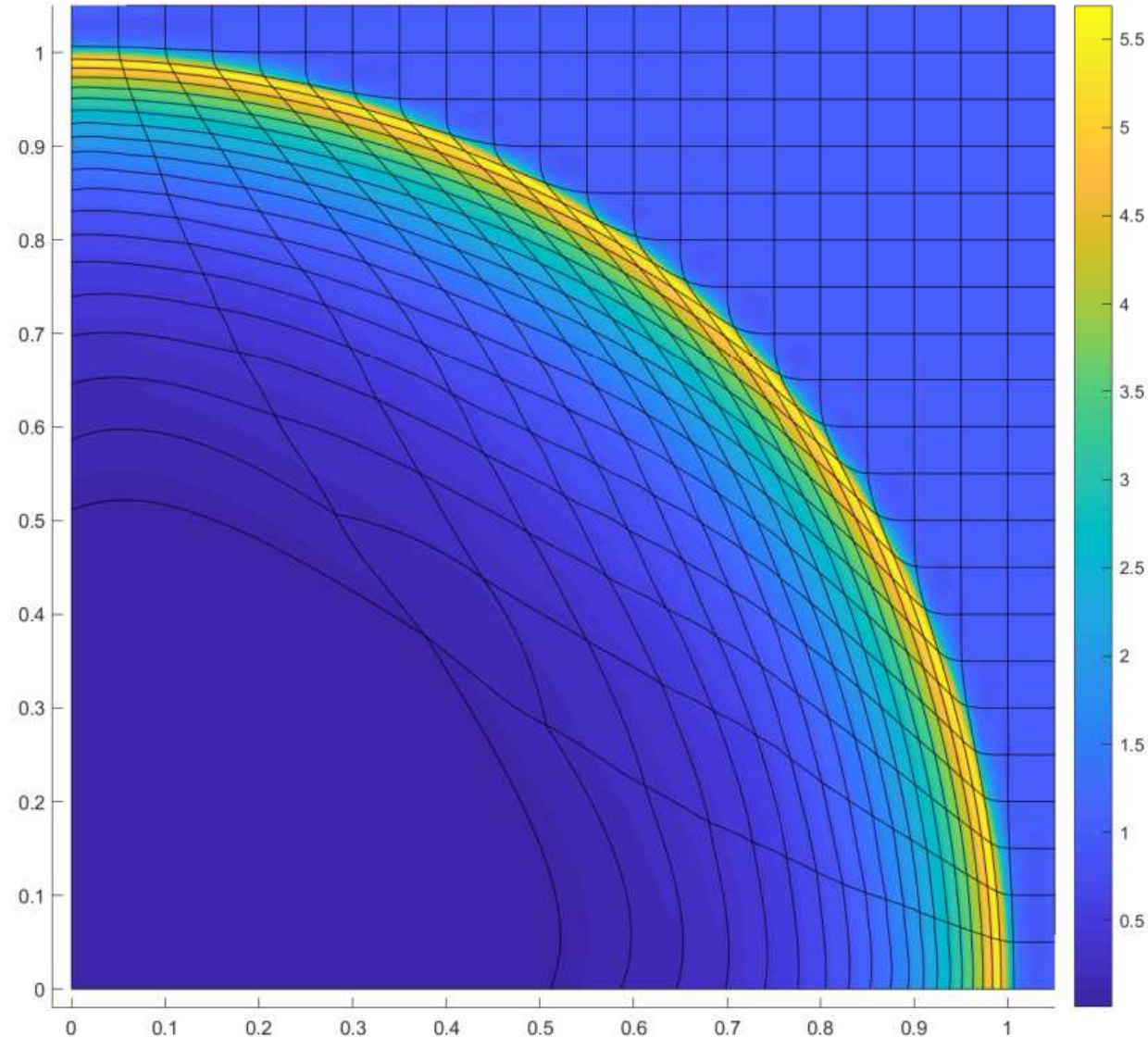}
    \caption{Density contour (Original method).}
    \label{Original_Sedov_Q3}
  \end{subfigure}

  \vspace{0.5em}

  \begin{subfigure}[b]{0.32\textwidth}
    \includegraphics[width=\textwidth]{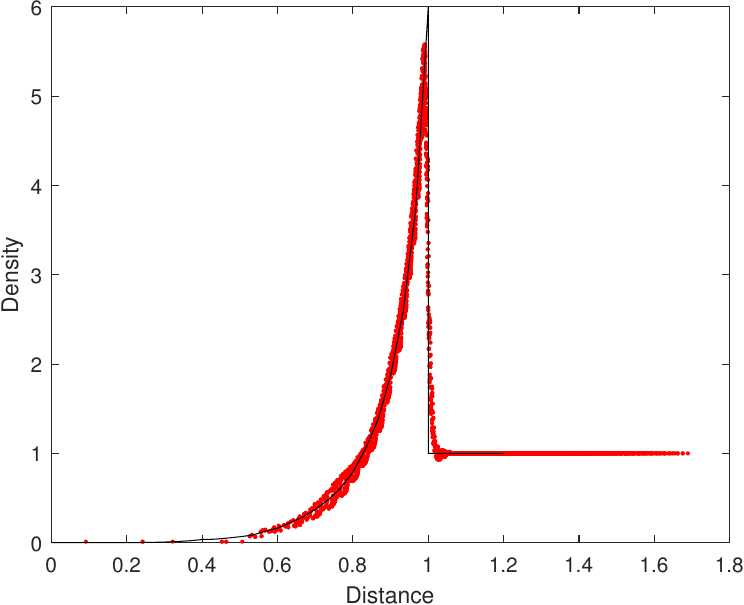}
    \caption{Pointwise density (without HG).}
    \label{Sedov_Q3_point_nohg}
  \end{subfigure}
  \hfill
  \begin{subfigure}[b]{0.32\textwidth}
    \includegraphics[width=\textwidth]{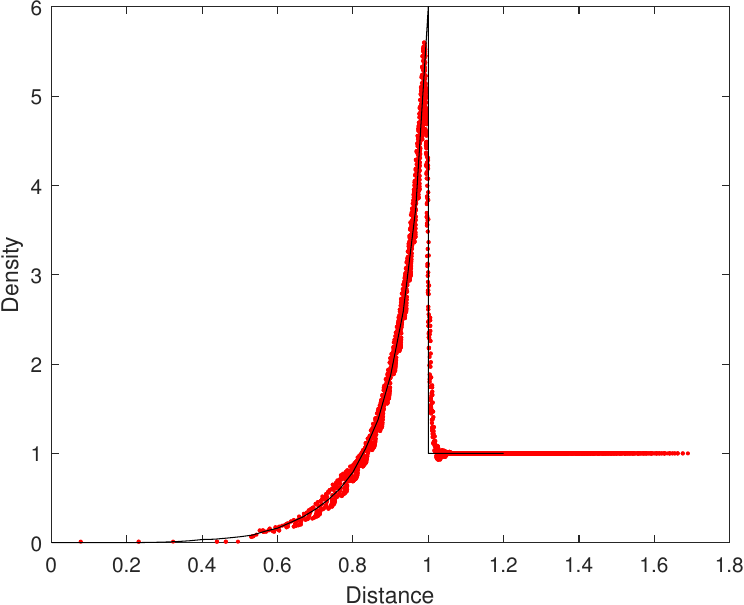}
    \caption{Pointwise density (with HG).}
    \label{Sedov_Q3_point}
  \end{subfigure}
  \hfill
  \begin{subfigure}[b]{0.32\textwidth}
    \includegraphics[width=\textwidth]{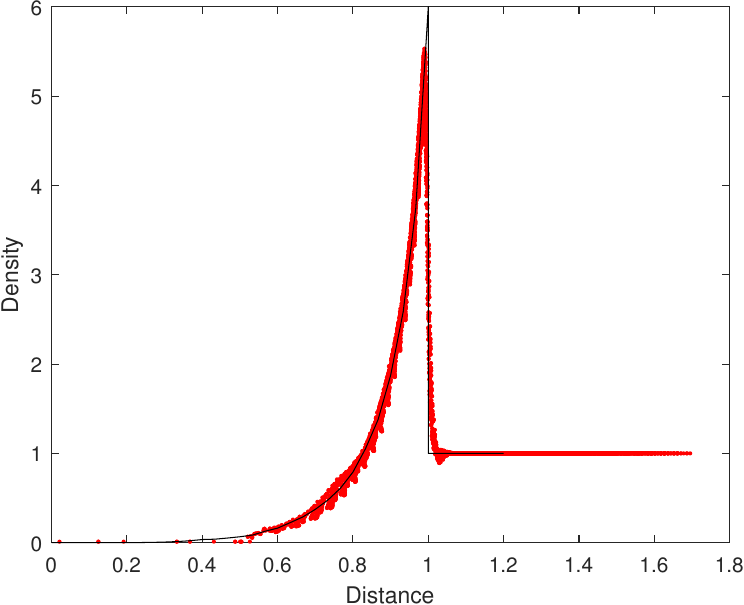}
    \caption{Pointwise density (Original method).}
    \label{Original_Sedov_Q3_point}
  \end{subfigure}
  \caption{Density contours and 1D pointwise density profiles along the radius for the Sedov problem using the $Q^3-Q^2$ element pair on a uniform mesh at $t = 1$.}
  \label{Sedov_Q3_Q2}
\end{figure}

Figure~\ref{Sedov_Q3_Q2} shows the corresponding results using the $Q^3-Q^2$ element pair on the uniform mesh. The hourglass parameter is again set to $0.1$. Without hourglass control, a severely distorted grid is observed (Figure~\ref{Sedov_Q3_nohg}), whereas activating the control effectively corrects this distortion (Figure~\ref{Sedov_Q3_hg}). Interestingly, in this higher-order $Q^3-Q^2$ scenario, the Original method does not suffer from the same severe stiffness observed in the $Q^2-Q^1$ case. Consequently, both the Original method and our proposed framework with hourglass control maintain excellent grid quality, and their pointwise density profiles closely capture the analytical solution.

\begin{figure}[htbp]
  \centering
  \begin{subfigure}[b]{0.32\textwidth}
    \includegraphics[width=\textwidth]{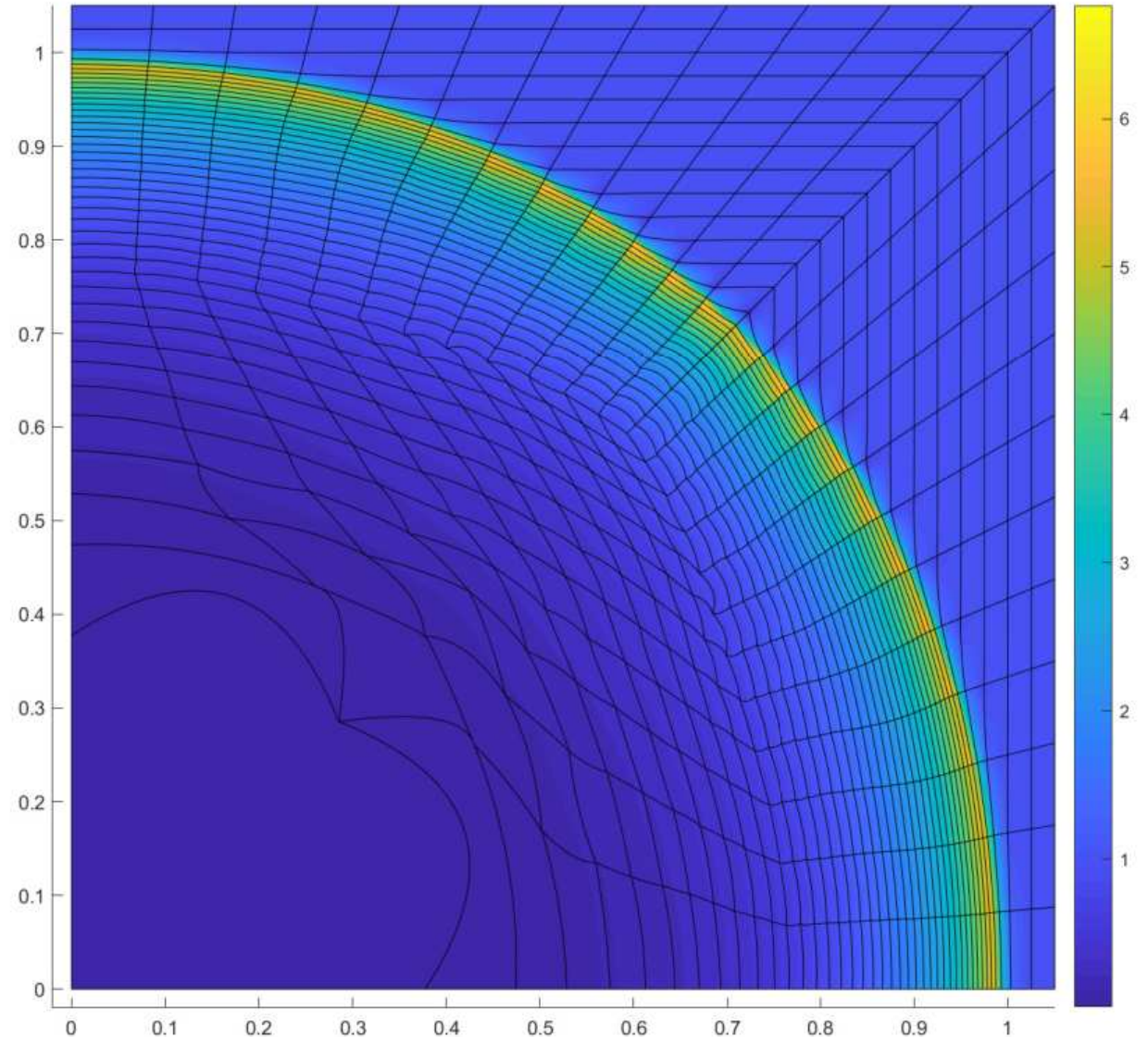}
    \caption{Density contour (without HG).}
    \label{Sedov_Q2_butterfly_nohg}
  \end{subfigure}
  \hfill
  \begin{subfigure}[b]{0.32\textwidth}
    \includegraphics[width=\textwidth]{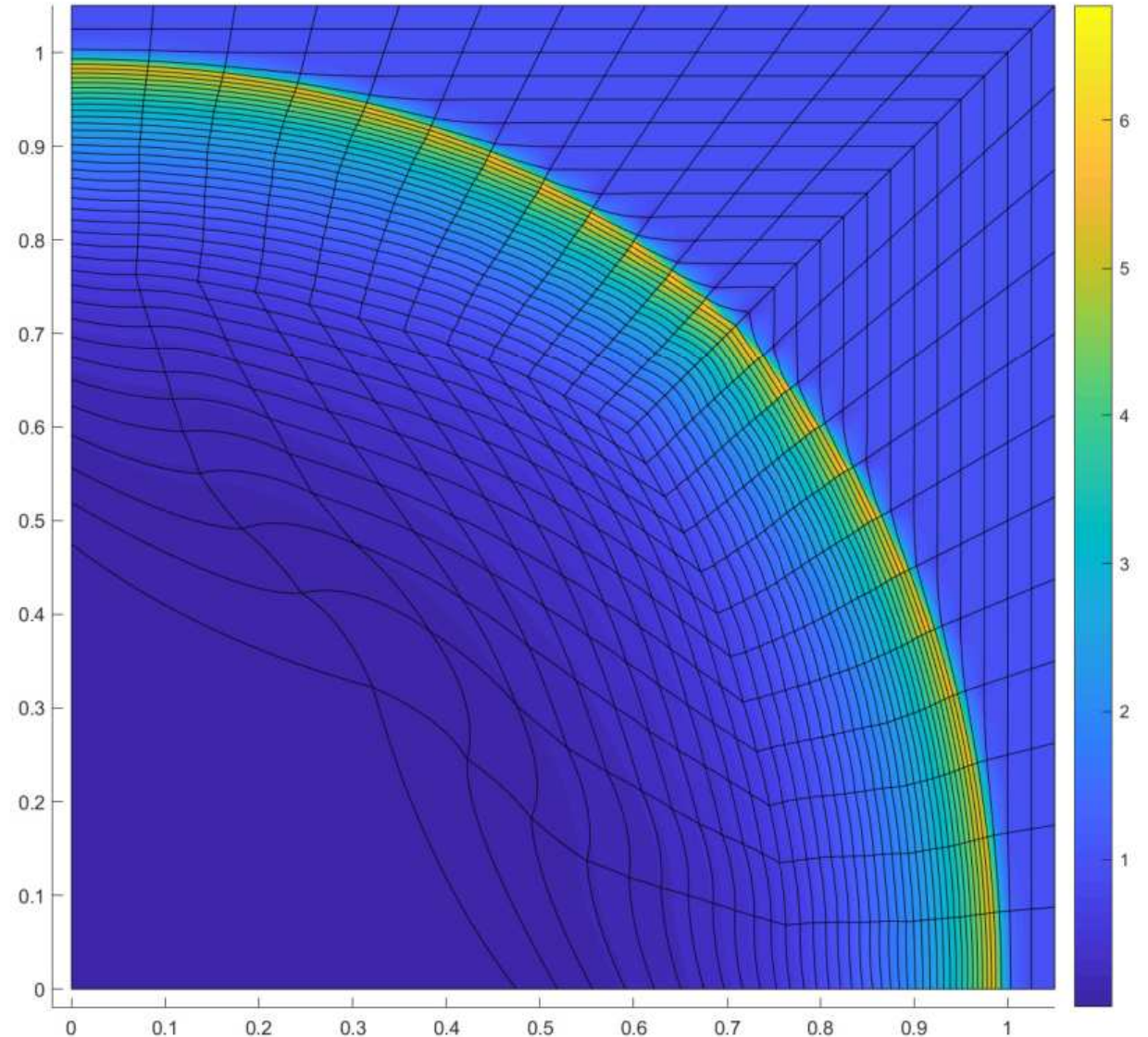}
    \caption{Density contour (with HG).}
    \label{Sedov_Q2_butterfly_hg}
  \end{subfigure}
  \hfill
  \begin{subfigure}[b]{0.32\textwidth}
    \includegraphics[width=\textwidth]{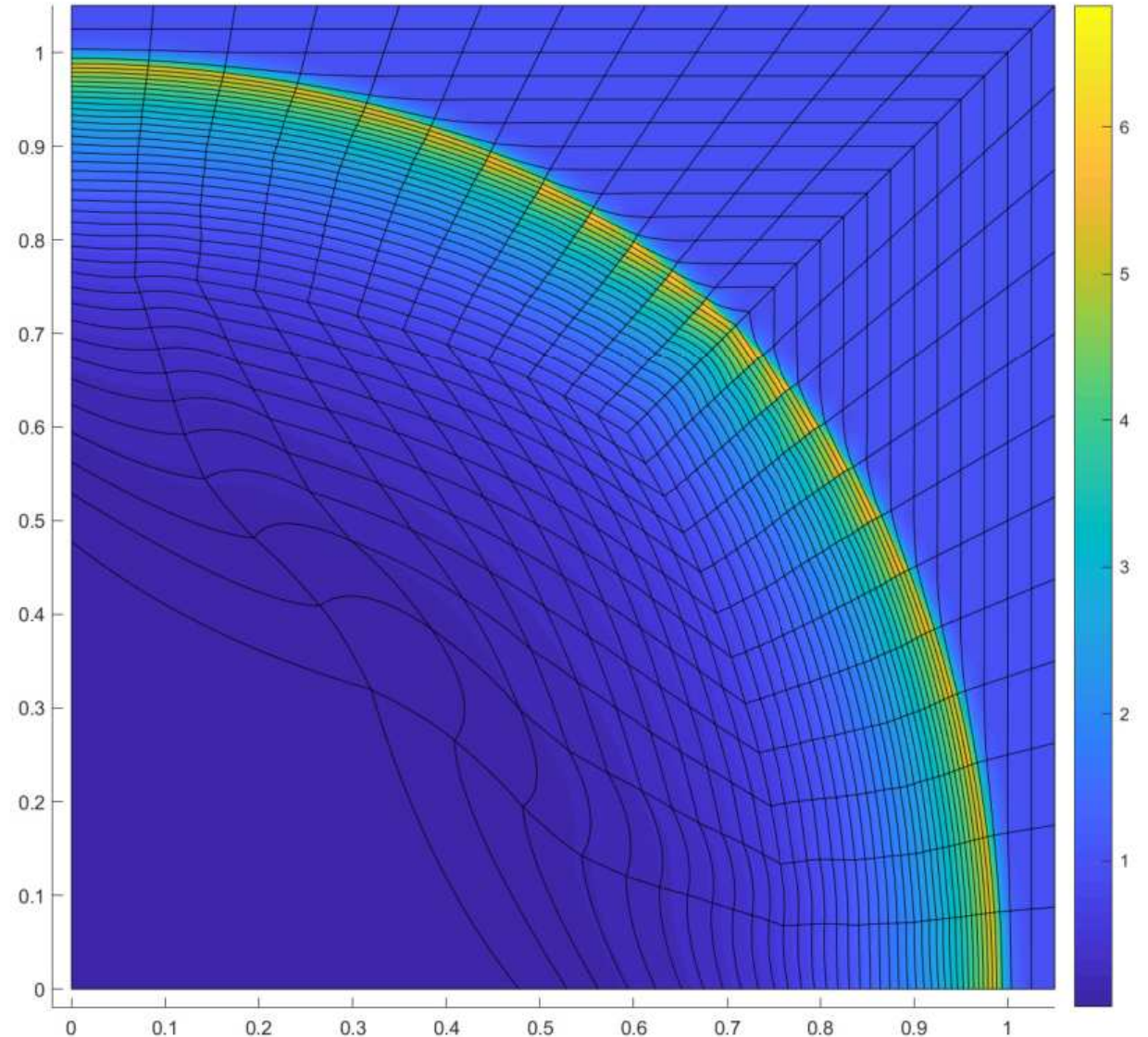}
    \caption{Density contour (Original method).}
    \label{Original_Sedov_Q2_butterfly}
  \end{subfigure}

  \vspace{0.5em}

  \begin{subfigure}[b]{0.32\textwidth}
    \includegraphics[width=\textwidth]{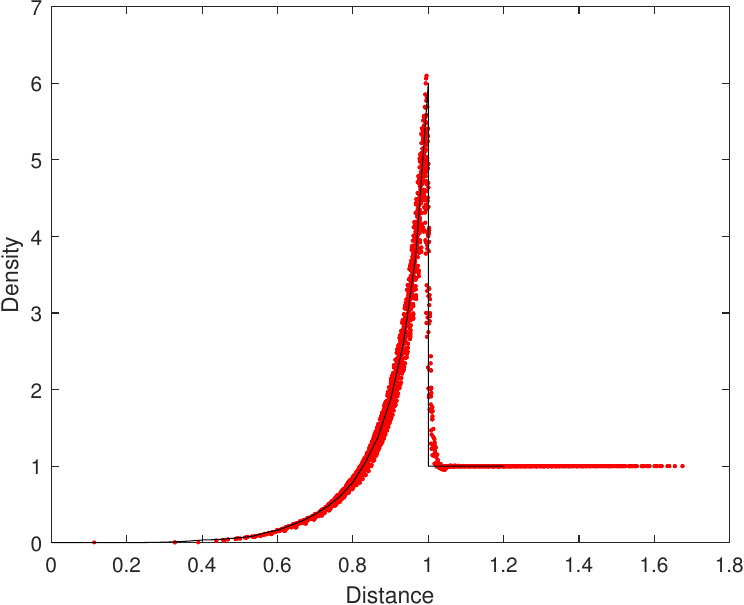}
    \caption{Pointwise density (without HG).}
    \label{Sedov_Q2_butterfly_point_nohg}
  \end{subfigure}
  \hfill
  \begin{subfigure}[b]{0.32\textwidth}
    \includegraphics[width=\textwidth]{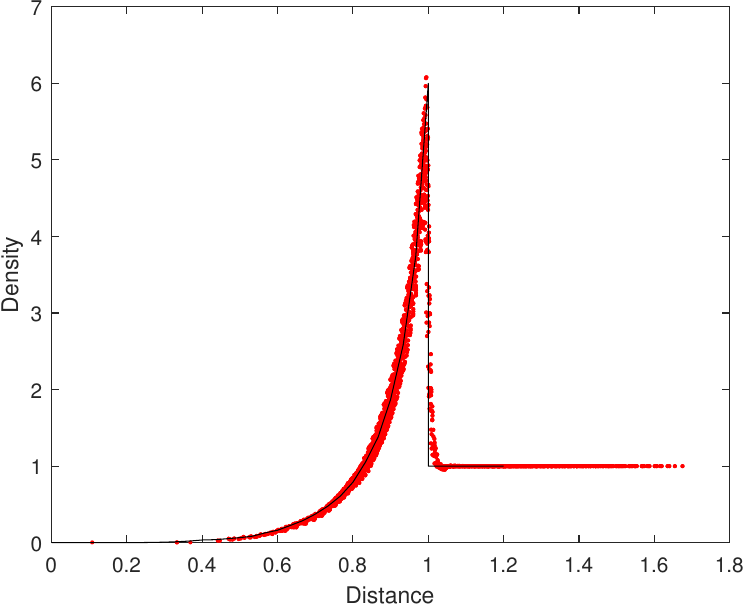}
    \caption{Pointwise density (with HG).}
    \label{Sedov_Q2_butterfly_point}
  \end{subfigure}
  \hfill
  \begin{subfigure}[b]{0.32\textwidth}
    \includegraphics[width=\textwidth]{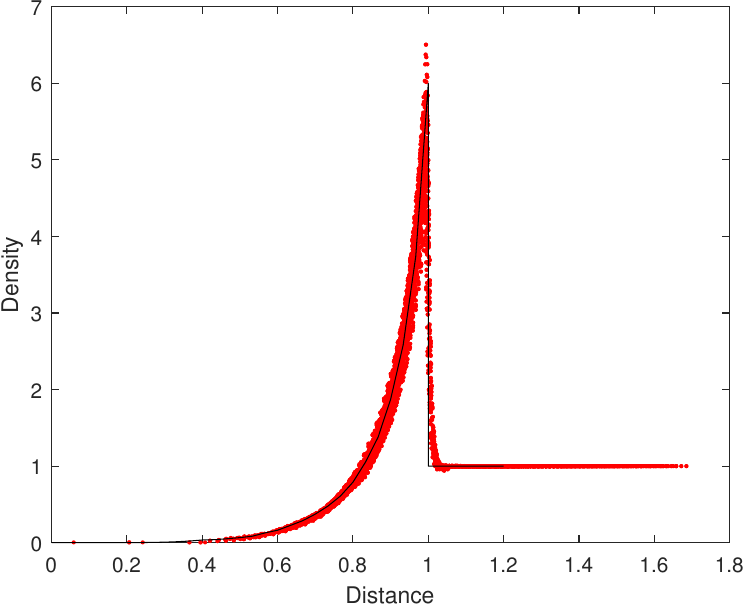}
    \caption{Pointwise density (Original method).}
    \label{Original_Sedov_Q2_butterfly_point}
  \end{subfigure}
  \caption{Density contours and 1D pointwise density profiles along the radius for the Sedov problem using the $Q^2-Q^1$ element pair on a butterfly mesh at $t = 1$.}
  \label{Sedov_butterfly_Q2_Q1}
\end{figure}

Figure~\ref{Sedov_butterfly_Q2_Q1} displays the numerical results for the $Q^2-Q^1$ element pair on a butterfly mesh. Because the butterfly mesh design inherently exacerbates hourglass deformation at the block interfaces, we increased the hourglass parameter to $0.6$. The results without hourglass control exhibit significant mesh tangling, while applying the control maintains a much more regular structure. However, this structure is noticeably stiffer than the uniform mesh case. We observe that the grid near the origin in our framework (Figure~\ref{Sedov_Q2_butterfly_hg}) is visually similar to that of the Original method (Figure~\ref{Original_Sedov_Q2_butterfly}). In this challenging configuration, it becomes obvious that excessive grid stiffness negatively impacts the pointwise density profile: the Original method produces a noticeable density overshoot near the shock peak (Figure~\ref{Original_Sedov_Q2_butterfly_point}), whereas our tunable framework successfully avoids this artifact.

\begin{figure}[htbp]
  \centering
  \begin{subfigure}[b]{0.32\textwidth}
    \includegraphics[width=\textwidth]{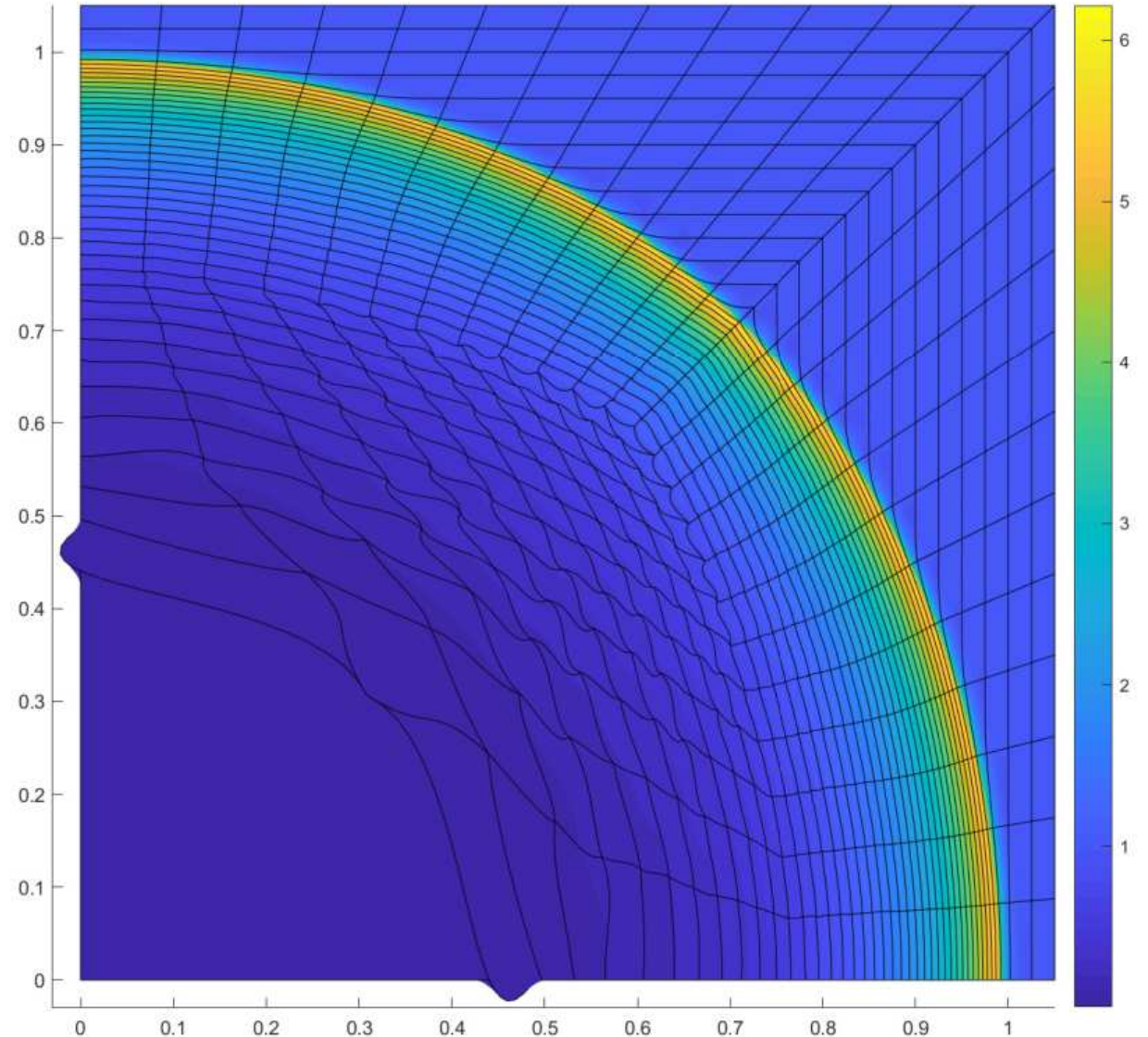}
    \caption{Density contour (without HG).}
    \label{Sedov_Q3_butterfly_nohg}
  \end{subfigure}
  \hfill
  \begin{subfigure}[b]{0.32\textwidth}
    \includegraphics[width=\textwidth]{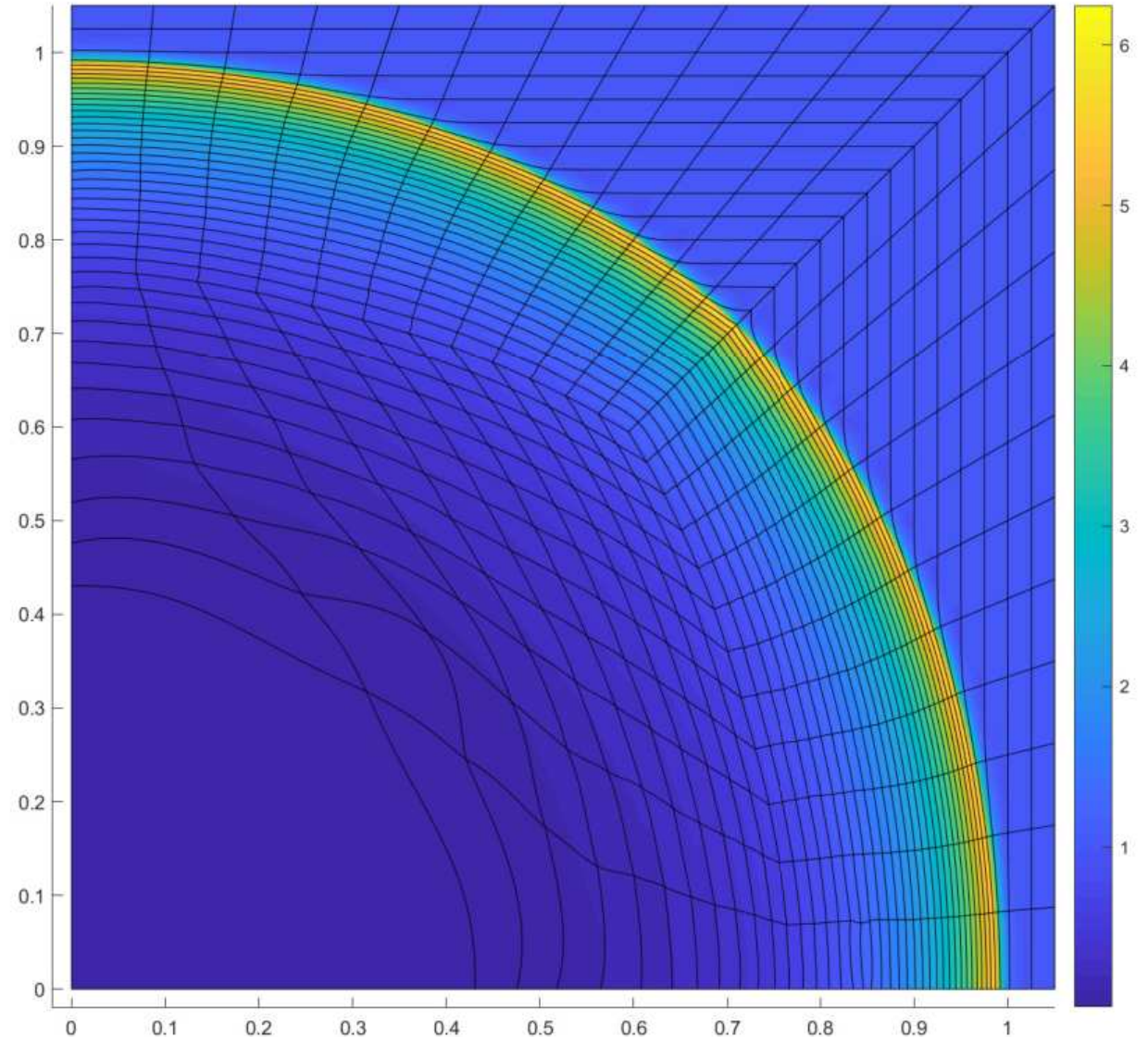}
    \caption{Density contour (with HG).}
    \label{Sedov_Q3_butterfly_hg}
  \end{subfigure}
  \hfill
  \begin{subfigure}[b]{0.32\textwidth}
    \includegraphics[width=\textwidth]{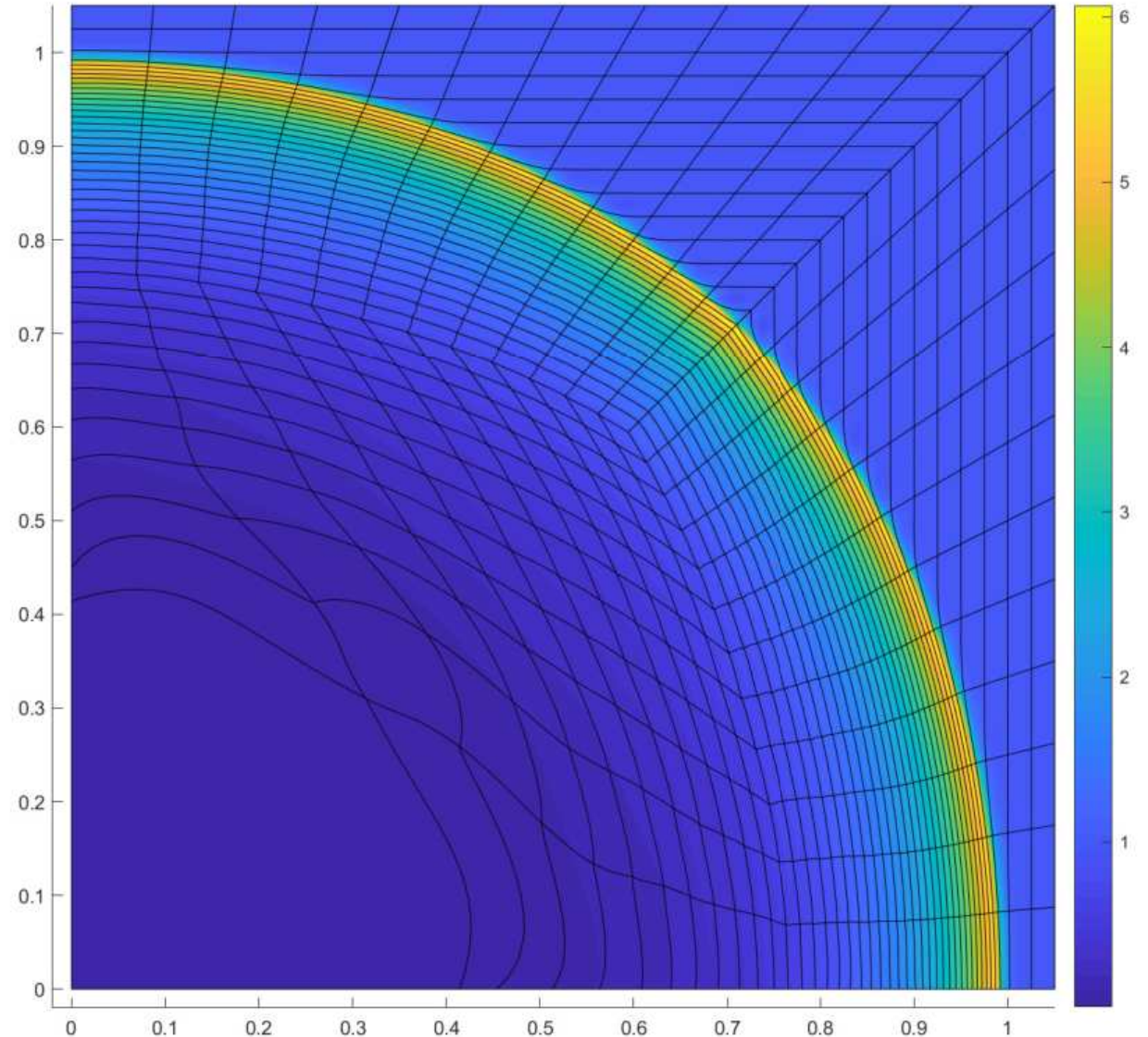}
    \caption{Density contour (Original method).}
    \label{Original_Sedov_Q3_butterfly}
  \end{subfigure}

  \vspace{0.5em}

  \begin{subfigure}[b]{0.32\textwidth}
    \includegraphics[width=\textwidth]{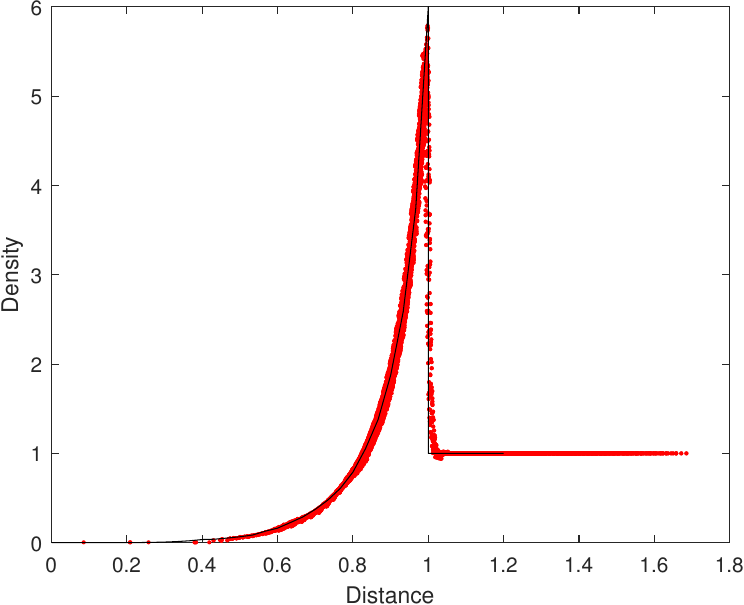}
    \caption{Pointwise density (without HG).}
    \label{Sedov_Q3_butterfly_point_nohg}
  \end{subfigure}
  \hfill
  \begin{subfigure}[b]{0.32\textwidth}
    \includegraphics[width=\textwidth]{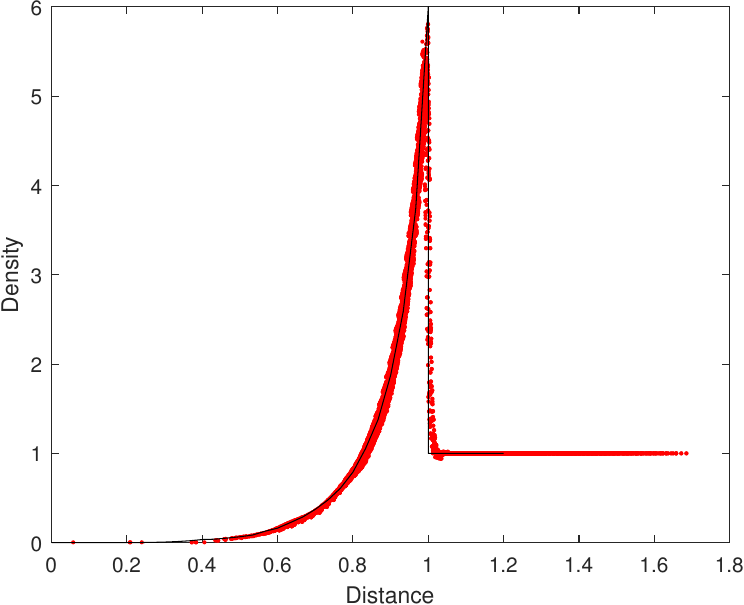}
    \caption{Pointwise density (with HG).}
    \label{Sedov_Q3_butterfly_point}
  \end{subfigure}
  \hfill
  \begin{subfigure}[b]{0.32\textwidth}
    \includegraphics[width=\textwidth]{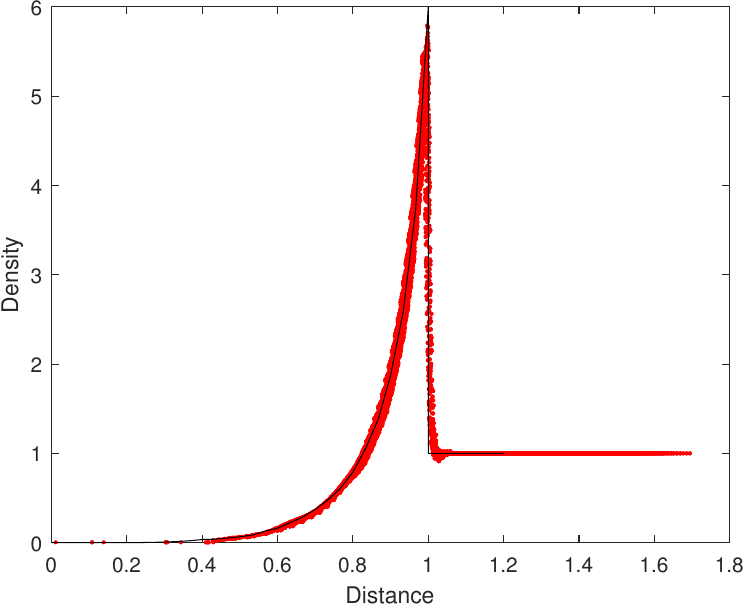}
    \caption{Pointwise density (Original method).}
    \label{Original_Sedov_Q3_butterfly_point}
  \end{subfigure}
  \caption{Density contours and 1D pointwise density profiles along the radius for the Sedov problem using the $Q^3-Q^2$ element pair on a butterfly mesh at $t = 1$.}
  \label{Sedov_butterfly_Q3_Q2}
\end{figure}

Figure~\ref{Sedov_butterfly_Q3_Q2} presents the butterfly mesh results using the $Q^3-Q^2$ element pair, with the hourglass parameter again set to $0.6$. The results without hourglass control exhibit severe mesh breakdown, which is entirely rectified once the control is activated. In this higher-order case, the numerical stiffness of the Original method is less detrimental, making its results largely comparable to our framework with hourglass control. Although Figure~\ref{Original_Sedov_Q3_butterfly} still displays slight stiffness near the origin, the maximum density capturing capabilities of all three strategies remain closely aligned.

\begin{figure}[htbp]
  \centering
  \begin{subfigure}[b]{0.32\textwidth}
    \includegraphics[width=\textwidth]{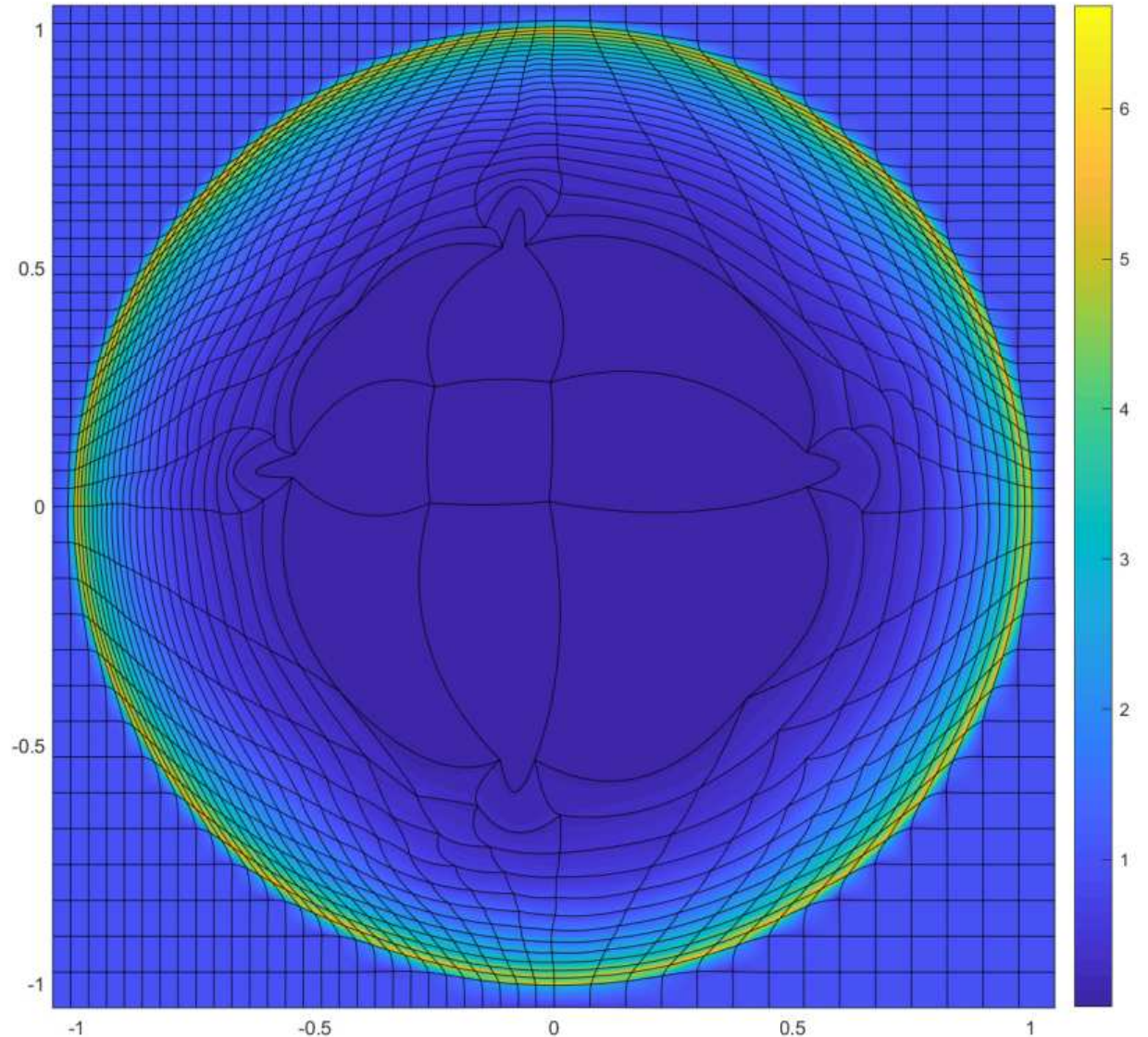}
    \caption{Density contour (without HG).}
    \label{Sedov_Q2_asym_nohg}
  \end{subfigure}
  \hfill
  \begin{subfigure}[b]{0.32\textwidth}
    \includegraphics[width=\textwidth]{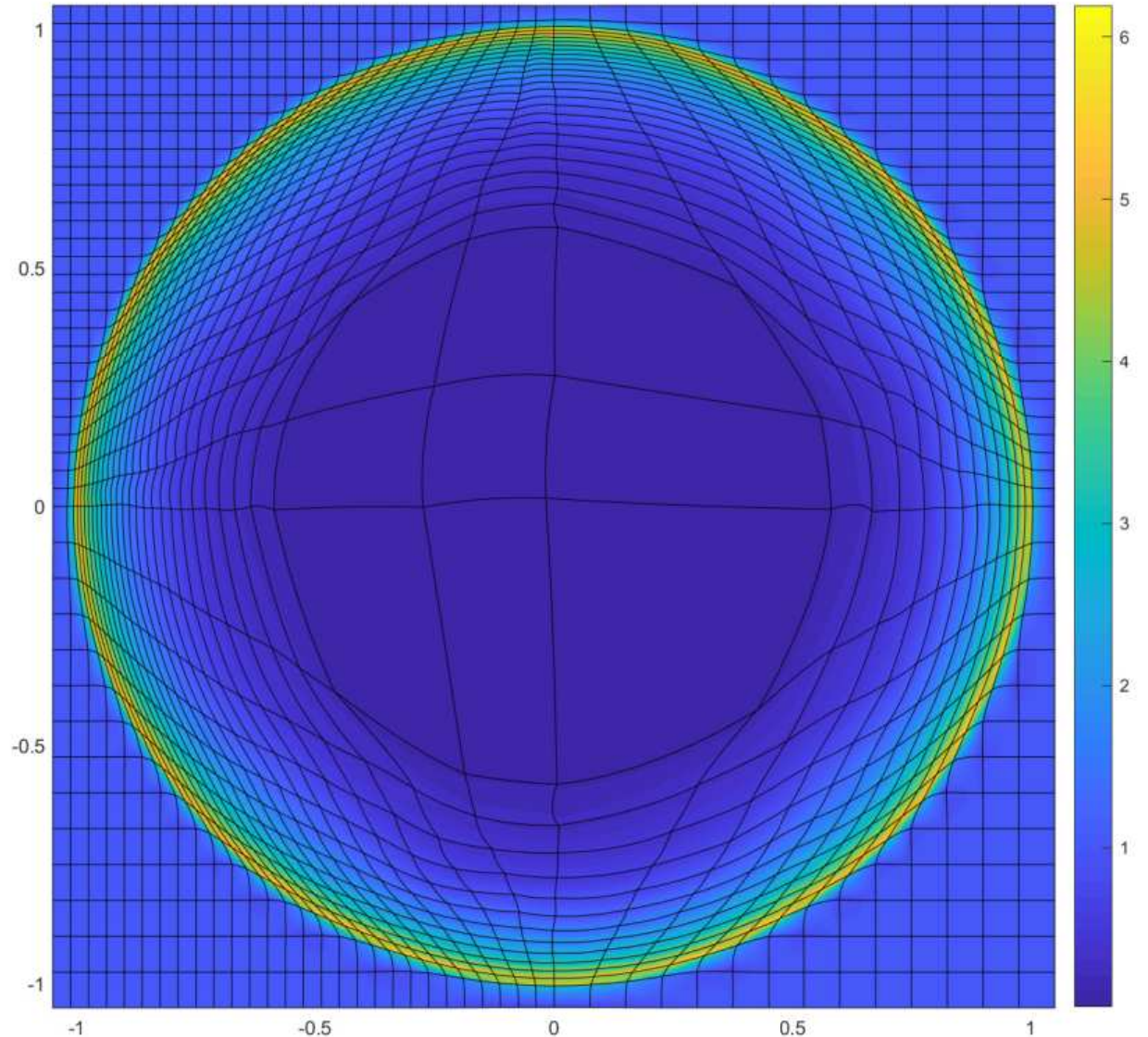}
    \caption{Density contour (with HG).}
    \label{Sedov_Q2_asym_hg}
  \end{subfigure}
  \hfill
  \begin{subfigure}[b]{0.32\textwidth}
    \includegraphics[width=\textwidth]{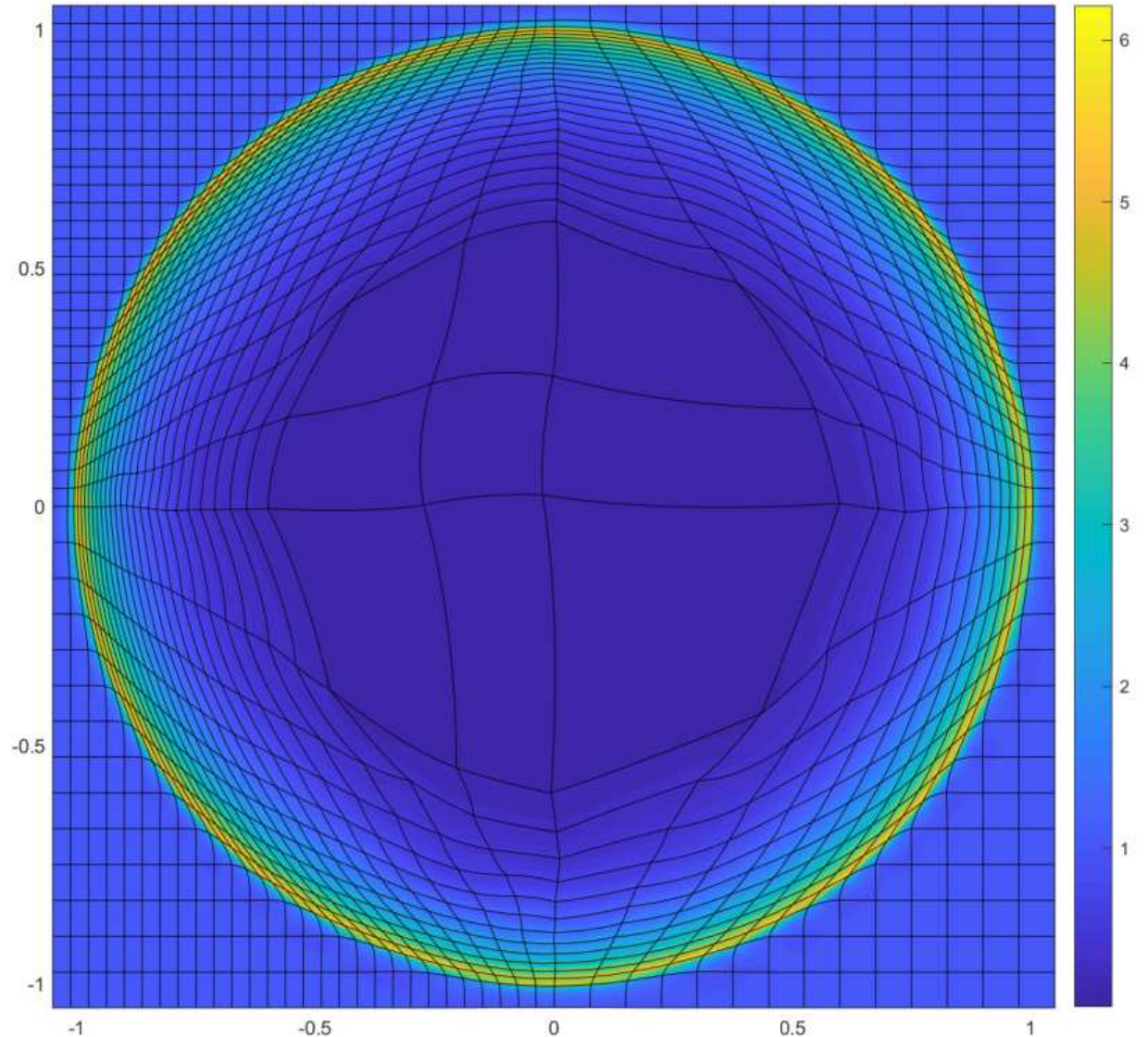}
    \caption{Density contour (Original method).}
    \label{Original_Sedov_Q2_asym}
  \end{subfigure}

  \vspace{0.5em}

  \begin{subfigure}[b]{0.32\textwidth}
    \includegraphics[width=\textwidth]{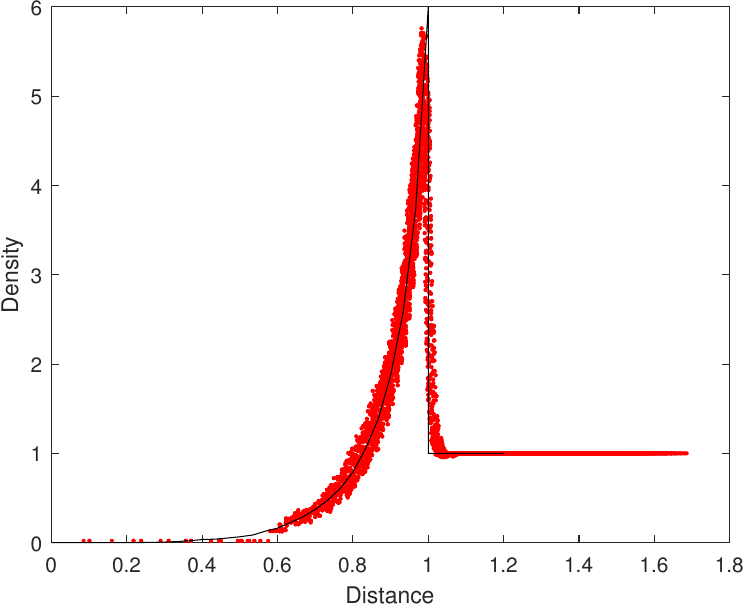}
    \caption{Pointwise density (without HG).}
    \label{Sedov_Q2_asym_point_nohg}
  \end{subfigure}
  \hfill
  \begin{subfigure}[b]{0.32\textwidth}
    \includegraphics[width=\textwidth]{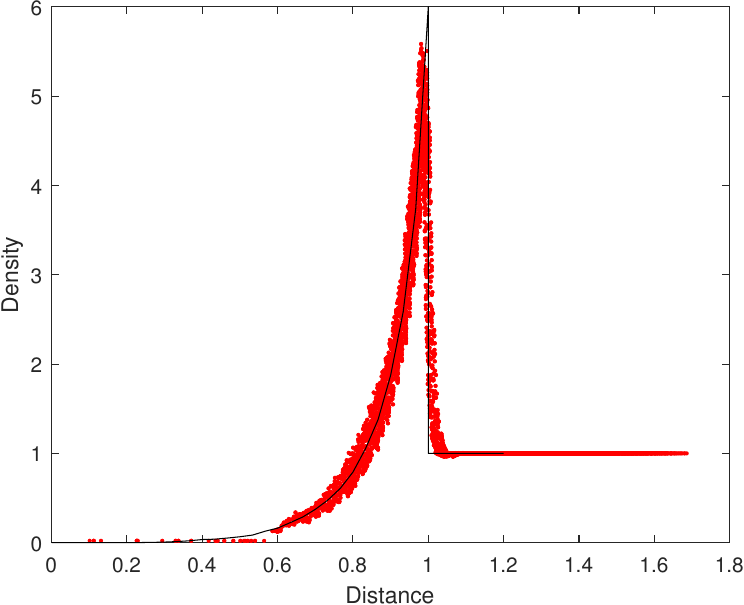}
    \caption{Pointwise density (with HG).}
    \label{Sedov_Q2_asym_point}
  \end{subfigure}
  \hfill
  \begin{subfigure}[b]{0.32\textwidth}
    \includegraphics[width=\textwidth]{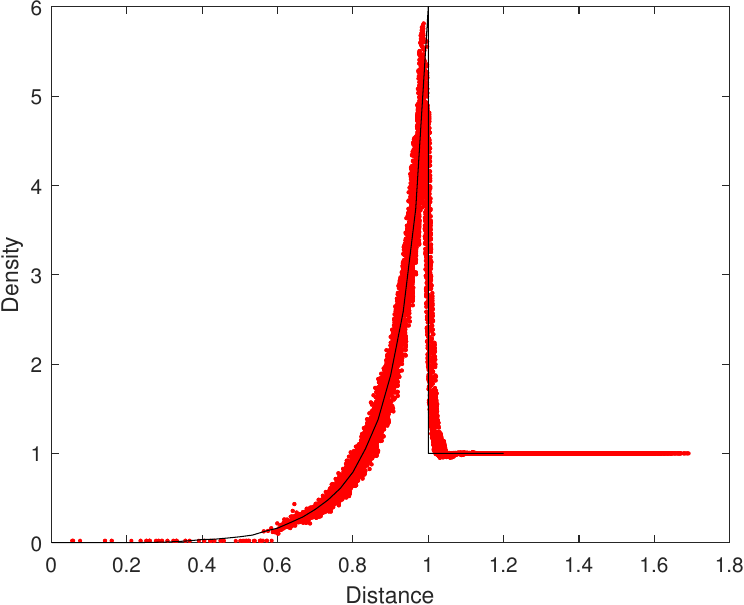}
    \caption{Pointwise density (Original method).}
    \label{Original_Sedov_Q2_asym_point}
  \end{subfigure}
  \caption{Density contours and 1D pointwise density profiles along the radius for the Sedov problem using the $Q^2-Q^1$ element pair on an asymmetric mesh at $t = 1$.}
  \label{Sedov_asym_Q2_Q1}
\end{figure}

Figure~\ref{Sedov_asym_Q2_Q1} shows the numerical results of the Sedov problem using the $Q^2-Q^1$ element pair on an asymmetric mesh at $t = 1$. As with the uniform and butterfly meshes, omitting hourglass control yields extreme grid distortion. By applying a tuned hourglass control parameter of $0.4$, this distortion is effectively suppressed, and the overall solution quality improves dramatically. However, the Original method again suffers from inherent numerical stiffness, resulting in an unphysical overshoot in the maximum density compared to the other two approaches. While Figure~\ref{Sedov_Q2_asym_hg} still displays some minor grid imperfections—such as slight node displacement along the $x$ and $y$ axes—the Original method is similarly incapable of maintaining perfect nodal symmetry along these axes due to its excessive stiffness.

\begin{figure}[htbp]
  \centering
  \begin{subfigure}[b]{0.32\textwidth}
    \includegraphics[width=\textwidth]{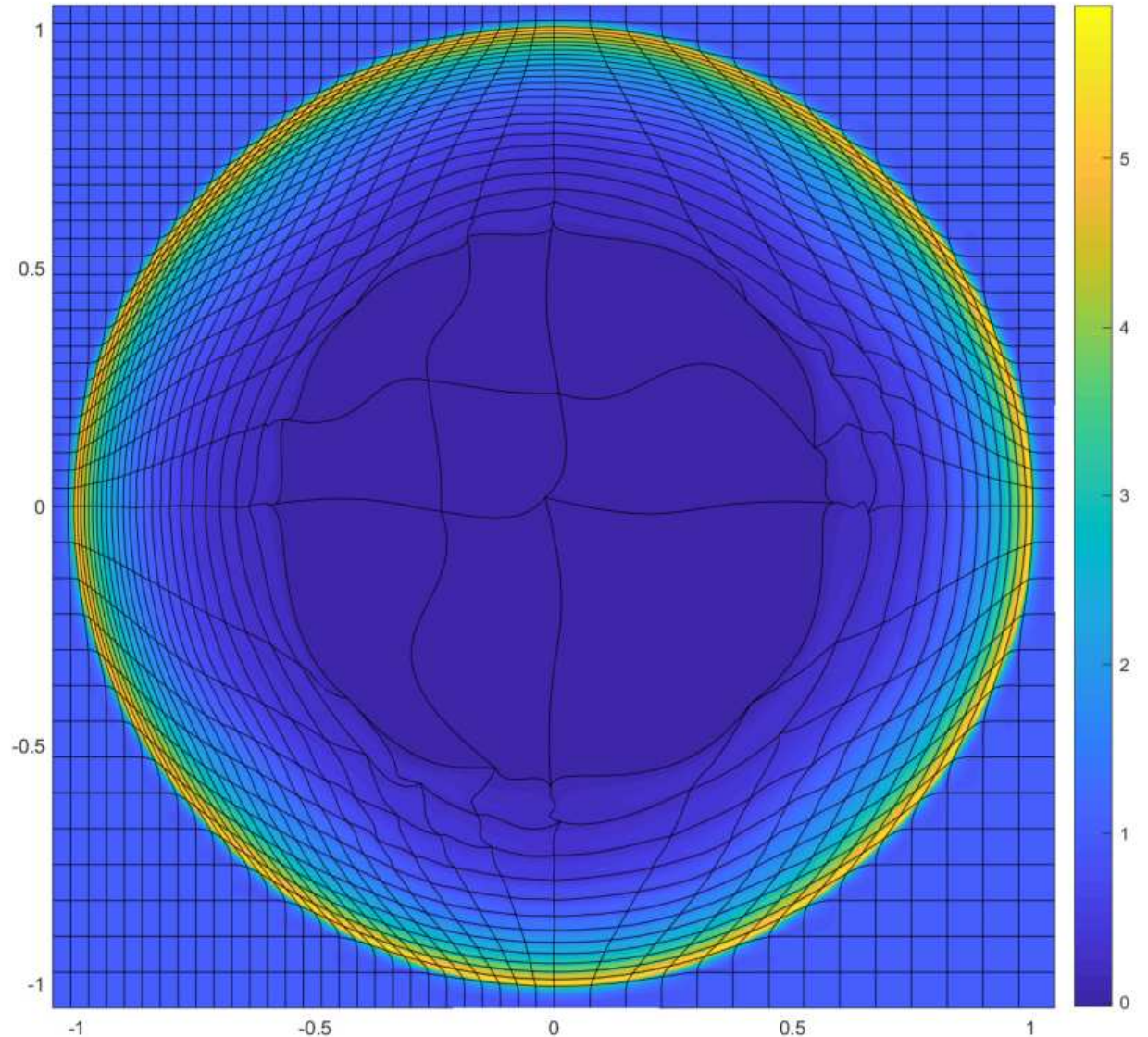}
    \caption{Density contour (without HG).}
    \label{Sedov_Q3_asym_nohg}
  \end{subfigure}
  \hfill
  \begin{subfigure}[b]{0.32\textwidth}
    \includegraphics[width=\textwidth]{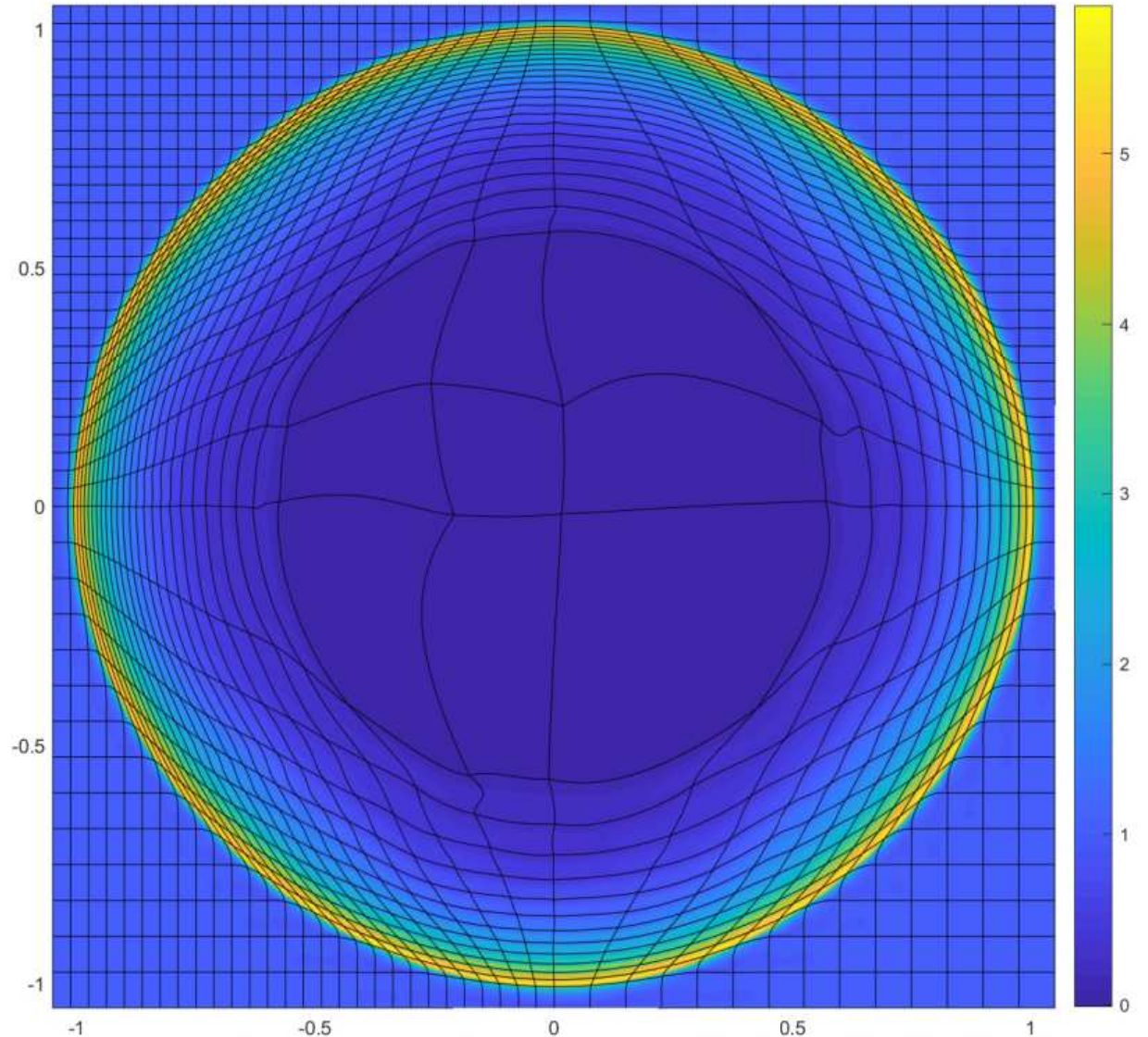}
    \caption{Density contour (with HG).}
    \label{Sedov_Q3_asym_hg}
  \end{subfigure}
  \hfill
  \begin{subfigure}[b]{0.32\textwidth}
    \includegraphics[width=\textwidth]{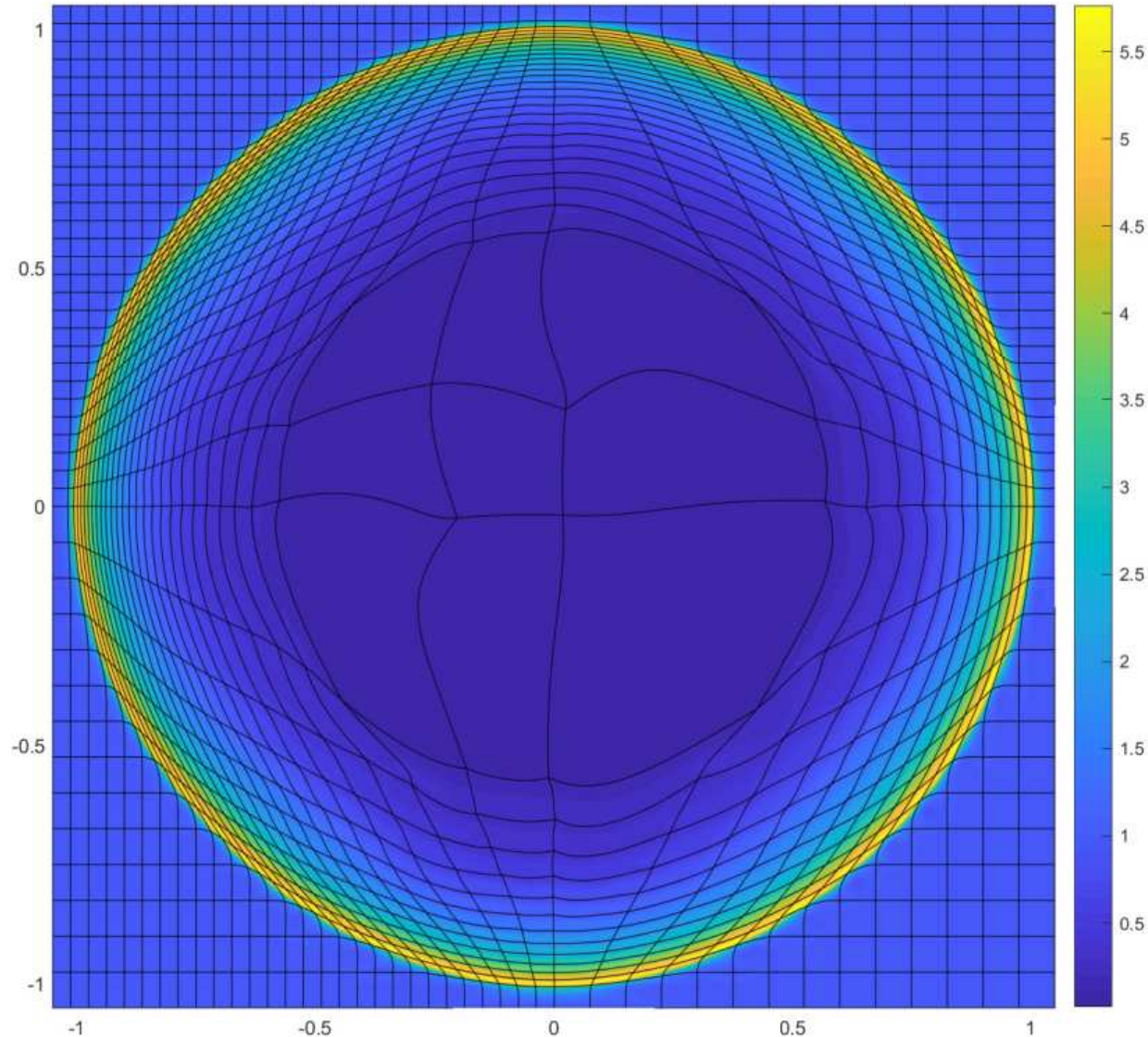}
    \caption{Density contour (Original method).}
    \label{Original_Sedov_Q3_asym}
  \end{subfigure}

  \vspace{0.5em}

  \begin{subfigure}[b]{0.32\textwidth}
    \includegraphics[width=\textwidth]{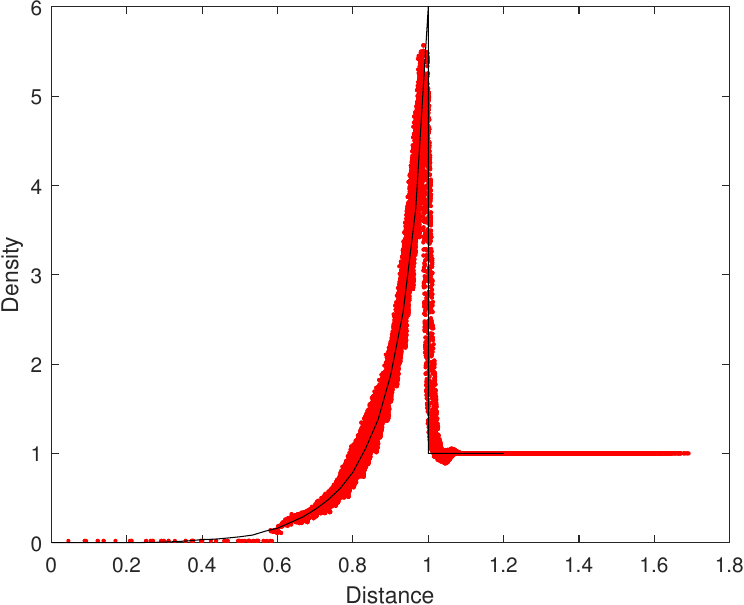}
    \caption{Pointwise density (without HG).}
    \label{Sedov_Q3_asym_point_nohg}
  \end{subfigure}
  \hfill
  \begin{subfigure}[b]{0.32\textwidth}
    \includegraphics[width=\textwidth]{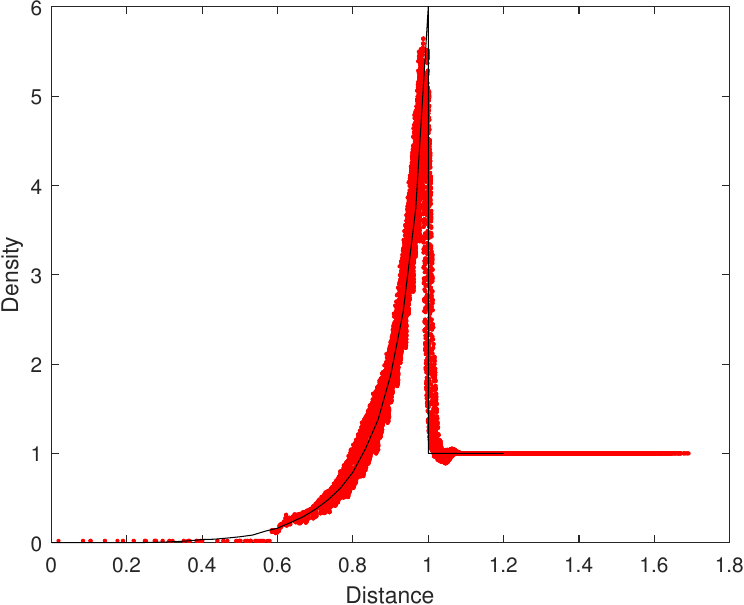}
    \caption{Pointwise density (with HG).}
    \label{Sedov_Q3_asym_point}
  \end{subfigure}
  \hfill
  \begin{subfigure}[b]{0.32\textwidth}
    \includegraphics[width=\textwidth]{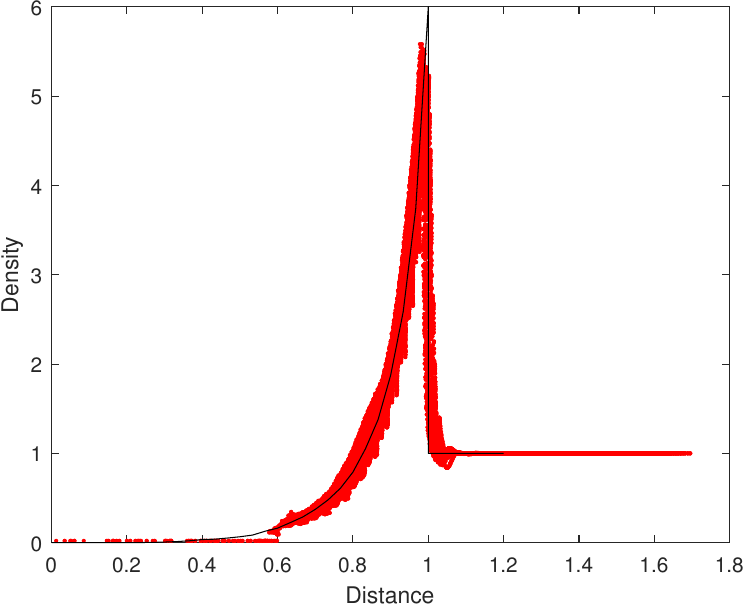}
    \caption{Pointwise density (Original method).}
    \label{Original_Sedov_Q3_asym_point}
  \end{subfigure}
  \caption{Density contours and 1D pointwise density profiles along the radius for the Sedov problem using the $Q^3-Q^2$ element pair on an asymmetric mesh at $t = 1$.}
  \label{Sedov_asym_Q3_Q2}
\end{figure}

For the $Q^3-Q^2$ element pair on the asymmetric mesh, Figure~\ref{Sedov_asym_Q3_Q2} presents the numerical results at $t = 1$. Consistent with the $Q^2-Q^1$ case, the absence of hourglass control leads to severe grid distortion, which is well-managed by setting the hourglass parameter to $0.4$. In this configuration, the performance of the proposed framework and the Original method is once again comparable. The Original method still retains a degree of grid stiffness, whereas our proposed framework allows for slightly more flexibility, leaving only minor, acceptable grid distortions.

\subsubsection{Discussion on stabilization and numerical stiffness}
Based on the extensive numerical experiments conducted for the Noh and Sedov problems, we offer a few concluding remarks regarding the proposed framework and the Original method. Generally, the numerical performance of the two approaches is comparable, as both effectively control grid distortion and maintain solution quality. However, a significant drawback of the Original method is that it lacks a tunable mechanism to control the degree of hourglass suppression. As demonstrated, this rigid suppression often renders the numerical solution excessively stiff, which can artificially distort the density field (e.g., producing overshoots near the shock front) and degrade overall grid quality in skewed domains. Conversely, our proposed framework provides a fully tunable hourglass parameter, enabling users to strike an optimal balance between grid regularity and solution fidelity. 

It is worth noting that in our framework, the pressure variation is estimated from the density variation and an interpolated sound speed. While this approach guarantees thermodynamically consistent behavior, it may not perfectly capture the actual pressure variation in all extreme scenarios, potentially affecting grid evolution. Because the ideal gas law represents a simple, linear Equation of State (EOS), both the Original method and our framework maintain reasonable accuracy in these tests. However, for materials governed by complex EOS, the direct interpolation of internal energy utilized in the Original method may introduce additional unphysical errors. Therefore, while our proposed framework successfully offers the flexibility to avoid numerical stiffness while preserving thermodynamic consistency, careful parameter tuning remains important when extending these simulations to complex equations of state.

\subsubsection{Dukowicz-Meltz problem} 
Having established the robustness of the proposed hourglass control and the shock-capturing fidelity of the framework, we now evaluate the method on complex interaction problems involving strong fluid deformation. This problem employs a piston-driven setup involving two distinct regions of an ideal gas. The left region is a right-angled trapezoid with a vertical left boundary, while the right region is a slanted parallelogram sharing an interface inclined at $60^{\circ}$. A piston on the far-left boundary moves inward with a constant velocity of $1.48$. The initial conditions for the left region are set to $\gamma = 1.4$, $\rho = 1$, and $e = 2.5$; the right region is initialized with $\gamma = 1.4$, $\rho = 1.5$, and $e = 2.5$.

\begin{figure}[htbp]
  \centering
  \begin{subfigure}[c]{0.49\textwidth}
    \includegraphics[width=0.90\textwidth]{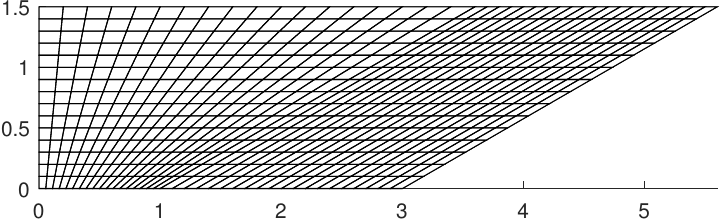}
    \caption{Initial computational mesh.}
    \label{DM_grid}
  \end{subfigure}
  \hfill
  \begin{subfigure}[c]{0.245\textwidth}
    \includegraphics[width=\textwidth]{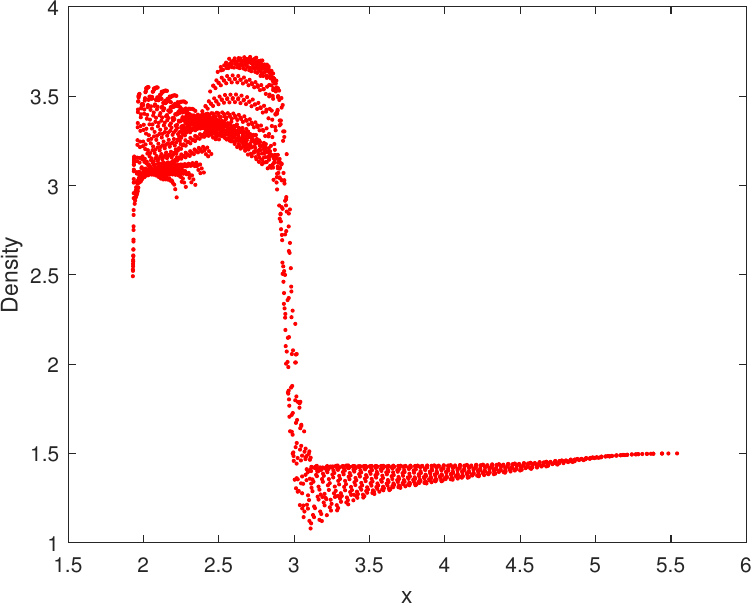}
    \caption{Pointwise density ($Q^2-Q^1$).}
    \label{DM_Q2_Q1_point}
  \end{subfigure}
  \hfill
  \begin{subfigure}[c]{0.245\textwidth}
    \includegraphics[width=\textwidth]{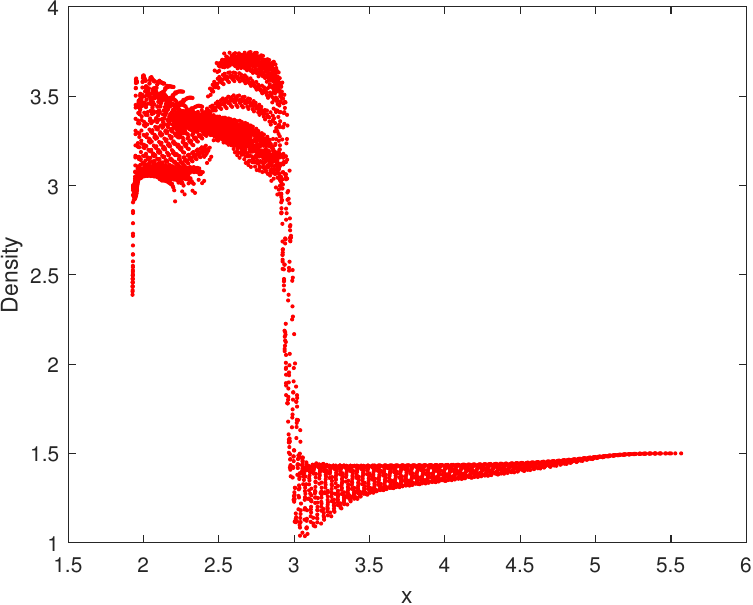}
    \caption{Pointwise density ($Q^3-Q^2$).}
    \label{DM_Q3_Q2_point}
  \end{subfigure}  
  \caption{Initial computational mesh and 1D pointwise density distributions along the $x$-axis for the Dukowicz-Meltz problem at $t = 1.3$.}
\label{DM}
\end{figure}

The initial $38 \times 15$ computational mesh is detailed in Figure~\ref{DM_grid}. The left trapezoidal region is partitioned into $18 \times 15$ non-uniform elements, whereas the right parallelogram region is uniformly partitioned into $20 \times 15$ elements. The terminal time of the simulation is set to $t = 1.3$. The artificial viscosity coefficients are chosen as $c_1=0.5$ and $c_2=0.5$, and the hourglass control parameter is fixed at $1$.

\begin{figure}[htbp]
  \centering
  \begin{subfigure}[b]{0.485\textwidth}
    \includegraphics[width=\textwidth]{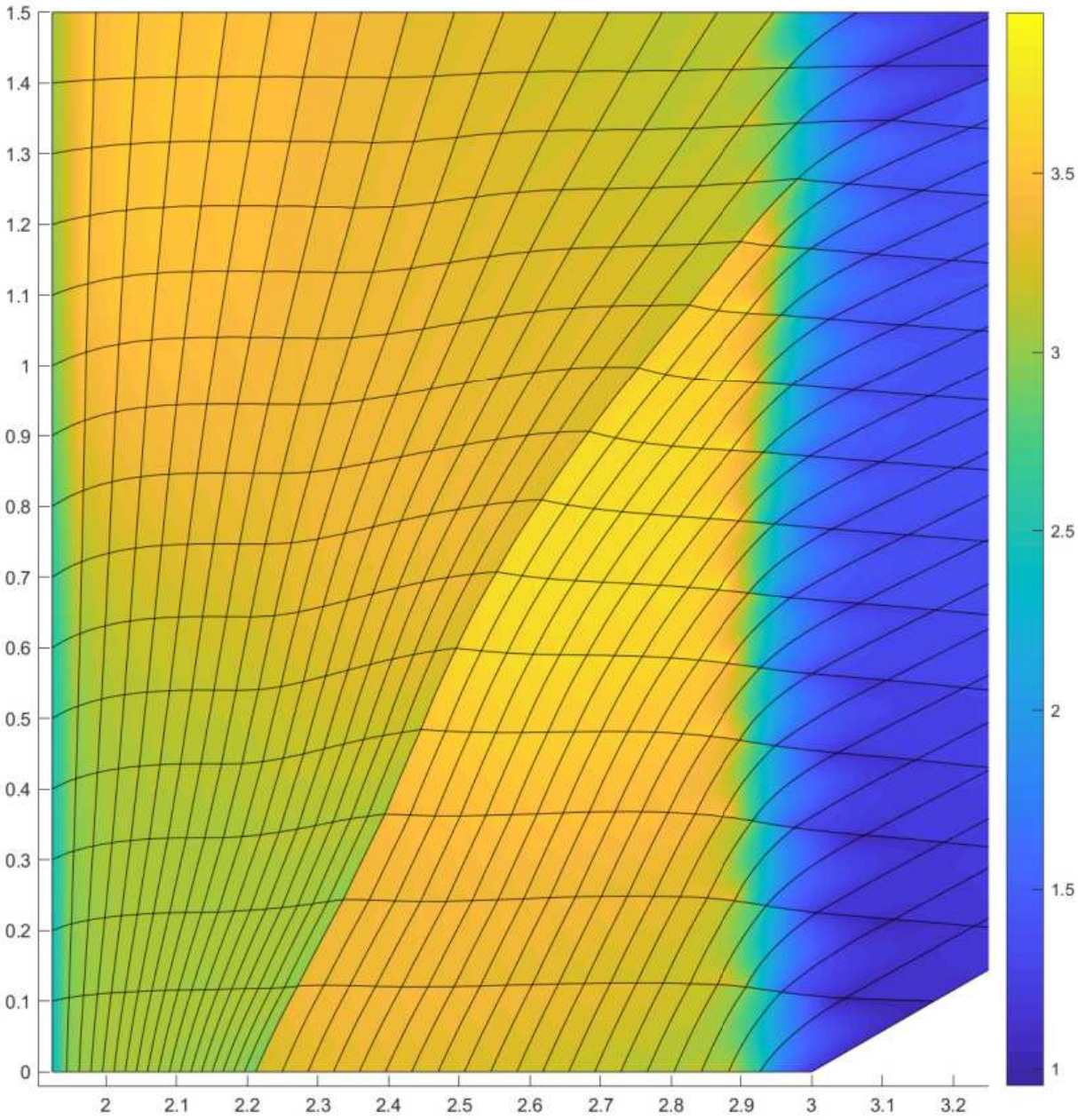}
    \caption{$Q^2-Q^1$ element pair.}
    \label{DM_Q2}
  \end{subfigure}
  \hfill
  \begin{subfigure}[b]{0.485\textwidth}
    \includegraphics[width=\textwidth]{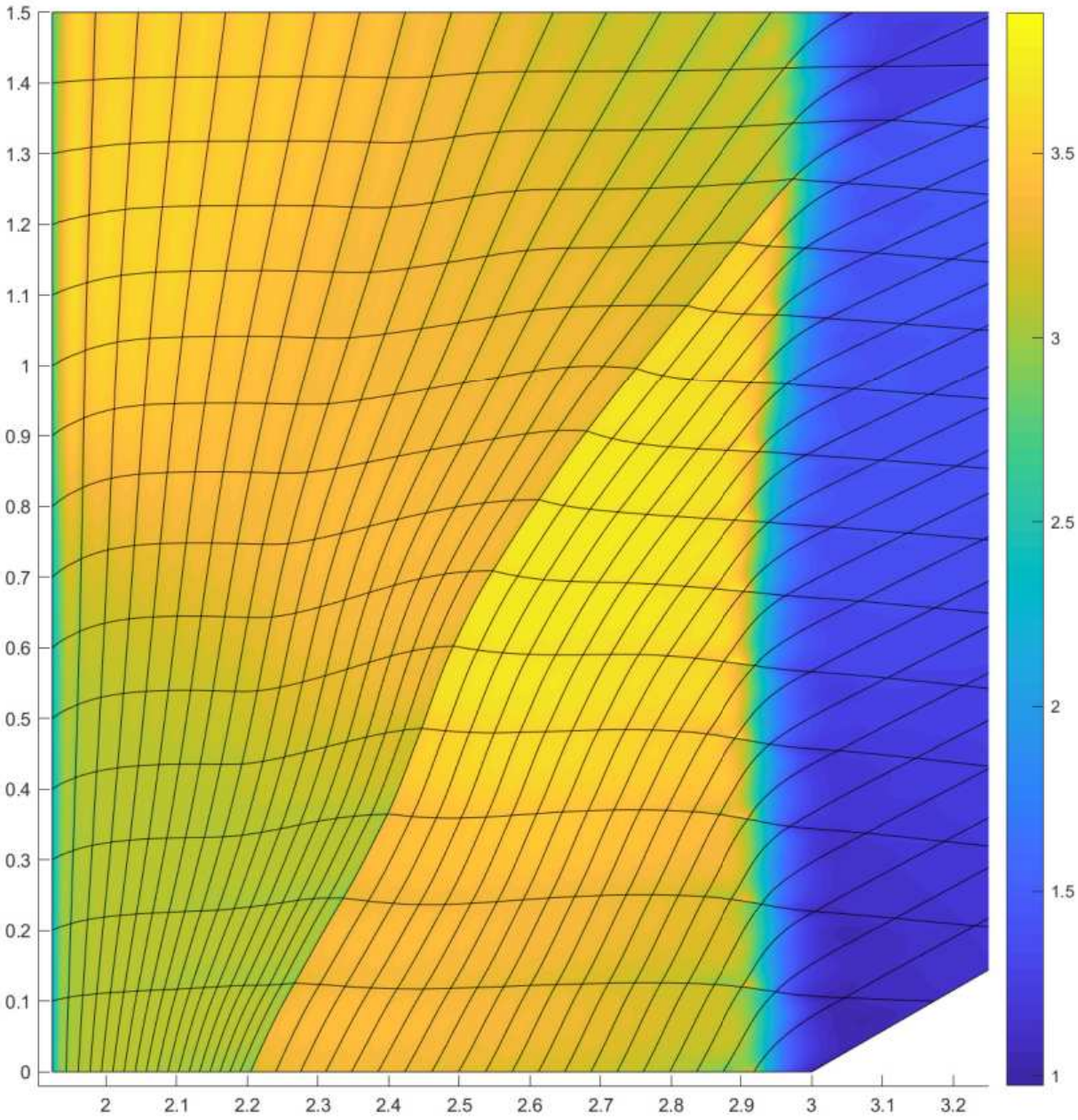}
    \caption{$Q^3-Q^2$ element pair.}
    \label{DM_Q3}
  \end{subfigure}
  \caption{2D density contours and deformed mesh configurations for the Dukowicz-Meltz problem at $t = 1.3$.}
  \label{DM_2D}
\end{figure}

Figures~\ref{DM_Q2_Q1_point} and \ref{DM_Q2} present the 1D pointwise density profile and the 2D density contour with the grid distribution for the $Q^2-Q^1$ element pair, respectively. Similarly, Figures~\ref{DM_Q3_Q2_point} and \ref{DM_Q3} display the corresponding numerical results for the higher-order $Q^3-Q^2$ element pair. As observed in the 2D contour plots, both numerical schemes successfully maintain high grid quality throughout the fluid deformation. The mesh remains robust without severe tangling, and unphysical grid distortion is effectively suppressed near the slanted contact discontinuities and the shock front. Furthermore, the 1D scatter plots demonstrate a sharp and accurate capture of the shock wave located near $x = 3$. The remarkably tight clustering of the pointwise density data indicates that spurious numerical oscillations are well mitigated by the artificial viscosity and hourglass control mechanisms, highlighting the stability and accuracy of both element pairs. 

\subsubsection{Triple-point shock interaction}
Finally, to demonstrate the inherent superiority of the high-order $Q^3-Q^2$ kinematics in resolving complex vortex dynamics, we simulate the Triple Point problem. The computational domain spans $[0,7] \times [0,3]$ and is divided into three distinct regions filled with ideal gases at varying initial states. Their rapid interaction generates an intricate shock system accompanied by a physical roll-up vortex. The left region ($[0,1] \times [0,3]$) is initialized with $\gamma_1 = 1.5$, $\rho_1 = 1$, and $p_1 = 1$. The lower-right region ($[1,7] \times [0,1.5]$) starts with $\gamma_2 = 1.4$, $\rho_2 = 1$, and $p_2 = 0.1$. The upper-right region ($[1,7] \times [1.5,3]$) has $\gamma_3 = 1.6$, $\rho_3 = 0.125$, and $p_3 = 0.1$. The terminal time of the simulation is set to $t = 2.5$. The artificial viscosity coefficients are chosen as $c_1=0.5$ and $c_2=0.5$, and the hourglass control parameter is fixed at $1$.

\begin{figure}[htbp]
  \centering
  \begin{subfigure}[b]{0.95\textwidth}
    \includegraphics[width=\textwidth]{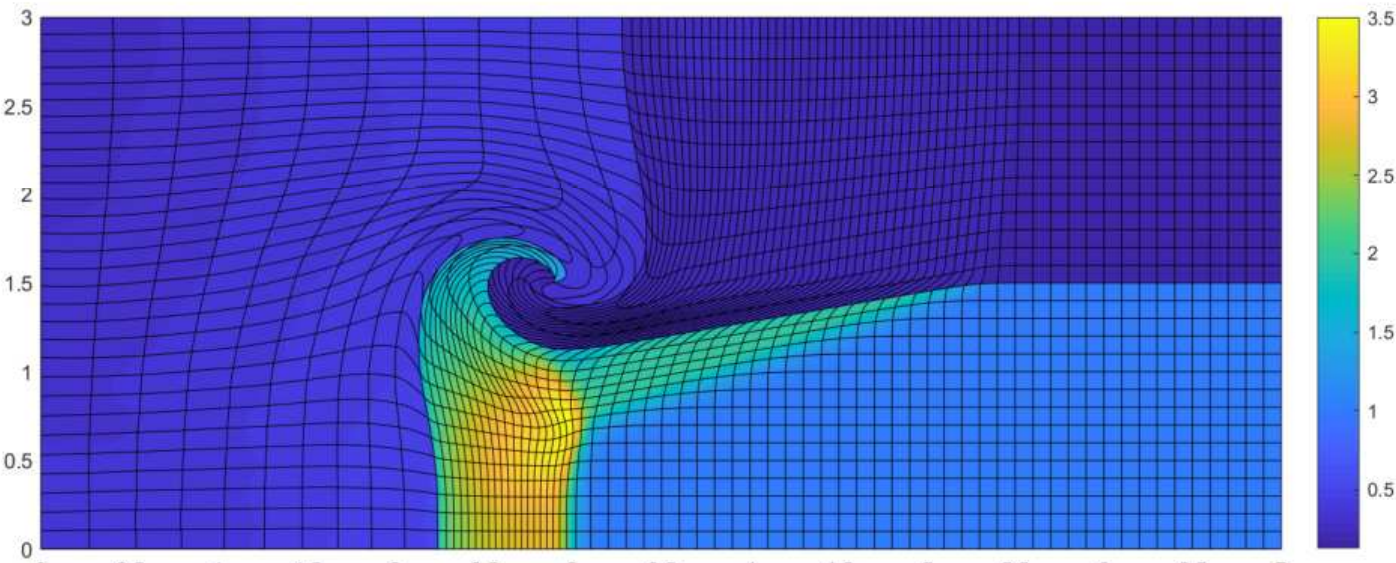}
    \caption{$Q^2-Q^1$ element pair.}
    \label{triple_Q2}
  \end{subfigure}
  
  \vspace{0.5em}
  
  \begin{subfigure}[b]{0.95\textwidth}
    \includegraphics[width=\textwidth]{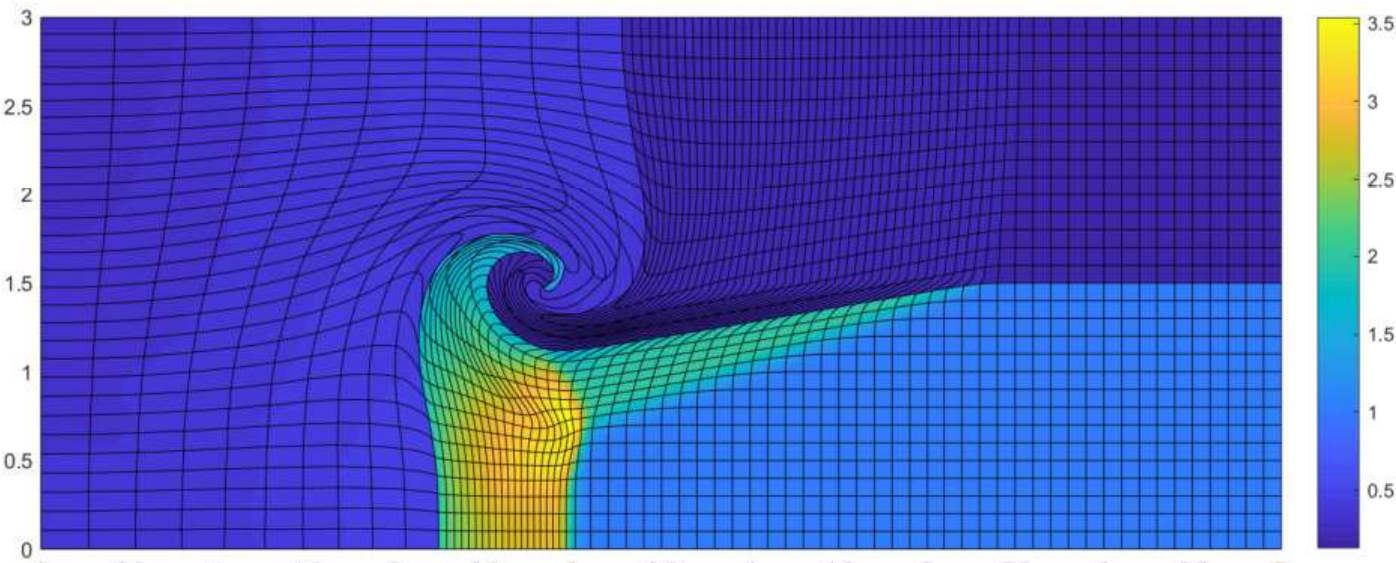}
    \caption{$Q^3-Q^2$ element pair.}
    \label{triple_Q3}
  \end{subfigure}
  \caption{Macroscopic 2D density contours and deformed mesh configurations for the Triple Point problem at $t = 2.5$.}
  \label{triple}
\end{figure}

Figures~\ref{triple_Q2} and \ref{triple_Q3} display the macroscopic density fields and the corresponding deformed meshes across the entire computational domain at $t=2.5$ using the $Q^2-Q^1$ and $Q^3-Q^2$ element pairs, respectively. Both numerical schemes robustly capture the propagation of the primary shock wave and the complex contact discontinuities. Despite the severe shear deformations inherent in this problem, the mesh moves seamlessly with the fluid flow without suffering from fatal entanglement, demonstrating exceptional geometric adaptability. 

\begin{figure}[htbp]
  \centering
  \begin{subfigure}[b]{0.49\textwidth}
    \includegraphics[width=\textwidth]{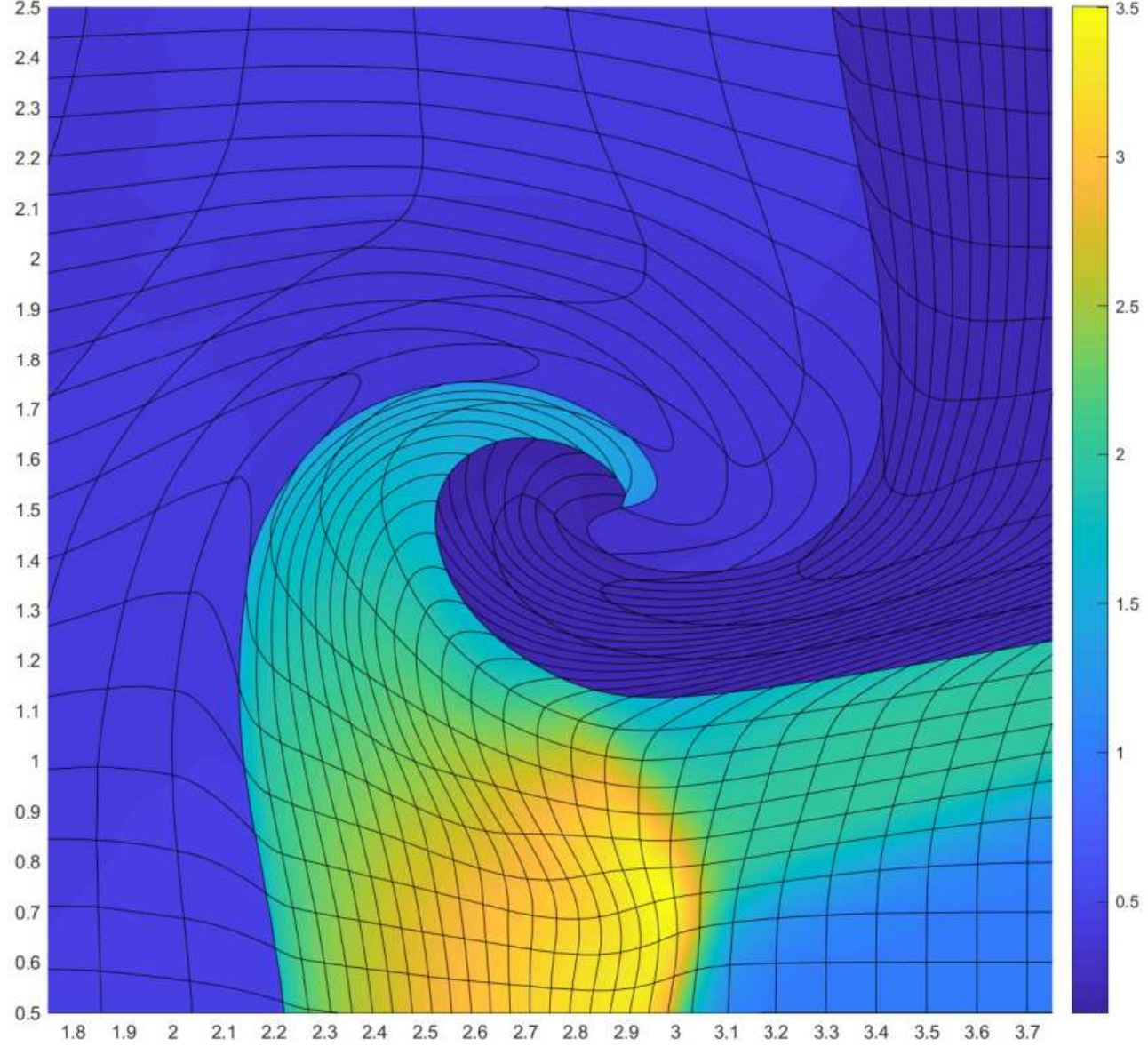}
    \caption{$Q^2-Q^1$ element pair.}
    \label{triple_Q2_zoom}
  \end{subfigure}
  \hfill
  \begin{subfigure}[b]{0.49\textwidth}
    \includegraphics[width=\textwidth]{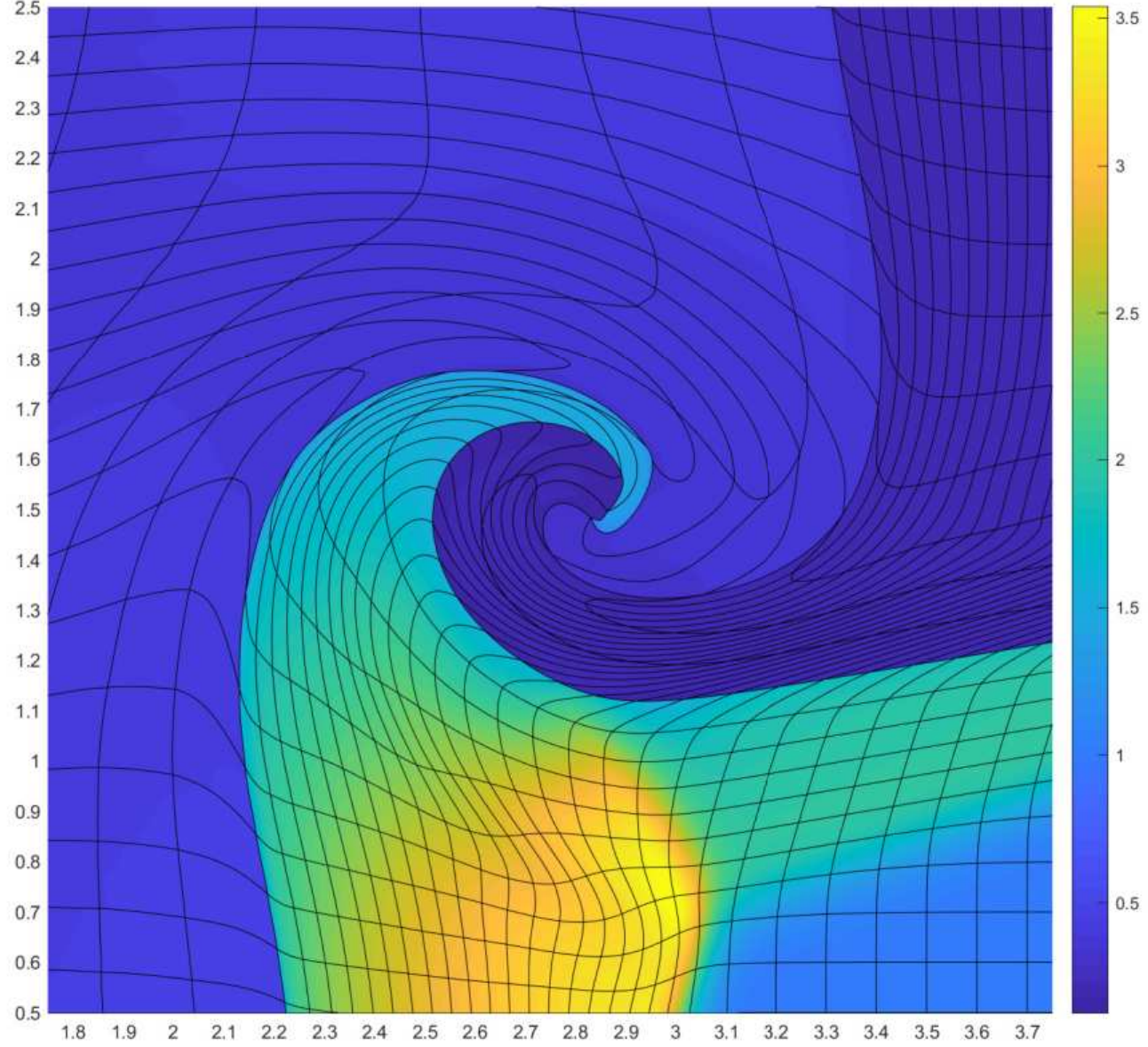}
    \caption{$Q^3-Q^2$ element pair.}
    \label{triple_Q3_zoom}
  \end{subfigure}
  \caption{Zoomed-in views of the density contours highlighting the vortex roll-up induced by the Kelvin-Helmholtz instability at $t = 2.5$.}
  \label{triple_zoom}
\end{figure}

Figure~\ref{triple_zoom} isolates a zoomed-in perspective of the vortex generation zone driven by the Kelvin-Helmholtz instability. The mesh lines trace the rotational flow field in both cases. Notably, while the $Q^2-Q^1$ elements capture the general topological features of the roll-up, the $Q^3-Q^2$ solution naturally resolves a significantly tighter, more pronounced, and deeply coiled twisting structure at the vortex core. This highlights that the enriched kinematic space and higher resolution of the $Q^3-Q^2$ elements provide a superior capability for capturing fine-scale flow structures and complex fluid instabilities.

\section{Conclusion}\label{sec:conclusion}
In this paper, we introduced a high-order staggered Lagrangian hydrodynamics (SGH) framework, extending established high-order curvilinear finite element methods~\cite{Dobrev2012High} into a strictly thermodynamically consistent context. We demonstrated that the classical compatible hydrodynamics algorithm is seamlessly recovered from this unified framework simply by restricting it to the $Q^1-P^0$ finite element space pair. A key theoretical contribution of this work is the rigorous clarification of the relationship between selected Gauss quadrature rules and the discrete mass conservation law, particularly concerning the exact determination of the density degrees of freedom and the elimination of equation-of-state (EOS) update ambiguities. To facilitate robust and highly efficient code implementation, we derived explicit matrix formulations for the Jacobian determinants and the shape function derivatives, significantly reducing computational overhead in high-order scenarios. 

A primary advantage of the proposed SGH architecture is its strategic collocation of kinematic and thermodynamic degrees of freedom. By pairing specific quadrature rules with these nodal definitions, the framework naturally yields perfectly diagonal mass and energy matrices. This mass-lumping strategy completely eliminates the need to assemble and solve expensive global linear algebraic systems during the explicit momentum and energy updates. Furthermore, by introducing a high-order generalization of the subzonal pressure method through an enriched pressure field, we successfully derived a restorative anti-hourglass force that stabilizes high-order spaces against unphysical zero-energy modes. 

Finally, a comprehensive suite of numerical benchmarks confirmed the robustness and accuracy of the proposed method. Smooth test cases validated the optimal formal convergence rates and superior computational efficiency of the $Q^m-Q^{m-1}$ spaces, while severe shock-driven problems demonstrated that our framework—when equipped with tensor artificial viscosity and a tunable anti-hourglass control parameter—maintains strict energy conservation and excellent grid quality under extreme deformations. Crucially, compared to legacy formulations that suffer from excessive numerical stiffness and mesh rigidity near shock fronts, our tunable framework successfully achieves an optimal balance between mesh regularity and solution fidelity.

\section*{Acknowledgment}
The authors would like to thank the anonymous referees. They have very
constructively helped to improve the original version of this paper. 

The research is supported by the National Natural Science Foundation of China (Grant No. 12571414) and Science Challenge Project (Grant No. TZ2025016).

\appendix
\section{Detailed IMEX Runge-Kutta Stages}\label{app:imex_rk}

This appendix provides the complete stage-by-stage formulations for the $s$-stage IMEX RK schemes employed in Section~\ref{sec:numerical_results}. As introduced previously, these IMEX RK schemes~\cite{Sandu2021Conservative} are specifically designed to guarantee exact total energy conservation. 

While the second-order RK2-average scheme can be implemented directly without defining an auxiliary velocity vector, the general $s$-stage IMEX RK scheme for $s \ge 3$ necessitates the definition of an auxiliary velocity vector $\vec{v} \equiv \vec{u}$. The parameters of the IMEX RK scheme are classified into implicit and explicit categories, denoted respectively by the implicit weights $b_1, b_2, \dots, b_s$ and the lower-triangular explicit coefficient matrix entries:
\begin{equation}
\begin{bmatrix}
  0 & 0 & 0 & 0 \\
  a_{2,1}^{\{E\}} & 0 & 0 & 0 \\
  a_{3,1}^{\{E\}} & a_{3,2}^{\{E\}} & 0 & 0 \\
  a_{4,1}^{\{E\}} & a_{4,2}^{\{E\}} & a_{4,3}^{\{E\}} & 0
\end{bmatrix}  
\end{equation}
for the 3-stage method, alongside the corresponding weights $b_1 \equiv 0, b_2, b_3, b_4$.

We describe the $s$-stage IMEX RK integration from time $t^{n}$ to $t^{n+1}$, denoting the intermediate velocity, position, and internal energy for the $i$-th stage as $\vec{v}_{i}, \vec{u}_{i}, \vec{x}_{i}, e_{i}$ for $i=1, 2, \dots, s$. The 3-stage, third-order IMEX RK scheme unfolds as follows. The first stage simply assigns the $t^n$ values to the intermediate variables:
\begin{equation*}
\begin{bmatrix}
\vec{v}_{1}\\
\vec{u}_{1}\\
\vec{x}_{1}\\
e_{1}
\end{bmatrix}=
\begin{bmatrix}
\vec{u}^{n}\\
\vec{u}^{n}\\
\vec{x}^{n}\\
e^{n}
\end{bmatrix}.
\end{equation*}
We then evaluate the force $\vec{f}_1$ using these initial intermediate variables and update the second stage as:
\begin{equation*}
\begin{bmatrix}
  \vec{u}_{2}\\
  \vec{x}_{2}\\
  e_{2}
\end{bmatrix}=
\begin{bmatrix}
  \vec{u}^n\\
  \vec{x}^n\\
  e^n
\end{bmatrix}
+
a_{2,1}^{\{E\}} \Delta t 
\begin{bmatrix}
  \frac{\vec{f}_1}{M_K}\\
  \vec{v}_1\\
  \frac{\vec{f}_1\cdot \vec{v}_1 }{M_T}
\end{bmatrix}.
\end{equation*}
Next, we calculate the force $\vec{f}_2$ using the updated intermediate variables $\vec{u}_{2},\vec{x}_{2},e_{2}$ and explicitly update the auxiliary velocity:
\begin{equation*}
  \vec{v}_2=\vec{v}_1+\frac{1}{2}\Delta t b_2 \frac{\vec{f}_2}{M_K}.
\end{equation*}

The third stage is updated as:
\begin{equation*}
\begin{bmatrix}
  \vec{u}_{3}\\
  \vec{x}_{3}\\
  e_{3}
\end{bmatrix}=
\begin{bmatrix}
  \vec{u}^n\\
  \vec{x}^n\\
  e^n
\end{bmatrix}
+
a_{3,1}^{\{E\}} \Delta t 
\begin{bmatrix}
  \frac{\vec{f}_1}{M_K}\\
  \vec{v}_1\\
  \frac{\vec{f}_1\cdot \vec{v}_1 }{M_T}
\end{bmatrix}
+
a_{3,2}^{\{E\}} \Delta t 
\begin{bmatrix}
  \frac{\vec{f}_2}{M_K}\\
  \vec{v}_2\\
  \frac{\vec{f}_2\cdot \vec{v}_2 }{M_T}
\end{bmatrix}.
\end{equation*}
We then calculate the force $\vec{f}_3$ using the variables $\vec{u}_{3},\vec{x}_{3},e_{3}$, updating the auxiliary velocity:
\begin{equation*}
  \vec{v}_3=\vec{v}_2+\frac{1}{2}\Delta t \left( b_2 \frac{\vec{f}_2}{M_K} +b_3 \frac{\vec{f}_3}{M_K}\right).
\end{equation*}

The fourth stage is updated as:
\begin{equation*}
\begin{bmatrix}
  \vec{u}_{4}\\
  \vec{x}_{4}\\
  e_{4}
\end{bmatrix}=
\begin{bmatrix}
  \vec{u}^n\\
  \vec{x}^n\\
  e^n
\end{bmatrix}
+
a_{4,1}^{\{E\}} \Delta t 
\begin{bmatrix}
  \frac{\vec{f}_1}{M_K}\\
  \vec{v}_1\\
  \frac{\vec{f}_1\cdot \vec{v}_1 }{M_T}
\end{bmatrix}
+
a_{4,2}^{\{E\}} \Delta t 
\begin{bmatrix}
  \frac{\vec{f}_2}{M_K}\\
  \vec{v}_2\\
  \frac{\vec{f}_2\cdot \vec{v}_2 }{M_T}
\end{bmatrix}
+
a_{4,3}^{\{E\}} \Delta t 
\begin{bmatrix}
  \frac{\vec{f}_3}{M_K}\\
  \vec{v}_3\\
  \frac{\vec{f}_3\cdot \vec{v}_3 }{M_T}
\end{bmatrix}.
\end{equation*}
Finally, we calculate the intermediate force $\vec{f}_4$ using $\vec{u}_{4},\vec{x}_{4},e_{4}$, and update the auxiliary velocity:
\begin{equation*}
  \vec{v}_4=\vec{v}_3+\frac{1}{2}\Delta t \left( b_3 \frac{\vec{f}_3}{M_K} +b_4 \frac{\vec{f}_4}{M_K}\right).
\end{equation*}

The final macroscopic update from $t^n$ to $t^{n+1}$ is computed as:
\begin{equation}
\begin{bmatrix}
  \vec{x}_{n+1}\\
  e_{n+1}
\end{bmatrix}=
\begin{bmatrix}
  \vec{x}^n\\
  e^n
\end{bmatrix}
+
b_2 \Delta t 
\begin{bmatrix}
  \vec{v}_2\\
  \frac{\vec{f}_2\cdot \vec{v}_2 }{M_T}
\end{bmatrix}
+
b_3 \Delta t 
\begin{bmatrix}
  \vec{v}_3\\
  \frac{\vec{f}_3\cdot \vec{v}_3 }{M_T}
\end{bmatrix}
+
b_4 \Delta t 
\begin{bmatrix}
  \vec{v}_4\\
  \frac{\vec{f}_4\cdot \vec{v}_4 }{M_T}
\end{bmatrix}.
\end{equation}
The corresponding velocity $\vec{u}^{n+1}$ update is straightforwardly given by:
\begin{equation}
\vec{u}^{n+1}=\vec{v}_4+\frac{1}{2}\Delta t b_4 \frac{\vec{f}_4}{M_K}.
\end{equation}

For the 4-stage, fourth-order IMEX RK scheme, the implicit weights are $b_1 \equiv 0, b_2, b_3, b_4, b_5$, and the explicit coefficient matrix is expanded to a $5 \times 5$ lower-triangular matrix:
\begin{equation}
\begin{bmatrix}
  0 & 0 & 0 & 0 & 0 \\
  a_{2,1}^{\{E\}} & 0 & 0 & 0 & 0 \\
  a_{3,1}^{\{E\}} & a_{3,2}^{\{E\}} & 0 & 0 & 0 \\
  a_{4,1}^{\{E\}} & a_{4,2}^{\{E\}} & a_{4,3}^{\{E\}} & 0 & 0 \\
  a_{5,1}^{\{E\}} & a_{5,2}^{\{E\}} & a_{5,3}^{\{E\}} & a_{5,4}^{\{E\}} & 0
\end{bmatrix}.  
\end{equation}

The execution of the 4-stage, fourth-order IMEX RK scheme proceeds identically in structure. The first stage assigns the $t^n$ values:
\begin{equation*}
\begin{bmatrix}
\vec{v}_{1}\\
\vec{u}_{1}\\
\vec{x}_{1}\\
e_{1}
\end{bmatrix}=
\begin{bmatrix}
\vec{u}^{n}\\
\vec{u}^{n}\\
\vec{x}^{n}\\
e^{n}
\end{bmatrix}.
\end{equation*}
We compute $\vec{f}_1$ and update the second stage:
\begin{equation*}
\begin{bmatrix}
  \vec{u}_{2}\\
  \vec{x}_{2}\\
  e_{2}
\end{bmatrix}=
\begin{bmatrix}
  \vec{u}^n\\
  \vec{x}^n\\
  e^n
\end{bmatrix}
+
a_{2,1}^{\{E\}} \Delta t 
\begin{bmatrix}
  \frac{\vec{f}_1}{M_K}\\
  \vec{v}_1\\
  \frac{\vec{f}_1\cdot \vec{v}_1 }{M_T}
\end{bmatrix}.
\end{equation*}
We compute $\vec{f}_2$ and update the auxiliary velocity:
\begin{equation*}
  \vec{v}_2=\vec{v}_1+\frac{1}{2}\Delta t b_2 \frac{\vec{f}_2}{M_K}.
\end{equation*}

The third stage is:
\begin{equation*}
\begin{bmatrix}
  \vec{u}_{3}\\
  \vec{x}_{3}\\
  e_{3}
\end{bmatrix}=
\begin{bmatrix}
  \vec{u}^n\\
  \vec{x}^n\\
  e^n
\end{bmatrix}
+
a_{3,1}^{\{E\}} \Delta t 
\begin{bmatrix}
  \frac{\vec{f}_1}{M_K}\\
  \vec{v}_1\\
  \frac{\vec{f}_1\cdot \vec{v}_1 }{M_T}
\end{bmatrix}
+
a_{3,2}^{\{E\}} \Delta t 
\begin{bmatrix}
  \frac{\vec{f}_2}{M_K}\\
  \vec{v}_2\\
  \frac{\vec{f}_2\cdot \vec{v}_2 }{M_T}
\end{bmatrix}.
\end{equation*}
We compute $\vec{f}_3$ and update the auxiliary velocity:
\begin{equation*}
  \vec{v}_3=\vec{v}_2+\frac{1}{2}\Delta t \left( b_2 \frac{\vec{f}_2}{M_K} +b_3 \frac{\vec{f}_3}{M_K}\right).
\end{equation*}

The fourth stage is:
\begin{equation*}
\begin{bmatrix}
  \vec{u}_{4}\\
  \vec{x}_{4}\\
  e_{4}
\end{bmatrix}=
\begin{bmatrix}
  \vec{u}^n\\
  \vec{x}^n\\
  e^n
\end{bmatrix}
+
a_{4,1}^{\{E\}} \Delta t 
\begin{bmatrix}
  \frac{\vec{f}_1}{M_K}\\
  \vec{v}_1\\
  \frac{\vec{f}_1\cdot \vec{v}_1 }{M_T}
\end{bmatrix}
+
a_{4,2}^{\{E\}} \Delta t 
\begin{bmatrix}
  \frac{\vec{f}_2}{M_K}\\
  \vec{v}_2\\
  \frac{\vec{f}_2\cdot \vec{v}_2 }{M_T}
\end{bmatrix}
+
a_{4,3}^{\{E\}} \Delta t 
\begin{bmatrix}
  \frac{\vec{f}_3}{M_K}\\
  \vec{v}_3\\
  \frac{\vec{f}_3\cdot \vec{v}_3 }{M_T}
\end{bmatrix}.
\end{equation*}
We compute $\vec{f}_4$ and update the auxiliary velocity:
\begin{equation*}
  \vec{v}_4=\vec{v}_3+\frac{1}{2}\Delta t \left( b_3 \frac{\vec{f}_3}{M_K} +b_4 \frac{\vec{f}_4}{M_K}\right).
\end{equation*}

The fifth stage is updated as follows:
\begin{equation*}
\begin{bmatrix}
  \vec{u}_{5}\\
  \vec{x}_{5}\\
  e_{5}
\end{bmatrix}=
\begin{bmatrix}
  \vec{u}^n\\
  \vec{x}^n\\
  e^n
\end{bmatrix}
+
\sum_{j=1}^{4} a_{5,j}^{\{E\}} \Delta t 
\begin{bmatrix}
  \frac{\vec{f}_j}{M_K}\\
  \vec{v}_j\\
  \frac{\vec{f}_j\cdot \vec{v}_j }{M_T}
\end{bmatrix}.
\end{equation*}
We compute the final intermediate force $\vec{f}_5$ and update the auxiliary velocity:
\begin{equation*}
  \vec{v}_5=\vec{v}_4+\frac{1}{2}\Delta t \left( b_4 \frac{\vec{f}_4}{M_K} +b_5 \frac{\vec{f}_5}{M_K}\right).
\end{equation*}

The final macroscopic update from $t^n$ to $t^{n+1}$ is computed as:
\begin{equation}
\begin{bmatrix}
  \vec{x}_{n+1}\\
  e_{n+1}
\end{bmatrix}=
\begin{bmatrix}
  \vec{x}^n\\
  e^n
\end{bmatrix}
+
\sum_{j=2}^{5} b_j \Delta t 
\begin{bmatrix}
  \vec{v}_j\\
  \frac{\vec{f}_j\cdot \vec{v}_j }{M_T}
\end{bmatrix},
\end{equation}
and the final velocity $\vec{u}^{n+1}$ update is explicitly given by:
\begin{equation}
\vec{u}^{n+1}=\vec{v}_5+\frac{1}{2}\Delta t b_5 \frac{\vec{f}_5}{M_K}.
\end{equation}

\bibliographystyle{abbrv}
\bibliography{ref}
\end{document}